\documentclass[10pt,a4paper]{article}
\usepackage[utf8]{inputenc}

\usepackage{mathrsfs,amsthm,xcolor,verbatim,bbm,amsmath,amsfonts,amssymb,nicefrac,enumitem,hyperref,bm,mathtools,xparse,etoolbox}

\numberwithin{equation}{section}
\usepackage{thmtools}
\usepackage[margin=1in]{geometry}
\usepackage[capitalise,nameinlink]{cleveref} 
\usepackage{tikz}
\usetikzlibrary{matrix,chains,positioning,decorations.pathreplacing,arrows,shapes,math,arrows.meta,calc, shapes.geometric}
\usepackage[most]{tcolorbox}

\makeatletter
\AtBeginEnvironment{tcolorbox}{%
   \def\@mpfn{footnote}%
}
\makeatother


\tikzset{
    neuron/.style={circle, draw, minimum size=0.5cm, text centered},
    input neuron/.style={neuron, fill=yellow!30},
    hidden neuron/.style={neuron, fill=blue!10, draw=blue, thick},
    output neuron/.style={neuron, fill=green!20},
    conn/.style={-stealth, thick},
    ellipsis/.style={font=\footnotesize},
}
\usepackage{capt-of}
\usepackage{stmaryrd} 

\crefname{enumi}{item}{items}
\crefname{equation}{}{}
\crefname{subsection}{Subsection}{Subsections}
\crefname{figure}{Figure}{Figures}

\usepackage{xurl}

\hypersetup{
    colorlinks,
    linkcolor={red!80!black},
    citecolor={green},
    urlcolor={blue!80!black}
}

\theoremstyle{plain}
\newtheorem{theorem}{Theorem}[section]

\newtheorem{prop}[theorem]{Proposition}

\newtheorem{definition}[theorem]{Definition}



\theoremstyle{definition}

\DeclareMathAlphabet{\mathpzc}{OT1}{pzc}{m}{it}

\DeclareFontEncoding{LS1}{}{}
\DeclareFontSubstitution{LS1}{stix}{m}{n}
\DeclareMathAlphabet{\mathscr}{LS1}{stixscr}{m}{n}
\DeclarePairedDelimiter{\pr}{(}{)}
\DeclarePairedDelimiter{\PR}{[}{]}

\newcommand{\muonprocess}{\MUON\ process\cfadd{defifinition: MUON process}}
\newcommand{\defmuon}[3]{$\alpha$-$\gamma$-$#2$-$\nscoe$-$#3$-\muonprocess\ $#1$ for $\smalll$}
\newcommand{\defmuonfull}[4]{$\alpha$-$\gamma$-$#4$-$\bfa$-$M$-\muonprocess\ $#1$ for $\smalll$ on $( \Omega, \cF, \P )$ with respect to $#2$ with data $#3$}

\newcommand{\E}{\mathbb{E}}
\renewcommand{\P}{\mathbb{P}}

\newcommand{\R}{\mathbb{R}}
\newcommand{\N}{\mathbb{N}}

\newcommand{\bbG}{\mathbb{G}}

\newcommand{\smalll}{\mathbb{L}}

\newcommand{\Tr}{\operatorname{Tr}}
\newcommand{\prb}[1]{\pr[\big]{ #1 }}
\newcommand{\PRb}[1]{\PR[\big]{ #1 }}

\makeatletter
\newcommand{\bvert}{%
  \mathord{%
    \mathchoice
      {\makebox[0.40em][c]{\scalebox{4.0}[1]{$\m@th\displaystyle\vert$}}}
      {\makebox[0.40em][c]{\scalebox{4.0}[1]{$\m@th\textstyle\vert$}}}
      {\makebox[0.38em][c]{\scalebox{4.0}[1]{$\m@th\scriptstyle\vert$}}}
      {\makebox[0.36em][c]{\scalebox{3.5}[1]{$\m@th\scriptscriptstyle\vert$}}}%
  }%
}
\makeatother

\newcommand{\nnorm}[1]{\mathopen{\bvert}#1\mathclose{\bvert}}

\renewcommand{\d}{ \mathrm{d}}

\renewcommand{\c}[1]{\mathfrak{c}^{#1}}

\newcommand{\NewtonSchulzAlgorithm}[2]{\Phi^{#1}_{#2}}
\newcommand{\TNewtonSchulzAlgorithm}[3]{\operatorname{NS}_{#1}^{#2}(#3)}
\newcommand{\TbigNewtonSchulzAlgorithm}[3]{\operatorname{NS}_{#1}^{#2}\bigl(#3\bigr)}
\newcommand{\nscoe}{b}

\newcommand{\cA}{\mathcal{A}}
\newcommand{\cB}{\mathcal{B}}

\newcommand{\cF}{\mathcal{F}}
\newcommand{\cG}{\mathcal{G}}

\newcommand{\cK}{\mathcal{K}}
\newcommand{\cL}{\mathcal{L}}

\newcommand{\cO}{\mathcal{O}}

\newcommand{\cU}{\mathcal{U}}
\newcommand{\cV}{\mathcal{V}}

\newcommand{\bfa}{\mathbf{a}}

\newcommand{\bfd}{\mathbf{d}}
\newcommand{\bfe}{\mathbf{e}}

\newcommand{\bfi}{\mathbf{i}}

\newcommand{\bfm}{\mathbf{m}}

\newcommand{\bfx}{\mathbf{x}}

\newcommand{\sgn}{\Phi}

\newcommand{\bfA}{\mathbf{A}}
\newcommand{\bfB}{\mathbf{B}}
\newcommand{\bfC}{\mathbf{C}}

\newcommand{\bfK}{\mathbf{K}}

\newcommand{\bfM}{\mathbf{M}}
\newcommand{\bfN}{\mathbf{N}}

\newcommand{\bfX}{\mathbf{X}}

\newcommand{\scrC}{\mathscr{C}}
\newcommand{\scrc}{\mathscr{c}}

\newcommand{\g}{\|\mathfrak{g}\|_{Lip}}

\newcommand{\fC}{\mathfrak{C}}

\newcommand{\fL}{\mathfrak{L}}
\newcommand{\fM}{\mathfrak{M}}
\newcommand{\fN}{\mathfrak{N}}

\newcommand{\bbX}{\mathbb{X}}

\newcommand{\bbC}{\mathbb{C}}

\newcommand{\fc}{\mathfrak{c}}
\newcommand{\fd}{\mathfrak{d}}

\newcommand{\fn}{\mathfrak{n}}

\newcommand{\imagine}{\operatorname{Im}}
\newcommand{\imag}{\mathbf{i}}

\newcommand{\dimX}{\mathscr{d}}

\renewcommand{\emptyset}{\varnothing}

\newcommand{\setX}{U}

\newcommand\restr[2]{{
  \left.\kern-\nulldelimiterspace 
  #1 
  \vphantom{|} 
  \right|_{#2} 
  }}

\DeclarePairedDelimiter{\spro}{\langle}{\rangle}

\newcommand{\qandq}{\quad \text{and} \quad }
\newcommand{\qqandqq}{\qquad\text{and}\qquad}

\newcommand{\bbM}{\mathbb{M}}

\makeatletter
\newcommand{\doublefootnotemark}{%
  \refstepcounter{footnote}%
  \Hy@raisedlink{\hyper@anchorstart{Hfootnote.\theHfootnote}\hyper@anchorend}%
  \edef\firstfootnote{\thefootnote}%
  \edef\firstHfootnote{\theHfootnote}%
  \refstepcounter{footnote}%
  \Hy@raisedlink{\hyper@anchorstart{Hfootnote.\theHfootnote}\hyper@anchorend}%
  \edef\secondfootnote{\thefootnote}%
  \edef\secondHfootnote{\theHfootnote}%
  \textsuperscript{\normalfont
    \hyperlink{Hfootnote.\firstHfootnote}{\firstfootnote},%
    \hyperlink{Hfootnote.\secondHfootnote}{\secondfootnote}}%
}
\makeatother

\usepackage{cleveref}

\ExplSyntaxOn

\seq_new:N \g_cflist_loaded
\seq_new:N \g_cflist_pending

\NewDocumentCommand{\cfadd} { m } {
  \seq_if_in:NnF \g_cflist_loaded { #1 } {
    \seq_if_in:NnF \g_cflist_pending { #1 } {
      \seq_gput_right:Nn \g_cflist_pending { #1 }
    }
  }
}

\NewDocumentCommand{\cfconsiderloaded} { m } {
  \seq_gput_right:Nn \g_cflist_loaded {#1}
}

\NewDocumentCommand{\cfremove} { m } {
  \seq_gremove_all:Nn \g_cflist_pending { #1 }
}

\NewDocumentCommand{\cfload} { o } {
  \seq_if_empty:NTF \g_cflist_pending {
    \IfValueTF{#1}{\ignorespaces}{\unskip}
  } {
    (cf.\ \cref{\seq_use:Nn \g_cflist_pending {,}})\IfValueTF{#1}{#1~}{\unskip}
    \seq_gconcat:NNN \g_cflist_loaded \g_cflist_loaded \g_cflist_pending
    \seq_gclear:N \g_cflist_pending
    \IfValueT{#1}{\ignorespaces}
  }
}

\NewDocumentCommand{\cfclear} {} {
  \seq_gclear:N \g_cflist_loaded
  \seq_gclear:N \g_cflist_pending
}

\NewDocumentCommand{\cfout} { o } {
  \seq_if_empty:NTF \g_cflist_pending {\unskip\IfValueT{#1}{\ignorespaces}} {
    (cf.\ \cref{\seq_use:Nn \g_cflist_pending {,}})\IfValueTF{#1}{#1~}{\unskip}
    \seq_gclear:N \g_cflist_pending
    \IfValueT{#1}{\ignorespaces}
  }
}

\NewDocumentCommand{\ifnocf} { m } {
  \seq_if_empty:NT \g_cflist_pending { #1 }
}

\ExplSyntaxOff

\ExplSyntaxOn

\bool_new:N \g_noteobserve

\NewDocumentCommand{\setnote}{}{
  \bool_gset_true:N \g_noteobserve
}

\NewDocumentCommand{\setobserve}{}{
  \bool_gset_false:N \g_noteobserve
}

\NewDocumentCommand{\nobs}{ o }{
  \IfValueT{#1}{
    \str_if_eq:noTF {note} {#1} {
      \bool_gset_true:N \g_noteobserve
    } {
      \str_if_eq:noTF {Note} {#1} {
        \bool_gset_true:N \g_noteobserve
      } {
        \bool_gset_false:N \g_noteobserve
      }
    }
  }
  \bool_if:nTF { \g_noteobserve } {
    \bool_gset_false:N \g_noteobserve
    note
  } {
    \bool_gset_true:N \g_noteobserve
    observe
  }
  \IfValueF{#1}{~}
}

\NewDocumentCommand{\Nobs}{ o }{
  \IfValueT{#1}{
    \str_if_eq:noTF {note} {#1} {
      \bool_gset_true:N \g_noteobserve
    } {
      \str_if_eq:noTF {Note} {#1} {
        \bool_gset_true:N \g_noteobserve
      } {
        \bool_gset_false:N \g_noteobserve
      }
    }
  }
  \bool_if:nTF { \g_noteobserve } {
    \bool_gset_false:N \g_noteobserve
    Note
  } {
    \bool_gset_true:N \g_noteobserve
    Observe
  }
  \IfValueF{#1}{~}
}

\ExplSyntaxOff

\ExplSyntaxOn

\bool_new:N \g_hencetherefore

\NewDocumentCommand{\hence}{ o }{
  \IfValueT{#1}{
    \str_if_eq:noTF {hence} {#1} {
      \bool_gset_true:N \g_hencetherefore
    } {
      \str_if_eq:noTF {Hence} {#1} {
        \bool_gset_true:N \g_hencetherefore
      } {
        \bool_gset_false:N \g_hencetherefore
      }
    }
  }
  \bool_if:nTF { \g_hencetherefore } {
    \bool_gset_false:N \g_hencetherefore
    hence
  } {
    \bool_gset_true:N \g_hencetherefore
    therefore
  }
  \IfValueF{#1}{~}
}

\NewDocumentCommand{\Hence}{ o }{
  \IfValueT{#1}{
    \str_if_eq:noTF {hence} {#1} {
      \bool_gset_true:N \g_hencetherefore
    } {
      \str_if_eq:noTF {Hence} {#1} {
        \bool_gset_true:N \g_hencetherefore
      } {
        \bool_gset_false:N \g_hencetherefore
      }
    }
  }
  \bool_if:nTF { \g_hencetherefore } {
    \bool_gset_false:N \g_hencetherefore
    Hence,~we~obtain
  } {
    \bool_gset_true:N \g_hencetherefore
    Therefore,~we~obtain
  }
  \IfValueF{#1}{~}
}

\ExplSyntaxOff

\ExplSyntaxOn

\seq_const_from_clist:Nn \g_prove_mru {
  establish,
  demonstrate,
  prove,
  show,
  imply,
  ensure
}

\prop_new:N \l__verbs
\prop_put:Nnn \l__verbs {show} {shows}
\prop_put:Nnn \l__verbs {imply} {implies}
\prop_put:Nnn \l__verbs {demonstrate} {demonstrates}
\prop_put:Nnn \l__verbs {prove} {proves}
\prop_put:Nnn \l__verbs {establish} {establishes}
\prop_put:Nnn \l__verbs {ensure} {ensures}
\prop_put:Nnn \l__verbs {assure} {assures}

\tl_new:N \g_wordtmp
\seq_new:N \l_mytmps

\cs_generate_variant:Nn \str_if_in:nnTF { nVTF }
\cs_generate_variant:Nn \str_if_in:nnTF { xVTF }

\NewDocumentCommand{\prove}{ o }{
  \IfValueTF{#1}{
    \seq_clear:N \l_mytmps
    \seq_map_inline:Nn \g_prove_mru {
      \str_if_eq:nnTF {##1} {ensure} {
        \str_set:Nn \l_temps {n}
      } {
        \str_set:Nx \l_temps {\str_head_ignore_spaces:n {##1}}
      }
      \str_if_in:xVTF {#1} \l_temps {
        \seq_put_right:Nn \l_mytmps {##1}
      } { }
    }
    \seq_get_right:NN \l_mytmps \g_wordtmp
  } {
    \seq_get_right:NN \g_prove_mru \g_wordtmp
  }
  \tl_use:N \g_wordtmp
  \IfValueTF{#1}{}{~}
  \seq_gput_left:NV \g_prove_mru \g_wordtmp
  \seq_gremove_duplicates:N \g_prove_mru
}

\NewDocumentCommand{\proves}{ o }{
  \IfValueTF{#1}{
    \seq_clear:N \l_mytmps
    \seq_map_inline:Nn \g_prove_mru {
      \str_if_eq:nnTF {##1} {ensure} {
        \str_set:Nn \l_temps {n}
      } {
        \str_set:Nx \l_temps {\str_head_ignore_spaces:n {##1}}
      }
      \str_if_in:xVTF {#1} \l_temps {
        \seq_put_right:Nn \l_mytmps {##1}
      } { }
    }
    \seq_get_right:NN \l_mytmps \g_wordtmp
  } {
    \seq_get_right:NN \g_prove_mru \g_wordtmp
  }
  \str_set:NV \l_tmpa_str \g_wordtmp
  \prop_get:NVN \l__verbs \l_tmpa_str \l_tmpa_tl
  \tl_use:N \l_tmpa_tl
  \IfValueTF{#1}{}{~}
  \seq_gput_left:NV \g_prove_mru \g_wordtmp
  \seq_gremove_duplicates:N \g_prove_mru
}

\newcommand{\llabel}[1]{\savelabel{#1}\label{\loc.#1}\ignorespaces}

\clist_new:N \l_localreflist
\clist_new:N \l_reflist

\NewDocumentCommand{\lref} { m } {
  \clist_set:No \l_localreflist {#1}
  \clist_clear:N \l_reflist
  \clist_map_inline:Nn \l_localreflist { \clist_put_right:Nn \l_reflist {\loc.##1} }
  \cref{\l_reflist}
}

\NewDocumentCommand{\Lref} { m } {
  \clist_set:No \l_localreflist {#1}
  \clist_clear:N \l_reflist
  \clist_map_inline:Nn \l_localreflist { \clist_put_right:Nn \l_reflist {\loc.##1} }
  \Cref{\l_reflist}
}

\NewDocumentCommand{\itref}{ m m }{
  \clist_set:No \l_localreflist {#2}
  \clist_clear:N \l_reflist
  \clist_map_inline:Nn \l_localreflist { \clist_put_right:Nn \l_reflist {#1.##1} }
  \cref{\l_reflist}~in~\cref{#1}
}

\seq_new:N \l_enum_seq
\int_new:N \l_num_items

\bool_new:N \g_commaused_bool

\providecommand{\comma}{}

\cs_new:Nn \enum_it:nn {
  \int_case:nnF {\l_num_items - #1} {
    {0} {
      \renewcommand{\comma}{}
      #2\space
    }
    {1} {
      \bool_gset_false:N \g_commaused_bool
      \renewcommand{\comma}{,~\bool_gset_true:N \g_commaused_bool}
      #2
      \bool_if:NTF \g_commaused_bool {} {,~}
      and~
    }
  } {
    \bool_gset_false:N \g_commaused_bool
    \renewcommand{\comma}{,~\bool_gset_true:N \g_commaused_bool}
    #2
    \bool_if:NTF \g_commaused_bool {} {,~}
  }
}

\cs_new:Nn \enum_it_U:nn {
  \int_case:nnF {\l_num_items - #1} {
    {0} {
      \renewcommand{\comma}{}
      #2
      \space
    }
    {1} {
      \bool_gset_false:N \g_commaused_bool
      \renewcommand{\comma}{,~\bool_gset_true:N \g_commaused_bool}
      #2
      \bool_if:NTF \g_commaused_bool {} {,~}
      and~
    }
  } {
    \bool_gset_false:N \g_commaused_bool
    \renewcommand{\comma}{,~\bool_gset_true:N \g_commaused_bool}
    \int_compare:nTF {#1=1} {
      \text_titlecase_first:n {#2}
    } {
      #2
    }
    \bool_if:NTF \g_commaused_bool {} {,~}
  }
}

\cs_generate_variant:Nn \tl_if_eq:nnTF {onTF}

\cs_new:Nn \enum:nnnn {
  \seq_set_split:Nnn \l_enum_seq ; {#1}
  \seq_remove_all:Nn \l_enum_seq { }
  \seq_remove_all:Nn \l_enum_seq {#2}
  \seq_log:N \l_enum_seq
  \int_set:Nn \l_num_items {\seq_count:N \l_enum_seq}
  \int_log:N \l_num_items
  \int_case:nnF {\l_num_items} {
    { 0 } { 0 }
    { 1 } {
      \IfBooleanTF{#4} {
        \tl_set:Nn \l_text_case_exclude_arg_tl {\cref}
        \bool_set_false:N \l_text_titlecase_check_letter_bool
        \text_titlecase_first:n {\seq_use:Nn \l_enum_seq {}}
      } {
        \seq_use:Nn \l_enum_seq {}
      }
      \space
      \tl_if_eq:onTF{#3}{-}{}{
        \bool_if:NTF \l_plural_bool {
          \prove[#3]~
        } {
          \proves[#3]~
        }
      }
    }
    { 2 } {
      \IfBooleanTF{#4} {
        \tl_set:Nn \l_text_case_exclude_arg_tl {\cref}
        \bool_set_false:N \l_text_titlecase_check_letter_bool
        \text_titlecase_first:n {\seq_use:Nn \l_enum_seq {~and~}}
      } {
        \seq_use:Nn \l_enum_seq {~and~}
      }
      \space
      \tl_if_eq:onTF{#3}{-}{}{
        \prove[#3]~
      }
    }
  } {
    \IfBooleanTF{#4} {
      \tl_set:Nn \l_text_case_exclude_arg_tl {\cref}
      \seq_indexed_map_function:NN \l_enum_seq \enum_it_U:nn
    } {
      \seq_indexed_map_function:NN \l_enum_seq \enum_it:nn
    }
    \tl_if_eq:onTF{#3}{-}{}{
      \prove[#3]~
    }
  }
}

\cs_generate_variant:Nn \enum:nnnn {nxnn}
\cs_generate_variant:Nn \enum:nnnn {nxxn}

\NewDocumentCommand{\enum}{O{} m O{-} s}{
  \IfBooleanTF{#4}{
    \enum:nxnn {#2} {#1} {sindep} \BooleanFalse
  } {
    \enum:nxxn {#2} {#1} {#3} \BooleanFalse
  }
}

\NewDocumentCommand{\dott}{}{\ifnocf{.}\space}

\bool_new:N \g_arg_start_bool
\bool_gset_true:N \g_arg_start_bool

\NewDocumentCommand{\startnewargseq}{}{\bool_gset_true:N \g_arg_start_bool \tl_set:Nn \g_label_tl {}}

\cs_generate_variant:Nn \seq_if_in:NnTF {NxTF}
\cs_generate_variant:Nn \seq_remove_all:Nn {Nx}

\int_new:N \l_random_int

\cs_generate_variant:Nn \tl_if_head_eq_catcode:nNTF {oNTF}
\cs_generate_variant:Nn \tl_if_head_eq_catcode:nNTF {VNTF}

\cs_generate_variant:Nn \tl_if_head_eq_catcode:nNTF {eNTF}
\cs_generate_variant:Nn \tl_log:n {o}
\cs_generate_variant:Nn \tl_log:n {f}
\cs_generate_variant:Nn \tl_log:n {x}
\cs_generate_variant:Nn \tl_log:n {e}

\cs_generate_variant:Nn \tl_if_in:nnTF {onTF}
\cs_generate_variant:Nn \tl_if_in:NnTF {NeTF}

\cs_generate_variant:Nn \tl_if_head_eq_meaning:nNTF {VNTF}

\seq_const_from_clist:Nn \g_arg_mru_this {
  Ahpr,
  Tapr,
  Ctapr,
  H
}

\seq_const_from_clist:Nn \g_arg_mru_nothis {
  Ia,
  Nwc,
  N,
  Itns,
  Fm,
  Itnswc,
  Mo
}

\seq_new:N \l_arg_seq
\tl_new:N \l_cons_tl
\tl_new:N \l_dummy_tl

\bool_new:N \g_debug_bool
\bool_gset_false:N \g_debug_bool

\bool_new:N \l_insidearg_bool

\bool_new:N \g_firstargletter_bool

\sys_gset_rand_seed:n {0903}

\bool_new:N \l_plural_bool
\tl_new:N \l_arg_verbs_tl

\NewDocumentCommand{\argument}{mom}{
\color{blue}
  \bool_set_false:N \l_plural_bool
  \tl_set:Nn \l_arg_verbs_tl {sindep}
  \keys_define:nn { benno/argument } {
    plural .value_forbidden:n = true,
    plural .code:n = {\bool_set_true:N \l_plural_bool},
    verbs .value_required:n = false,
    verbs .tl_set:N = \l_arg_verbs_tl,
  }
  \IfValueT{#2}{
    \keys_set:nn { benno/argument } {#2}
  }
  \bool_log:N \l_plural_bool
  \bool_gset_true:N \l_insidearg_bool
  \seq_set_split:Nnn \l_arg_seq ; {#1}
  \seq_remove_all:Nn \l_arg_seq { }
  \seq_log:N \l_arg_seq
  \tl_set:Nn \l_cons_tl {#3}
  \tl_trim_spaces:N \l_cons_tl
  \seq_if_in:NxTF \l_arg_seq {\lref{\g_label_tl}} {
    \seq_remove_all:Nx \l_arg_seq {\lref{\g_label_tl}}
    \seq_get_left:NNTF \l_arg_seq \l_dummy_tl {
      \tl_trim_spaces:N \l_dummy_tl
      \bool_gset_false:N \g_firstargletter_bool
      \tl_if_head_eq_catcode:VNTF \l_dummy_tl a {
        \bool_gset_true:N \g_firstargletter_bool
      } {
        \tl_if_head_eq_meaning:VNTF \l_dummy_tl {\cref} {
          \tl_set:Nx \l_tmpa_tl {\tl_tail:N \l_dummy_tl}
          \tl_set:Nx \l_tmpb_tl {\tl_head:N \l_tmpa_tl}
          \bool_gset_true:N \g_firstargletter_bool
          \tl_if_in:NeTF \l_tmpb_tl {lem\c_colon_str} {} {
            \tl_if_in:NeTF \l_tmpb_tl {thm\c_colon_str} {} {
              \tl_if_in:NeTF \l_tmpb_tl {prop\c_colon_str} {} {
                \tl_if_in:NeTF \l_tmpb_tl {cor\c_colon_str} {} {
                  \bool_gset_false:N \g_firstargletter_bool
                }
              }
            }
          }
        } {
        }
      }
      \bool_if:NTF \g_firstargletter_bool {
        \seq_set_eq:NN \l_tmpa_seq \g_arg_mru_this
        \seq_remove_all:Nn \l_tmpa_seq {H}
        \seq_get_right:NN \l_tmpa_seq \l_tmpa_tl
        \int_case:nnF {\seq_count:N \l_arg_seq} {
          {1} {
            \str_case:VnF {\l_tmpa_tl} {
              {Ahpr} {
                \bool_if:NT \g_debug_bool {C1.1}
                \seq_gput_left:Nn \g_arg_mru_this {Ahpr}
                \seq_gremove_duplicates:N \g_arg_mru_this
                \enum:nxnn {#1} {\lref{\g_label_tl}} {-} {\BooleanTrue}
                \hence~
                \bool_if:NTF \l_plural_bool {
                  \prove[\l_arg_verbs_tl]~\ignorespaces #3
                } {
                  \proves[\l_arg_verbs_tl]~\ignorespaces #3
                }
              }
              {Tapr} {
                \bool_if:NT \g_debug_bool {C1.2}
                \seq_gput_left:Nn \g_arg_mru_this {Tapr}
                \seq_gremove_duplicates:N \g_arg_mru_this
                \enum[\lref{\g_label_tl}]{
                  This;
                  #1
                }[\l_arg_verbs_tl]\ignorespaces #3
              }
              {Ctapr} {
                \bool_if:NT \g_debug_bool {C1.3}
                \seq_gput_left:Nn \g_arg_mru_this {Ctapr}
                \seq_gremove_duplicates:N \g_arg_mru_this
                Combining~
                \enum[\lref{\g_label_tl}]{
                  this;
                  #1
                } \proves[\l_arg_verbs_tl]~\ignorespaces #3
              }
            } {}
          }
        } {
          \str_case:VnF {\l_tmpa_tl} {
             {Ahpr} {
              \bool_if:NT \g_debug_bool {C2.1}
              \seq_gput_left:Nn \g_arg_mru_this {Ahpr}
              \seq_gremove_duplicates:N \g_arg_mru_this
              \enum:nxnn {#1} {\lref{\g_label_tl}} {-} {\BooleanTrue}
              \hence~
              \prove[\l_arg_verbs_tl]~\ignorespaces #3
            }
            {Tapr} {
              \bool_if:NT \g_debug_bool {C2.2}
              \seq_gput_left:Nn \g_arg_mru_this {Tapr}
              \seq_gremove_duplicates:N \g_arg_mru_this
              \enum[\lref{\g_label_tl}]{
                This;
                #1
              }[\l_arg_verbs_tl]\ignorespaces #3
            }
            {Ctapr} {
              \int_case:nn {\int_rand:nn {0} {1}} {
                {0} {
                  \bool_if:NT \g_debug_bool {C2.3}
                  \seq_gput_left:Nn \g_arg_mru_this {Ctapr}
                  \seq_gremove_duplicates:N \g_arg_mru_this
                  Combining~
                  \enum[\lref{\g_label_tl}]{
                    this;
                    #1
                  } \proves[\l_arg_verbs_tl]~\ignorespaces #3
                }
                {1} {
                  \bool_if:NT \g_debug_bool {C2.4}
                  \seq_gput_left:Nn \g_arg_mru_this {Ctapr}
                  \seq_gremove_duplicates:N \g_arg_mru_this
                  Combining~
                  \enum:nxnn {#1} {\lref{\g_label_tl}} {-} {\BooleanFalse}
                  \hence~
                  \proves[\l_arg_verbs_tl]~\ignorespaces #3
                }
              }
            }
          } {}
        }
      } {
        \seq_set_eq:NN \l_tmpa_seq \g_arg_mru_this
        \seq_remove_all:Nn \l_tmpa_seq {H}
        \seq_remove_all:Nn \l_tmpa_seq {Ahpr}
        \seq_get_right:NN \l_tmpa_seq \l_tmpa_tl
        \int_case:nnF {\seq_count:N \l_arg_seq} {
          {1} {
            \str_case:VnF {\l_tmpa_tl} {
              {Tapr} {
                \bool_if:NT \g_debug_bool {C3.1}
                \seq_gput_left:Nn \g_arg_mru_this {Tapr}
                \seq_gremove_duplicates:N \g_arg_mru_this
                \enum[\lref{\g_label_tl}]{
                  This;
                  #1
                }[\l_arg_verbs_tl]\ignorespaces #3
              }
              {Ctapr} {
                \bool_if:NT \g_debug_bool {C3.2}
                \seq_gput_left:Nn \g_arg_mru_this {Ctapr}
                \seq_gremove_duplicates:N \g_arg_mru_this
                Combining~
                \enum[\lref{\g_label_tl}]{
                  this;
                  #1
                } \proves[\l_arg_verbs_tl]~\ignorespaces #3
              }
            } {}
          }
        } {
          \str_case:VnF {\l_tmpa_tl} {
            {Tapr} {
              \bool_if:NT \g_debug_bool {C4.1}
              \seq_gput_left:Nn \g_arg_mru_this {Tapr}
              \seq_gremove_duplicates:N \g_arg_mru_this
              \enum[\lref{\g_label_tl}]{
                This;
                #1
              }[\l_arg_verbs_tl]\ignorespaces #3		
            }
            {Ctapr} {
              \int_case:nn {\int_rand:nn {0} {1}} {
                {0} {
                  \bool_if:NT \g_debug_bool {C4.2}
                  \seq_gput_left:Nn \g_arg_mru_this {Ctapr}
                  \seq_gremove_duplicates:N \g_arg_mru_this
                  Combining~
                  \enum[\lref{\g_label_tl}]{
                    this;
                    #1
                  } \proves[\l_arg_verbs_tl]~\ignorespaces #3		
                }
                {1} {
                  \bool_if:NT \g_debug_bool {C4.3}
                  \seq_gput_left:Nn \g_arg_mru_this {Ctapr}
                  \seq_gremove_duplicates:N \g_arg_mru_this
                  Combining~
                  \enum:nxnn {#1} {\lref{\g_label_tl}} {-} {\BooleanFalse}
                  \hence~
                  \proves[\l_arg_verbs_tl]~\ignorespaces #3    
                }
              }
            }
          } {}
        }
      }
    } {
      \tl_if_head_eq_catcode:oNTF \l_cons_tl a {
        \seq_set_eq:NN \l_tmpa_seq \g_arg_mru_this
        \seq_remove_all:Nn \l_tmpa_seq {Ctapr}
        \seq_remove_all:Nn \l_tmpa_seq {Ahpr}
        \seq_get_right:NN \l_tmpa_seq \l_tmpa_tl
        \str_case:VnF {\l_tmpa_tl} {
          {H} {
            \bool_if:NT \g_debug_bool {C5.1}
            \seq_gput_left:Nn \g_arg_mru_this {H}
            \seq_gremove_duplicates:N \g_arg_mru_this
            Hence,~we~obtain~\ignorespaces #3
          }
          {Tapr} {
            \bool_if:NT \g_debug_bool {C5.2}
            \seq_gput_left:Nn \g_arg_mru_this {Tapr}
            \seq_gremove_duplicates:N \g_arg_mru_this
            This~\proves[\l_arg_verbs_tl]~\ignorespaces #3
          }
        } {}
      } {
        \bool_if:NT \g_debug_bool {C6.1}
        \seq_gput_left:Nn \g_arg_mru_this {Tapr}
        \seq_gremove_duplicates:N \g_arg_mru_this
        This~\proves[\l_arg_verbs_tl]~\ignorespaces #3
      }
    } 
  } {
    \int_compare:nNnTF {\seq_count:N \l_arg_seq} = {0} {
      \bool_if:NTF \g_arg_start_bool {
        \bool_if:NT \g_debug_bool {C7.1}
        \Nobs\unskip
        #3
      } {
        \bool_if:NT \g_debug_bool {C7.2}
        \Moreover~
        #3
      }
    } {
      \bool_if:NTF \g_arg_start_bool {
        \bool_if:NT \g_debug_bool {C8.1}
        \tl_log:N \l_arg_verbs_tl
        \Nobs~that~
        \enum{
          #1
        }[\l_arg_verbs_tl]\ignorespaces #3
      } {
        \int_compare:nNnTF {\seq_count:N \l_arg_seq} = {1} {
          \seq_set_eq:NN \l_tmpa_seq \g_arg_mru_nothis
          \seq_remove_all:Nn \l_tmpa_seq {Nwc}
          \seq_remove_all:Nn \l_tmpa_seq {Itnswc}
          \seq_get_right:NN \l_tmpa_seq \l_tmpa_tl
        } {
          \seq_get_right:NN \g_arg_mru_nothis \l_tmpa_tl
        }
        \str_case:VnF {\l_tmpa_tl} {
          {Mo} {
            \bool_if:NT \g_debug_bool {C9.1}
            \seq_gput_left:Nn \g_arg_mru_nothis {Mo}
            \seq_gremove_duplicates:N \g_arg_mru_nothis
            Moreover,~\nobs~that~
            \enum{
              #1
            }[\l_arg_verbs_tl]\ignorespaces #3		
          }
          {Fm} {
            \bool_if:NT \g_debug_bool {C9.2}
            \seq_gput_left:Nn \g_arg_mru_nothis {Fm}
            \seq_gremove_duplicates:N \g_arg_mru_nothis
            Furthermore,~\nobs~that~
            \enum{
              #1
            }[\l_arg_verbs_tl]\ignorespaces #3		
          }
          {Ia} {
            \bool_if:NT \g_debug_bool {C9.3}
            \seq_gput_left:Nn \g_arg_mru_nothis {Ia}
            \seq_gremove_duplicates:N \g_arg_mru_nothis
            In~addition,~\nobs~that~
            \enum{
              #1
            }[\l_arg_verbs_tl]\ignorespaces #3		
          }
          {N} {
            \bool_if:NT \g_debug_bool {C9.4}
            \seq_gput_left:Nn \g_arg_mru_nothis {N}
            \seq_gremove_duplicates:N \g_arg_mru_nothis
            Next,~\nobs~that~
            \enum{
              #1
            }[\l_arg_verbs_tl]\ignorespaces #3		
          }
          {Itns} {
            \bool_if:NT \g_debug_bool {C9.5}
            \seq_gput_left:Nn \g_arg_mru_nothis {Itnswc}
            \seq_gput_left:Nn \g_arg_mru_nothis {Itns}
            \seq_gremove_duplicates:N \g_arg_mru_nothis
            In~the~next~step~we~\nobs~that~
            \enum{
              #1
            }[\l_arg_verbs_tl]\ignorespaces #3		
          }
          {Nwc} {
            \bool_if:NT \g_debug_bool {C9.6}
            \seq_gput_left:Nn \g_arg_mru_nothis {Nwc}
            \seq_gremove_duplicates:N \g_arg_mru_nothis
            Next~we~combine~
            \enum{
              #1
            }to~obtain~\ignorespaces #3
          }
          {Itnswc} {
            \bool_if:NT \g_debug_bool {C9.7}
            \seq_gput_left:Nn \g_arg_mru_nothis {Itns}
            \seq_gput_left:Nn \g_arg_mru_nothis {Itnswc}
            \seq_gremove_duplicates:N \g_arg_mru_nothis
            In~the~next~step~we~combine~
            \enum{
              #1
            }to~obtain~\ignorespaces #3
          }
        } {}
      }
    }
  }
  \bool_gset_false:N \g_arg_start_bool
  \bool_gset_false:N \l_insidearg_bool
  \cfload[.]
  \color{black}
}

\tl_new:N \g_label_tl
\tl_gset:Nn \g_label_tl { }

\NewDocumentCommand{\savelabel}{m}{
  \bool_if:NTF \l_insidearg_bool {
    \tl_gset:Nn \g_label_tl {#1}
  } {
    \tl_gset:Nn \g_label_tl { }
  }
}

\ExplSyntaxOff

\ExplSyntaxOn

\NewDocumentEnvironment {athm} {m m o} {
\str_if_eq:noTF {example} {#1} {
  \bool_gset_true:N \g_example_bool
} {
  \bool_gset_false:N \g_example_bool
}
\cfclear
\IfNoValueTF{#3}{
\begin{#1}\label{#2}\global\def\loc{#2}
}{
\begin{#1}[#3]\label{#2}\global\def\loc{#2}
}
}{
\end{#1}
}

\NewDocumentEnvironment {adef} {m} {
\begin{definition}\label{#1}\global\def\loc{#1}
}{
\end{definition}
}

\NewDocumentEnvironment{aproof} {} {
\bool_if:NTF \g_example_bool {
  \bool_gset_true:N \g_arg_start_bool
  \begin{proof}[Proof~for~\cref{\loc}]
} {
  \bool_gset_true:N \g_arg_start_bool
  \begin{proof}[Proof~of~\cref{\loc}]
}
\bool_gset_false:N \g_finishproof_bool
}{
\bool_if:NTF \g_finishproof_bool {}
{\finishproofthus}
\end{proof}
}

\NewDocumentCommand{\finishproofthus} {} {
  \bool_gset_true:N \g_finishproof_bool 
  \bool_if:NTF \g_example_bool {
    The~proof~for~\cref{\loc}~is~thus~complete.
  } {
    The~proof~of~\cref{\loc}~is~thus~complete.
  }
}
\NewDocumentCommand{\finishproofthis} {} {
  \bool_gset_true:N \g_finishproof_bool 
  \bool_if:NTF \g_example_bool {
    This~completes~the~proof~for~\cref{\loc}.
  } {
    This~completes~the~proof~of~\cref{\loc}.
  }
}

\ExplSyntaxOff

\NewDocumentEnvironment{cproof}{m}
{\begin{proof}[Proof of \cref{#1}]}%
{\noindent The proof of \cref{#1} is thus complete.
\end{proof}}

\NewDocumentEnvironment{cproof2}{m}
{\begin{proof}[Proof of \cref{#1}]}%
{\noindent This completes the proof of \cref{#1}.
\end{proof}}

\hypersetup{
    colorlinks,
    linkcolor={blue!80!black},
    citecolor={green},
    urlcolor={blue!80!black}
}

\makeatletter
\ExplSyntaxOn
\seq_new:N \g_abbrs
\prop_new:N \g_abbr_counts
\tl_new:N \l_abbr_count_tl

\ExplSyntaxOn

\bool_new:N \g_forexample

\NewDocumentCommand{\eg}{ o }{
	\IfValueT{#1}{
		\str_if_eq:noTF {fe} {#1} {
			\bool_gset_true:N \g_forexample
		} {\bool_gset_false:N \g_forexample}
	}
	\bool_if:nTF { \g_forexample } {
		\bool_gset_false:N \g_forexample
		for~example
	}{
		\bool_gset_true:N \g_forexample
		for~instance
	}
}

\NewDocumentCommand{\abbr}{m m O{#1} m m O{#4} m}{
	\expandafter\newcommand\csname#3\endcsname[1][]{
		\seq_if_in:NnTF \g_abbrs {#1} {
			\prop_get:NnN \g_abbr_counts {#1} \l_abbr_count_tl
			\prop_gput:Nnx \g_abbr_counts {#1} {\int_eval:n {\l_abbr_count_tl + 1}}
			\hyperref[#1]{#7}
		} {
			\seq_gput_left:Nn \g_abbrs {#1}
			\prop_gput:Nnn \g_abbr_counts {#1} {1}
			\expandafter\gdef\csname#1@def\endcsname{#2}
			\phantomsection\label{#1}
			\str_if_eq:nnTF{##1}{}{\emph{#2}}{##1}~(\hyperref[#1]{#7})
		}
	}
	\expandafter\newcommand\csname#6\endcsname[1][]{
		\seq_if_in:NnTF \g_abbrs {#1} {
			\prop_get:NnN \g_abbr_counts {#1} \l_abbr_count_tl
			\prop_gput:Nnx \g_abbr_counts {#1} {\int_eval:n {\l_abbr_count_tl + 1}}
			\hyperref[#1]{#4}
		} {
			\expandafter\gdef\csname#1@def\endcsname{#5}
			\seq_gput_left:Nn \g_abbrs {#1}
			\prop_gput:Nnn \g_abbr_counts {#1} {1}
			\phantomsection\label{#1}
			\str_if_eq:nnTF{##1}{}{\emph{#5}}{##1}~(\hyperref[#1]{#4})
		}
	}
}

\ExplSyntaxOff
\makeatother
\abbr{GD}{gradient descent}{GDs}{gradient descents}{GD}
\abbr{FW}{Frank-Wolfe}{FWs}{Frank-Wolfe}{FW}
\abbr{ANN}{artificial neural network}{ANNs}{artificial neural networks}{ANN}
\abbr{RMSprop}{root mean square propagation}{RMSprops}{the root mean square propagation}{RMSprop}
\abbr{DNN}{deep neural network}{DNNs}{deep neural networks}{DNN}
\abbr{SGD}{Stochastic gradient descent}{SGDs}{Stochastic gradient descent}{SGD}
\abbr{iid}{independent and identically distributed}{i.i.d}{independent and identically distributed}{i.i.d.}
\abbr{SOP}{stochastic optimization problem}{SOPs}{stochastic optimization problems}{SOP}
\abbr{DL}{Deep learning}{DLs}{Deep learning}{DL}
\abbr{Adam}{adaptive moment estimation}{Adams}{adaptive moment estimation \SGD}{Adam}
\abbr{NS}{Newton-Schulz}{NSs}{Newton-Schulz}{NS}
\abbr{AdamW}{\Adam\ with decoupled weight decay}{AdamWs}{\Adam\ with decoupled weight decay}{AdamW}
\abbr{MUON}{momentum orthogonalized by Newton-Schulz}{MUONs}{momentum orthogonalized by Newton-Schulz}{MUON}
\abbr{AI}{artificial intelligence}{AIs}{artificial intelligence}{AI}
\abbr{LLM}{large language model}{LLMs}{large language models}{LLM}
\abbr{PDE}{partial differential equation}{PDEs}{partial differential equations}{PDE}
\abbr{as}{almost surely}{ass}{almost surely}{a.s.}
\abbr{pas}{$\P$-almost surely}{pass}{almost surely}{$\P$-a.s.}
\abbr{ReLUs}{rectified linear unit}{ReLU}{rectified linear unit}
\makeatother

\title{
On MUON optimization: From non-convergence\\
to an error analysis with Polar Express and the\\
Newton-Schulz polynomial from implementations
}
\author{Thang Do$^{1,2}$, Steffen Dereich$^{3}$, and Arnulf Jentzen$^{4,5}$
	\bigskip
	\\
    \small{$^1$ School of Data Science, The Chinese University of Hong Kong, Shenzhen}
	\vspace{-0.1cm}\\
	\small{ (CUHK-Shenzhen), China, e-mail: \texttt{minhthangdo@link.cuhk.edu.cn}}
 \smallskip
	\\
     \small{$^2$ Department of Probability and Statistics, Institute of Mathematics,}
	\vspace{-0.1cm}\\
	\small{Vietnam Academy of Science and Technology, Vietnam, e-mail: \texttt{dmthang@math.ac.vn}}
	\smallskip
	\\
	\small{$^3$ Institute for Mathematical Stochastics, Faculty of Mathematics and Computer Science,}\vspace{-0.1cm}\\
\small{University of M\"unster, Germany, e-mail: \texttt{steffen.dereich@uni-muenster.de}}
\smallskip
\\
	\small{$^4$ School of Data Science and School of Artificial Intelligence, The Chinese University}
	\vspace{-0.1cm}\\
	\small{of Hong Kong, Shenzhen (CUHK-Shenzhen), China, e-mail: \texttt{ajentzen@cuhk.edu.cn}}
	\smallskip
	\\
 \small{$^5$ Applied Mathematics: Institute for Analysis and Numerics, Faculty of Mathematics and}
	\vspace{-0.1cm}\\
	\small{Computer Science, University of M{\"u}nster, Germany, e-mail: \texttt{ajentzen@uni-muenster.de}}
	\smallskip
	\\
}

\date{\today}
\allowdisplaybreaks 
\begin{document}
\maketitle
\begin{abstract}
   \emph{Stochastic gradient descent} (SGD) optimization methods are the standard instruments for the training of \emph{deep neural networks} (DNNs). In many relevant \emph{artificial intelligence} (AI) systems – such as popular \emph{large language models} (LLMs) – not the standard SGD scheme is used as the optimization method but instead suitable accelerated variants of SGD are employed. One of the most popular methods of such accelerated SGD variants is the \emph{momentum orthogonalized by Newton-Schulz} (MUON) optimizer proposed by Jordan et al. in 2024. The MUON optimizer exploits the special matrix structure of the weight parameters in the training of the DNNs and, in its original form, employs five Newton-Schultz (NS) matrix steps in each MUON iteration. 
   
   In this work we propose and study a generalized variant of the MUON optimizer involving an arbitrary number of generalized NS steps with polynomials of possibly arbitrary high degree. The considered optimizer covers MUON with the original NS polynomial as well as MUON combined with the recently proposed Polar Express method as special cases. For a simple class of stochastic optimization problems (SOPs) we show for almost every mini-batch size that MUON \emph{fails to converge} to the solution of the SOP as the number of gradient steps converges to infinity. We also establish an error analysis for MUON with the generalized NS steps that provides convergence rates in terms of the number of gradient steps and in terms of the size of the mini-batch. We illustrate our general error analysis for MUON in the case of several concrete examples including quadratic stochastic optimization problems (SOPs) as well as $\ell_2$-regularized logistic regression for binary classification.
\end{abstract}
\pagebreak
\tableofcontents
\pagebreak
\section{Introduction}

\SGD\ optimization methods are the standard tools to train \DNNs\ \cite{Ruder2016AnOO}. In many relevant \AI\ systems -- such as popular \LLMs\ -- not the standard \SGD\ scheme is used as the optimization method but instead suitable sophisticated, accelerated variants of \SGD\ are employed such as the famous \Adam\ optimizer \cite{KingmaBa2024_Adam} and its variant with decoupled weight decay referred to as \AdamW\ optimizer \cite{Loshchilov2017DecoupledWD}. For about one decade \Adam\ and \AdamW, respectively, have been the method of choice among the sophisticated optimizers for training large \AI\ systems.

However, in 2024 another highly competitive optimization method, referred to as \MUON\ optimizer, has been proposed and significantly influenced the state-of-the-art in the training of \AI\ systems \cite{jordan2024muon}. In particular, according to the intelligence score \cite{ArtificialAnalysis} the most intelligent \LLMs, in which the choice of the optimizer has been publicly disclosed, are currently all trained with \MUON, variants of \MUON, or the \AdamW\ optimizer, whereby \MUON\ and its variants are significantly overrepresented compared to \AdamW\ \cite[Figure 1]{DeDoArPhi2025}. 

\subsection{Generalized Newton-Schulz (NS) iterations}
The \MUON\ optimizer employs the \NS\ method (cf., \eg, Appendix A in \cite{BernLaker2024}, (5.22) in \cite{MR2396439}, and (6) and (8) in \cite{MR295547}) with the polynomial proposed in \cite{jordan2024muon}. In this work we study \MUON\ combined with generalized \NS\ iterations. We first formulate the generalized \NS\ iterations that we consider in \cref{definition: NS} below and, thereafter, we employ \cref{definition: NS} in the formulation of our error analysis for \MUON\ in \cref{main theorem} below. The natural numbers $n, m \in \N = \{1, 2, 3, \dots\}$ in \cref{definition: NS} represent the dimensionalities of the matrices to which we apply the generalized \NS\ iterations.

 \begin{tcolorbox}[colback=white!95!gray,
                  colframe=black,
                  boxrule=0.5pt,
                  sharp corners,
                  enhanced,
                  breakable,
                 ]
\begin{definition}[\textcolor{red}{Generalized \NS\ method}]\label{definition: NS}
   Let $\nscoe=(\nscoe_{i,j})_{(i,j)\in(\N_0)^2}\colon (\N_0)^2\to\R$ satisfy $\min\{\nscoe_{0,0},\nscoe_{0,1}\}>0$ and $\#(\nscoe^{ - 1 }( \R\backslash\{0\} )) <\infty$ and let $n,m\in \N$, $A\in \R^{n\times m}$. Then we denote by
$
	(\TNewtonSchulzAlgorithm{\nscoe}{k}{A})_{k\in \N_0}\allowbreak\subseteq \R^{n \times m}
$ the sequence of matrices
 which satisfies\footnotemark\! for all
	$k \in \N$, $ M \in \R^{n\times m}$ with $\textcolor{magenta}{M=\TNewtonSchulzAlgorithm{\nscoe}{k-1}{A}} $
that
$\textcolor{magenta}{\TNewtonSchulzAlgorithm{\nscoe}{0}{A} = (\nscoe_{0,0}\|A\| + \nscoe_{0,1})^{-1} A}$
and
\begin{equation}\label{T_B_D}
\begin{split}
\textcolor{magenta}{\TNewtonSchulzAlgorithm{\nscoe}{k}{A}= \textstyle\sum_{ i = 0 }^{ \infty } \nscoe_{ k, i } ( M M^{\top} )^{ i} M}
.
\end{split}
\end{equation}
\end{definition}
\end{tcolorbox}
\footnotetext{Note that for all $v, w \in \N$, $A \in \R^{ v \times w }$ it holds that $A^{\top}$ is the transpose of $A$.}
Note that \cref{T_B_D} ensures that for every $k, n, m \in \N$, $\varepsilon \in (0,\infty)$, $A \in \R^{ n \times m }$, every $\nscoe=(\nscoe_{i,j})_{(i,j)\in(\N_0)^2}\colon\allowbreak (\N_0)^2\to\R$, and every $P \colon \R^{ n \times m } \to \R^{ n \times m }$ with the property that for all $M \in \R^{ n \times m }$ it holds that $P( M ) = \textstyle\sum_{ i = 0 }^{ \infty } \nscoe_{ k, i } ( M M^{\top} )^{ i} M$ we have that $
  \TNewtonSchulzAlgorithm{\nscoe}{k}{A} = P(\TNewtonSchulzAlgorithm{\nscoe}{k-1}{A} )$ .

\subsection{First main result: Error analysis for MUON with generalized NS iterations}

Using the \NS\ method from \cref{definition: NS} above we now present in the following theorem, \cref{main theorem}, an error analysis result for the \MUON\ optimizer. The \iid\ random variables $X_{ n,m} \colon \Omega \to \R^\dimX$, $( n,m ) \in \N^2$, in \cref{main theorem} represent the (input-output) data in the considered \SOP.
    
\begin{samepage}
 \begin{tcolorbox}[colback=white!95!gray,
                  colframe=black,
                  boxrule=0.5pt,
                  sharp corners,
                  enhanced,
                  breakable,
                 ]
\begin{athm}{theorem}{main theorem}[\textcolor{red}{\MUON\ error analysis}]
      Let $(\Omega,\cF,\P)$ be a probability space, let $\delta,\dimX,K\in \N$, $d_1,d_2,\dots,d_\delta,\fd_1,\fd_2,\dots,\fd_\delta\in \N$, $\alpha\in (0,1)$,  $\kappa\in (0,\infty)$, let $\nscoe=(\nscoe_{i,j})_{(i,j)\in (\N_0)^2}\colon(\N_0)^2\to \R$ satisfy $\textcolor{magenta}{\min\{\nscoe_{0,0},\nscoe_{0,1}\}>0}$ and $\textcolor{magenta}{\#(\nscoe^{ - 1 }( \R\backslash\{0\} )) <\infty}$,  assume for all $i\in \N\cap[0,K]$, $x\in \R$ that $\textcolor{magenta}{\sum_{j=0}^\infty\nscoe_{i,j}x^{2j}>0}$,
 let $X_{n,m}\colon \Omega\to\R^{\dimX}$, $(n,m)\in \N^2$, be bounded \iid\ random variables,  let $\cA=\times_{i=1}^\delta\R^{d_i\times\fd_i}$, $\xi\in \cA$, $\smalll=(\smalll(\theta,x))_{(\theta,x)\in \cA\times \R^{\dimX}}\in \allowbreak C^{1,0}(\cA\times \R^{\dimX},\R)$, assume\footnotemark\! for all $x\in \R^\dimX$ that $\cA\ni\theta\mapsto \smalll(\theta,x)-\kappa\|\theta\|^2\in \R$ is convex,
      assume that $\nabla_\theta\smalll$ is Lipschitz continuous, let $\vartheta\in \cA$ satisfy 
      \begin{equation}
       \textstyle  \textcolor{magenta}{ \E[\smalll(\vartheta,X_{1,1})]=\inf_{\theta\in \cA}\E[\smalll(\theta,X_{1,1})]},
      \end{equation} let $(\gamma_n)_{n\in \N}\subseteq (0,\infty)$ be non-increasing, and assume $\limsup_{n\to\infty} (\gamma_n+(\gamma_n)^{-2}(\gamma_{n}-\gamma_{n+1}))=0$.
       Then there exists $\fC\in \R$ such that for every $M\in \N$, every $\Theta=(\Theta^{1},\dots,\Theta^{\delta})\colon \N_0\times\Omega\to\cA$, and every $\bfm=(\bfm^{1},\dots,\bfm^{\delta})\colon \N_0\times\Omega\to\cA$ with the property that for all $n\in \N$, $i\in \{1,2,\dots,\delta\}$ it holds that 
         \begin{equation}\llabel{def: bfm}
            \textcolor{magenta}{\bfm_0=0},\qquad \textcolor{magenta}{\bfm_n= \alpha \bfm_{n-1}+(1-\alpha)\bigl[\textstyle \frac 1M \sum_{m=1}^M(\nabla_\theta\smalll)(\Theta_{n-1},X_{n,m})\bigr]},
            \end{equation}
            \begin{equation}\llabel{def: Theta}
            \textcolor{magenta}{\Theta_0=\xi},\qquad\text{and}\qquad  \textcolor{magenta}{\Theta_n^{i}=\Theta_{n-1}^{i}-\gamma_n
            \TNewtonSchulzAlgorithm{\nscoe}{K}{\bfm_n^i}}
         \end{equation}
         (cf.\ \cref{definition: NS})
     we have for all $n \in \N$ that
\begin{equation} \label{conclude: main theorem}
     \textcolor{magenta}{\E\bigl[\|\Theta_n-\vartheta\|^2\bigr]\leq \fC (M^{-1}+\gamma_{n})}.
\end{equation}
\end{athm}
\end{tcolorbox}
\end{samepage}
\footnotetext{Note that for all $N\in \N$, $v_1, v_2, \dots, v_{N }, w_1, w_2, \dots, w_{ N } \in \N$, $\theta=( ( \theta^n_{ i, j } )_{ (i, j) \in \{1,2,\dots,v_n\}\times\{1,2,\dots,w_n\} } )_{ n \in \{1,2,\dots,N\} } \in \times_{ n = 1 }^{N} \R^{ v_n \times w_n }$ it holds that $\| \theta \| = \bigl[\sum_{n=1}^N\sum_{i=1}^{v_n}\sum_{j=1}^{w_n}|\theta_{i,j}^n|^2\bigr]^{1/2}$ (standard norm, Frobenius norm, Hilbert-Schmidt norm).}

\cref{main theorem} is an immediate consequence of the more general result in \cref{muon stochastic convergence} in \cref{subsec: without stays bound} below (applied with $p\curvearrowleft 2$, $U\curvearrowleft \{x\in \R^\dimX\colon\|x\|\leq \inf\{ \rho \in \R \colon \P( \| X_{ 1, 1 } \| \leq \rho ) = 1 \}\}$ in the notation of \cref{muon stochastic convergence}), which is one of the main results of this work. \cref{muon stochastic convergence} treats a more general class of \MUON\ optimization methods and also allows the optimization process to start in a random variable instead of a deterministic vector $\xi \in \R^d$.

We emphasize that the class of \MUON\ optimization methods with generalized \NS\ iterations described through different choices of the function $\nscoe \colon ( \N_0 )^2 \to \R$ in \cref{main theorem} is so general that it covers the original \MUON\ method with the original \MUON\ polynomial \cite{jordan2024muon} as implemented in the {\sc Pytorch} library \cite{pytorch_muon} (see \cref{lem: abc verify} in \cref{subsec: jordan muon} below) and \MUON\ combined with the Polar Express method \cite{ar2505.16932} as implemented in the {\sc modded-nanoGPT} code \cite{Jordanmuoncode} (see \cref{lem: abc verify 2} in \cref{subsec: polar express} below) as special cases.

We also note that the assumption in \cref{main theorem} that $\limsup_{n\to\infty} (\gamma_n+(\gamma_n)^{-2}(\gamma_{n}-\gamma_{n+1}))=0$ is equivalent to the conditions that $\limsup_{ n \to \infty } \gamma_n = 0$ and $\limsup_{n\to\infty} ((\gamma_n)^{-2}(\gamma_{n}-\gamma_{n+1}))=0$. The condition $\limsup_{n\to\infty} ((\gamma_n)^{-2}(\gamma_{n}-\gamma_{n+1}))=0$ ensures that the sequence of learning rates $( \gamma_n )_{ n \in \N }$ does not converge to quickly to zero (cf., \eg,  \cite[Lemma 2.13]{DeDoArPhi2025}). These assumptions cover, \eg, polynomially decaying learning rates of the form $\gamma_n = c n^{ - r }$, $n \in \N$, for arbitrary $c \in (0,\infty)$, $r \in (0,1)$ (cf., \eg, \cite[Lemma 2.14]{DeDoArPhi2025}) and they have also been frequently used in error analyses for \Adam\ and other \SGD\ optimization methods (cf., \eg, \cite[Assumption 2.2]{MR1167814} and \cite[Theorem 1.1]{DereichAdamconvergence2024}).
\subsection{Second main result: Non-convergence of MUON}

It is well-known that under suitable assumptions (cf.\ \cref{main theorem} above) we have that the mean square error of the standard \SGD\ method after $n \in \N$ gradient steps can be bounded from above by a constant multiplied by the learning rate $\gamma_n$ at step $n$. In \cref{conclude: main theorem} in \cref{main theorem} we only bound the mean square error of \MUON\ with generalized \NS\ iterations from above by a constant $\fC$ multiplied by the sum of the learning rate $\gamma_n$ and the reciprocal $M^{ - 1 }$ of the mini-batch size $M$. In the next result, \cref{main theorem 2}, we prove that the stronger bound $\fC \gamma_n$ from standard \SGD\ can in general not be established in the case of \MUON. Specifically, in \cref{main theorem 2} we reveal for a simple quadratic \SOP\ that for almost every mini-batch size $M \in \N$ we have that \MUON\ (with generalized NS iterations) \emph{does not converge} to the solution of \SOP.
\begin{samepage}
\begin{tcolorbox}[colback=white!95!gray,
                  colframe=black,
                  boxrule=0.5pt,
                  sharp corners,
                  enhanced,
                  breakable,
                 ]
\begin{athm}{theorem}{main theorem 2}[\textcolor{red}{Non-convergence of \MUON}]
     Let $(\Omega,\cF,\P)$ be a probability space, let $K\in\N$, $\alpha\in (0,1)$,  $\kappa\in (0,\infty)$, $\xi\in \R$, let $\nscoe=(\nscoe_{i,j})_{(i,j)\in (\N_0)^2}\colon(\N_0)^2\to \R$ satisfy $\textcolor{magenta}{\min\{\nscoe_{0,0},\nscoe_{0,1}\}>0}$ and $\textcolor{magenta}{\#(\nscoe^{ - 1 }( \R\backslash\{0\} )) <\infty}$,  assume for all $i\in \N\cap[0,K]$, $x\in \R$ that $\textcolor{magenta}{\sum_{j=0}^\infty\nscoe_{i,j}x^{2j}>0}$, let $\smalll=(\smalll(\theta,x))_{(\theta,x)\in \R\times \R} \allowbreak \colon \R\times \R\to\R$ satisfy for all  $\theta,x\in \R$ that $\smalll(\theta,x)=|\theta-x|^2$, let $X_{n,m}\colon \Omega\to\R$, $(n,m)\in \N^2$, be bounded \iid\ random variables, assume $\E[(X_{1,1}-\E[X_{1,1}])^3]\neq 0$, let $(\gamma_n)_{n\in \N}\subseteq (0,\infty)$ be non-increasing, assume $\limsup_{ n \to \infty } ( ( \gamma_{n+1} )^{ - 1 } \gamma_n ) < \alpha^{-1} < \sum_{ n=1 }^{ \infty } \gamma_n = \infty$, and for every $M\in \N$ let $\Theta^M \colon\N_0\times\Omega\to\R$ and $\bfm^M\colon \N_0\times\Omega\to\R$ satisfy for all $n\in \N$ that 
        \begin{equation}\label{eq1: main theorem 2}
           \textcolor{magenta}{ \bfm_0^M=0},\qquad \textcolor{magenta}{\bfm_n^M= \alpha \bfm_{n-1}^M+(1-\alpha)\bigl[\textstyle \frac 1M \sum_{m=1}^M(\nabla_\theta\smalll)(\Theta_{n-1}^M,X_{n,m})\bigr]},
            \end{equation}
            \begin{equation}\label{eq2: main theorem 2}
           \textcolor{magenta}{\Theta_0^M= \xi},\qquad \text{and}\qquad \textcolor{magenta}{\Theta_n^{M}=\Theta_{n-1}^{M}-\gamma_n\TNewtonSchulzAlgorithm{\nscoe}{K}{\bfm_n^{M}}}.
         \end{equation}
         Then it holds for almost all $M\in \N$ that 
         \begin{equation}\label{eq3: main theorem 2}
           \textcolor{magenta}{\textstyle\limsup_{ n \to \infty }\E\bigl[ \min\{ 1, | \Theta_n^M - \E[X_{1,1}]|\} \bigr] > 0}.
           \end{equation}
  (cf.\ \cref{definition: NS}).
\end{athm}
\end{tcolorbox}
\end{samepage}

\cref{main theorem 2} is a direct consequence of the more general result in  \cref{MUON non-convergence 1d big batch} in \cref{sec: non convergence} below, which is one of the main results of this work. \cref{MUON non-convergence 1d big batch} treats a more general class of \MUON-type optimization methods and also allows the optimization process to start in a random variable instead of a deterministic vector $\xi \in \R^d$. In \cref{main theorem 2} we intend to solve the optimization problem 
\begin{equation}\label{ZZZ}
  \operatorname{argmin}_{ \theta \in \R } \cL( \theta )
\end{equation} 
where $\cL \colon \R \to \R$ is the function which satisfies for all $\theta \in \R$ that $\cL( \theta ) = \E[ \smalll( \theta, X_{ 1, 1 } ) ] = \E[ |\theta - X_{ 1, 1 } |^2 ]$ (the objective function, the function we intend to minimize). We observe that the objective function $\cL \colon \R \to \R$ has a unique global mimimizer at $\theta = \E[ X_{ 1, 1 } ]$. In \cref{eq3: main theorem 2} in \cref{main theorem 2} we then reveal for almost every mini-batch size $M \in \N$ that we have that the considered generalized \MUON\ opimization process $( \Theta^M_n )_{ n \in \N_0 }$ \emph{fails} to converge to the solution of the \SOP\ as the number of gradient steps $n$ converges to infinity. 

However, from \cref{conclude: main theorem} in \cref{main theorem} we see that the mean square error of the generalized \MUON\ method can still be made arbitrarily small if, both, the number of gradient steps n and the mini-batch size M are sufficiently large. \cref{main theorem 2} only demonstrates that an additional error, beside the standard error term $\fC \sqrt{ \gamma_n }$, must be included in the error estimates for \MUON. The error of standard \SGD\ and momentum \SGD\ can, however, be bounded above only by $\fC \sqrt{ \gamma_n }$; cf., \cite[Theorem ...]{DaAr2026strong}, \cite[Theorem 5.8]{Guro2023handbook}, \cite[Theorem 3.7]{Jentzen_2020}, \cite[Proposition 3.3]{MR4055054}, and the references therein.

\subsection{An example: MUON for quadratic stochastic optimization problems (SOPs)}
A key feature of \cref{main theorem} and the more general result in \cref{muon stochastic convergence}, respectively, is that they are not merely abstract convergence results for the \MUON\ optimizer, but that they also applicable to several concretely specified examples of \SOPs. For example, \cref{main theorem} applies to this class of quadratic \SOPs
\begin{equation}\label{def: example}
 \textstyle \operatorname{argmin}_{ \theta = ( \theta_1, \dots, \theta_{ \delta } ) \in \cA } \E\bigl[ \sum_{ i = 1 }^{ \delta } \|M_i (\theta_i) - p_i( X_{1,1} ) \|^2 \bigr]
\end{equation} 
where $\cA=\times_{i=1}^\delta\R^{d_i\times \fd_i}$, where $M_1\colon\R^{d_1\times\fd_1}\to\R^{d_1\times\fd_1}$, $M_2\colon\R^{d_2\times\fd_2}\to\R^{d_2\times\fd_2}$, $\dots$, $M_\delta\colon\R^{d_\delta\times\fd_\delta}\to\R^{d_\delta\times\fd_\delta}$ are invertible linear functions, where $p_1 \colon \R^\dimX\to \R^{d_1\times\fd_1} $, $p_2 \colon \R^\dimX\to\R^{d_2\times\fd_2} $, $\dots$, $p_{ \delta }\colon \R^\dimX\to \R^{d_\delta\times\fd_\delta} $ are Lipschitz continuous functions, where $X_{ n, m } \colon \Omega \to \R^\dimX$, $(n,m) \in \N^2$, are bounded \iid\ random variables, 
where $( \Omega, \cF, \P )$ is a probability space, and where $d_1, d_2, \dots, d_{ \delta }, \fd_1,\fd_2,\dots,\fd_\delta \in \N$,  $\delta, \mathscr{d} \in \N$.

We refer to \cref{lem: loss verify} and \cref{main theorem 2 quadratic 2} in \cref{subsec: quadradic muon} below for an error analysis for \MUON\ applied to a more general class of example \SOPs\ (including \cref{def: example} as a special case). \cref{main theorem 2 quadratic 2} is based on an application of our more general error analysis result for the \MUON\ optimizer in \cref{muon stochastic convergence}.

\cref{muon stochastic convergence} is also applicable to several other concretely specified optimization problems -- including $\ell_2$-regularized logistic regression for binary classification (\cref{subsec: logistic}) -- and we refer to \cref{sec: example} below for applications of \cref{muon stochastic convergence} to some example optimization problems.

\subsection{Literature review}\label{subsec: literature review}

There is nowadays an extensive literature on the {\bf theoretical analysis of \SGD\ optimization methods},  
including works on the classical standard \SGD\ method (cf., \eg, \cite{Bottou1998}), 
articles on accelerated \SGD\ optimizers such as Nesterov \cite{MR701288} (cf., \eg, \cite{pmlr-v28-sutskever13}) 
and momentum \SGD\ \cite{Polyak1964SomeMO} (cf., \eg, \cite{ar1808.10396}), 
papers on adaptive \SGD\ methods such as \RMSprop\ \cite{Hinton24_RMSprop} (cf., \eg, \cite{ZouShen2019}) 
and stochastic sign descent type methods \cite{sign1} (cf., \eg, \cite{sign2,sign4,sign3,Mhersignsgd2019,sign5}), 
as well as also a large number of works on accelerated and adaptive optimizers, 
such as the famous \Adam\ \cite{KingmaBa2024_Adam} and \AdamW\ \cite{Loshchilov2017DecoupledWD} optimizers (cf., \eg, \cite{Defossez2022,DeDoArPhi2025,DereichAdamconvergence2024,DeArAdsharp2025,li2023convergenceadamrelaxedassumptions,ReddiKale2019,ZhangChen2022} and the references therein). 

The \MUON\ optimizer, which is based on momentum \SGD\ and approximate matrix normalization 
(approximate projection to the class of semi-orthogonal matrices; cf.\ \cref{def: Pi: main theorem stochastic} in \cref{main theorem stochastic} below), 
can be regarded as an optimizer that belongs to this last class of accelerated (due to the involved momentum) 
and adaptive (due to the approximate matrix normalization) \SGD\ optimization methods, 
but also appends additional features as it respects the matrix structures of the parameters to optimize. 
Since its introduction in December 2024 \cite{jordan2024muon}, 
there have already appeared a series of works that study \MUON-type methods 
analytically; see, \eg, \cite{muon5,muon2,muon7,muon1,muon6,muon8,muon9,muon4,muon3,Weimuonconvergence2025,ar2510.19933}.

Some of these works analyze \MUON-type methods in non-convex settings. For example, \cite{muon7}, \cite[Theorem 3.1, Theorem 3.2, and Proposition 3.3]{muon4}, \cite[Theorem 3]{muon3}, \cite[Theorem 4.1]{Weimuonconvergence2025}, and \cite[Theorem 2]{ar2510.19933}
 neither assume strong convexity nor convexity but impose suitable $L$-smoothness assumptions (global Lipschitz continuity (with respect to different topologies) of the gradient of the objective function) and bounded variance conditions (global boundedness of the variance of the unbiased estimator 
of the gradient of the objective function) to establish for every $ T \in \N $ upper bounds 
for the first absolute moment of the gradient of the objective function 
composed with the considered \MUON\ optimization process evaluated at an on $\{0,1,2, \dots, T\}$ 
uniformly distributed random variable (ergodic first-moment gradient bounds). 

In \cite[Theorem 2.1]{muon1}, \cite[Theorem 3.1]{muon8}, \cite[Theorem 4.1]{Weimuonconvergence2025},
 \cite[Theorem 3.1 and Theorem 3.2]{muon4}, and \cite[Algorithm 2]{muon3} the authors study the \emph{idealized} \MUON\ optimizer (see \cite[Section ”The design of Muon”]{jordan2024muon}
and, for example, \cite[Definition 6.12.2]{ArBePhi2024} and \cref{main theorem stochastic} below), which corresponds to the situation where 
the number of \NS\ steps converges to infinity and where appropriate parameters $a$, $b$, and $c$ 
in the \NS\ method in \cref{definition: NS} are chosen \cite[Definition 2]{muon7}, and the work \cite{muon7}
studies \MUON\ with a finite number of Newton-Schulz steps with a polynomial that ensures 
convergence to the idealized \MUON\ optimizer (cf.\ also \cite[Proposition 3.3]{muon4}) but not with the polynomial that has been originally proposed for \MUON\ \cite{jordan2024muon} and that is used in the {\sc Pytorch} library \cite{pytorch_muon}.

The work \cite[Theorem 3]{muon3} establishes such ergodic first-moment gradient bounds for stochastic
\FW\ methods \cite{MarWolfe1956} and shows that \MUON\ can be regarded as an instance of the considered stochastic \FW\ schemes \cite[Theorem 1]{muon3}.

In the work \cite[Section 4.2]{muon5} (cf., \eg, also \cite[Section 4.2.2]{muon9}) 
convergence rates for certain idealized \MUON-type methods are shown under the assumption that the popular \emph{Polyak–Łojasiewicz inequality} holds and under the assumption that the momentum decay parameter is not kept constant but increases to one during the training. This adjustment to increase the momentum decay parameter 
to one during the training allows them to establish convergence of \MUON\ as the number of \MUON\ steps 
go to infinity \cite[Section 4.2]{muon5}. 
In implementations the momentum decay parameter is, however, typically kept constant at the value 
$\alpha = 0.95$ \cite{pytorch_muon,Jordanmuoncode}. 
We also note that the findings in \cite[Section 4.2]{muon5} do not contradict to the non-convergence result
for \MUON\ in \cref{main theorem 2} above as the momentum decay parameter $\alpha \in (0,1)$ is kept constant during the training 
in \cref{main theorem 2} above (and as \cref{main theorem 2} considers \MUON\ with the \NS\ steps instead of idealized \MUON-type methods).

We also emphasize the instructive work \cite{ar2510.19933} which studies a generalized \MUON\ method that \emph{only incorporates approximations of the polar factor but not necessarily the exact polar factor} in the idealized \MUON\ method. In \cite{ar2510.19933} it is assumed that the approximation errors (with respect to the exact polar factor) can be bounded by real numbers $\delta_k \geq 0$, $k \in \N_0$; see \cite[Assumption 1]{ar2510.19933}. A crucial assumption in the error analysis of \cite{{ar2510.19933}} is then to assume that $\sup_{k\in\N_0}\delta_k < 1$ is strictly smaller than 1; see \cite[Theorem 2 and Theorem 4]{{ar2510.19933}}. The work \cite{{ar2510.19933}} does not verify this assumption in the case of concrete optimization methods (cf.\ \cite[Appendix A]{ar2510.19933}) and in some settings this assumption excludes the \NS\ iterations with the polynomial originally introduced within the \MUON\ optimizer \cite{jordan2024muon} as implemented in the {\sc Pytorch} library \cite{pytorch_muon} and also excludes the Polar Express method \cite{ar2505.16932} as implemented in the {\sc modded-nanoGPT} code \cite{Jordanmuoncode}; see \cref{lem: counter example 1} and \cref{lem: counter example 3} in \cref{subsec: counter example} below.

To the best of our knowledge, \cref{main theorem} above and its generalization in \cref{muon stochastic convergence} below, respectively, are the first results in the scientific literature that provide an error analysis for \MUON\ with an arbitrary finite number of \NS\ steps with the \NS\ polynomial that has been originally proposed for \MUON\ \cite{jordan2024muon} and that is used in the {\sc Pytorch} library \cite{pytorch_muon} (corresponding to the choice that for all $i,j\in\N_0$ it holds that $\nscoe_{ i, j } = 3.4445  \mathbbm{1}_{  \{1,2,\dots,5\}\times \{0\}  }(  i, j  ) -4.7750   \mathbbm{1}_{  \{1,2,\dots,5\} \times \{1\} }(  i, j  )+  2.0315\mathbbm{1}_{ \{1,2,\dots,5\}\times \{2\}  }(  i, j ) +\mathbbm{1}_{ \{(0,0)\}  }(  i, j )+10^{-7}\mathbbm{1}_{ \{(0,1)\}}(  i, j )$ in \cref{definition: NS}; see \cref{lem: abc verify} in \cref{subsec: jordan muon} below) and also for \MUON\ combined with the Polar Express method \cite{ar2505.16932} as implemented in the {\sc modded-nanoGPT} code \cite{Jordanmuoncode} (see \cref{lem: abc verify 2} in \cref{subsec: polar express} below).

In \cref{main theorem 2} above and \cref{MUON non-convergence 1d big batch} below we also establish {\bf non-convergence results for \MUON}. 
There are also a series of of lower error bounds, counterexamples, as well as 
non-convergence and divergence results for gradient based optimization methods in the literature.

In particular, we refer, \eg, to \cite{CheriditoJentzenRossmanek2021,HannibalJentzenThang2024,DoArAd2025,ArAd2024,LSSK2020}
for non-convergence results for \SGD\ optimization methods 
in the training of \ANNs\ (some of them, showing that the probability to converge to a global minimizer  
is strictly positive or even large, close to one), 
we refer, \eg, to \cite{gallon2022blowphenomenagradientdescent,ar2505.09572,Lyu2020Gradient,vardi2022on} for divergence results for gradient based optimization processes in the training of \ANNs\ (which provide sufficient conditions that the norm of the \ANN\ parameter vector 
diverges during the training to infinity) (cf., \eg, also \cite[Section 1.2]{MR4243432}), 
we refer, \eg, to \cite{DeDoArPhi2025,ar2607.03519,ReddiKale2019} for non-convergence results specifically for the \Adam\ optimizer 
studying the \emph{internal error of \Adam} (showing in certain examples that \Adam\ does not converge 
to the solution (the unique minimizer) of the \SOP\ under consideration), 
we refer, \eg, to \cite[Section 3]{sign3} for non-convergence results for stochastic sign descent type optimization methods,  
and we refer, \eg, to \cite{DeRoAr2024nonconvergence} for non-convergence results for \Adam\ and other \SGD\ optimization processes  where the learning rates (the step sizes) do not converge to zero. 

We also refer to \cite{ar2602.11948,ar2605.08980} for certain non-convergence results for \MUON\ type methods. 
The work \cite[Section 3]{ar2605.08980} shows 
non-convergence of \MUON\ to global minimum points for certain specific deterministic optimization problems 
with a non-smooth objective function. 
The non-smoothness of the objective function in \cite{ar2605.08980} 
is essential and can not be avoided as \MUON\ with asymptotically vanishing learning rates 
does converge to the global minimum point for smooth strongly convex deterministic optimization problems; cf.\ \cref{deterministic muon convergence} and \cref{main theorem stochastic} below.

The work \cite[Section 3]{ar2602.11948} studies a certain deterministic optimization problem in which the dynamics of the idealized \MUON\ optimizer with a fixed step size becomes confined to discrete lattices and fails to converge. The non-convergence of the learning rates to zero is essential and can not be avoided as \MUON\ with asymptotically vanishing learning rates applied to smooth strongly convex deterministic optimization problems does converge; cf.\ \cref{deterministic muon convergence} and \cref{main theorem stochastic} below. 

We also refer to \cite[Figure 3]{ar2505.16932} for a numerical study that indicates that the \NS\ method with the polynomial as originally proposed in \cite{jordan2024muon} and as typically used in \MUON\ implementations fails to converge to
the polar factor.

To the best of our knowledge, \cref{main theorem 2} above and its generalization in \cref{MUON non-convergence 1d big batch} below, respectively, are the first non-convergence results for the \MUON\ optimizer 
in the literature that show that \MUON\ (with the \NS\ polynomial as originally proposed in \cite{jordan2024muon} and as typically used in implementations) does, in general, not converge to the solution of the \SOP\ (the minimizer of the \SOP) as the number of gradient steps $n$ converges to infinity.

A crucial assumption in \cref{main theorem 2} and \cref{MUON non-convergence 1d big batch} is that the expectation of the third central moment of the data does not vanish $\E[ ( X_{ 1, 1 } - \E[ X
_{1,1}] )^3 ] \neq 0$. We recall that if the optimization problem is deterministic (corresponding to the case where the data $X_{ 1, 1 }$ is one-point distributed, in which we have that $\E[ (X_{1,1} - \E[ X_{1,1} ] )^3 ] = \E[ ( \E[ X_{1,1} ] - \E[ X_{1,1} ] )^3 ] = 0$), then convergence of \MUON\ does hold (see \cref{deterministic muon convergence} below) and the non-convergence conclusion in \cref{eq3: main theorem 2} in \cref{main theorem 2} can not hold anymore.

\subsection{Organization of this article}

The remainder of this work is organized in the following way. In \cref{sec: priori bound} we establish pathwise a priori bounds, which are also of independent interest, for a general class of optimization methods including \MUON\ as a special case; see \cref{theo: MUON bound pre} and \cref{theo: MUON bound}. 

In \cref{sec: convergence of stochastic muon} we first establish in \cref{muon bound convergence stochastic} an error bound for the \MUON\ optimizer under the assumption that the \MUON\ optimization process is almost surely bounded by a fixed constant and, thereafter, we combine \cref{muon bound convergence stochastic} with the a priori bound in \cref{theo: MUON bound} of \cref{sec: priori bound} to establish in \cref{muon stochastic convergence} an error estimate for \MUON\ without assuming that the \MUON\ optimization process stays bounded. \cref{main theorem} above in this introduction section follows from an application of \cref{muon stochastic convergence} (see \cref{subsec: without stays bound} for details).

Furthermore, in \cref{subsec: idealized muon} we also develop in \cref{main theorem stochastic} an error analysis for the idealized \MUON\ optimizer (see \cite[Section "The design of Muon"]{jordan2024muon} and, \eg, \cite[Definition 6.12.2]{ArBePhi2024}), which corresponds to the situation where the number of \NS\ steps converges to infinity and where appropriate parameters $a$, $b$, and $c$ in the \NS\ method in \cref{definition: NS} are chosen.

In \cref{sec: example} we illustrate the conclusion of \cref{muon stochastic convergence} in the situation of several example optimization problems, namely, logistic regression for binary classification, quadratic \SOPs, and deterministic optimization problems.

In \cref{sec: practical} we show that the framework in \cref{main theorem} and \cref{muon stochastic convergence} is general enough to cover \MUON\ combined with five \NS\ steps with the original \MUON\ polynomial \cite{jordan2024muon} (see \cref{lem: abc verify} in \cref{subsec: jordan muon} below) as well as \MUON\ combined with the Polar Express method \cite{ar2505.16932} (see \cref{lem: abc verify 2} in \cref{subsec: polar express} below).

In \cref{sec: non convergence} we show for \MUON\ optimization methods applied to a simple class of one-dimensional quadratic \SOPs\ that for almost every mini-batch size $M \in \N$ we have that the optimizer fails to converge to the solution of the \SOP\ (the unique global minimizer of the \SOP) as the number of gradient steps $n$ goes to infinity; see \cref{MUON non-convergence 1d big batch}. \cref{main theorem 2} above in this introductory section is an immediate consequence of \cref{MUON non-convergence 1d big batch}.

\subsection{Use of large language models}

\cref{lem: stochastic recurrence''}, \cref{lem: stochastic recurrence 2}, \cref{lem: stochastic bfm analysis}, \cref{theo: stochastic convergence}, \cref{lem: measuability}, \cref{Casero}, \cref{vector field}, \cref{lem: weak converge}, \cref{lem: tg1},  \cref{lem: integral}, \cref{lem: tg3}, and \cref{lem: tg2} in this work have been developed based on discussions with {\sc GPT 5.5 Pro} or {\sc GPT 5.6 Pro}. Specifically, the core ideas for the proofs of \cref{lem: stochastic recurrence''}, \cref{lem: stochastic recurrence 2}, \cref{lem: stochastic bfm analysis}, \cref{theo: stochastic convergence}, \cref{lem: measuability}, \cref{Casero}, \cref{vector field}, \cref{lem: weak converge}, \cref{lem: tg1}, and \cref{lem: tg2} and the core ideas for the statements of \cref{lem: stochastic recurrence''}, \cref{lem: stochastic recurrence 2}, \cref{lem: stochastic bfm analysis}, \cref{Casero}, \cref{lem: weak converge}, \cref{lem: tg1}, \cref{lem: integral}, and \cref{lem: tg3} are due to {\sc GPT 5.5 Pro} or {\sc GPT 5.6 Pro}. In addition, {\sc GPT 5.5 Pro} and {\sc GPT 5.6 Pro} have supported us in creating the literature review in Subsection 1.5 above and have helped us to reveal several inaccuracies and errors within this manuscript. The entire article (including all statements and all arguments in this article) has been carefully written by the authors and the authors take full responsibility for each of the sentences/statements made in the paper.

\section{A priori bounds for MUON and closely related optimization methods}\label{sec: priori bound}

In this section we establish in \cref{theo: MUON bound pre} and \cref{theo: MUON bound} pathwise a priori bounds for a general class of optimization methods including \MUON\ as a special case. In \cref{sec: convergence of stochastic muon} we employ the a priori bound in \cref{theo: MUON bound} to develop our error analysis for the \MUON\ optimizer (see \cref{muon stochastic convergence} below).
\subsection{Dynamics and approximations for momentum processes}\label{subsec: dynamic for momentum}

In the following elementary and well-known lemma, \cref{momentum:representation}, we represent a momentum process at time $n \in \N_0$ 
in terms of the underlying gradients at times $0, 1, \dots, n - 1$. \cref{momentum:representation} can, \eg, be found in \cite[Lemma 2.7]{DeRoArAd2025} in a slightly modified/specialized form.

\begin{samepage}
\begin{tcolorbox}[colback=white!95!gray,
                  colframe=black,
                  boxrule=0.5pt,
                  sharp corners,
                  enhanced,
                  breakable,
                 ]
\begin{athm}{lemma}{momentum:representation}
			Let $\cA$ be an $\R$-vector space and 
			let $\bbM\colon\N_0\to\cA$, $g\colon \N\to\cA$, and $\alpha\colon\N\to\R$ satisfy for all $n\in \N$ that
			\begin{equation}
				\label{setup:basic:geom:mom}
				\bbM_n=\alpha_n \bbM_{n-1}+ g_n.
			\end{equation}
			Then it holds for all $n\in\N_0$ that
			$
			\bbM_n
			=\bigl[\prod_{j=1}^{n}\alpha_j\bigr]\bbM_0
			+\sum_{k=1}^{n}\bigl[\prod_{j=k+1}^{n}\alpha_j\bigr]g_k
			$.
		\end{athm}
        \end{tcolorbox}
		\end{samepage}
		\begin{aproof}
			\argument{
				\eqref{setup:basic:geom:mom}}{that
				\begin{equation}
					\llabel{eq:base:momentum:geom}
					\begin{split}
						\bbM_0
						\textstyle
						=\bbM_0
						\qqandqq
						\bbM_1
						=\textstyle\alpha_1 \bbM_0+ g_1
						.
					\end{split}
				\end{equation}
			}
			\argument{\eqref{setup:basic:geom:mom}}{that for all $n\in\N$ with $\bbM_n=\bigl[\prod_{j=1}^n\alpha_j\bigr]\bbM_0+\sum_{k=1}^n\bigl[\prod_{j=k+1}^{n}\alpha_j\bigr]g_k$ it holds that
				\begin{equation}
					\llabel{eq:base:momentum:geom:2}
					\begin{split}
						&\bbM_{n+1}
						=\textstyle\alpha_{n+1} \bbM_n+g_{n+1}
						=\textstyle\alpha_{n+1} \prb{\PRb{\prod_{j=1}^n\alpha_j}\bbM_0+\sum_{k=1}^n\PRb{\prod_{j=k+1}^n\alpha_j}g_k}+g_{n+1}
						\\&=\textstyle\PRb{\prod_{j=1}^{n+1}\alpha_j}\bbM_0
						+\sum_{k=1}^{n}\PRb{\prod_{j=k+1}^{n+1}\alpha_j}g_k+g_{n+1}
						=\textstyle
						\PRb{\prod_{j=1}^{n+1}\alpha_j}\bbM_0
						+\sum_{k=1}^{n+1}\PRb{\prod_{j=k+1}^{n+1}\alpha_j}g_k
						.
					\end{split}
				\end{equation}
			}
			\argument{\lref{eq:base:momentum:geom:2};
				\lref{eq:base:momentum:geom};
				induction}{that for all $n\in\N_0$ it holds that
				$
				\bbM_n=\PRb{\prod_{j=1}^{n}\alpha_j}\bbM_0
				+\sum_{k=1}^{n}\PRb{\prod_{j=k+1}^{n}\alpha_j}g_k
				$}.
		\end{aproof}

        In the following elementary result, \cref{lem: bfm analysis 1}, we establish for every $n \in \N$ an upper bound for the distance between the auxiliary averaged gradient process in the momentum method at time $n$ and the gradient of the objective function (here denoted by the function $\cG \colon \cV \to \cV$) evaluated at the optimization process in the momentum method at time $n-1$. \cref{lem: bfm analysis 1} and its proof are strongly inspired by, \eg, \cite[Lemma 2.8]{ShoAr2026} (in the special situation where $\delta = 1$, $\forall \, n \in \N \colon \alpha_n = \alpha_1$, and $\forall \, n \in \N \colon \mu_n = 0$ in \cite[Lemma 2.8]{ShoAr2026}). 

\begin{tcolorbox}[colback=white!95!gray,
                  colframe=black,
                  boxrule=0.5pt,
                  sharp corners,
                  enhanced,
                  breakable,
                 ]
\begin{athm}{lemma}{lem: bfm analysis 1}
     Let $\alpha\in (0,1)$, $\cK\in \R$, let $(\cV,\nnorm{\cdot})$ be a normed $\R$-vector space, and let $\cG\colon \cV\to\cV$, $\gamma\colon \N\to(0,\infty)$, $\bfm\colon \N_0\to\cV$, and $\Theta\colon \N_0\to\cV$ satisfy for all $n\in \N$ that
     \begin{equation}\llabel{def: bfm}
        \nnorm{\cG(\Theta_n)-\cG(\Theta_{n-1})}\leq \cK\nnorm{\Theta_n-\Theta_{n-1}},\qquad \bfm_{n}= \alpha\bfm_{n-1}+ (1-\alpha)\cG(\Theta_{n-1}),
     \end{equation}
     \begin{equation}\llabel{assume}
        \qqandqq\textstyle\max\bigl\{ \sum_{i=1}^n(\alpha^{n-i}\gamma_i), \nnorm{\Theta_n-\Theta_{n-1}}\bigr\}\leq \cK\gamma_n.
     \end{equation} Then it holds for all $n\in \N$ that
     \begin{equation}\llabel{conclude}
\nnorm{\bfm_{n}-\cG(\Theta_{n-1})}\leq \bigl(\cK(\gamma_1)^{-1}\nnorm{\bfm_0-\cG(\Theta_0)}+\cK^3+\cK^2\bigr)\gamma_n.
     \end{equation}
    \end{athm}
    \end{tcolorbox}
    \begin{aproof}
          \argument{\lref{def: bfm};\lref{assume};}{that for all $n\in \N$ it holds that
        \begin{equation}\llabel{eq2.1}
            \nnorm{\cG(\Theta_n)-\cG(\Theta_{n-1})}\leq \cK\nnorm{\Theta_n-\Theta_{n-1}}\leq \cK^2\gamma_n.
        \end{equation}}
        \argument{\cite[Lemma 2.5]{ShoAr2026};}{that for all $n\in \N$ it holds that
        \begin{equation}\llabel{eq1}
           \textstyle \bfm_n-\cG(\Theta_n)=\alpha^n[\bfm_0-\cG(\Theta_0)]-\sum_{i=1}^n[\alpha^{n-i}(\cG(\Theta_i)-\cG(\Theta_{i-1}))].
        \end{equation}}
        \argument{\lref{eq1};\lref{eq2.1}}{that for all $n\in \N$ it holds that
        \begin{equation}\llabel{eq2}
            \begin{split}
                \nnorm{\bfm_n-\cG(\Theta_n)}&\leq \alpha^n\nnorm{\bfm_0-\cG(\Theta_0)}+\sum_{i=1}^n\bigl[\alpha^{n-i}\nnorm{\cG(\Theta_i)-\cG(\Theta_{i-1})}\bigr]\\
                &\leq \alpha^n\nnorm{\bfm_0-\cG(\Theta_0)}+\sum_{i=1}^n\alpha^{n-i}\cK^2\gamma_i.
            \end{split}
        \end{equation} }
        \argument{\lref{eq2};\lref{assume};}{that for all $n\in \N$ it holds that
        \begin{equation}\llabel{eq4}
            \nnorm{\bfm_n-\cG(\Theta_n)}\leq \alpha^n\nnorm{\bfm_0-\cG(\Theta_0)}+\cK^3\gamma_n.
        \end{equation}}
        \argument{\lref{assume};}{that for all $n\in \N$ it holds that \llabel{eq5}
            $\alpha^{n-1}\gamma_1\leq \cK\gamma_n$\dott}
        \argument{\lref{eq5};the fact that $\alpha<1$}{that for all $n\in \N$ it holds that
        \begin{equation}\llabel{eq6}
           \alpha^n\leq \cK\gamma_n\alpha(\gamma_1)^{-1}\leq \cK(\gamma_1)^{-1}\gamma_n.
        \end{equation}}
        \argument{\lref{eq6};\lref{eq4}}{that for all $n\in \N$ it holds that
        \begin{equation}\llabel{eq7}
            \nnorm{\bfm_n-\cG(\Theta_n)}\leq \alpha^n\nnorm{\bfm_0-\cG(\Theta_0)}+\cK^3\gamma_n\leq \cK(\gamma_1)^{-1}\gamma_n\nnorm{\bfm_0-\cG(\Theta_0)}+\cK^3\gamma_n.
        \end{equation}}
        \argument{\lref{eq7};\lref{eq2.1}}{that for all $n\in \N$ it holds that
        \begin{equation}\llabel{eq8}
        \begin{split}
            \nnorm{\bfm_n-\cG(\Theta_{n-1})}&\leq \nnorm{\bfm_n-\cG(\Theta_n)}+\nnorm{\cG(\Theta_n)-\cG(\Theta_{n-1})}\\
            &\leq \cK\gamma_n(\gamma_1)^{-1}\nnorm{\bfm_0-\cG(\Theta_0)}+\cK^3\gamma_n+\cK^2\gamma_n\\
            &=\bigl(\cK(\gamma_1)^{-1}\nnorm{\bfm_0-\cG(\Theta_0)}+\cK^3+\cK^2\bigr)\gamma_n.
        \end{split}
        \end{equation}}
        \argument{\lref{eq8};}{\lref{conclude}\dott}
    \end{aproof}
   \subsection{Quantitative pathwise a priori bounds for MUON type optimization processes}

   In this subsection we establish in \cref{theo: MUON bound 0} below a quantitative a priori estimate for a general class of \MUON-type optimization processes which includes the \MUON\ optimization process as a special case. In \cref{subsec: pathwise bound} below we apply \cref{theo: MUON bound 0} to derive in \cref{theo: MUON bound pre} below qualitative a priori estimates for such optimization processes.
    \begin{samepage}
   \begin{tcolorbox}[colback=white!95!gray,
                  colframe=black,
                  boxrule=0.5pt,
                  sharp corners,
                  enhanced,
                  breakable,
                 ]
\begin{athm}{lemma}{lem: recurrence}
      Let $c,\kappa,\varepsilon\in (0,\infty)$ and let $\gamma\colon \N\to(0,\infty)$ and $e\colon \N_0\to[0,\infty)$ satisfy for all $n\in \N$ that
      \begin{equation}\llabel{assume}
            e_n\leq  e_{n-1}-\kappa\gamma_ne_{n-1}(\sqrt{e_{n-1}}+\varepsilon)^{-1}+c\gamma_n.
      \end{equation}
      Then $\sup_{n\in \N_0}e_n\leq \max\bigl\{e_0,\max\{\varepsilon^2,4c^2\kappa^{-2}\}+c\bigl[\textstyle \sup_{n\in\N}\gamma_n\bigr]\bigr\}$.
\end{athm}
\end{tcolorbox}
\end{samepage}
\begin{aproof}
Throughout this proof let $\tau\in [0,\infty]$ satisfy $\tau=\textstyle \sup_{n\in\N}\gamma_n$.
    In the following we prove that for all $n\in \N_0$ it holds that
    \begin{equation}\llabel{need to prove}
        e_n\leq \max\{e_0,\max\{\varepsilon^2,4c^2\kappa^{-2}\}+c\tau\}.
    \end{equation}
    We prove \lref{need to prove} by induction on $n\in \N_0$. For the base case $n=0$ note that
    \begin{equation}\llabel{ind: eq1}
        e_0\leq  \max\{e_0,\max\{\varepsilon^2,4c^2\kappa^{-2}\}+c\tau\}.
    \end{equation}
    This establishes \lref{need to prove} in the base case $n=0$. For the induction step we assume that there exists $n\in \N$ which satisfies that
    \begin{equation}\llabel{ind: assume}
        e_{n-1}\leq \max\{e_0,\max\{\varepsilon^2,4c^2\kappa^{-2}\}+c\tau\}.
    \end{equation}
    In the following we prove that
    \begin{equation}\llabel{ind: need to prove}
        e_n\leq \max\{e_0,\max\{\varepsilon^2,4c^2\kappa^{-2}\}+c\tau\}.
    \end{equation}
    In our proof of \lref{ind: need to prove} we distinguish between the case $e_{n-1}\geq  \max\{\varepsilon^2,4c^2\kappa^{-2}\}$ and the case $e_{n-1}\leq \max\{\varepsilon^2,4c^2\kappa^{-2}\}$. We first prove \lref{ind: need to prove} in the case 
    \begin{equation}\llabel{case 1}
        e_{n-1}\geq  \max\{\varepsilon^2,4c^2\kappa^{-2}\}.
    \end{equation}
    \startnewargseq
    \argument{\lref{case 1};}{that
\begin{equation}\llabel{eq1}
   \kappa e_{n-1}\geq 2c\sqrt{e_{n-1}}\geq c\sqrt{e_{n-1}}+c\varepsilon.
\end{equation}}
\argument{\lref{eq1};}{that
\begin{equation}\llabel{eqq1.1}
    \frac{\kappa\gamma_ne_{n-1}}{\sqrt{e_{n-1}}+\varepsilon}\geq c\gamma_n.
\end{equation}}
\argument{\lref{eqq1.1};\lref{assume};\lref{ind: assume}}{that
\begin{equation}
    e_{n}\leq e_{n-1}\leq \max\{e_0,\max\{\varepsilon^2,4c^2\kappa^{-2}\}+c\tau\}.
\end{equation}}
This proves \lref{ind: need to prove} in the case $e_{n-1}\geq \max\{\varepsilon^2,4c^2\kappa^{-2}\}$. We now prove \lref{ind: need to prove} in the case
\begin{equation}\llabel{case 2}
e_{n-1}\leq \max\{\varepsilon^2,4c^2\kappa^{-2}\}.
\end{equation}
\startnewargseq
\argument{\lref{assume};\lref{case 2};;the fact that $\textstyle \sup_{n\in\N}\gamma_n=\tau$}{that
\begin{equation}
    e_n\leq e_{n-1}+c\gamma_n\leq e_{n-1}+c\tau\leq \max\{\varepsilon^2,4c^2\kappa^{-2}\}+c\tau\leq \max\{e_0,\max\{\varepsilon^2,4c^2\kappa^{-2}\}+c\tau\}.
\end{equation}}
This establishes \lref{ind: need to prove} in the case $e_{n-1}\leq \max\{\varepsilon^2,4c^2\kappa^{-2}\}$. Note that \lref{ind: eq1}, \lref{ind: need to prove}, and induction prove \lref{need to prove}.
\end{aproof}

In \cref{conclude: MUON bound 0} in \cref{theo: MUON bound 0} below we present the promised qualitative a priori estimate for a general class of \MUON-type optimization processes. In \cref{theo: MUON bound 0} the optimization process is denoted by $( \Theta_n )_{ n \in \N_0 }$ and takes values in the set $\cA = \times_{i=1}^\delta\R^{d_i\times\fd_i}$.

    \begin{tcolorbox}[colback=white!95!gray,
                  colframe=black,
                  boxrule=0.5pt,
                  sharp corners,
                  enhanced,
                  breakable,
                 ]
\begin{athm}{prop}{theo: MUON bound 0}
      Let $\delta,M\in \N$, $d_1,d_2,\dots,d_\delta,\fd_1,\fd_2,\dots,\fd_\delta\in \N$, $\alpha\in (0,1)$,  $\varepsilon,\kappa,\cK\in (0,\infty)$, let $\setX$ be a set, let $\cA=\times_{i=1}^\delta\R^{d_i\times\fd_i}$, let $\smalll=(\smalll(\theta,x))_{(\theta,x)\in \cA\times \setX}\colon \cA\times \setX\to\R$ satisfy for all $x\in \setX$ that $(\cA\ni\theta\mapsto \smalll(\theta,x)\in \R)\in C^1(\cA,\R)$, assume\footnotemark\! for all  $\theta,\vartheta\in \cA$, $x,y\in \setX$ that $ \spro{(\nabla_{\theta}\smalll)(\theta,x)-(\nabla_{\theta}\smalll)(\vartheta,x),\theta-\vartheta}\geq \kappa\|\theta-\vartheta\|^2$,
      \begin{equation}\llabel{assume}
          \|(\nabla_\theta\smalll)(\theta,x)-(\nabla_\theta\smalll)(\vartheta,x)\|\leq \cK\|\theta-\vartheta\|,\qqandqq \|(\nabla_\theta\smalll)(\theta,x)-(\nabla_\theta\smalll)(\theta,y)\|\leq  \cK,
      \end{equation}
       for every $n,m\in \N$ let $X_{n,m}\in\setX$, let $(\gamma_n)_{n\in \N}\subseteq (0,\infty)$ be non-increasing, assume for all $n\in \N$ that $\sum_{i=1}^n(\alpha^{n-i}\gamma_i)\leq \cK\gamma_n$, for every $i\in \{1,2,\dots,\delta\}$ let $\sgn_i\colon\R^{d_i\times\fd_i}\to\R^{d_i\times\fd_i}$ satisfy for all  $\theta\in \R^{d_i\times\fd_i}$ that
      \begin{equation}\llabel{def: Pi}
        \|\sgn_i(\theta)\|\leq\cK\qqandqq  \spro{\theta,\sgn_i(\theta)}\geq \kappa\|\theta\|^2(\|\theta\|+\varepsilon)^{-1},
      \end{equation}
let $\Theta=(\Theta^1,\dots,\Theta^\delta)\colon \N_0\to\cA$ and $\bfm=(\bfm^1,\dots,\bfm^\delta)\colon \N_0\to\cA$ satisfy for all $n\in \N$, $i\in \{1,2,\dots,\delta\}$ that 
         \begin{equation}\llabel{def: bfm}
            \bfm_0=0,\qquad \bfm_n= \alpha \bfm_{n-1}+(1-\alpha)\bigl[\textstyle \frac 1M \sum_{m=1}^M(\nabla_\theta\smalll)(\Theta_{n-1},X_{n,m})\bigr],
            \end{equation}
            \begin{equation}\llabel{def: Theta}
            \text{and}\qquad \Theta_n^i=\Theta_{n-1}^i-\gamma_n\sgn_i(\bfm_n^i),
         \end{equation}
          and let $\vartheta\in \cA$, $\bfx\in \setX$, $\varrho,v,\bfA\in \R$ satisfy $\smalll(\vartheta,\bfx)=\inf_{\theta\in \cA}\smalll(\theta,\bfx)$,
         \begin{equation}\llabel{def: vartheta}
             \textstyle \varrho=\cK\delta\|(\nabla_{\theta}\smalll)(\vartheta,\bfx)\|+\cK^2\delta\|\Theta_0-\vartheta\|+\cK^3\delta^3\gamma_1+\cK^2\delta^2\gamma_1+\cK,\qquad v=\delta ( 3 \kappa \varrho + \cK \varrho + \cK^3 \gamma_1 ) ,
         \end{equation}
         \begin{equation}\llabel{def: bfA}
         \begin{split}
            \text{and}\qquad\textstyle  \bfA=\max\{\kappa^{-2}\max\{\varepsilon^2\cK,4\kappa^{-2}\delta^2\cK v^2\}+\gamma_1 v,
              \smalll(\Theta_0,\bfx)-\smalll(\vartheta,\bfx)\},
      \end{split}
         \end{equation}
     Then 
     \begin{equation}\label{conclude: MUON bound 0}
        \textstyle\sup_{n\in \N_0} \|\Theta_n\|\leq[2\bfA\kappa^{-1}]^{1/2}+\|\vartheta\|.
     \end{equation}
\end{athm}
\end{tcolorbox}
\footnotetext{Note that for all $N\in \N$, $v_1, v_2, \dots, v_{N }, w_1, w_2, \dots, w_{ N } \in \N$, $\theta=( ( \theta^n_{ i, j } )_{ (i, j) \in \{1,2,\dots,v_n\}\times\{1,2,\dots,w_n\} } )_{ n \in \{1,2,\dots,N\} } \in \times_{ n = 1 }^{N} \R^{ v_n \times w_n }$, $\vartheta=( ( \vartheta^n_{ i, j } )_{ (i, j) \in \{1,2,\dots,v_n\}\times\{1,2,\dots,w_n\} } )_{ n \in \{1,2,\dots,N\} } \in \times_{ n = 1 }^{N} \R^{ v_n \times w_n }$ it holds that $\spro{\theta,\vartheta} = \sum_{n=1}^N\sum_{i=1}^{v_n}\sum_{j=1}^{w_n}\theta_{i,j}^n\vartheta_{i,j}^n$.}
\begin{aproof}
Throughout this proof let $\bbM=(\bbM^1,\dots,\bbM^\delta)\colon \N_0\to\cA$ satisfy for all $n\in \N$ that
\begin{equation}\llabel{def: bbM}
    \bbM_0=0\qqandqq \bbM_n=\alpha\bbM_{n-1}+(1-\alpha)(\nabla_\theta\smalll)(\Theta_{n-1},\bfx)
\end{equation}
and let $\cL=(\cL(\theta,x))_{(\theta,x)=((\theta_1,\dots,\theta_\delta),x)\in \cA\times \setX} \colon\cA\times \setX\to\R$ satisfy for all $\theta\in \cA$, $x\in \setX$ that $\cL(\theta,x)=\smalll(\theta,x)$.
 \argument{\lref{assume};the fundamental theorem of calculus;the Cauchy-Schwarz inequality}{that for all $\theta_1,\theta_2\in \cA$ it holds that
\begin{equation}\llabel{eqp1}
\begin{split}
    &\cL(\theta_1,\bfx)-\cL(\theta_2,\bfx)
    =\int_{0}^1\spro{(\nabla_\theta\cL)(\theta_2+t(\theta_1-\theta_2),\bfx),\theta_1-\theta_2}\,\d t\\
    &=\spro{(\nabla_\theta\cL)(\theta_2,\bfx),\theta_1-\theta_2}+\int_{0}^1\spro{(\nabla_\theta\cL)(\theta_2+t(\theta_1-\theta_2),\bfx)-(\nabla_\theta\cL)(\theta_2,\bfx),\theta_1-\theta_2}\,\d t\\
    &\leq \spro{(\nabla_\theta\cL)(\theta_2,\bfx),\theta_1-\theta_2}+\int_{0}^1\|(\nabla_\theta\cL)(\theta_2+t(\theta_1-\theta_2),\bfx)-(\nabla_\theta\cL)(\theta_2,\bfx)\|\|\theta_1-\theta_2\|\,\d t\\
    &\leq \spro{(\nabla_\theta\cL)(\theta_2,\bfx),\theta_1-\theta_2}+\int_{0}^1\cK\|t(\theta_1-\theta_2)\|\|\theta_1-\theta_2\|\,\d t\\
    &\leq \spro{(\nabla_\theta\cL)(\theta_2,\bfx),\theta_1-\theta_2}+\cK\|\theta_1-\theta_2\|^2.
    \end{split}
\end{equation}}
\argument{\lref{def: Theta};\lref{eqp1}}{that for all $n\in \N$ it holds that
    \begin{equation}\llabel{eq1}
    \begin{split}
    \cL(\Theta_{n},\bfx)-\cL(\Theta_{n-1},\bfx)&\leq \spro{(\nabla_\theta\cL)(\Theta_{n-1},\bfx),\Theta_{n}-\Theta_{n-1}}+\cK\|\Theta_n-\Theta_{n-1}\|^2\\
    &\leq\textstyle -\gamma_n\bigl[\sum_{i=1}^\delta\spro{(\nabla_{\theta_i}\cL)(\Theta_{n-1},\bfx),\sgn_i(\bfm_n^i)}\bigr]+\cK(\gamma_n)^2\bigl[\sum_{i=1}^\delta\|\sgn_i(\bfm_n^i)\|^2\bigr].
    \end{split}
    \end{equation}}
    \argument{\lref{eq1};\lref{def: Pi}}{for all $n\in \N$ that
    \begin{equation}\llabel{eq2}
    \begin{split}
        \cL(\Theta_{n},\bfx)&\leq \textstyle\cL(\Theta_{n-1},\bfx)-\gamma_n\bigl[\sum_{i=1}^\delta\spro{(\nabla_{\theta_i}\cL)(\Theta_{n-1},\bfx),\sgn_i(\bfm_n^i)}\bigr]+\cK^3\delta(\gamma_n)^2.
        \end{split}
    \end{equation}}
      \argument{\lref{def: Pi};\lref{def: Theta}}{that for all $n\in \N$, $i\in \{1,2,\dots,\delta\}$ it holds that
    \begin{equation}\llabel{eq3t1}
        \|\Theta_n^i-\Theta_{n-1}^i\|=\gamma_n\|\sgn_i(\bfm_n^i)\|\leq \cK\gamma_n.
    \end{equation}}
    \argument{\lref{eq3t1};the fact that for all $\theta=(\theta_1,\dotsm,\theta_\delta)\in \cA$ it holds that $\|\theta\|\leq \sum_{i=1}^\delta \|\theta_i\|$}{that for all $n\in\N$ it holds that
    \begin{equation}\llabel{eq3}
          \|\Theta_n-\Theta_{n-1}\|\textstyle\leq \sum_{i=1}^\delta  \|\Theta_n^i-\Theta_{n-1}^i\|\leq \cK\delta\gamma_n.
    \end{equation}}
    \argument{\lref{assume};the fact that 
    $\bfx\in \setX$}{that for all $n\in \N$ it holds that
    \begin{equation}\llabel{eq3''}
       \| (\nabla_{\theta}\cL)(\Theta_{n},\bfx)-(\nabla_{\theta}\cL)(\Theta_{n-1},\bfx)\|\leq\cK\|\Theta_{n}-\Theta_{n-1}\|\leq \cK\delta\|\Theta_{n}-\Theta_{n-1}\|.
    \end{equation}}
    \argument{\lref{assume};\lref{def: vartheta};\lref{def: bbM};\lref{eq3};\lref{eq3''};the fact that for all $n\in \N$ it holds that $\sum_{i=1}^n(\alpha^{n-i}\gamma_i)\leq \cK\delta\gamma_n$;\cref{lem: bfm analysis 1};}{that for all $n\in \N$ it holds that
    \begin{equation}\llabel{eqq4}
    \begin{split}
       &\|\bbM_{n}-(\nabla_{\theta}\cL)(\Theta_{n-1},\bfx)\|\leq (\cK\delta(\gamma_1)^{-1}\|\bbM_0-(\nabla_{\theta}\cL)(\Theta_0,\bfx)\|+\cK^3\delta^3+\cK^2\delta^2)\gamma_n\\
       &=(\cK\delta(\gamma_1)^{-1}\|(\nabla_{\theta}\cL)(\Theta_0,\bfx)\|+\cK^3\delta^3+\cK^2\delta^2)\gamma_n\\
       &\leq (\cK\delta(\gamma_1)^{-1}\|(\nabla_{\theta}\cL)(\vartheta,\bfx)\|+\cK\delta(\gamma_1)^{-1}\|(\nabla_{\theta}\cL)(\Theta_0,\bfx)-(\nabla_{\theta}\cL)(\vartheta,\bfx)\|+\cK^3\delta^3+\cK^2\delta^2)\gamma_n\\
       &\leq (\cK\delta(\gamma_1)^{-1}\|(\nabla_{\theta}\cL)(\vartheta,\bfx)\|+\cK^2\delta(\gamma_1)^{-1}\|\Theta_0-\vartheta\|+\cK^3\delta^3+\cK^2\delta^2)\gamma_n .
       \end{split}
    \end{equation}}
     \argument{\lref{def: bfm};\lref{def: bbM};\cref{momentum:representation}}{that for all $n\in \N_0$ it holds that
    \begin{equation}\llabel{eqq4.1}
        \bfm_{n}=\textstyle(1-\alpha)\biggl[\sum\limits_{k=1}^n\alpha^{n-k}\textstyle \frac 1M \bigl[\textstyle\sum_{m=1}^M(\nabla_\theta\cL)(\Theta_{k-1},X_{k,m})\bigr]\biggr]
        \end{equation}
        \begin{equation}\llabel{eqq4.2}
            \text{and}\qquad \textstyle \bbM_{n}=(1-\alpha)\bigl[\textstyle\sum_{k=1}^n\alpha^{n-k}(\nabla_\theta\cL)(\Theta_{k-1},\bfx)\bigr].
    \end{equation}}
    \argument{\lref{eqq4.2};\lref{assume};}{that for all $n\in \N$ it holds that
    \begin{equation}\llabel{eqq4.3}
    \begin{split}
        \|\bfm_n-\bbM_n\|&=\biggl\|\textstyle(1-\alpha)\biggl[\sum\limits_{k=1}^n\alpha^{n-k}\bigl(\textstyle \frac 1M \bigl[\textstyle\sum_{m=1}^M(\nabla_\theta\cL)(\Theta_{k-1},X_{k,m})\bigr]-(\nabla_\theta\cL)(\Theta_{k-1},\bfx)\bigr)\biggr]\biggr\|\\
        &\leq \textstyle(1-\alpha)\biggl[\sum\limits_{k=1}^n\alpha^{n-k}\textstyle \frac 1M \bigl[\textstyle\sum_{m=1}^M\|(\nabla_\theta\cL)(\Theta_{k-1},X_{k,m})-(\nabla_\theta\cL)(\Theta_{k-1},\bfx)\|\bigr]\biggr]\\
        &\leq \textstyle(1-\alpha)\biggl[\sum\limits_{k=1}^n\alpha^{n-k}\textstyle \frac 1M \bigl[\textstyle\sum_{m=1}^M\cK\bigr]\biggr]\leq \cK.
        \end{split}
    \end{equation}}
    \argument{\lref{def: vartheta};\lref{eqq4.3};\lref{eqq4}; the fact that for all $n\in \N$ it holds that $\gamma_n\leq \gamma_1$}{that for all $n\in \N$ it holds that
    \begin{equation}\llabel{eq4}
    \begin{split}
       & \|\bfm_{n}-(\nabla_{\theta}\cL)(\Theta_{n-1},\bfx)\|\leq \|\bbM_{n}-(\nabla_{\theta}\cL)(\Theta_{n-1},\bfx)\|+\|\bfm_n-\bbM_n\|\\
        &\leq (\cK\delta(\gamma_1)^{-1}\|(\nabla_{\theta}\cL)(\vartheta,\bfx)\|+\cK^2\delta(\gamma_1)^{-1}\|\Theta_0-\vartheta\|+\cK^3\delta^3+\cK^2\delta^2)\gamma_n +\cK\\
        &\leq  (\cK\delta(\gamma_1)^{-1}\|(\nabla_{\theta}\cL)(\vartheta,\bfx)\|+\cK^2\delta(\gamma_1)^{-1}\|\Theta_0-\vartheta\|+\cK^3\delta^3+\cK^2\delta^2)\gamma_1 +\cK\\
        &=\cK\delta\|(\nabla_{\theta}\cL)(\vartheta,\bfx)\|+\cK^2\delta\|\Theta_0-\vartheta\|+\cK^3\delta^3\gamma_1+\cK^2\delta^2\gamma_1+\cK\\
        &=\varrho
        .
        \end{split}
    \end{equation}}
     \argument{the fact that $\varepsilon>0$;the triangle inequality}{that for all $x,y\in [0,\infty)$ it holds that
    \begin{equation}\llabel{eq4.0}
    \begin{split}
        \textstyle\bigl|\frac{x^2}{x+\varepsilon}-\frac{y^2}{y+\varepsilon}\bigr|&\textstyle=\bigl|\frac{x^2y+\varepsilon x^2-y^2x-\varepsilon y^2}{(x+\varepsilon)(y+\varepsilon)}\bigr|\leq \bigl|\frac{xy(x-y)}{(x+\varepsilon)(y+\varepsilon)}\bigr|+\bigl|\frac{\varepsilon(x-y)(x+y)}{(x+\varepsilon)(y+\varepsilon)}\bigr|\\
        &\leq \textstyle \bigl|\frac{xy(x-y)}{(x+\varepsilon)(y+\varepsilon)}\bigr|+\bigl|\frac{\varepsilon(x-y)x}{(x+\varepsilon)(y+\varepsilon)}\bigr|+\bigl|\frac{\varepsilon(x-y)y}{(x+\varepsilon)(y+\varepsilon)}\bigr|\leq |x-y|+|x-y|+|x-y|\\
        &=\textstyle3|x-y|.
        \end{split}
    \end{equation}}
    \argument{\lref{eq4.0};\lref{def: Pi};\lref{eq4};}{that for all $n\in \N$, $i\in \{1,2,\dots,\delta\}$ it holds that
    \begin{equation}\llabel{eq4.1}
    \begin{split}
        \spro{\bfm_n^i,\sgn_i(\bfm_n^i)}\geq \frac{\kappa \|\bfm_n^i\|^2}{\|\bfm_n^i\|+\varepsilon}&\geq  \frac{\kappa \|(\nabla_{\theta_i}\cL)(\Theta_{n-1},\bfx)\|^2}{\|(\nabla_{\theta_i}\cL)(\Theta_{n-1},\bfx)\|+\varepsilon}-3\kappa\bigl|\|\bfm_n^i\|-\|(\nabla_{\theta_i}\cL)(\Theta_{n-1},\bfx)\|\bigr|\\
        &\geq \frac{\kappa \|(\nabla_{\theta_i}\cL)(\Theta_{n-1},\bfx)\|^2}{\|(\nabla_{\theta_i}\cL)(\Theta_{n-1},\bfx)\|+\varepsilon}-3\kappa\|\bfm_n^i-(\nabla_{\theta_i}\cL)(\Theta_{n-1},\bfx)\|\\
        &\geq \frac{\kappa \|(\nabla_{\theta_i}\cL)(\Theta_{n-1},\bfx)\|^2}{\|(\nabla_{\theta_i}\cL)(\Theta_{n-1},\bfx)\|+\varepsilon}-3\kappa\varrho .
        \end{split}
    \end{equation}}
    \argument{\lref{eq4.1};\lref{def: Pi};\lref{eq4};the Cauchy-Schwarz inequality}{that for all $n\in \N$, $i\in \{1,2,\dots,\delta\}$ it holds that
    \begin{equation}\llabel{arg1}
    \begin{split}
       & \spro{(\nabla_{\theta_i}\cL)(\Theta_{n-1},\bfx),\sgn_i(\bfm_n^i)}=\spro{\bfm_n^i,\sgn_i(\bfm_n^i)}+\spro{(\nabla_{\theta_i}\cL)(\Theta_{n-1},\bfx)-\bfm_n^i,\sgn_i(\bfm_n^i)}\\
        &\geq \frac{\kappa \|(\nabla_{\theta_i}\cL)(\Theta_{n-1},\bfx)\|^2}{\|(\nabla_{\theta_i}\cL)(\Theta_{n-1},\bfx)\|+\varepsilon}-3\kappa\varrho-\|(\nabla_{\theta_i}\cL)(\Theta_{n-1},\bfx)-\bfm_n^i\|\|\sgn_i(\bfm_n^i)\|\\
        &\geq \frac{\kappa \|(\nabla_{\theta_i}\cL)(\Theta_{n-1},\bfx)\|^2}{\|(\nabla_{\theta_i}\cL)(\Theta_{n-1},\bfx)\|+\varepsilon}-3\kappa\varrho-\cK\varrho.
        \end{split}
    \end{equation}}
    \argument{\lref{arg1};\lref{eq2}; the fact that for all $n\in \N$ it holds that $\gamma_n\leq \gamma_1$}{that for all $n\in \N$ it holds that
    \begin{equation}\llabel{eq7.1}
    \begin{split}
       \cL(\Theta_n,\bfx)&\leq  \cL(\Theta_{n-1},\bfx)-\gamma_n\sum_{i=1}^\delta\biggl( \frac{\kappa \|(\nabla_{\theta_i}\cL)(\Theta_{n-1},\bfx)\|^2}{\|(\nabla_{\theta_i}\cL)(\Theta_{n-1},\bfx)\|+\varepsilon}-3\kappa\varrho-\cK\varrho\biggr)+\cK^3\delta(\gamma_n)^2\\
       &\leq   \cL(\Theta_{n-1},\bfx)-\gamma_n\sum_{i=1}^\delta\biggl( \frac{\kappa \|(\nabla_{\theta_i}\cL)(\Theta_{n-1},\bfx)\|^2}{\|(\nabla_{\theta_i}\cL)(\Theta_{n-1},\bfx)\|+\varepsilon}-3\kappa\varrho-\cK\varrho\biggr)+\cK^3\delta\gamma_1\gamma_n\\
       &\leq  \cL(\Theta_{n-1},\bfx)-\gamma_n\sum_{i=1}^\delta\biggl( \frac{\kappa \|(\nabla_{\theta_i}\cL)(\Theta_{n-1},\bfx)\|^2}{\|(\nabla_{\theta_i}\cL)(\Theta_{n-1},\bfx)\|+\varepsilon}-3\kappa\varrho-\cK\varrho-\cK^3\gamma_1\biggr).
       \end{split}
    \end{equation}}
    \argument{\lref{eq7.1};the Bunyakovsky inequality;the fact that for all $n,m\in \N$, $x=(x_1,\dots,x_n)\in (\R^m)^n$ it holds that $ \|x\|\leq \sum_{i=1}^n \|x_i\|\leq n\|x\|$}{that for all it holds $n\in \N$ that
    \begin{equation}\llabel{eq8}
    \begin{split}
        &\cL(\Theta_n,\bfx)-\cL(\vartheta,\bfx)\\
        &\leq \cL(\Theta_{n-1},\bfx)-\cL(\vartheta,\bfx)-\gamma_n\sum_{i=1}^\delta\biggl( \frac{\kappa \|(\nabla_{\theta_i}\cL)(\Theta_{n-1},\bfx)\|^2}{\|(\nabla_{\theta_i}\cL)(\Theta_{n-1},\bfx)\|+\varepsilon}-3\kappa\varrho-\cK\varrho-\cK^3\gamma_1\biggr)\\
        &\leq \cL(\Theta_{n-1},\bfx)-\cL(\vartheta,\bfx)-\gamma_n\kappa\biggl( \frac{\bigl[\sum_{i=1}^\delta \|(\nabla_{\theta_i}\cL)(\Theta_{n-1},\bfx)\|\bigr]^2}{\sum_{i=1}^\delta\|(\nabla_{\theta_i}\cL)(\Theta_{n-1},\bfx)\|+\varepsilon\delta}\biggr)+(3\kappa\varrho+\cK\varrho+\cK^3\gamma_1)\gamma_n\delta\\
        &\leq \cL(\Theta_{n-1},\bfx)-\cL(\vartheta,\bfx)-\gamma_n\kappa\biggl( \frac{ \|(\nabla_{\theta}\cL)(\Theta_{n-1},\bfx)\|^2}{\delta\|(\nabla_{\theta}\cL)(\Theta_{n-1},\bfx)\|+\varepsilon\delta}\biggr)+(3\kappa\varrho+\cK\varrho+\cK^3\gamma_1)\gamma_n\delta.
        \end{split}
    \end{equation}}
    \argument{the fact that $\smalll(\vartheta,\bfx)=\inf_{\theta\in \cA}\smalll(\theta,\bfx)$;the fact that for all $\theta\in \cA$ it holds that
          $\spro{\theta-\vartheta,(\nabla_\theta\cL)(\theta,\bfx)}\geq \kappa\|\theta-\vartheta\|^2$;\unskip, \eg, \cite[Lemma 5.7.29]{ArBePhi2024}}{that for all $\theta\in\cA$ it holds that
          \begin{equation}\llabel{eq7.2.1}
         \cL(\theta,\bfx)-\cL(\vartheta,\bfx)\geq \frac{\kappa}{2}\|\theta-\vartheta\|^2\geq 0.
          \end{equation}}
          \argument{the fact that $\smalll(\vartheta,\bfx)=\inf_{\theta\in \cA}\smalll(\theta,\bfx)$;the fact that $(\cA\ni\theta\mapsto\cL(\theta,\bfx)\in \R)\in C^{1}(\cA,\R)$}{that \llabel{arg2} $(\nabla_\theta\cL)(\vartheta,\bfx)=0$\dott}
\argument{\lref{arg2};\lref{eqp1};\lref{eq7.2.1};}{that for all $n\in \N_0$ it holds that
\begin{equation}\llabel{eq7.2.3}
    0\leq \cL(\Theta_n,\bfx)-\cL(\vartheta,\bfx)\leq \cK\|\Theta_n-\vartheta\|^2.
\end{equation}}
          \argument{the fact that for all $\theta\in \cA$ it holds that
          $\spro{\theta-\vartheta,(\nabla_\theta\cL)(\theta,\bfx)}\geq \kappa\|\theta-\vartheta\|^2$;the Cauchy-Schwarz inequality}{that for all $\theta\in \cA$ it holds that
          \begin{equation}\llabel{eq7.2.4}
         \kappa\|\theta-\vartheta\|^2\leq \|\theta-\vartheta\|\|(\nabla_\theta\cL)(\theta,\bfx)\|.
          \end{equation}}
          \argument{\lref{eq7.2.4};}{for all $\theta\in \cA$ that
          \begin{equation}\llabel{eq7.2.5}
              \|(\nabla_\theta\cL)(\theta,\bfx)\|\geq\kappa\|\theta-\vartheta\|.
          \end{equation}}
           \argument{\lref{eq7.2.5};\lref{eq7.2.3};}{that for all $n\in \N_0$ it holds that
          \begin{equation}\llabel{eq7.2.6}
              \|(\nabla_\theta\cL)(\Theta_n,\bfx)\|\geq \kappa\cK^{-1/2}\sqrt{\cL(\Theta_n,\bfx)-\cL(\vartheta,\bfx)}.
          \end{equation}}
   \argument{\lref{eq7.2.6};\lref{eq8};the fact that $[0,\infty)\ni x\mapsto x^2(x+\varepsilon)^{-1}\in \R$ is increasing}{that for all $n\in \N$ it holds that
   \begin{equation}\llabel{eq9}
   \begin{split}
         &\cL(\Theta_n,\bfx)-\cL(\vartheta,\bfx)\\&\leq\textstyle    \cL(\Theta_{n-1},\bfx)-\cL(\vartheta,\bfx)-\gamma_n\Bigl( \frac{\kappa(\kappa\cK^{-1/2})^2(\cL(\Theta_{n-1},\bfx)-\cL(\vartheta,\bfx))}{\delta\kappa\cK^{-1/2}\sqrt{\cL(\Theta_{n-1},\bfx)-\cL(\vartheta,\bfx)}+\delta\varepsilon}\Bigr)+(3\kappa\varrho+\cK\varrho+\cK^3\gamma_1)\gamma_n\delta\\
         &\textstyle=\cL(\Theta_{n-1},\bfx)-\cL(\vartheta,\bfx)- \gamma_n\Bigl(\frac{\kappa^2\cK^{-1/2} (\cL(\Theta_{n-1},\bfx)-\cL(\vartheta,\bfx))}{\delta\sqrt{\cL(\Theta_{n-1},\bfx)-\cL(\vartheta,\bfx)}+\delta\varepsilon\kappa^{-1}\cK^{1/2}}\Bigr)+(3\kappa\varrho+\cK\varrho+\cK^3\gamma_1)\gamma_n\delta\\
         &= \textstyle\cL(\Theta_{n-1},\bfx)-\cL(\vartheta,\bfx)- \gamma_n\Bigl(\frac{\kappa^2\cK^{-1/2}\delta^{-1} (\cL(\Theta_{n-1},\bfx)-\cL(\vartheta,\bfx))}{\sqrt{\cL(\Theta_{n-1},\bfx)-\cL(\vartheta,\bfx)}+\varepsilon\kappa^{-1}\cK^{1/2}}\Bigr)+(3\kappa\varrho+\cK\varrho+\cK^3\gamma_1)\gamma_n\delta.
         \end{split}
   \end{equation}}
   \argument{\lref{eq9};\lref{def: bfA};the fact that $\gamma_1=\sup_{n\in \N}\gamma_n$;\cref{lem: recurrence}}{that for all $n\in \N_0$ it holds that
   \begin{equation}\llabel{eq10}
   \begin{split}
       \cL(\Theta_{n},\bfx)-\cL(\vartheta,\bfx)
      &\leq \max\{\cL(\Theta_0,\bfx)-\cL(\vartheta,\bfx),\max\{\varepsilon^2\kappa^{-2}\cK,4(3\kappa\varrho+\cK\varrho+\cK^3\gamma_1)^2\kappa^{-4}\cK\delta^4\}\\
      &\quad+(3\kappa\varrho+\cK\varrho+\cK^3\gamma_1)\delta\gamma_1\}\\
      &=\bfA.
       \end{split}
   \end{equation}}
   \argument{\lref{eq10};\lref{eq7.2.1}}{for all $n\in \N_0$ that
\begin{equation}\llabel{eq11}
    \|\Theta_n-\vartheta\|^2\leq2\kappa^{-1}(\cL(\Theta_n,\bfx)-\cL(\vartheta,\bfx))\leq2\bfA\kappa^{-1}.
\end{equation}}
\argument{\lref{eq11};}{that for all $n\in \N_0$ it holds that
\begin{equation}\llabel{eq12}
    \|\Theta_n\|\leq \|\Theta_n-\vartheta\|+\|\vartheta\|\leq [2\bfA\kappa^{-1}]^{1/2}+\|\vartheta\|.
\end{equation}}
\argument{\lref{eq12};}{\cref{conclude: MUON bound 0}\dott}
\end{aproof}
\subsection{Qualitative pathwise a priori bounds for MUON type optimization processes}\label{subsec: pathwise bound}

In the next result, \cref{theo: MUON bound pre}, we employ the quantitative a priori estimate for \MUON-type optimization processes in \cref{theo: MUON bound 0} to establish in \cref{conclude: theo: MUON bound 0} a qualitative a priori estimate for such optimization processes.

  \begin{samepage}
    \begin{tcolorbox}[colback=white!95!gray,
                  colframe=black,
                  boxrule=0.5pt,
                  sharp corners,
                  enhanced,
                  breakable,
                 ]
\begin{athm}{theorem}{theo: MUON bound pre}
      Let $\delta\in \N$, $d_1,d_2,\dots,d_\delta,\fd_1,\fd_2,\dots,\fd_\delta\in \N$, $\alpha\in (0,1)$,  $\varepsilon,\kappa,\cK\in (0,\infty)$, let $\setX$ be a set, let $\cA=\times_{i=1}^\delta\R^{d_i\times\fd_i}$, let $\smalll=(\smalll(\theta,x))_{(\theta,x)\in \cA\times \setX}\colon \cA\times \setX\to\R$ satisfy for all $x\in \setX$ that $(\cA\ni\theta\mapsto \smalll(\theta,x)\in \R)\in C^1(\cA,\R)$, assume for all  $\theta,\vartheta\in \cA$, $x,y\in \setX$ that $ \spro{(\nabla_{\theta}\smalll)(\theta,x)-(\nabla_{\theta}\smalll)(\vartheta,x),\theta-\vartheta}\geq \kappa\|\theta-\vartheta\|^2$, $\|(\nabla_\theta\smalll)(\theta,x)-(\nabla_\theta\smalll)(\vartheta,x)\|\leq \cK\|\theta-\vartheta\|$, and $ \|(\nabla_\theta\smalll)(\theta,x)-(\nabla_\theta\smalll)(\theta,y)\|\leq  \cK$, let $(\gamma_n)_{n\in \N}\subseteq (0,\infty)$ be non-increasing, assume $\limsup_{n\to\infty}((\gamma_{n+1})^{-1}\alpha\gamma_n)<1$, and for every $i\in \{1,2,\dots,\delta\}$ let $\sgn_i\colon\R^{d_i\times\fd_i}\to\R^{d_i\times\fd_i}$ satisfy for all  $\theta\in \R^{d_i\times\fd_i}$ that
      \begin{equation}\llabel{def: Pi}
        \textcolor{magenta}{\|\sgn_i(\theta)\|\leq\cK}\qqandqq  \textcolor{magenta}{\spro{\sgn_i(\theta),\theta}\geq \kappa\|\theta\|^2(\|\theta\|+\varepsilon)^{-1}}.
      \end{equation}
     Then there exists $\fC\in \R$ such that for every $M\in \N$, every $\Theta=(\Theta^{1},\dots,\Theta^{\delta})\colon \N_0\to\cA$, every $\bfm=(\bfm^{1},\dots,\bfm^{\delta})\colon \N_0\to\cA$, and every $X=(X_{n,m})_{(n,m)\in \N^2}\colon \N^2\to\setX$ with the property that
     for all $n\in \N$, $i\in \{1,2,\dots,\delta\}$ it holds that 
         \begin{equation}\llabel{def: bfm}
            \textcolor{magenta}{\bfm_0=0},\qquad \textcolor{magenta}{\bfm_n= \alpha \bfm_{n-1}+(1-\alpha)\bigl[\textstyle \frac 1M \sum_{m=1}^M(\nabla_\theta\smalll)(\Theta_{n-1},X_{n,m})\bigr]},
            \end{equation}
            \begin{equation}\llabel{def: Theta}
            \textcolor{magenta}{ \|\Theta_0\|\leq \cK},\qquad\text{and}\qquad  \textcolor{magenta}{\Theta_n^{i}=\Theta_{n-1}^{i}-\gamma_n\sgn_i(\bfm_n^{i})}
         \end{equation}
         we have that
\begin{equation} \label{conclude: theo: MUON bound 0}
   \textcolor{magenta}{\textstyle\sup_{n\in \N_0}\|\Theta_n\|\leq \fC}.
\end{equation}
\end{athm}
\end{tcolorbox}
\end{samepage}
\begin{aproof}
Throughout this proof let $\bfx\in \setX$.
\argument{the fact that for all $\theta,\vartheta\in \cA$ it holds that $ \spro{(\nabla_{\theta}\smalll)(\theta,\bfx)-(\nabla_{\theta}\smalll)(\vartheta,\bfx),\theta-\vartheta}\geq \kappa\|\theta-\vartheta\|^2$;\cite[Corollary 5.7.22]{ArBePhi2024}}{that there exists $\vartheta\in \cA$ which satisfies
\begin{equation}\llabel{def: vartheta}
\smalll(\vartheta,\bfx)=\inf_{\theta\in \cA}\smalll(\theta,\bfx).    
\end{equation}}
\startnewargseq
   \argument{the assumption that $\limsup_{n\to\infty}((\gamma_{n+1})^{-1}\alpha\gamma_n)<1$;\cite[Lemma 2.6]{ShoAr2026};}{that there exists $\fC\in (0,\infty)$ such that for all $n\in\N$ it holds that
        \begin{equation}\llabel{eq3}
            \sum_{i=1}^n(\alpha^{n-i}\gamma_i)\leq \fC\gamma_n.
        \end{equation}}
        \argument{\lref{def: vartheta};\lref{eq3};\cref{theo: MUON bound 0}}{\cref{conclude: theo: MUON bound 0}\dott}
\end{aproof}
  \begin{samepage}

  In \cref{sec: convergence of stochastic muon} we develop an error analysis for the \MUON\ optimizer when applied to \SOPs. In the next result, \cref{theo: MUON bound}, we slightly reformulate the purely deterministic result in \cref{theo: MUON bound pre} to a stochastic statement involving a probability space so that it is more directly applicable to the \SOPs\ considered in \cref{sec: convergence of stochastic muon}.
    \begin{tcolorbox}[colback=white!95!gray,
                  colframe=black,
                  boxrule=0.5pt,
                  sharp corners,
                  enhanced,
                  breakable,
                 ]
\begin{athm}{cor}{theo: MUON bound}
      Let $(\Omega,\cF,\P)$ be a probability space, let $\delta\in \N$, $d_1,d_2,\dots,d_\delta,\fd_1,\fd_2,\dots,\fd_\delta\in \N$, $\alpha\in (0,1)$,  $\varepsilon,\kappa,\cK\in (0,\infty)$, let $(\setX,\cU)$ be a measurable space, let $\cA=\times_{i=1}^\delta\R^{d_i\times\fd_i}$, let $\smalll=(\smalll(\theta,x))_{(\theta,x)\in \cA\times \setX}\colon \cA\times \setX\to\R$ be measurable, assume for all $x\in \setX$ that $(\cA\ni\theta\mapsto \smalll(\theta,x)\in \R)\in C^1(\cA,\R)$, assume for all  $\theta,\vartheta\in \cA$, $x,y\in \setX$ that $ \spro{(\nabla_{\theta}\smalll)(\theta,x)-(\nabla_{\theta}\smalll)(\vartheta,x),\theta-\vartheta}\geq \kappa\|\theta-\vartheta\|^2$, $\|(\nabla_\theta\smalll)(\theta,x)-(\nabla_\theta\smalll)(\vartheta,x)\|\leq \cK\|\theta-\vartheta\|$, and $ \|(\nabla_\theta\smalll)(\theta,x)-(\nabla_\theta\smalll)(\theta,y)\|\leq  \cK$, let $X_{n,m}\colon \Omega\to\setX$, $(n,m)\in \N^2$, be random variables, let $(\gamma_n)_{n\in \N}\subseteq (0,\infty)$ be non-increasing, assume $\limsup_{n\to\infty}((\gamma_{n+1})^{-1}\alpha\gamma_n)<1$, and for every $i\in \{1,2,\dots,\delta\}$ let $\sgn_i\colon\R^{d_i\times\fd_i}\to\R^{d_i\times\fd_i}$ satisfy for all  $\theta\in \R^{d_i\times\fd_i}$ that
      \begin{equation}\llabel{def: Pi}
       \textcolor{magenta}{ \|\sgn_i(\theta)\|\leq\cK}\qqandqq  \textcolor{magenta}{\spro{\sgn_i(\theta),\theta}\geq \kappa\|\theta\|^2(\|\theta\|+\varepsilon)^{-1}}.
      \end{equation}
     Then there exists $\fC\in \R$ such that for every $M\in \N$, every stochastic process $\bfm=(\bfm^{1},\dots,\bfm^{\delta})\colon \N_0\times\Omega\to\cA$, and every stochastic process $\Theta=(\Theta^{1},\dots,\Theta^{\delta})\colon \N_0\times\Omega\to\cA$ with the property that $\P(\|\Theta_0\|\leq \cK)=1$ and with the property that
     for all $n\in \N$, $i\in \{1,2,\dots,\delta\}$ it holds that 
         \begin{equation}\llabel{def: bfm}
            \textcolor{magenta}{\bfm_0=0},\qquad \textcolor{magenta}{\bfm_n= \alpha \bfm_{n-1}+(1-\alpha)\bigl[\textstyle \frac 1M \sum_{m=1}^M(\nabla_\theta\smalll)(\Theta_{n-1},X_{n,m})\bigr]},
            \end{equation}
            \begin{equation}\llabel{def: Theta}
            \text{and}\qquad  \textcolor{magenta}{\Theta_n^{i}=\Theta_{n-1}^{i}-\gamma_n\sgn_i(\bfm_n^{i})}
         \end{equation}
         we have that
\begin{equation} \llabel{conclude}
   \textcolor{magenta}{\P\bigl(\textstyle\sup_{n\in \N_0}\|\Theta_n\|\leq \fC\bigr)=1}.
\end{equation}
\end{athm}
\end{tcolorbox}
\end{samepage}
\begin{aproof}
\argument{\cref{theo: MUON bound pre}}{\lref{conclude}\dott}
\end{aproof}
  \section{Error analysis for MUON and closely related optimization methods}\label{sec: convergence of stochastic muon}

   In this section we establish in \cref{muon bound convergence stochastic} in \cref{subsec: muon under boundness} below an error bound for the \MUON\ optimizer under the assumption that the \MUON\ optimization process is almost surely bounded by a fixed constant and, thereafter, we combine \cref{muon bound convergence stochastic} with the a priori bound in \cref{theo: MUON bound} of \cref{sec: priori bound} to establish in \cref{muon stochastic convergence} in \cref{subsec: without stays bound} below an error estimate for \MUON\ without assuming that the \MUON\ optimization process stays bounded. \cref{main theorem} above in this introduction section follows from an application of \cref{muon stochastic convergence}.

Our proof of \cref{muon bound convergence stochastic} is based on an application of an error analysis for a general class of \MUON-type optimization methods (which includes \MUON\ as a special case) that we establish in \cref{theo: stochastic convergence} in \cref{subsec: general analysis} below. In our proof of \cref{theo: stochastic convergence} we first derive in \cref{theo: stochastic convergence.eq7} a recursive real inequality for the weak optimization errors $(\E[|\cL(\Theta_n^{\xi,M})-\cL(\vartheta)|^p])^{1/p}$ for $n \in \N$ and, 
  thereafter, we employ the estimate from \cref{lem: stochastic recurrence 2} in this subsection below for solutions of the derived recursive inequality in \cref{theo: stochastic convergence.eq7} (cf.\ \cref{assume: stochastic recurrence 2} in \cref{lem: stochastic recurrence 2} below).

  Our proof of \cref{lem: stochastic recurrence 2} is based on the elementary estimates for the recursive inequalities in \cref{lem: stochastic recurrence'} and \cref{lem: stochastic recurrence''}. Only for completeness we include here detailed proofs for \cref{lem: stochastic recurrence'}, \cref{lem: stochastic recurrence''}, and \cref{lem: stochastic recurrence 2}.
  \subsection{Upper estimates for solutions of some recursive inequalities for the analysis of MUON}
  \begin{tcolorbox}[colback=white!95!gray,
                  colframe=black,
                  boxrule=0.5pt,
                  sharp corners,
                  enhanced,
                  breakable,
                 ]
\begin{athm}{lemma}{lem: stochastic recurrence'}
      Let $c\in \R$, $\tau\in (0,\infty)$, $N\in \N$ and let $\gamma\colon \N\to(0,\infty)$
    and $e\colon \N_0\to[0,\infty)$ satisfy for all $n\in \N\cap(N,\infty)$ that
      \begin{equation}\llabel{assume}
            e_n\leq  e_{n-1}-\tau\gamma_{n-1}e_{n-1}+c\gamma_{n-1}.
      \end{equation}
      Then 
      \begin{equation}\llabel{conclude}
         \textstyle \sup_{n\in \N}e_n\leq \max\{e_1,e_2,\dots,e_N,\tau^{-1}|c|+\sup_{n\in \N}(\gamma_n|c|)\}.
      \end{equation}
\end{athm}
\end{tcolorbox}
\begin{aproof}
Throughout this proof let  $\Gamma\in [0,\infty]$, $\fC\in \R$ satisfy 
    \begin{equation}\llabel{def: fC}
    \textstyle  \Gamma=\sup_{n\in \N}(\gamma_n|c|)\qqandqq \fC=\textstyle \max\{e_1,e_2,\dots,e_N\}.
    \end{equation} 
    In the following we prove that for all $n\in \N\cap[N,\infty)$ it holds that
    \begin{equation}\llabel{need to prove}
        e_n\leq \max\{\fC,|c|\tau^{-1}+\Gamma\}.
    \end{equation}
    We prove \lref{need to prove} by induction on $n\in \N\cap[N,\infty)$. For the base case $n=N$ note that \lref{def: fC} shows that
    \begin{equation}\llabel{ind: eq1}
        e_N\leq \fC\leq \max\{\fC,|c|\tau^{-1}+\Gamma\}.
    \end{equation}
    This establishes \lref{need to prove} in the base case $n=N$. For the induction step we assume that there exists $n\in \N\cap(N,\infty)$ which satisfies that
    \begin{equation}\llabel{ind: assume}
        e_{n-1}\leq \max\{\fC,|c|\tau^{-1}+\Gamma\}.
    \end{equation}
    In the following we prove that
    \begin{equation}\llabel{ind: need to prove}
        e_n\leq \max\{\fC,|c|\tau^{-1}+\Gamma\}.
    \end{equation}
    In our proof of \lref{ind: need to prove} we distinguish between the case $e_{n-1}\geq |c|\tau^{-1}$ and the case $e_{n-1}\leq |c|\tau^{-1}$. We first prove \lref{ind: need to prove} in the case 
    \begin{equation}\llabel{case 1}
        e_{n-1}\geq |c|\tau^{-1}.
    \end{equation}
    \startnewargseq
\argument{\lref{case 1};}{that \llabel{eqq1.1}
    $\tau\gamma_{n-1}e_{n-1}\geq |c|\gamma_{n-1}$\dott}
\argument{\lref{eqq1.1};\lref{assume};\lref{ind: assume}}{that
\begin{equation}
    e_{n}\leq e_{ n - 1 } - \tau \gamma_{ n - 1 } e_{ n - 1 } + c \gamma_{ n - 1 } \leq e_{ n - 1 } + |c| \gamma_{ n - 1 } - \tau \gamma_{ n - 1 } e_{ n - 1 } \leq e_{ n - 1 }\leq \max\{\fC,|c|\tau^{-1}+\Gamma\}.
\end{equation}}
This proves \lref{ind: need to prove} in the case $e_{n-1}\geq |c|\tau^{-1}$. We now prove \lref{ind: need to prove} in the case
\begin{equation}\llabel{case 2}
e_{n-1}\leq |c|\tau^{-1}.
\end{equation}
\startnewargseq
\argument{\lref{assume};\lref{case 2};the fact that for all $k\in \N$ it holds that $\gamma_k|c|\leq \Gamma$}{that
\begin{equation}
    e_n\leq e_{n-1}+|c|\gamma_{n-1}\leq e_{n-1}+\Gamma\leq |c|\tau^{-1}+\Gamma\leq \max\{\fC,|c|\tau^{-1}+\Gamma\}.
\end{equation}}
This establishes \lref{ind: need to prove} in the case $e_{n-1}\leq |c|\tau^{-1}$. Note that \lref{ind: eq1}, \lref{ind: need to prove}, and induction prove \lref{need to prove}.
\startnewargseq
\argument{\lref{def: fC};\lref{need to prove};}{\lref{conclude}\dott}
\end{aproof}
\begin{samepage}
 \begin{tcolorbox}[colback=white!95!gray,
                  colframe=black,
                  boxrule=0.5pt,
                  sharp corners,
                  enhanced,
                  breakable,
                 ]
\begin{athm}{lemma}{lem: stochastic recurrence''}
      Let $c\in [0,\infty)$, $\tau\in (0,\infty)$, $N\in \N$, let $(\gamma_n)_{n\in \N}\subseteq (0,\infty)$ be non-increasing, let $\fn \in \N \cap [N,\infty)$ satisfy $\sup_{ n \in \N \cap ( \fn, \infty) } ( 2 (\gamma_{ n })^{-1}-2(\gamma_{n-1})^{-1}) \leq \tau $,
   and let $(e_n)_{n\in \N}\subseteq[0,\infty)$ satisfy for all $n\in \N\cap(N,\infty)$ that
      \begin{equation}\llabel{assume}
            e_n\leq  e_{n-1}-\tau\gamma_{n}e_{n-1}+c(\gamma_{n})^2.
      \end{equation}
      Then 
      \begin{equation}\llabel{conclude}
          \textstyle\sup_{n\in \N}((\gamma_n)^{-1}e_n)\leq \max\{2c\tau^{-1}+c\gamma_1,\max_{n\in \{1,2,\dots,\fn\}}((\gamma_n)^{-1}e_n)\}.
      \end{equation}
\end{athm}
\end{tcolorbox}
\end{samepage}
\begin{aproof}
Throughout this proof for every $n\in \N$ let $x_n\in \R$ satisfy
\begin{equation}\llabel{def: e}
    x_n=(\gamma_n)^{-1}e_n
\end{equation}
and  let $\fC\in [0,\infty)$ satisfy
    \begin{equation}\llabel{def: fC}
       \fC=\textstyle \sup_{n\in \{1,2,\dots,\fn\}}((\gamma_n)^{-1}e_n).
    \end{equation}
     \argument{\lref{assume};\lref{def: e}}{that for all $n\in \N\cap(N,\infty)$ it holds that
    \begin{equation}\llabel{eq8}
    \begin{split}
        x_n&\leq \frac{\gamma_{n-1} x_{n-1}}{\gamma_{n}}-\tau\gamma_{n-1}x_{n-1}+c\gamma_n\\
        &= x_{n-1}+\frac{(\gamma_{n-1}-\gamma_n) x_{n-1}}{\gamma_{n}}-\tau\gamma_{n-1}x_{n-1}+c\gamma_n.
        \end{split}
    \end{equation}}
    \argument{\lref{eq8};the fact that for all $n\in \N\cap(\fn,\infty)$ it holds that
        $\gamma_{n-1}-\gamma_n\leq \frac{\tau\gamma_n\gamma_{n-1}}{2}$;the fact that for all $n\in \N$ it holds that $\gamma_n\leq \gamma_{n-1}$}{that for all $n\in \N\cap(\fn,\infty)$ it holds that
    \begin{equation}\llabel{eq10}
    \begin{split}
        x_n&\leq
     x_{n-1}+\frac{\tau\gamma_{n-1} x_{n-1}}{2}-\tau\gamma_{n-1}x_{n-1}+c\gamma_{n-1}= x_{n-1}-\frac{\tau\gamma_{n-1} x_{n-1}}{2}+c\gamma_{n-1}.
     \end{split}
    \end{equation}}
    \argument{\lref{def: fC};\lref{def: e}}{that
    \begin{equation}\llabel{eq12}
    \textstyle    \fC\geq \sup_{n\in \{1,2,\dots,\fn\}}x_n.
    \end{equation}}
  \argument{\lref{eq12};\lref{def: e};\lref{eq10};the fact that $\gamma_1=\sup_{n\in \N} \gamma_n$;\cref{lem: stochastic recurrence'}}{that 
  \begin{equation}\llabel{eq11}
      \textstyle\sup_{n\in \N}
      ((\gamma_{n})^{-1}e_{n})= \sup_{n\in \N }x_n\leq  \max\{\fC,2c\tau^{-1}+c\gamma_1\}.
  \end{equation}}
\end{aproof}
\begin{samepage}
 \begin{tcolorbox}[colback=white!95!gray,
                  colframe=black,
                  boxrule=0.5pt,
                  sharp corners,
                  enhanced,
                  breakable,
                 ]
\begin{athm}{prop}{lem: stochastic recurrence 2}
      Let $c,\chi,\rho\in (0,\infty)$, let $(\gamma_n)_{n\in \N}\subseteq (0,\infty)$ be non-increasing, and assume $\limsup_{n\to\infty} (\gamma_n+(\gamma_n)^{-2}(\gamma_{n}-\gamma_{n+1}))=0$.
      Then there exists $\fC\in \R$ such that for every $\tau\in [0,\infty)$ and every $e\colon \N_0\to[0,\infty)$ with the property that for all $n\in \N$ it holds that $e_0\leq \chi$ and
      \begin{equation}\label{assume: stochastic recurrence 2}
            e_n\leq e_{n-1}(1-\rho\gamma_n)+c\tau\gamma_n+c(\gamma_n)^2
      \end{equation}
      we have  for all $n\in \N_0$  that $e_n\leq \fC\tau+\fC\gamma_{n+1}$.
\end{athm}
\end{tcolorbox}
\end{samepage}
\begin{aproof}
Throughout this proof
      let $N\in \N$ satisfy
\begin{equation}\llabel{eq0}
    \sup_{n\in n\cap[N,\infty)}\rho\gamma_n\leq 1
\end{equation}
 and for every $n\in \N$ let $S_n\in \R$ satisfy
 \begin{equation}\llabel{def: e}
     \textstyle S_n=\sum\limits_{k=N}^n\biggl[
(\gamma_k)^2
\biggl(\prod\limits_{j=k+1}^n (1-\rho\gamma_j)\biggr)\biggr].
 \end{equation}
\argument{the discrete Gronwall inequality;}{that for every $\tau\in [0,\infty)$ and every $e\colon \N_0\to[0,\infty)$ with the property that for all $n\in \N$ it holds that 
            $e_n\leq e_{n-1}(1-\rho\gamma_n)+c\tau\gamma_n+c(\gamma_n)^2$
      we have  for all $n\in \N\cap[N,\infty)$ that
\begin{equation}\llabel{eq1}
\begin{split}
\textstyle
   e_n &\leq\textstyle  e_{N-1} \biggl[
\prod\limits_{j=N}^n (1-\rho\gamma_j)\biggr]+
\biggl[\sum\limits_{k=N}^n
(c\tau\gamma_k+c(\gamma_k)^2)
\biggl(\prod\limits_{j=k+1}^n (1-\rho\gamma_j)\biggr)\biggr]\\
&\textstyle =e_{N-1} \biggl[
\prod\limits_{j=N}^n (1-\rho\gamma_j)\biggr]+
c\tau\biggl[\sum\limits_{k=N}^n\gamma_k
\biggl(\prod\limits_{j=k+1}^n (1-\rho\gamma_j)\biggr)\biggr]+c\biggl[\sum\limits_{k=N}^n
(\gamma_k)^2
\biggl(\prod\limits_{j=k+1}^n (1-\rho\gamma_j)\biggr)\biggr].
\end{split}
\end{equation}}
 \argument{the fact that $\limsup_{n\to\infty} \bigl(\frac{\gamma_{n}-\gamma_{n+1}}{(\gamma_{n})^2}\bigr)\allowbreak<\infty$;the fact that $\limsup_{n\to\infty}\gamma_n=0$;}{that
    \begin{equation}\llabel{eq2}
        \textstyle\limsup_{n\to\infty} \bigl(\frac{\gamma_{n}-\gamma_{n+1}}{\gamma_{n}}\bigr)=0.
    \end{equation}}
    \argument{\lref{eq2};}{that
    \begin{equation}\llabel{eq3}
        \textstyle\lim_{n\to\infty} \bigl(\frac{\gamma_{n}}{\gamma_{n+1}}\bigr)=1.
    \end{equation}}
\argument{\lref{eq3};the fact that $\limsup_{n\to\infty} \bigl(\frac{\gamma_{n}-\gamma_{n+1}}{(\gamma_{n})^2}\bigr)=0$;}{that
    \begin{equation}\llabel{eq4}
        \textstyle\limsup_{n\to\infty} \bigl(\frac{\gamma_{n}-\gamma_{n+1}}{\gamma_n\gamma_{n+1}}\bigr)=0.
    \end{equation}}
    \argument{\lref{eq4};}{that there exists $M\in \N\cap(N,\infty)$ such that for all $n\in \N\cap[M,\infty)$ it holds that
    \begin{equation}\llabel{eq5}
        \gamma_{n-1}-\gamma_n\leq \frac{\rho}{2}\gamma_n\gamma_{n-1}.
    \end{equation}}
    \argument{\lref{eq5};}{that there exists $M\in \N\cap(N,\infty)$ such that for all $n\in \N\cap[M,\infty)$ it holds that
    \begin{equation}\llabel{eq6}
        1-\rho\gamma_n\leq 1-\frac{\rho\gamma_n}{2}\leq \frac{\gamma_n}{\gamma_{n-1}}.
    \end{equation}}
    \argument{\lref{eq6};\lref{eq0}}{that there exist $\fC\in (0,\infty)$, $M\in \N\cap(N,\infty)$ such that for all $n\in \N\cap[M,\infty)$ it holds that
    \begin{equation}\llabel{eq7}
       \begin{split}
           \textstyle\prod\limits_{j=N}^n (1-\rho\gamma_j)&\textstyle=\biggl[ \prod\limits_{j=N}^{M-1} (1-\rho\gamma_j)\biggr]\biggl[ \prod\limits_{j=M}^n (1-\rho\gamma_j)\biggr]\leq \biggl[ \prod\limits_{j=N}^{M-1} (1-\rho\gamma_j)\biggr]\biggl[ \prod\limits_{j=M}^n \frac{\gamma_j}{\gamma_{j-1}}\biggr]\\
          &= \textstyle\biggl[ \prod\limits_{j=N}^{M-1} (1-\rho\gamma_j)\biggr](\gamma_{M-1})^{-1}\gamma_n=\fC\gamma_n.
       \end{split}
    \end{equation}}
    \argument{\lref{eq7};}{that there exists $\fC\in (0,\infty)$ such that for all $n\in \N\cap[N,\infty)$ it holds that
    \begin{equation}\llabel{evidence 1}
        \textstyle\prod\limits_{j=N}^n (1-\rho\gamma_j)\leq \fC\gamma_n.
    \end{equation}}
      \argument{the fact that for all $n\in \N$, $k\in \{1,2,\dots,n\}$ it holds that
    \begin{equation}
    \begin{split}
       \textstyle \gamma_k\bigg[\prod\limits_{j=k+1}^n (1-\rho\gamma_j)\biggr]&\textstyle=\frac {1}{\rho}\bigg[\prod\limits_{j=k+1}^n (1-\rho\gamma_j)\biggr]-\frac {1}{\rho}(1-\rho\gamma_k)\bigg[\prod\limits_{j=k+1}^n (1-\rho\gamma_j)\biggr]\\
       &\textstyle=\frac {1}{\rho}\biggl[\prod\limits_{j=k+1}^n (1-\rho\gamma_j)\biggr]-\frac{1}{\rho}\biggl[\prod\limits_{j=k}^n (1-\rho\gamma_j)\biggr]
       \end{split}
    \end{equation}}{that there exists $\fC\in (0,\infty)$ such that for all $n\in \N\cap[N,\infty)$ it holds that
    \begin{equation}\llabel{evidence 2}
        \textstyle\sum\limits_{k=N}^n\biggl[\gamma_k
\biggl(\prod\limits_{j=k+1}^n (1-\rho\gamma_j)\biggr)\biggr]=\frac {1}{\rho}-\frac{1}{\rho}\biggl[\prod\limits_{j=N}^n (1-\rho\gamma_j)\biggr]\leq \frac {1}{\rho}.
    \end{equation}}
    \argument{\lref{def: e};}{that for all $n\in \N\cap[N+1,\infty)$ it holds that
    \begin{equation}\llabel{eq8}
        \begin{split}
        S_n&=\textstyle    \sum\limits_{k=N}^n\biggl[
(\gamma_k)^2
\biggl(\prod\limits_{j=k+1}^n (1-\rho\gamma_j)\biggr)\biggr]\textstyle=\biggl[\sum\limits_{k=N}^{n-1}
(\gamma_k)^2
\biggl(\prod\limits_{j=k+1}^{n} (1-\rho\gamma_j)\biggr)\biggr]+(\gamma_n)^2\\
&=\textstyle(1-\rho\gamma_n)\biggl[\sum\limits_{k=N}^{n-1}
(\gamma_k)^2
\biggl(\prod\limits_{j=k+1}^{n-1} (1-\rho\gamma_j)\biggr)\biggr]+(\gamma_n)^2=(1-\rho\gamma_n)S_{n-1}+(\gamma_n)^2\\
&=S_{n-1}-\rho\gamma_nS_{n-1}+(\gamma_n)^2.
        \end{split}
    \end{equation}}
    \argument{\lref{def: e};\lref{eq5};\lref{eq8};\cref{lem: stochastic recurrence''}}{that there exists $\fC\in (0,\infty)$ such that for all $n\in \N$ it holds that
    \begin{equation}\llabel{evidence 3}
      \textstyle    \sum\limits_{k=N}^n\biggl[
(\gamma_k)^2
\biggl(\prod\limits_{j=k+1}^n (1-\rho\gamma_j)\biggr)\biggr] = S_n\leq \fC\gamma_n.
    \end{equation}}
    \argument{\lref{evidence 3};\lref{eq1};\lref{evidence 1};\lref{evidence 2}}{that there exists $\fC\in (0,\infty)$ such that for every $\tau\in [0,\infty)$ and every $e\colon \N_0\to[0,\infty)$ with the property that for all $n\in \N$ it holds that 
            $e_n\leq e_{n-1}(1-\rho\gamma_n)+c\tau\gamma_n+c(\gamma_n)^2$
      we have for all $n\in \N\cap[N+1,\infty)$ that
    \begin{equation}\llabel{eq10}
        e_n \leq \fC\tau+(\fC+\fC e_{N-1})\gamma_n.
    \end{equation}}
    \argument{the fact that for all $n\in \N$ it holds that $\gamma_n \leq \gamma_1$}{that for every $\tau\in [0,\infty)$ and every $e\colon \N_0\to[0,\infty)$ with the property that for all $n\in \N$ it holds that $e_0\leq \chi$ and 
            $e_n\leq e_{n-1}(1-\rho\gamma_n)+c\tau\gamma_n+c(\gamma_n)^2$
      we have for all $n\in \{1,2,\dots,N\}$ that
    \begin{equation}\llabel{eq11}
        e_n\leq e_{n-1}+c\tau\gamma_1+c(\gamma_1)^2.
    \end{equation}}
    \argument{\lref{eq11};induction}{that there exist $\fC_1,\fC_2\in (0,\infty)$ such that for every $\tau\in [0,\infty)$ and every $e\colon \N_0\to[0,\infty)$ with the property that for all $n\in \N$ it holds that $e_0\leq \chi$ and 
            $e_n\leq e_{n-1}(1-\rho\gamma_n)+c\tau\gamma_n+c(\gamma_n)^2$ we have for all $n\in \{1,2,\dots,N\}$ that
    \begin{equation}\llabel{eq12}
    \begin{split}
       e_n&\leq e_0+nc\tau\gamma_1+nc(\gamma_1)^2\leq \chi+Nc\tau\gamma_1+Nc(\gamma_1)^2\\
       &=\fC_1+\fC_1\tau\leq\fC_2\tau+\fC_2\gamma_N \leq \fC_2\tau+\fC_2\gamma_n.
      \end{split}
    \end{equation}}
    \argument{\lref{eq12};\lref{eq10}}{that there exists $\fC\in (0,\infty)$ such that for every $\tau\in [0,\infty)$ and every $e\colon \N_0\to[0,\infty)$ with the property that for all $n\in \N$ it holds that $e_0\leq \chi$ and 
            $e_n\leq e_{n-1}(1-\rho\gamma_n)+c\tau\gamma_n+c(\gamma_n)^2$ we have for all $n\in \N$ that
            \begin{equation}\llabel{arg2}
                e_n\leq \fC\tau+\fC\gamma_n.
            \end{equation}}
               \argument{the fact that $\limsup_{n\to\infty}\gamma_n=0$; the fact that $\limsup_{n\to\infty} \bigl(\frac{\gamma_{n}-\gamma_{n+1}}{(\gamma_{n})^2}\bigr)<\infty$;}{that
    \begin{equation}\llabel{eq13}
        \textstyle\limsup_{n\to\infty} \bigl(\frac{\gamma_{n}-\gamma_{n+1}}{\gamma_{n}}\bigr)=0.
    \end{equation}}
    \argument{\lref{eq13};}{that
    \begin{equation}\llabel{eq14}
        \textstyle\limsup_{n\to\infty} \bigl(\frac{\gamma_{n}}{\gamma_{n+1}}\bigr)=1.
    \end{equation}}
    \argument{\lref{eq14};}{that there exists $\fC\in (0,\infty)$ such that for all $n\in \N$ it holds that
    \begin{equation}\llabel{eq15}
        \gamma_n\leq \fC\gamma_{n+1}.
    \end{equation}}
    \argument{\lref{eq15};\lref{arg2};}{that there exist $\fC_1,\fC_2,\fC_3\in (0,\infty)$ such that for every $\tau\in [0,\infty)$ and every $e\colon \N_0\to[0,\infty)$ with the property that for all $n\in \N$ it holds that $e_0\leq \chi$ and 
            $e_n\leq e_{n-1}(1-\rho\gamma_n)+c\tau\gamma_n+c(\gamma_n)^2$ we have for all $n\in \N$ that
            \begin{equation}\llabel{arg3}
                e_n\leq \fC_1\tau+\fC_1\gamma_n\leq \fC_1\tau+\fC_1\fC_2\gamma_{n+1}\leq \fC_3\tau+\fC_3\gamma_{n+1}.
            \end{equation}}
            \argument{\lref{arg3};}{that there exists $\fC\in (0,\infty)$ such that for every $\tau\in [0,\infty)$ and every $e\colon \N_0\to[0,\infty)$ with the property that for all $n\in \N$ it holds that $e_0\leq \chi$ and 
            $e_n\leq e_{n-1}(1-\rho\gamma_n)+c\tau\gamma_n+c(\gamma_n)^2$ we have for all $n\in \N_0$ that
            \begin{equation}\llabel{arg4}
                e_n\leq \fC\tau+\fC\gamma_{n+1}.
            \end{equation}}
\end{aproof}
 \subsection{Error estimates for MUON-type optimizers under boundedness assumptions}\label{subsec: general analysis}
 \begin{samepage}
 \begin{tcolorbox}[colback=white!95!gray,
                  colframe=black,
                  boxrule=0.5pt,
                  sharp corners,
                  enhanced,
                  breakable,
                 ]
\begin{athm}{lemma}{lem: stochastic bfm analysis pre}
    Let $(\Omega,\cF,\P)$ be a probability space, let $(\setX,\cU)$ be a measurable space, let $(\cA, \langle \cdot,\cdot \rangle)$ be a Hilbert space, let $M\in \N$, $\alpha\in (0,1)$, $\fc,\cK\in (0,\infty)$, $p\in[1,\infty)$, let $X_{n,m}\colon \Omega\to\setX$, $(n,m)\in \N^2$, be \iid\ random variables, let $\bbG\colon \cA\times \setX\to\cA$ be measurable, assume for all $\theta,\vartheta\in \cA$, $x\in \setX$ with $\max\{\|\theta\|,\|\vartheta\|\}<\fc$ that
    \begin{equation}\llabel{eq1}
       \|\bbG(\theta,x)\|\leq \cK\qqandqq\|\bbG(\theta,x)-\bbG(\vartheta,x)\|\leq \cK\|\theta-\vartheta\|,
    \end{equation}
    let $\gamma\colon \N\to(0,\infty)$ satisfy for all $n\in\N$ that
        \begin{equation}\llabel{eqq3}
          \textstyle  \sum_{i=1}^n(\alpha^{n-i}\gamma_i)\leq \cK\gamma_n,
        \end{equation}
         let $\Theta\colon \N_0\times\Omega\to\cA$ and $\bfm\colon \N_0\times\Omega\to\cA$ be stochastic processes, assume for all $n\in\N$ that
    \begin{equation}\llabel{def: Theta}
        \textstyle\P(\|\Theta_{n-1}\|< \fc)=\P(\|\Theta_n-\Theta_{n-1}\|\leq \cK\gamma_{n})=1,
     \end{equation}
    assume for all $n\in \N$ that
         \begin{equation}\llabel{def: bfm}
            \bfm_0=0\qqandqq \bfm_n= \alpha \bfm_{n-1}+(1-\alpha)\bigl[\textstyle \frac 1M \sum_{m=1}^M\bbG(\Theta_{n-1},X_{n,m})\bigr],
         \end{equation}
         assume for all $n\in \N$ that $\Theta_{n-1}$ and $(X_{n,m})_{m\in \N}$ are independent,
         let $\cG\colon\cA\to\cA$ be measurable, and let $\fC\in [0,\infty)$ satisfy for all $\theta\in \cA$ with $\|\theta\|<\fc$ that $\cG(\theta)=\E[\bbG(\theta,X_{1,1})]$ and 
         \begin{equation}\llabel{def: fC}
             \E\bigl[\|\E[\bbG(\theta,X_{1,1})]-\textstyle \frac 1M \sum_{m=1}^M\bbG(\theta,X_{1,m})\|^p\bigr]\leq \fC^p M^{-p/2}.
         \end{equation}
      Then it holds for all $n\in \N$ that
     \begin{equation}\llabel{conclude}
       \textstyle\bigl(\E\bigl[\|\bfm_{n}-\cG(\Theta_{n-1})\|^p\bigr]\bigr)^{1/p}\leq \fC M^{-1/2}+(\cK^2(\gamma_1)^{-1}+\cK^3+\cK^2)\gamma_n.
     \end{equation}
    \end{athm}
    \end{tcolorbox}
    \end{samepage}
    \begin{aproof}
    Throughout this proof let $\bbM\colon\N_0\times\Omega\to\cA$ satisfy for all $n\in \N$ that
    \begin{equation}\llabel{def: bbM}
         \bbM_0=0\qqandqq \bbM_n= \alpha \bbM_{n-1}+(1-\alpha)\cG(\Theta_{n-1}).
    \end{equation}
    \argument{\lref{def: fC};the fact that for all $\theta\in \cA$ it holds that $\bbG(\theta,X_{n,m})$, $(n,m)\in \N^2$, are \iid;}{that for all $n\in\N$, $\theta\in \cA$ with $\|\theta\|<\fc$ it holds that
    \begin{equation}\llabel{eq2}
    \begin{split}
        \E\Bigl[\bigl\|\bigl(\textstyle \frac 1M \sum_{m=1}^M\bbG(\theta,X_{n,m})\bigr)-\E[\bbG(\theta,X_{1,1})]\bigr\|^p\Bigr]&\textstyle=\E\Bigl[\bigl\|\bigl(\textstyle \frac 1M \sum_{m=1}^M\bbG(\theta,X_{1,m})\bigr)-\E[\bbG(\theta,X_{1,1})]\bigr\|^p\Bigr]\\
        &\leq \fC^p M^{-p/2}.
        \end{split}
    \end{equation}}
    \argument{\lref{eq2};the assumption that assume for all $n\in \N$ that $\Theta_{n-1}$ and $(X_{n,m})_{m\in \N}$ are independent;the assumption that for all $\theta\in \cA$ with $\|\theta\|< \fc$ it holds that $\cG(\theta)=\E[\bbG(\theta,X_{1,1})]$; the fact that for all $n\in \N$ it holds that $\P(\|\Theta_{n-1}\|< \fc)=1$;the factorization lemma, \eg, \cite[Proposition 3.15]{DeRoAr2024nonconvergence}}{that for all $n\in \N$ it holds that
    \begin{equation}\llabel{eq3}
    \begin{split}
       & \textstyle\E\Bigl[\bigl\|\bigl(\textstyle \frac 1M \sum_{m=1}^M\bbG(\Theta_{n-1},X_{n,m})\bigr)-\cG(\Theta_{n-1})\bigr\|^p\Bigr]\\
       &=  \textstyle\E\Bigl[\E\Bigl[\bigl\|\bigl(\textstyle \frac 1M \sum_{m=1}^M\bbG(\Theta_{n-1},X_{n,m})\bigr)-\cG(\Theta_{n-1})\bigr\|^p|\sigma(\Theta_{n-1})\Bigr]\Bigr]\\
       &\leq \textstyle\sup_{\theta\in \cA,\,\|\theta\|< \fc}\E\Bigl[\bigl\|\bigl(\textstyle \frac 1M \sum_{m=1}^M\bbG(\theta,X_{n,m})\bigr)-\cG(\theta)\bigr\|^p\Bigr]\\
       &=\textstyle\sup_{\theta\in\cA,\,\|\theta\|< \fc}\E\Bigl[\bigl\|\bigl(\textstyle \frac 1M \sum_{m=1}^M\bbG(\theta,X_{n,m})\bigr)-\E[\bbG(\theta,X_{1,1})]\bigr\|^p\Bigr]\\
       &\leq \fC^p M^{-p/2}.
       \end{split}
    \end{equation}}
    \argument{\lref{def: bfm};\lref{def: bbM};\cref{momentum:representation}}{that for all $n\in \N_0$ it holds that
    \begin{equation}\llabel{eq4}
        \bfm_{n}=\textstyle(1-\alpha)\biggl[\sum\limits_{k=1}^n\alpha^{n-k}\textstyle \frac 1M \bigl[\textstyle\sum_{m=1}^M\bbG(\Theta_{k-1},X_{k,m})\bigr]\biggr]
        \end{equation}
        \begin{equation}\llabel{eq5}
            \text{and}\qquad \textstyle \bbM_{n}=(1-\alpha)\bigl[\textstyle\sum_{k=1}^n\alpha^{n-k}\cG(\Theta_{k-1})\bigr].
    \end{equation}}
    \argument{\lref{eq5};\lref{eq3};the triangle inequality}{that for all $n\in \N_0$ it holds that
    \begin{equation}\llabel{evidence 1}
    \begin{split}
       &\bigl( \E\bigl[\|\bfm_n-\bbM_n\|^p\bigr]\bigr)^{1/p}\\&=\textstyle\biggl(\E\biggl[\Bigl\|(1-\alpha)\sum\limits_{k=1}^n\alpha^{n-k}\textstyle\bigl( \frac 1M \bigl[\textstyle\sum_{m=1}^M\bbG(\Theta_{k-1},X_{k,m})\bigr]-\cG(\Theta_{k-1})\bigr)\Bigr\|^p\biggr]\biggr)^{1/p}\\
       &\leq (1-\alpha)\biggl[\textstyle\sum\limits_{k=1}^n\alpha^{n-k}\Bigl(\E\Bigl[\bigl\|\textstyle\bigl( \frac 1M \bigl[\textstyle\sum_{m=1}^M\bbG(\Theta_{k-1},X_{k,m})\bigr]-\cG(\Theta_{k-1})\bigr)\bigr\|^p\Bigr]\Bigr)^{1/p}\biggr]\\
       &\leq  (1-\alpha)\biggl[\textstyle\sum\limits_{k=1}^n\alpha^{n-k}\fC M^{-1/2}\biggr]\leq \fC M^{-1/2}.
       \end{split}
    \end{equation}}
    \argument{\lref{eq1};\lref{def: Theta};\lref{def: fC};}{that for all $\theta_1,\theta_2\in \cA$ with $\max\{\|\theta_1\|,\|\theta_2\|\}< \fc$ it holds that
    \begin{equation}\llabel{eq6}
    \begin{split}
      \| \cG(\theta_1)-\cG(\theta_2)\| &=\bigl\|\E[\bbG(\theta_1,X_{1,1})]-\E[\bbG(\theta_2,X_{1,1})]\bigr\|\\
      &\leq \E\bigl[\|\bbG(\theta_1,X_{1,1})-\bbG(\theta_2,X_{1,1})\|\bigr]\leq \cK\|\theta_1-\theta_2\|.
      \end{split}
    \end{equation}}
    \argument{\lref{eq6};\lref{def: Theta};}{that for all $n\in \N$ it holds that
    \begin{equation}\llabel{eq7}
        \P(\|\cG(\Theta_n)-\cG(\Theta_{n-1})\|\leq \cK\|\Theta_n-\Theta_{n-1}\|)=1.
    \end{equation}}
        \argument{\lref{eqq3};\lref{def: Theta};\lref{def: bbM};\lref{eq7};\cref{lem: bfm analysis 1}}{that for all $n\in \N$ it holds that
        \begin{equation}\llabel{evidence 2.0}
            \P\bigl(\|\bbM_n-\cG(\Theta_{n-1})\|\leq (\cK(\gamma_1)^{-1}\|\bbM_0-\cG(\Theta_0)\|+\cK^3+\cK^2)\gamma_n\bigr)=1.
        \end{equation}}
        \argument{\lref{eq1};\lref{def: Theta};\lref{def: fC}}{that for all $\theta\in \cA$ with $\|\theta\|< \fc$ it holds that
        \begin{equation}\llabel{evidence 2.1}
            \|\cG(\theta)\|=\bigl\|\E[\bbG(\theta,X_{1,1})]\bigr\|\leq \E\bigl[\|\bbG(\theta,X_{1,1})\|\bigr] \leq \cK.
        \end{equation}}
        \argument{\lref{evidence 2.1};\lref{def: Theta}}{that
        \begin{equation}\llabel{evidence 2.2}
            \P(\|\cG(\Theta_0)\|\leq \cK)=1.
        \end{equation}}
        \argument{\lref{evidence 2.2};\lref{def: bbM};\lref{evidence 2.0};}{that for all $n\in \N$ it holds that
        \begin{equation}\llabel{evidence 2}
        \begin{split}
         &\P\bigl(\|\bbM_n-\cG(\Theta_{n-1})\|\leq (\cK^2(\gamma_1)^{-1}+\cK^3+\cK^2)\gamma_n\bigr)\\
        &\geq \P\bigl(\|\bbM_n-\cG(\Theta_{n-1})\|\leq (\cK(\gamma_1)^{-1}\|\cG(\Theta_0)\|+\cK^3+\cK^2)\gamma_n\bigr)\\
             &=\P\bigl(\|\bbM_n-\cG(\Theta_{n-1})\|\leq (\cK(\gamma_1)^{-1}\|\bbM_0-\cG(\Theta_0)\|+\cK^3+\cK^2)\gamma_n\bigr)\\
             &=1.
             \end{split}
        \end{equation}}
        \argument{\lref{def: bbM};\lref{evidence 1};\lref{evidence 2};the triangle inequality}{that  for all $n\in \N$ it holds that
        \begin{equation}\llabel{eq8}
        \begin{split}
            \textstyle\bigl(\E\bigl[\|\bfm_{n}-\cG(\Theta_{n-1})\|^p\bigr]\bigr)^{1/p}&\leq \bigl(\E\bigl[\|\bfm_{n}-\bbM_n\|^p\bigr]\bigr)^{1/p}+\bigl(\E\bigl[\|\bbM_n-\cG(\Theta_{n-1})\|^p\bigr]\bigr)^{1/p}\\
            &\leq \fC  M^{-1/2}+(\cK^2(\gamma_1)^{-1}+\cK^3+\cK^2)\gamma_n.
            \end{split}
        \end{equation}}
    \end{aproof}
    \begin{samepage}
    \begin{tcolorbox}[colback=white!95!gray,
                  colframe=black,
                  boxrule=0.5pt,
                  sharp corners,
                  enhanced,
                  breakable,
                 ]
\begin{athm}{cor}{lem: stochastic bfm analysis}
    Let $(\Omega,\cF,\P)$ be a probability space, let $(\setX,\cU)$ be a measurable space, let $(\cA, \langle \cdot,\cdot \rangle)$ be a finite-dimensional Hilbert space, let $\alpha\in (0,1)$, $p,\fc,\cK\in (0,\infty)$, let $X_{n,m}\colon \Omega\to\setX$, $(n,m)\in \N^2$, be \iid\ random variables, let $\smalll=(\smalll(\theta,x))_{(\theta,x)\in \cA\times \setX}\colon \cA\times \setX\to\R$ be measurable, assume for all $x\in \setX$ that $(\cA\ni\theta\mapsto \smalll(\theta,x)\in \R)\in C^1(\cA,\R)$, assume for all  $\theta,\vartheta\in \cA$, $x\in \setX$ with $\max\{\|\theta\|,\|\vartheta\|\}< \fc$ that 
    \begin{equation}\llabel{assume}
        \|(\nabla_\theta\smalll)(\theta,x)-(\nabla_\theta\smalll)(\vartheta,x)\|\leq \cK\|\theta-\vartheta\|\qqandqq\|(\nabla_\theta\smalll)(\theta,x)\|\leq  \cK,
    \end{equation}  let $\gamma\colon \N\to(0,\infty)$ satisfy $\limsup_{n\to\infty}((\gamma_{n+1})^{-1}\alpha\gamma_n)<1$, 
         and let $\cL\colon\cA\to\R$ satisfy for all $\theta\in \cA$ with $\|\theta\|< \fc$ that $\cL(\theta)=\E[\smalll(\theta,X_{1,1})]$.
      Then there exists $\fC\in \R$ such that for every $M\in \N$, every stochastic process $\bfm\colon \N_0\times\Omega\to\cA$, and every stochastic process $\Theta\colon \N_0\times\Omega\to\cA$ with the property that for all $n\in\N$ it holds that
    \begin{equation}\llabel{def: Theta}
         \textstyle\P(\|\Theta_{n-1}\|< \fc)=\P(\|\Theta_n-\Theta_{n-1}\|\leq \cK\gamma_{n})=1,
     \end{equation}
         \begin{equation}\llabel{def: bfm}
            \bfm_0=0,\qqandqq \bfm_n= \alpha \bfm_{n-1}+(1-\alpha)\bigl[\textstyle \frac 1M \sum_{m=1}^M(\nabla_\theta\smalll)(\Theta_{n-1},X_{n,m})\bigr]
         \end{equation}
         and with the property that for all $n\in \N$ it holds that $\Theta_{n-1}$ and $(X_{n,m})_{m\in \N}$ are independent we have for all $n\in \N$ that
     \begin{equation}\llabel{conclude}
       \textstyle\bigl(\E\bigl[\|\bfm_{n}-(\nabla\cL)(\Theta_{n-1})\|^p\bigr]\bigr)^{1/p}\leq \fC(M^{-1/2}+\gamma_n).
     \end{equation}
    \end{athm}
    \end{tcolorbox}
    \end{samepage}
       \begin{aproof}
       Throughout this proof let $q\in \R$ satisfy $q=\max\{p,1\}$.
       \argument{\lref{assume};}{that for all $n,m\in \N$, $\theta\in \cA$ with $\|\theta\|< \fc$ it holds that
    \begin{equation}\llabel{eq1}
        \P(\|(\nabla_\theta\smalll)(\theta,X_{n,m})\|\leq \cK)=1.
    \end{equation}}
    \argument{\lref{eq1};the fact that for all $n\in \N$, $\theta\in \cA$ it holds that $(\nabla_\theta\smalll)(\theta,X_{n,m})$, $m\in \N$, are \iid;the fact that $q\geq 1$;the Marcinkiewicz–Zygmund inequality, \eg, \cite[Theorem 2 in Subsection 10.3]{MR1476912}}{that there exist $\fC_1,\fC_2\in (0,\infty)$ such that for all $M\in\N$, $\theta\in \cA$ with $\|\theta\|<\fc$ it holds that
    \begin{equation}\llabel{eq2}
    \begin{split}
        &\E\Bigl[\bigl\|\bigl(\textstyle \frac 1M \sum_{m=1}^M(\nabla_\theta\smalll)(\theta,X_{1,m})\bigr)-\E[(\nabla_\theta\smalll)(\theta,X_{1,1})]\bigr\|^q\Bigr]\textstyle \\
&\textstyle\leq\fC_1\E\Bigl[\bigl(\sum_{m=1}^M\bigl\|\frac{(\nabla_\theta\smalll)(\theta,X_{1,m})}{M}-\frac{\E[(\nabla_\theta\smalll)(\theta,X_{1,m})]}{M}\bigr\|^2\bigr)^{q/2}\Bigr]\\
        &\textstyle\leq \fC_1\E\Bigl[\bigl(\sum_{m=1}^M\bigl(2\bigl\|\frac{(\nabla_\theta\smalll)(\theta,X_{1,m})}{M}\bigr\|^2+2\bigl\|\frac{\E[(\nabla_\theta\smalll)(\theta,X_{1,m})]}{M}\bigr\|^2\bigr)\bigr)^{q/2}\Bigr]\\
        &\textstyle\leq \fC_1\E\Bigl[\bigl(\sum_{m=1}^M\bigl(2\bigl|\frac{\cK}{M}\bigr|^2+2\bigl|\frac{\cK}{M}\bigr|^2\bigr)\bigr)^{q/2}\Bigr]\leq\fC_2M^{-q/2}.
        \end{split}
    \end{equation}}
     \argument{the assumption that $\limsup_{n\to\infty}((\gamma_{n+1})^{-1}\alpha\gamma_n)<1$;\cite[Lemma 2.6]{ShoAr2026};}{that there exists $\fC\in (0,\infty)$ such that for all $n\in\N$ it holds that
        \begin{equation}\llabel{eqq3}
            \sum_{i=1}^n(\alpha^{n-i}\gamma_i)\leq \fC\gamma_n.
        \end{equation}}
        \argument{\lref{eqq3};\lref{eq2};the fact that $q\geq 1$;\cref{lem: stochastic bfm analysis pre}}{that there exists $\fC\in \R$ such that for every $M\in \N$, every stochastic process $\bfm\colon \N_0\times\Omega\to\cA$, and every stochastic process $\Theta\colon \N_0\times\Omega\to\cA$ with the property that
    \begin{equation}
       \textstyle\P(\|\Theta_{n-1}\|< \fc)=\P(\|\Theta_n-\Theta_{n-1}\|\leq \cK\gamma_{n})=1,
     \end{equation}
    with the property that for all $n\in \N$ it holds that
         \begin{equation}
            \bfm_0=0\qqandqq \bfm_n= \alpha \bfm_{n-1}+(1-\alpha)\bigl[\textstyle \frac 1M \sum_{m=1}^M(\nabla_\theta\smalll)(\Theta_{n-1},X_{n,m})\bigr],
         \end{equation}
         and with the property that for all $n\in \N$ it holds that $\Theta_{n-1}$ and $(X_{n,m})_{m\in \N}$ are independent we have for all $n\in \N$ that
     \begin{equation}\llabel{preconclude}
       \textstyle\bigl(\E\bigl[\|\bfm_{n}-(\nabla\cL)(\Theta_{n-1})\|^q\bigr]\bigr)^{1/q}\leq \fC(M^{-1/2}+\gamma_n).
     \end{equation}}
       \argument{\lref{preconclude}}{\lref{conclude}\dott}
    \end{aproof}
    \begin{samepage}
    \begin{tcolorbox}[colback=white!95!gray,
                  colframe=black,
                  boxrule=0.5pt,
                  sharp corners,
                  enhanced,
                  breakable,
                 ]
\begin{athm}{theorem}{theo: stochastic convergence}[\textcolor{red}{Error analysis for general \MUON-type processes}]
      Let $(\Omega,\cF,\P)$ be a probability space, let $\delta\in \N$, $d_1,d_2,\dots,d_\delta,\fd_1,\fd_2,\dots,\fd_\delta\in \N$, $\alpha\in (0,1)$, $r\in [0,1]$,  $\varepsilon,\kappa,\cK,\fc,p\in (0,\infty)$, let $(\setX,\cU)$ be a measurable space, let $X_{n,m}\colon \Omega\to\setX$, $(n,m)\in \N^2$, be \iid\ random variables, let $\cA=\times_{i=1}^\delta\R^{d_i\times\fd_i}$, let $\smalll=(\smalll(\theta,x))_{(\theta,x)\in \cA\times \setX}\colon \cA\times \setX\to\R$ be measurable, assume for all $x\in \setX$ that $(\cA\ni\theta\mapsto \smalll(\theta,x)\in \R)\in C^1(\cA,\R)$, assume for all $\theta,\vartheta\in \cA$, $x\in \setX$ with $\max\{\|\theta\|,\|\vartheta\|\}< \fc$ that
      \begin{equation}\llabel{assume: locally Lipschitz}
          \textcolor{magenta}{\|(\nabla_\theta\smalll)(\theta,x)-(\nabla_\theta\smalll)(\vartheta,x)\|\leq \cK\|\theta-\vartheta\|}\qqandqq \textcolor{magenta}{\|(\nabla_\theta\smalll)(\theta,x)\|\leq  \cK},
      \end{equation}
      let $\cL\colon \cA\to\R$ satisfy for all $\theta\in \cA$ with $\|\theta\|< \fc$ that $\cL(\theta)=\E[\smalll(\theta,X_{1,1})]$, assume that $\cL$ is strongly convex, let $\vartheta\in \cA$ satisfy $\|\vartheta\|<\fc$ and $\cL(\vartheta)=\inf_{\theta\in \cA}\cL(\theta)$, let $(\gamma_n)_{n\in \N}\subseteq (0,\infty)$ be non-increasing, assume $\limsup_{n\to\infty} (\gamma_n+(\gamma_n)^{-2}(\gamma_{n}-\gamma_{n+1}))=0$, and for every $i\in \{1,2,\dots,\delta\}$ let $\sgn_i\colon\R^{d_i\times\fd_i}\to\R^{d_i\times\fd_i}$ be measurable and assume for all $\theta\in \R^{d_i\times\fd_i}$ that
      \begin{equation}\label{def: Pi: stochastic convergence}
       \textcolor{magenta}{ \|\sgn_i(\theta)\|\leq\cK\|\theta\|^r}\qqandqq  \textcolor{magenta}{\spro{\sgn_i(\theta),\theta}\geq \kappa\|\theta\|^2(\|\theta\|+\varepsilon)^{-1}}.
      \end{equation}
     Then there exists $\fC\in \R$ such that for every $M\in \N$, every stochastic process $\Theta=(\Theta^{1},\dots,\Theta^{\delta})\colon \N_0\times\Omega\to\cA$, and every stochastic process $\bfm=(\bfm^{1},\dots,\bfm^{\delta})\colon \N_0\times\Omega\to\cA$ with the property that
    for all $n\in \N$, $i\in \{1,2,\dots,\delta\}$ it holds that
         \begin{equation}
            \textcolor{magenta}{\bfm_0=0},\qquad \textcolor{magenta}{\bfm_n= \alpha \bfm_{n-1}+(1-\alpha)\bigl[\textstyle \frac 1M \sum_{m=1}^M(\nabla_\theta\smalll)(\Theta_{n-1},X_{n,m})\bigr]},
            \end{equation}
            \begin{equation}
           \textcolor{magenta}{\P(\|\Theta_{n-1}\|< \fc)=1},\qquad \text{and}\qquad  \textcolor{magenta}{\Theta_n^{i}=\Theta_{n-1}^{i}-\gamma_n\sgn_i(\bfm_n^{i})}
                \end{equation}
                and  with the property that $\Theta_0$ and $(X_{n,m})_{(n,m)\in \N^2}$ are independent we have for all $n \in \N_0$ that
\begin{equation} \llabel{conclude}
     \textcolor{magenta}{\bigl(\E\bigl[\|\Theta_n-\vartheta\|^p\bigr]\bigr)^{1/p}\leq \fC M^{-1/(4-2r)}+\fC\sqrt{\gamma_{n+1}}}.
\end{equation}
\end{athm}
\end{tcolorbox}
\end{samepage}
\begin{aproof}
Throughout this proof let $q\in\R$ satisfy $q=\max\{p,1\}$, let $\cG=(\cG_1,\dots,\cG_\delta)\colon \cA\to\cA$ satisfy for all $i\in \{1,2,\dots,\delta\}$, $\theta=(\theta_1,\dots,\theta_\delta)\in \cA$ with $\|\theta\|<\fc$ that
\begin{equation}\llabel{def: cG}
    \cG_i(\theta)=\nabla_{\theta_i}\cL(\theta),
\end{equation} for every $M\in \N$ and every function $\xi\colon \Omega\to\cA$ let $\Theta^{\xi,M}=(\Theta^{\xi,M,1},\dots,\Theta^{\xi,M,\delta})\colon \N_0\times\Omega\to\cA$ and $\bfm^{\xi,M}=(\bfm^{\xi,M,1},\dots,\bfm^{\xi,M,\delta})\colon \N_0\times\Omega\to\cA$ satisfy for all $n\in \N$, $i\in \{1,2,\dots,\delta\}$ that 
         \begin{equation}\llabel{def: bfm}
           \bfm_0^{\xi,M}=0,\qquad\bfm_n^{\xi,M}= \alpha \bfm_{n-1}^{\xi,M}+(1-\alpha)\bigl[\textstyle \frac 1M \sum_{m=1}^M(\nabla_\theta\smalll)(\Theta_{n-1}^{\xi,M},X_{n,m})\bigr],
            \end{equation}
            \begin{equation}\llabel{def: Theta}
            \Theta_0^{\xi,M}=\xi,\qquad\text{and}\qquad  \Theta_n^{\xi,M,i}=\Theta_{n-1}^{\xi,M,i}-\gamma_n\sgn_i(\bfm_n^{\xi,M,i}),
         \end{equation}
         and for every $M\in\N$ let 
\begin{align}
    \Xi_M&=\bigl\{\xi\colon \Omega\to \cA \,\text{be a random variable with the property that for all $n\in\N_0$ it holds that }\notag\\
    &\quad \textstyle\P\bigl(\|\Theta_n^{\xi,M}\|< \fc\bigr)=1\,\text{and with the property that}\, \xi\,\text{and}\, (X_{n,m})_{(n,m)\in \N^2}\,\text{are independent}  \bigr\}.\llabel{def: Xi}
\end{align}
\argument{\lref{assume: locally Lipschitz};the assumption that for all $x\in \setX$ it holds that $(\cA\ni\theta\mapsto \smalll(\theta,x)\in \R)\in C^1(\cA,\R)$; the fact that for all $\theta\in \cA$ with $\|\theta\|< \fc$ it holds that $\cL(\theta)=\E[\smalll(\theta,X_{1,1})]$;the fact that $X_{1,1}\in \setX$;the dominated convergence theorem}{that 
\begin{enumerate}[label=(\roman*)]
    \item \llabel{item 1} it holds that $(\{x\in \cA\colon \|x\|<\fc\}\ni \theta\mapsto \cL(\theta)\in \R)\in C^1(\{x\in \cA\colon \|x\|< \fc\},\R)$ and
    \item \llabel{item 2} it holds for all $\theta\in \cA$ with $\|\theta\|< \fc$ that $\cG(\theta)=(\nabla\cL)(\theta)=\E[(\nabla_\theta\smalll)(\theta,X_{1,1})]$.
\end{enumerate}}
 \argument{\lref{item 2};\lref{def: cG};the fundamental theorem of calculus;}{that for all $\theta_1,\theta_2\in\cA$ with $\max\{\|\theta_1\|,\|\theta_2\|\}<\fc$ it holds that
\begin{equation}\llabel{eqp1}
\begin{split}
    &\cL(\theta_1)-\cL(\theta_2)
    =\int_{0}^1\spro{(\nabla\cL)(\theta_2+t(\theta_1-\theta_2)),\theta_1-\theta_2}\,\d t\\
    &=\spro{\cG(\theta_2),\theta_1-\theta_2}+\int_{0}^1\spro{\cG(\theta_2+t(\theta_1-\theta_2))-\cG(\theta_2),\theta_1-\theta_2}\,\d t\\
    &\leq \spro{\cG(\theta_2),\theta_1-\theta_2}+\int_{0}^1\bigl\|\spro{\cG(\theta_2+t(\theta_1-\theta_2))-\cG(\theta_2),\theta_1-\theta_2}\bigr
    \|\,\d t.
    \end{split}
\end{equation}}
\argument{\lref{item 2};\lref{assume: locally Lipschitz};the fact that for all $\theta\in \cA$ with $\|\theta\|<\fc$ it holds that $\cL(\theta)=\E[\smalll(\theta,X_{1,1})]$; the fact that $X_{1,1}\in \setX$}{that for all $\theta_1,\theta_2\in \cA$ with $\max\{\|\theta_1\|,\|\theta_2\|\}< \fc$ it holds that
\begin{equation}\llabel{eqp2}
\begin{split}
    \|\cG(\theta_1)-\cG(\theta_2)\|&=\bigl\|\E[(\nabla_\theta\smalll)(\theta_1,X_{1,1})]-\E[(\nabla_\theta\smalll)(\theta_2,X_{1,1})]\bigr\|\\
    &\leq \E\bigl[\|(\nabla_\theta\smalll)(\theta_1,X_{1,1})-(\nabla_\theta\smalll)(\theta_2,X_{1,1})\|\bigr]\leq \cK\|\theta_1-\theta_2\|.
    \end{split}
\end{equation}}
\argument{\lref{eqp2};the Cauchy-Schwarz inequality; the fact that for all $M,n\in\N$, $t\in [0,1]$, $\xi\in \Xi_M$ it holds $\P$-a.s.\ that \begin{equation}
    \|\Theta_{n-1}^{\xi,M}+t(\Theta_n^{\xi,M}-\Theta_{n-1}^{\xi,M})\|\leq (1-t)\|\Theta_{n-1}^{\xi,M}\|+t\|\Theta_{n}^{\xi,M}\|< (1-t)\fc+t\fc=\fc
\end{equation}}{that there exists $\fC\in (0,\infty)$ such that for all $M,n\in \N$, $t\in [0,1]$, $\xi\in \Xi_M$ it holds $\P$-a.s.\ that
\begin{equation}\llabel{eqp3}
\begin{split}
    &\bigl|\spro{\cG(\Theta_{n-1}^{\xi,M}+t(\Theta_{n}^{\xi,M}-\Theta_{n-1}^{\xi,M}))-\cG(\Theta_{n-1}^{\xi,M}),\Theta_{n}^{\xi,M}-\Theta_{n-1}^{\xi,M}}\bigr|\\
    &\leq  \|\cG(\Theta_{n-1}^{\xi,M}+t(\Theta_{n}^{\xi,M}-\Theta_{n-1}^{\xi,M}))-\cG(\Theta_{n-1}^{\xi,M})\|\|\Theta_{n}^{\xi,M}-\Theta_{n-1}^{\xi,M}\|\\
    &\leq \fC\|\Theta_{n-1}^{\xi,M}+t(\Theta_{n}^{\xi,M}-\Theta_{n-1}^{\xi,M})-\Theta_{n-1}^{\xi,M}\|\|\Theta_{n}^{\xi,M}-\Theta_{n-1}^{\xi,M}\|=\fC t\|\Theta_n^{\xi,M}-\Theta_{n-1}^{\xi,M}\|^2\\
    &\leq \fC\|\Theta_n^{\xi,M}-\Theta_{n-1}^{\xi,M}\|^2.
    \end{split}
\end{equation}}
\argument{\lref{eqp3};\cref{def: Pi: stochastic convergence};\lref{def: Theta};\lref{eqp1};}{that there exists $\fC\in (0,\infty)$ such that for all $M,n\in \N$, $\xi\in \Xi_M$ it holds $\P$-a.s.\ that
    \begin{equation}\llabel{eq1}
    \begin{split}
    \cL(\Theta_{n}^{\xi,M})-\cL(\Theta_{n-1}^{\xi,M})&\leq \spro{\cG(\Theta_{n-1}^{\xi,M}),\Theta_{n}^{\xi,M}-\Theta_{n-1}^{\xi,M}}+\fC\|\Theta_n^{\xi,M}-\Theta_{n-1}^{\xi,M}\|^2\\
    &= \textstyle-\gamma_n\bigl[\sum_{i=1}^\delta\spro{\cG_i(\Theta_{n-1}^{\xi,M}),\sgn_i(\bfm_n^{\xi,M,i})}\bigr]+\fC\bigl[\sum_{i=1}^\delta(\gamma_n)^2\|\sgn_i(\bfm_n^{\xi,M,i})\|^2\bigr]\\
    &\leq \textstyle-\gamma_n\bigl[\sum_{i=1}^\delta\spro{\cG_i(\Theta_{n-1}^{\xi,M}),\sgn_i(\bfm_n^{\xi,M,i})}\bigr]+\fC\cK^2\bigl[\sum_{i=1}^\delta(\gamma_n)^2\|\bfm_n^{\xi,M,i}\|^{2r}\bigr].
    \end{split}
    \end{equation}}
     \argument{\lref{def: bfm};\unskip, \eg, \cref{momentum:representation}}{that for all $M\in \N$, $n\in \N_0$, $\xi\in \Xi_M$ it holds that
    \begin{equation}\llabel{eqq4}
        \bfm_{n}^{\xi,M}=\textstyle(1-\alpha)\biggl[\sum\limits_{k=1}^n\alpha^{n-k}\textstyle \frac 1M \bigl[\textstyle\sum_{m=1}^M(\nabla_\theta\smalll)(\Theta_{k-1}^{\xi,M},X_{k,m})\bigr]\biggr].
        \end{equation}}
    \argument{\lref{assume: locally Lipschitz};\lref{def: Xi};the fact that for all $n,m\in \N$ it holds that $X_{n,m}\in\setX$}{that for all $M,n,m\in \N$, $\xi\in \Xi_M$ it holds that
    \begin{equation}\llabel{def: cK'}
         \P(\|(\nabla_\theta\smalll)(\Theta_{n-1}^{\xi,M},X_{n,m})\|\leq \cK)=1.
    \end{equation}}
    \argument{\lref{def: cK'};\lref{eqq4}}{that there exists $\rho\in (0,\infty)$ such that for all $M\in \N$, $n\in \N_0$, $i\in \{1,2,\dots,\delta\}$, $\xi\in \Xi_M$ it holds $\P$-a.s.\ that
    \begin{equation}\llabel{eq4.1}
         \begin{split}
         \|\bfm_n^{\xi,M,i}\|\leq \|\bfm_n^{\xi,M}\|&\leq \textstyle(1-\alpha)\biggl[\sum\limits_{k=1}^n\alpha^{n-k}\textstyle \frac 1M \bigl[\textstyle\sum_{m=1}^M\|(\nabla_\theta\smalll)(\Theta_{k-1}^{\xi,M},X_{k,m})\|\bigr]\biggr]\\
         &\leq\textstyle(1-\alpha)\biggl[\sum\limits_{k=1}^n\alpha^{n-k}\textstyle \frac 1M \bigl[\textstyle\sum_{m=1}^M\cK\bigr]\biggr]\leq \cK.
         \end{split}
    \end{equation}}
    \argument{\cref{def: Pi: stochastic convergence};\lref{eq1};\lref{eq4.1}}{that there exists $\bfK\in (0,\infty)$ which satisfies that for all $M,n\in \N$, $\xi\in\Xi_M$ it holds $\P$-a.s.\ that
    \begin{equation}\llabel{eq2}
    \begin{split}
        \cL(\Theta_{n}^{\xi,M})&\leq \textstyle \cL(\Theta_{n-1}^{\xi,M})-\gamma_n\bigl[\sum_{i=1}^\delta\spro{\cG_i(\Theta_{n-1}^{\xi,M}),\sgn_i(\bfm_n^{\xi,M,i})}\bigr]+\bfK(\gamma_n)^2.
        \end{split}
    \end{equation}}
    \startnewargseq
    \argument{the fact that $\limsup_{n\to\infty}\gamma_n=0$; the fact that $\limsup_{n\to\infty} \bigl(\frac{\gamma_{n}-\gamma_{n+1}}{(\gamma_{n})^2}\bigr)<\infty$;}{that
    \begin{equation}\llabel{eq3'}
        \textstyle\limsup_{n\to\infty} \bigl(\frac{\gamma_{n}-\gamma_{n+1}}{\gamma_{n}}\bigr)=0.
    \end{equation}}
    \argument{\lref{eq3'};}{that
    \begin{equation}\llabel{eq3'.1}
        \textstyle\limsup_{n\to\infty} \bigl(\frac{\gamma_{n}}{\gamma_{n+1}}\bigr)=1.
    \end{equation}}
    \argument{\lref{eq3'.1};the fact that $0<\alpha<1$}{that
    \begin{equation}\llabel{eq3''}
        \limsup_{n\to\infty}(\allowbreak(\gamma_{n+1})^{-1}\alpha\gamma_n)= \alpha<1.
    \end{equation}}
      \argument{\cref{def: Pi: stochastic convergence};\lref{def: Theta};\lref{eq4.1}; the fact that for all  $x=(x_1,\dots,x_\delta)\in \cA$ it holds that $\|x\|\leq \sum_{i=1}^\delta \|x_i\|$}{that for all $M,n\in \N$, $\xi\in \Xi_M$ it holds that
    \begin{equation}\llabel{eq3}
        \|\Theta_n^{\xi,M}-\Theta_{n-1}^{\xi,M}\|\leq \sum_{i=1}^\delta\|\Theta_{n}^{\xi,M,i}-\Theta_{n-1}^{\xi,M,i}\|\leq \sum_{i=1}^\delta\gamma_n\|\sgn_i(\bfm_n^{\xi,M,i})\|\leq \sum_{i=1}^\delta\gamma_n\cK\|\bfm_n^{\xi,M,i}\|^r\leq  \cK^{r+1}\delta\gamma_n.
    \end{equation}}
        \argument{\lref{eqq4}; the fact that for all $n\in \N$, $i\in \{1,2,\dots,\delta\}$, $\xi\in \Xi_M$ it holds that $\Theta_{n}^{\xi,M,i}=\Theta_0^{\xi,M,i}-\sum_{k=1}^n\gamma_k\sgn_i(\bfm_k^{\xi,M,i})$;induction}{that for all $M,n\in \N$, $i\in \{1,2,\dots,\delta\}$, $\xi\in \Xi_M$ it holds that
        \begin{equation}\llabel{eqq6'}
            \sigma(\Theta_{n}^{\xi,M,i})\subseteq\sigma \bigl(\Theta_0^{\xi,M},(X_{k,m})_{(k,m)\in \{1,2,\dots,n\}\times\N}\bigr).
        \end{equation}}
        \argument{\lref{eqq6'};}{that for all $M,n\in \N$, $\xi\in \Xi_M$ it holds that
        \begin{equation}\llabel{eqq6}
            \sigma(\Theta_{n}^{\xi,M})\subseteq\sigma \bigl(\Theta_0^{\xi,M},(X_{k,m})_{(k,m)\in \{1,2,\dots,n\}\times\N}\bigr).
        \end{equation}}
        \argument{\lref{eqq6};the assumption that for all $M\in \N$, $\xi\in \Xi_M$ it holds that $\Theta_{0}^{\xi,M}$ and $(X_{n,m})_{(n,m)\in \N^2}$ are independent; the fact that $X_{n,m}$, $(n,m)\in\N^2$, are independent}{that \llabel{argg1} for all $M,n\in \N$, $\xi\in \Xi_M$ it holds that $\Theta_{n-1}^{\xi,M}$ and $(X_{n,m})_{m\in \N}$ are independent\dott}
    \argument{\lref{argg1};\lref{def: bfm};\lref{eq3''};\lref{eq3};\lref{item 2};\cref{lem: stochastic bfm analysis};}{that there exist $\fC\in \R$ such that for all $M,n\in \N$, $\xi\in \Xi_M$ it holds that
    \begin{equation}\llabel{eq4p}
       \textstyle\bigl(\E\bigl[\|\bfm_{n}^{\xi,M}-\cG(\Theta_{n-1}^{\xi,M})\|^{2q/(2-r)}\bigr]\bigr)^{(2-r)/(2q)}\leq \fC M^{-1/2}+\fC\gamma_n.
    \end{equation}}
     \argument{\lref{eq4p};the fact that $\frac{2}{2-r}\geq 1$; the fact that for all $n\in \N$ it holds that $\gamma_n\leq \gamma_1$}{that there exists $\fC_1,\fC_2\in (0,\infty)$ such that for all $M,n\in \N$, $\xi\in \Xi_M$ that
    \begin{equation}\llabel{eq4}
       \textstyle\bigl(\E\bigl[\|\bfm_{n}^{\xi,M}-\cG(\Theta_{n-1}^{\xi,M})\|^{2q/(2-r)}\bigr]\bigr)^{1/q}\leq \fC_1 M^{-1/(2-r)}+\fC_1(\gamma_n)^{2/(2-r)}\leq \fC_2 M^{-1/(2-r)}+\fC_2\gamma_n.
    \end{equation}}
    \argument{the Young inequality;}{that for all $x,y\in [0,\infty)$, $n,m\in [1,\infty)$ with $\nicefrac{1}{n}+\nicefrac{1}{m}=1$ it holds that
    \begin{equation}\llabel{eqqp1}
        xy\leq \frac{x^n}{n}+\frac{y^m}{m}.    
        \end{equation}}
    \argument{\lref{eqqp1};}{that for all $s\in(0,1]$, $\rho\in (0,\infty)$, $x,y\in [0,\infty)$ it holds that
    \begin{equation}\llabel{eqqp5'}
    \begin{split}
        xy^r&\leq \frac{2-s}{2}\biggl(\frac{r}{2\rho}+\mathbbm 1_{\{0\}}(s)\biggr)^{s/(2-s)} x^{2/(2-s)}+ \biggl(\frac{s}{2\rho}+\mathbbm 1_{\{0\}}(s)\biggr)^{-1}\frac {sy^2}{2}\\
        &\leq  \frac{2-s}{2}\biggl(\frac{s}{2\rho}+\mathbbm 1_{\{0\}}(s)\biggr)^{s/(2-s)} x^{2/(2-s)}+ \rho y^2.
        \end{split}
    \end{equation}}
    \argument{\lref{eqqp5'};the fact that for all $\rho\in (0,\infty)$, $x,y\in[0,\infty)$ it holds that $xy^0\leq x+\rho y^2$}{that for all $\rho\in (0,\infty)$, $x,y\in [0,\infty)$ it holds that
    \begin{equation}\llabel{eqqp5}
    \begin{split}
        xy^r&\leq \frac{2-r}{2}\biggl(\frac{r}{2\rho}+\mathbbm 1_{\{0\}}(r)\biggr)^{r/(2-r)} x^{2/(2-r)}+ \biggl(\frac{r}{2\rho}+\mathbbm 1_{\{0\}}(r)\biggr)^{-1}\frac {ry^2}{2}\\
        &\leq  \frac{2-r}{2}\biggl(\frac{r}{2\rho}+\mathbbm 1_{\{0\}}(r)\biggr)^{r/(2-r)} x^{2/(2-r)}+ \rho y^2.
        \end{split}
    \end{equation}}
    \argument{\lref{item 2};\lref{assume: locally Lipschitz}}{that for all $\theta=(\theta_1,\dots,\theta_\delta)\in\cA$ with $\|\theta\|<\fc$ it holds that
    \begin{equation}\llabel{eqqp6'}
         \|\cG_i(\theta)\|=\bigl\|\E[(\nabla_{\theta_i}\smalll)(\theta,X_{1,1})]\bigr\|\leq \cK.
    \end{equation}}
   \argument{\lref{eqqp6'};the fact that for all $M\in\N$, $n\in\N_0$ it holds that $\P(\|\Theta_n^{\xi,M}\|<\fc)=1$}{that for all $M,n\in \N$, $i\in \{1,2,\dots,\delta\}$, $\xi\in \Xi_M$ it holds that
    \begin{equation}\llabel{eqqp6}
        \P(\|\cG_i(\Theta_{n-1}^{\xi,M})\|\leq \cK)=1.
    \end{equation}}
    \argument{\cref{def: Pi: stochastic convergence};\lref{eq4.1};\lref{eqqp5};\lref{eqqp6};the fact that $\frac{2}{2-r}\leq 2$}{that there exist $\rho_1,\rho_2,\rho_3,\rho_4\in (0,\infty)$ such that for all $M,n\in \N$, $i\in \{1,2,\dots,\delta\}$, $\xi\in \Xi_M$ it holds $\P$-a.s.\ that
    \begin{equation}\llabel{arg1}
    \begin{split}
        &\spro{\cG_i(\Theta_{n-1}^{\xi,M}),\sgn_i(\bfm_n^{\xi,M,i})}=\spro{\bfm_n^{\xi,M,i},\sgn_i(\bfm_n^{\xi,M,i})}+\spro{\cG_i(\Theta_{n-1}^{\xi,M})-\bfm_n^{\xi,M,i},\sgn_i(\bfm_n^{\xi,M,i})}\\
        &\geq \frac{\kappa \|\bfm_n^{\xi,M,i}\|^2}{\|\bfm_n^{\xi,M,i}\|+\varepsilon}-\|\cG_i(\Theta_{n-1}^{\xi,M})-\bfm_n^{\xi,M,i}\|\|\sgn_i(\bfm_n^{\xi,M,i})\|\\
        &\geq 2\rho_1\|\bfm_n^{\xi,M,i}\|^2-\|\cG_i(\Theta_{n-1}^{\xi,M})-\bfm_n^{\xi,M,i}\|\|\sgn_i(\bfm_n^{\xi,M,i})\|\\
        &\geq 2\rho_1\|\bfm_n^{\xi,M,i}\|^2-\cK\|\cG_i(\Theta_{n-1}^{\xi,M})-\bfm_n^{\xi,M,i}\|\|\bfm_n^{\xi,M,i}\|^r\\
        &\geq 2\rho_1\|\bfm_n^{\xi,M,i}\|^2-\rho_2\|\cG_i(\Theta_{n-1}^{\xi,M})-\bfm_n^{\xi,M,i}\|^{2/(2-r)}-\rho_1\|\bfm_n^{\xi,M,i}\|^2\\
        &=  \rho_1\|\bfm_n^{\xi,M,i}\|^2-\rho_2\|\cG_i(\Theta_{n-1}^{\xi,M})-\bfm_n^{\xi,M,i}\|^{2/(2-r)}\\
        &\geq \frac{\rho_1}{2}\|\cG_i(\Theta_{n-1}^{\xi,M})\|^2-2\rho_1\|\bfm_n^{\xi,M,i}-\cG_i(\Theta_{n-1}^{\xi,M})\|^{2}-\rho_3\|\cG_i(\Theta_{n-1}^{\xi,M})-\bfm_n^{\xi,M,i}\|^{2/(2-r)}\\
        &\geq \frac{\rho_1}{2}\|\cG_i(\Theta_{n-1}^{\xi,M})\|^2-\rho_4\|\cG_i(\Theta_{n-1}^{\xi,M})-\bfm_n^{\xi,M,i}\|^{2/(2-r)}.
        \end{split}
    \end{equation}}
    \startnewargseq
    \argument{\lref{eq2};\lref{arg1};}{that there exist $\rho_1,\rho_2,\rho_3\in (0,\infty)$ such that for all $M,n\in \N$, $\xi\in \Xi_M$ it holds $\P$-a.s.\ that
    \begin{equation}\llabel{eq7.1}
    \begin{split}
      & \cL(\Theta_n^{\xi,M})\\&\textstyle\leq  \cL(\Theta_{n-1}^{\xi,M})-\gamma_n\sum\limits_{i=1}^\delta\bigl( \rho_1\|\cG_i(\Theta_{n-1}^{\xi,M})\|^2-\rho_2\|\cG_i(\Theta_{n-1}^{\xi,M})-\bfm_n^{\xi,M,i}\|^{2/(2-r)}\bigr)+\bfK(\gamma_n)^2\\
      &=\textstyle\cL(\Theta_{n-1}^{\xi,M})-\gamma_n\rho_1\biggl[\sum\limits_{i=1}^\delta  \|\cG_i(\Theta_{n-1}^{\xi,M})\|^2\biggr]+\gamma_n\rho_2\biggl[\sum\limits_{i=1}^\delta\|\bfm_n^{\xi,M,i}-\cG_i(\Theta_{n-1}^{\xi,M})\|^{2/(2-r)}\biggr]+\bfK(\gamma_n)^2\\
       &\leq   \textstyle \cL(\Theta_{n-1}^{\xi,M})-\gamma_n\rho_1 \|\cG(\Theta_{n-1}^{\xi,M})\|^2+\gamma_n\rho_3\|\bfm_n^{\xi,M}-\cG(\Theta_{n-1}^{\xi,M})\|^{2/(2-r)}+\bfK(\gamma_n)^2.
       \end{split}
    \end{equation}}
    \argument{\lref{eq7.1};}{that there exist $\rho_1,\rho_2\in (0,\infty)$ such that for all $M,n\in \N$, $\xi\in \Xi_M$ it holds that $\P$-a.s.\ that \begin{equation}\llabel{eq7.2}
    \begin{split}
       &\cL(\Theta_n^{\xi,M})-\cL(\vartheta)\\
       &\leq \textstyle  \cL(\Theta_{n-1}^{\xi,M})-\cL(\vartheta)-\gamma_n\rho_1 \|\cG(\Theta_{n-1}^{\xi,M})\|^2+\gamma_n\rho_2\|\bfm_n^{\xi,M}-\cG(\Theta_{n-1}^{\xi,M})\|^{2/(2-r)}+\bfK(\gamma_n)^2.
       \end{split}
    \end{equation}}
          \argument{\lref{item 1};the assumption that $\|\vartheta\|< \fc$ ;the assumption that $\cL(\vartheta)=\inf_{\theta\in \cA}\cL(\theta)$}{that  \begin{equation}\llabel{arg2}
              \cG(\vartheta)=0.
          \end{equation}}
          \argument{\lref{arg2};\lref{eqp1};\lref{eqp2}}{that there exists $\varrho\in (0,\infty)$ such that for all $\theta\in \cA$ with $\|\theta\|< \fc$ it holds that
\begin{equation}\llabel{eq7.2.2}
\begin{split}
    \cL(\theta)-\cL(\vartheta)
    &=\spro{\cG(\vartheta),\theta-\vartheta}+\int_{0}^1\spro{\cG(\vartheta+t(\theta-\vartheta))-\cG(\vartheta),\theta-\vartheta}\,\d t\\
    &=\int_{0}^1\spro{\cG(\vartheta+t(\theta-\vartheta))-\cG(\vartheta),\theta-\vartheta}\,\d t\\
    &\leq \int_{0}^1\|\cG(\vartheta+t(\theta-\vartheta))-\cG(\vartheta)\|\|\theta-\vartheta\|\,\d t\\
    &\leq \int_{0}^1\varrho t\|(\theta-\vartheta)\|^2\,\d t\leq \varrho\|\theta-\vartheta\|^2.
    \end{split}
\end{equation}}
\argument{\lref{eq7.2.2};the fact that for all $M\in\N$, $n\in\N_0$, $\xi\in\Xi_M$ it holds that it holds that $\P\bigl(\|\Theta_n^{\xi,M}\|< \fc\bigr)$;the fact that $\cL(\vartheta)=\inf_{\theta\in \cA}\cL(\theta)$}{that there exists $\varrho\in (0,\infty)$ such that for all $M\in \N$, $n\in \N_0$, $\xi\in \Xi_M$ it holds $\P$-a.s.\ that
\begin{equation}\llabel{eq7.2.3}
    0\leq \cL(\Theta_n^{\xi,M})-\cL(\vartheta)\leq \varrho\|\Theta_n^{\xi,M}-\vartheta\|^2.
\end{equation}}
          \argument{\lref{arg2};the fact that $\cL$ is strongly convex;the Cauchy-Schwarz inequality;\unskip, \eg, \cite[item (iii) in Proposition 5.7.23]{ArBePhi2024}}{that there exists $\varrho\in (0,\infty)$ such that for all $\theta\in \cA$ with $\|\theta\|<\fc$ it holds that
          \begin{equation}\llabel{eq7.2.4}
       \varrho\|\theta-\vartheta\|^2\leq  \spro{\cG(\theta)-\cG(\vartheta),\theta-\vartheta}=\spro{\cG(\theta),\theta-\vartheta}\leq  \|\theta-\vartheta\|\|\cG(\theta)\|.
          \end{equation}}
          \argument{\lref{eq7.2.4};}{that there exists $\varrho\in (0,\infty)$ such that for all $\theta\in \cA$ with $\|\theta\|<\fc$ it holds that
          \begin{equation}\llabel{eq7.2.5}
              \|\cG(\theta)\|\geq \varrho\|\theta-\vartheta\|.
          \end{equation}}
          \argument{\lref{def: Xi};\lref{eq7.2.5};\lref{eq7.2.3};}{that there exist $\varrho_1,\varrho_2\in (0,\infty)$ such that for all $M\in \N$, $n\in \N_0$, $\xi\in \Xi_M$ it holds $\P$-a.s.\ that
          \begin{equation}\llabel{eq7.2.6}
              \|\cG(\Theta_n^{\xi,M})\|^2\geq (\varrho_1)^2\|\Theta_n^{\xi,M}-\vartheta\|^2\geq \varrho_2(\cL(\Theta_n^{\xi,M})-\cL(\vartheta)).
          \end{equation}}
    \argument{\lref{eq7.2.6};\lref{eq7.2};}{that there exist $\rho_1,\rho_2\in (0,\infty)$ which satisfy that for all $M,n\in \N$, $\xi\in \Xi_M$ it holds $\P$-a.s.\ that
    \begin{equation}\llabel{eq7'}
    \begin{split}
        &\cL(\Theta_n^{\xi,M})-\cL(\vartheta)\\
        &\leq \cL(\Theta_{n-1}^{\xi,M})-\cL(\vartheta)-\rho_1\gamma_n(\cL(\Theta_{n-1}^{\xi,M})-\cL(\vartheta))+\gamma_n\rho_2\|\bfm_n^{\xi,M}-\cG(\Theta_{n-1}^{\xi,M})\|^{2/(2-r)}+\bfK(\gamma_n)^2\\
        &=(1-\rho_1\gamma_n)(\cL(\Theta_{n-1}^{\xi,M})-\cL(\vartheta))+\gamma_n\rho_2\|\bfm_n^{\xi,M}-\cG(\Theta_{n-1}^{\xi,M})\|^{2/(2-r)}+\bfK(\gamma_n)^2.
        \end{split}
    \end{equation}}
    \startnewargseq
    \argument{the fact that $\limsup_{n\to\infty}\gamma_n=0$}{that there exists $\fn\in \N$ which satisfies for all $n\in \N\cap[\fn,\infty)$ that
    \begin{equation}\llabel{def: fn}
        \rho_1\gamma_n\leq 1.
    \end{equation}}
    \startnewargseq
    \argument{\lref{eq7'};\lref{def: fn};the fact that for all $\theta\in \cA$ it holds that $\cL(\theta)-\cL(\vartheta)\geq 0$;the triangle inequality;}{that for all $M\in \N$, $n\in \N\cap[\fn,\infty)$, $\xi\in \Xi_M$ it holds that
    \begin{equation}\llabel{eq7''}
        \begin{split}
            \bigl(\E\bigl[|\cL(\Theta_n^{\xi,M})-\cL(\vartheta)|^q\bigr]\bigr)^{1/q}
            &\leq (1-\rho_1\gamma_n)\bigl(\E\bigl[|\cL(\Theta_{n-1}^{\xi,M})-\cL(\vartheta)|^q\bigr]\bigr)^{1/q}\\
            &+\gamma_n\rho_2\bigl(\E\bigl[\|\bfm_n^{\xi,M}-\cG(\Theta_{n-1}^{\xi,M})\|^{2q/(2-r)}\bigr]\bigr)^{1/q}+\bfK(\gamma_n)^2.
        \end{split}
    \end{equation}}
    \argument{\lref{eq7''};\lref{eq4};}{that there exists $\fC\in (0,\infty)$ such that for all $M\in \N$, $n\in \N\cap[\fn,\infty)$, $\xi\in \Xi_M$ it holds that
    \begin{equation}\llabel{eq7}
    \begin{split}
        &\bigl(\E\bigl[|\cL(\Theta_n^{\xi,M})-\cL(\vartheta)|^q\bigr]\bigr)^{1/q}\\&\leq (1-\rho_1\gamma_n)\bigl(\E\bigl[|\cL(\Theta_{n-1}^{\xi,M})-\cL(\vartheta)|^{q}\bigr]\bigr)^{1/q} +\fC\gamma_n M^{-1/(2-r)}+\fC(\gamma_n)^2+\bfK(\gamma_n)^2\\
        &= (1-\rho_1\gamma_n)\bigl(\E\bigl[|\cL(\Theta_{n-1}^{\xi,M})-\cL(\vartheta)|^{q}\bigr]\bigr)^{1/q} +\fC\gamma_n M^{-1/(2-r)}+(\fC+\bfK)(\gamma_n)^2.
        \end{split}
    \end{equation}}
   \argument{\lref{eq7};\cref{lem: stochastic recurrence 2}; the assumption that $\limsup_{n\to\infty} (\gamma_n+(\gamma_n)^{-2}(\gamma_{n}-\gamma_{n+1}))=0$; the fact that for all $M\in \N$, $\xi\in \Xi_M$ it holds that
   \begin{equation}
      \bigl(\E\bigl[|\cL(\Theta_0^{\xi,M})-\cL(\vartheta)|^q\bigr]\bigr)^{1/q}=\bigl(\E\bigl[|\cL(\xi)-\cL(\vartheta)|^q\bigr]\bigr)^{1/q}\leq  \sup_{\theta\in \cA,\,\|\theta\|< \fc}\bigl[|\cL(\theta)|+|\cL(\vartheta)|\bigr]
   \end{equation}}{that there exists $\fC\in (0,\infty)$ such that for all $M\in \N$, $n\in \N\cap[\fn,\infty)$, $\xi\in \Xi_M$ it holds that 
  \begin{equation}\llabel{eq11}
      \bigl(\E\bigl[|\cL(\Theta_n^{\xi,M})-\cL(\vartheta)|^q\bigr]\bigr)^{1/q}\leq \fC M^{-1/(2-r)}+\fC\gamma_n.
  \end{equation}}
  \argument{\lref{eq11};\lref{eq3'.1};}{that there exists $\fC\in (0,\infty)$ such that for all $M\in \N$, $n\in \N\cap[\fn,\infty)$, $\xi\in \Xi_M$ it holds that 
  \begin{equation}\llabel{eq12}
    \bigl(\E\bigl[|\cL(\Theta_n^{\xi,M})-\cL(\vartheta)|^q\bigr]\bigr)^{1/q}\leq \fC M^{-1/(2-r)}+\fC\gamma_{n+1}.
  \end{equation}}
  \argument{\lref{eq7.2.3};the fact that for all $\xi\in \Xi_M$, $M\in\N$, $n\in\N_0$ it holds that $\P(\|\Theta_n^{\xi,M}\|<\fc)$;\lref{item 1};the fact that for all $n\in \{0,1,\dots,\fn\}$ it holds that $\gamma_{n+1}\geq \gamma_{\fn+1}$}{that there exist $\fC_1,\fC_2\in (0,\infty)$ such that for all $M\in \N$, $n\in \{0,1,\dots,\fn\}$, $\xi\in \Xi_M$ it holds that
  \begin{equation}\llabel{eq13}
       \bigl(\E\bigl[|\cL(\Theta_n^{\xi,M})-\cL(\vartheta)|^q\bigr]\bigr)^{1/q}\leq \fC_1\leq \fC_2\gamma_{\fn+1}\leq \fC_2\gamma_{n+1}\leq   \fC _2M^{-1/(2-r)}+\fC_2\gamma_{n+1}.
  \end{equation}}
  \argument{\lref{eq13};\lref{eq12};}{that there exists $\fC\in (0,\infty)$ such that for all $M\in \N$, $n\in \N_0$, $\xi\in \Xi_M$ it holds that 
  \begin{equation}\llabel{eq14}
      \bigl(\E\bigl[|\cL(\Theta_n^{\xi,M})-\cL(\vartheta)|^q\bigr]\bigr)^{1/q}\leq \fC M^{-1/(2-r)}+\fC\gamma_{n+1}.
  \end{equation}}
   \argument{the fact that $\cL$ is strongly convex;the fact that $\cL(\vartheta)=\inf_{\theta\in \cA}\cL(\theta)$;\unskip, \eg, \cite[item (ii) in Proposition 5.7.23]{ArBePhi2024}}{that there exists $c\in (0,\infty)$ such that for all $\theta\in \cA$ with $\|\theta\|<\fc$ it holds that
          \begin{equation}\llabel{eq7.2.1}
         \cL(\theta)-\cL(\vartheta)\geq \spro{(\nabla\cL)(\vartheta),\theta-\vartheta}+c\|\theta-\vartheta\|^2\geq c\|\theta-\vartheta\|^2.
          \end{equation}}
  \argument{\lref{eq14};\lref{eq7.2.1};the fact that for all  $M\in\N$, $n\in\N_0$, $\xi\in \Xi_M$ it holds that $\P(\|\Theta_n^{\xi,M}\|<\fc)=1$}{that there exists $\fC\in (0,\infty)$ such that for all $M\in \N$, $n\in \N_0$, $\xi\in \Xi_M$ it holds that 
  \begin{equation}\llabel{eq15}
      \bigl(\E\bigl[\|\Theta_n^{\xi,M}-\vartheta\|^q\bigr]\bigr)^{1/q}\leq \fC M^{-1/(4-2r)}+\fC\sqrt{\gamma_{n+1}}.
  \end{equation}}
  \argument{\lref{eq15};}{\lref{conclude}\dott}
\end{aproof}
    \subsection{Error estimates for idealized MUON under boundedness assumptions}\label{subsec: idealized muon}

    In this subsection we establish in \cref{main theorem stochastic} below an error analysis for the idealized \MUON\ optimizer (cf.\ \cite[Section "The design of Muon"]{jordan2024muon} and, \eg, \cite[Definition 6.12.2]{ArBePhi2024}). Our proof of \cref{main theorem stochastic} is based on an application of \cref{theo: stochastic convergence} above and also employs the elementary and well-known facts on (orthogonal) matrices in \cref{lem: Pi property 1}, \cref{lem: Pi property 2}, and \cref{lem: Pi property 3}. Only for completeness we include here detailed proofs for \cref{lem: Pi property 1}, \cref{lem: Pi property 2}, and \cref{lem: Pi property 3}.
        \begin{tcolorbox}[colback=white!95!gray,
                  colframe=black,
                  boxrule=0.5pt,
                  sharp corners,
                  enhanced,
                  breakable,
                 ]
\begin{athm}{lemma}{lem: Pi property 1}
      Let $d,\fd\in \N$, $A,B\in\R^{d\times \fd}$. Then\footnotemark\!
          $\spro{A,B}=\Tr(A^{\top}B)$.
\end{athm}
\end{tcolorbox}
\footnotetext{Note that for all $n \in \N$, $v = ( v_{i,j})_{(i,j)\in \{1,2,\dots,n\}\times\{1,2,\dots,n\}} \in \R^{n\times n}$ it holds that $\Tr(v) = \sum_{ i = 1 }^nv_{i,i}$ (trace of a matrix).} 
\begin{aproof}
    Throughout this proof let $\bfA=( a_{i,j})_{(i,j)\in \{1,2,\dots,d\}\times\{1,2,\dots,\fd\}}$, $\bfB=( b_{i,j}\allowbreak\break)_{(i,j)\in \{1,2,\dots,d\}\times\{1,2,\dots,\fd\}}$, $\bfC=( c_{i,j})_{(i,j)\in \{1,2,\dots,\fd\}\times\{1,2,\dots,\fd\}}\in \R^{\fd\times\fd}$ satisfy $\bfA=A$, $\bfB=B$, and $\bfC=A^{\top}B$.
    \argument{the fact that $\bfC=A^{\top}B$}{that for all $i,j\in \{1,2,\dots,\fd\}$ it holds that \llabel{eq1}
        $c_{i,j}=\sum_{k=1}^d a_{k,i}b_{k,j}$\dott}
    \argument{\lref{eq1};}{that
    \begin{equation}\llabel{eq2}
        \Tr(A^{\top}B)=\sum_{i=1}^\fd c_{i,i}=\sum_{i=1}^\fd\sum_{k=1}^d a_{k,i}b_{k,i}=\spro{A,B}.
    \end{equation}}
\end{aproof}
    \begin{tcolorbox}[colback=white!95!gray,
                  colframe=black,
                  boxrule=0.5pt,
                  sharp corners,
                  enhanced,
                  breakable,
                 ]
\begin{athm}{lemma}{lem: Pi property 2}
     Let\footnotemark\! $d,\fd\in \N$, $A\in \{O\in\R^{d\times\fd}\colon (OO^{\top}=I_{d})\vee (O^{\top}O=I_{\fd})\}$.  
      Then
      \begin{equation}\llabel{conclude}
          \|A\|^2=\min\{d,\fd\}.
      \end{equation}
\end{athm}
\end{tcolorbox}
\footnotetext{Note that for all $n\in \N$, $v\in \R^n$ it holds that $I_n v=v$ (identity matrix in $\R^{ n \times n }$).} 
\begin{aproof}
    In our proof of \lref{conclude} we distinguish between the case $A^{\top}A=I_\fd$ and the case $AA^{\top}=I_d$. We first prove \lref{conclude} in the case 
    \begin{equation}\llabel{case 1}
        A^{\top}A=I_\fd.
    \end{equation}
    \argument{\lref{case 1};the fact that $A\in \R^{d\times\fd}$}{that \llabel{arg1} $\fd\leq d$\dott}
    \argument{\lref{arg1};\lref{case 1};\cref{lem: Pi property 1};}{that
    \llabel{case 1: eq1}
        $\|A\|^2=\spro{A,A}=\Tr(A^{\top}A)=\Tr(I_\fd)=\fd=\min\{d,\fd\}$\dott}
    This proves \lref{conclude} in the case $A^{\top}A=I_\fd$. We now prove \lref{conclude} in the case
    \begin{equation}\llabel{case 2}
        AA^{\top}=I_d.
    \end{equation}
    \startnewargseq
     \argument{\lref{case 2};the fact that $A\in \R^{d\times\fd}$}{that \llabel{arg2} $d\leq \fd$\dott}
      \argument{\lref{arg2};\lref{case 2};\cref{lem: Pi property 1};}{that
    \begin{equation}\llabel{case 2: eq1}
    \|A\|^2=\|A^{\top}\|^2=\spro{A^{\top},A^{\top}}=\Tr(AA^{\top})=\Tr(I_d)=d=\min\{d,\fd\}.
    \end{equation}}
    This establishes \lref{conclude} in the case $AA^{\top}=I_d$.
\end{aproof}
\begin{samepage}
    \begin{tcolorbox}[colback=white!95!gray,
                  colframe=black,
                  boxrule=0.5pt,
                  sharp corners,
                  enhanced,
                  breakable,
                 ]
\begin{athm}{lemma}{lem: Pi property 3}
      Let $d,\fd\in \N$, let $\cO=\{O\in\R^{d\times\fd}\colon (OO^{\top}=I_{d})\vee (O^{\top}O=I_{\fd})\}$, let $\sgn\colon\R^{d\times\fd}\to\cO$ satisfy for all $A\in \R^{d\times\fd}$ that
      \begin{equation}\llabel{def: Pi}
        \textstyle  \|\sgn(A)-A\|=\inf_{O\in \cO}\|O-A\|,
      \end{equation}
      and let $A\in\R^{d\times \fd}$. Then
          $\spro{A,\sgn(A)}\geq \|A\|$.
\end{athm}
\end{tcolorbox}
\end{samepage}
\begin{aproof}
\argument{\unskip, \eg, \cite[Theorem 7.3.5]{horn_johnson_2012};}{that there exist $U\in \R^{d\times d}$, $V\in \R^{\fd\times\fd}$, $\Sigma=(\sigma_{i,j})_{(i,j)\in \{1,2,\dots,d\}\times \{1,2,\dots,\fd\}}\in \R^{d\times \fd}$ which satisfy that
\begin{enumerate}[label=(\roman*)]
    \item \llabel{item 1} it holds that $UU^{\top}=U^{\top}U=I_d$ and $VV^{\top}=V^{\top}V=I_\fd$,
    \item \llabel{item 2} it holds for all $i\in \{1,2,\dots,d\}$, $j\in \{1,2,\dots,\fd\}$ with $i\neq j$ that $\sigma_{i,j}=0$,
    \item \llabel{item 3} it holds for all $i\in \{1,2,\dots, \min\{d,\fd\}\}$ that $\sigma_{i,i}\geq 0$, and
    \item \llabel{item 4} it holds that $A=U\Sigma V^{\top}$.
\end{enumerate}}
In the following let $B=(b_{i,j})_{(i,j)\in \{1,2,\dots,d\}\times\{1,2,\dots,\fd\}}\in \R^{d\times\fd}$ satisfy for all $i\in \{1,2,\dots,d\}$, $j\in \{1,2,\dots,\fd\}$ that
\begin{equation}\llabel{def: B}
    b_{i,j}=\mathbbm 1_{\{i=j,\, i\leq \min \{d,\fd\}\}}.
\end{equation}
\startnewargseq
\argument{\lref{def: B};}{that
\begin{equation}\llabel{eq1.1}
    \begin{cases}
        B^{\top}B=I_\fd\colon & d\geq \fd\\
        BB^{\top}=I_d\colon &d\leq \fd.
    \end{cases}
\end{equation}}
\argument{\lref{eq1.1};}{that 
\begin{equation}\llabel{eq2}
    \|B\|^2=\min\{d,\fd\}.
\end{equation}}
\argument{\lref{eq2};the fact that $\sgn(A)\in \cO$;\cref{lem: Pi property 2}}{that
\begin{equation}\llabel{eq3}
    \|\sgn(A)\|^2= \min\{d,\fd\}=\|B\|^2.
\end{equation}}
\argument{\lref{item 1};\cref{lem: Pi property 1}; the fact that for all $X,Y\in \R^{\fd\times\fd}$ it holds that $\Tr(XY)=\Tr(YX)$}{that 
\begin{equation}\llabel{eq4}
\begin{split}
\|UBV^{\top}\|^2&=\Tr((UBV^{\top})^{\top}UBV^{\top})=\Tr(VB^{\top}U^{\top}UBV^{\top})=\Tr(VB^{\top}BV^{\top})\\
&=\Tr(V^{\top}VB^{\top}B)=\Tr(B^{\top}B)=\|B\|^2.
\end{split}
\end{equation}}
\argument{\lref{eq4};\lref{eq3};}{that
\begin{equation}\llabel{eq5}
    \|\sgn(A)\|^2= \|UBV^{\top}\|^2.
\end{equation}}
\argument{\lref{eq1.1};\lref{item 1}}{that 
\begin{equation}\llabel{eq6}
     \begin{cases}
       (UBV^{\top})^{\top}UBV^{\top}=VB^{\top}BV^{\top}=VV^{\top}=I_\fd\colon & d\geq \fd\\
        UBV^{\top}(UBV^{\top})^{\top}=UBB^{\top}U^{\top}=UU^{\top}=I_d\colon &d\leq \fd.
    \end{cases}
\end{equation}}
\argument{\lref{eq6};}{that
\begin{equation}\llabel{eq7}
    UBV^{\top}\in \cO.
\end{equation}}
\argument{\lref{eq7};\lref{def: Pi}}{that
\begin{equation}\llabel{eq8}
    \|A-\sgn(A)\|^2\leq \|A-UBV^{\top}\|^2.
\end{equation}}
\argument{\lref{eq8};}{that
\begin{equation}\llabel{eq9}
    \|A\|^2-2\spro{A,\sgn(A)}+\|\sgn(A)\|^2\leq \|A\|^2-2\spro{A,UBV^{\top}}+\|UBV^{\top}\|^2.
\end{equation}}
\argument{\lref{eq9};\lref{eq5}}{that
\begin{equation}\llabel{eq10}
    \spro{A,\sgn(A)}\geq \spro{A,UBV^{\top}}.
\end{equation}}
\argument{\lref{eq10};\lref{item 1};\lref{item 4};\cref{lem: Pi property 1};the fact that for all $X,Y\in \R^{\fd\times\fd}$ it holds that $\Tr(XY)=\Tr(YX)$}{that
\begin{equation}\llabel{eq11}
\begin{split}
    \spro{A,\sgn(A)}&\geq \Tr((UBV^{\top})^{\top}A)=\Tr((UBV^{\top})^{\top}U\Sigma V^{\top})=\Tr(VB^{\top}U^{\top}U\Sigma V^{\top})\\
    &=\Tr(VB^{\top}\Sigma V^{\top})=\Tr(V^{\top}VB^{\top}\Sigma)=\Tr(B^{\top}\Sigma).
    \end{split}
\end{equation}}
\argument{\lref{eq11};\lref{def: B};\lref{item 2};\cref{lem: Pi property 1}}{that
\begin{equation}\llabel{eq12}
    \spro{A,\sgn(A)}\geq \spro{\Sigma,B}= \sum_{i=1}^{\min\{d,\fd\}}\sigma_{i,i}.
\end{equation}}
\argument{\lref{eq12};\lref{item 1};\lref{item 2};\lref{item 3};\lref{item 4};\cref{lem: Pi property 1}; the fact that for all $X,Y\in \R^{\fd\times \fd}$ it holds that $\Tr(XY)=\Tr(YX)$}{that 
\begin{equation}\llabel{eq13}
\begin{split}
\|A\|^2&=\|U\Sigma V^{\top}\|^2=\Tr((U\Sigma V^{\top})^{\top}U\Sigma V^{\top})=\Tr(V\Sigma^{\top}U^{\top}U\Sigma V^{\top})=\Tr(V\Sigma^{\top}\Sigma V^{\top})\\
&\textstyle =\Tr(V^{\top}V\Sigma^{\top}\Sigma)=\Tr(\Sigma^{\top}\Sigma)=\|\Sigma\|^2=\sum\limits_{i=1}^{\min\{d,\fd\}}(\sigma_{i,i})^2\leq \biggl[\sum\limits_{i=1}^{\min\{d,\fd\}}\sigma_{i,i}\biggr]^2\leq  |\spro{A,\sgn(A)}|^2.
\end{split}
\end{equation}}
\argument{\lref{eq13};the fact that $\spro{A,\sgn(A)}\geq 0$;}{that   $\spro{A,\sgn(A)}\geq \|A\|$\dott}
\end{aproof}
\begin{samepage}
    \begin{tcolorbox}[colback=white!95!gray,
                  colframe=black,
                  boxrule=0.5pt,
                  sharp corners,
                  enhanced,
                  breakable,
                 ]
\begin{athm}{theorem}{main theorem stochastic}[\textcolor{red}{Idealized \MUON\ error analysis}]
     Let $(\Omega,\cF,\P)$ be a probability space, let $\delta\in \N$, $\alpha\in (0,1)$,  $\cK,\fc,p\in (0,\infty)$, let $(\setX,\cU)$ be a measurable space, let $X_{n,m}\colon \Omega\to\setX$, $(n,m)\in \N^2$, be \iid\ random variables, for every $i\in \{1,2,\dots,\delta\}$ let $d_i,\fd_i\in \N$, let $\cO_i=\{O\in\R^{d_i\times\fd_i}\colon (OO^{\top}=I_{d_i})\vee (O^{\top}O=I_{\fd_i})\}$, let $\sgn_i\colon\R^{d_i\times\fd_i}\to\cO_i$ be measurable, and assume for all $\theta\in \R^{d_i\times\fd_i}$ that
      \begin{equation}\label{def: Pi: main theorem stochastic}
        \textcolor{magenta}{ \textstyle \|\sgn_i(\theta)-\theta\|=\inf_{O\in \cO_i}\|O-\theta\|},
      \end{equation}
      let $\cA=\times_{i=1}^\delta\R^{d_i\times\fd_i}$, let $\smalll=(\smalll(\theta,x))_{(\theta,x)\in \cA\times \setX}\colon \cA\times \setX\to\R$ be measurable, assume for all $x\in \setX$ that $(\cA\ni\theta\mapsto \smalll(\theta,x)\in \R)\in C^1(\cA,\R)$, assume for all $\theta,\vartheta\in \cA$, $x\in \setX$ with $\max\{\|\theta\|,\|\vartheta\|\}< \fc$ that
      \begin{equation}\llabel{assume: locally Lipschitz}
          \textcolor{magenta}{\|(\nabla_\theta\smalll)(\theta,x)-(\nabla_\theta\smalll)(\vartheta,x)\|\leq \cK\|\theta-\vartheta\|}\qqandqq \textcolor{magenta}{\|(\nabla_\theta\smalll)(\theta,x)\|\leq  \cK},
          \end{equation} let $\cL\colon \cA\to\R$ satisfy for all $\theta\in \cA$ with $\|\theta\|<\fc$ that $\cL(\theta)=\E[\smalll(\theta,X_{1,1})]$, assume that $\cL$ is strongly convex, let $\vartheta\in \cA$ satisfy $\|\vartheta\|< \fc$ and $\cL(\vartheta)=\inf_{\theta\in \cA}\cL(\theta)$, let $(\gamma_n)_{n\in \N}\subseteq (0,\infty)$ be non-increasing, and assume $\limsup_{n\to\infty} (\gamma_n+(\gamma_n)^{-2}(\gamma_{n}-\gamma_{n+1}))=0$. Then there exists $\fC\in \R$ such that for
        every $M\in \N$, every stochastic process $\bfm=(\bfm^{1},\dots,\bfm^{\delta})\colon \N_0\times\Omega\to\cA$, and every stochastic process $\Theta=(\Theta^{1},\dots,\Theta^{\delta})\colon \N_0\times\Omega\to\cA$ 
     with the property that for all $n\in \N$, $i\in \{1,2,\dots,\delta\}$ it holds that 
         \begin{equation}\llabel{def: bfm}
            \textcolor{magenta}{\bfm_0=0},\qquad \textcolor{magenta}{\bfm_n= \alpha \bfm_{n-1}+(1-\alpha)\bigl[\textstyle \frac 1M \sum_{m=1}^M(\nabla_\theta\smalll)(\Theta_{n-1},X_{n,m})\bigr]},
            \end{equation}
            \begin{equation}\llabel{def: Theta}
           \textcolor{magenta}{ \P(\|\Theta_{n-1}\|< \fc)=1},\qquad  \text{and}\qquad  \textcolor{magenta}{\Theta_n^{i}=\Theta_{n-1}^{i}-\gamma_n\sgn_i(\bfm_n^{i})}
         \end{equation}
         and with the property that $\Theta_{0}$ and $(X_{n,m})_{(n,m)\in \N^2}$ are independent we have for all $n\in \N_0$ that
\begin{equation} \llabel{conclude}
     \textcolor{magenta}{\bigl(\E\bigl[\|\Theta_n-\vartheta\|^p\bigr]\bigr)^{1/p}\leq \fC (M^{-1/4}+\sqrt{\gamma_{n+1}})}.
\end{equation}
\end{athm}
\end{tcolorbox}
\end{samepage}
\begin{aproof}
    \argument{\cref{lem: Pi property 2};the fact that for all $i\in \{1,2,\dots,\delta\}$, $\theta\in \R^{d_i\times\fd_i}$ it holds that $\sgn_i(\theta)\in \cO_i$}{that for all $i\in \{1,2,\dots,\delta\}$, $\theta\in \R^{d_i\times\fd_i}$ it holds that
    \begin{equation}\llabel{eq1.1}
        \|\sgn_i(\theta)\|^2= \max\{d_i,\fd_i\}.
    \end{equation}}
    \argument{\lref{eq1.1};}{for all $i\in \{1,2,\dots,\delta\}$, $\theta\in \R^{d_i\times\fd_i}$ that
    \begin{equation}\llabel{eq1.2}
        \|\sgn_i(\theta)\|\leq [\max\{d_i,\fd_i\}]^{1/2}\leq  \textstyle\bigl[\max_{i\in \{1,2,\dots,\delta\}}\max\{d_i,\fd_i\}\bigr]^{1/2}.
    \end{equation}}
    \argument{\cref{def: Pi: main theorem stochastic};\cref{lem: Pi property 3}}{that for all $i\in \{1,2,\dots,\delta\}$, $\theta\in \R^{d_i\times\fd_i}$ it holds that
    \begin{equation}\llabel{eq4.1}
        \spro{\sgn_i(\theta),\theta}\geq \|\theta\|\geq \frac{\|\theta\|^2}{\|\theta\|+1}.
    \end{equation}}
   \argument{\lref{eq4.1};\lref{eq1.2};\cref{theo: stochastic convergence}}{\lref{conclude}\dott}
\end{aproof}

In the next elementary and well-known lemma, \cref{lem: measuability}, we show that the functions $\Phi_i \colon \R^{d_i\times\fd_i} \to \cO_i$, $i \in \{1,2,\dots,\delta\}$, in \cref{def: Pi: main theorem stochastic} do indeed exist. Only for completeness we include here a proof for \cref{lem: measuability}.
 \begin{tcolorbox}[colback=white!95!gray,
                  colframe=black,
                  boxrule=0.5pt,
                  sharp corners,
                  enhanced,
                  breakable,
                 ]
                 \begin{athm}{lemma}{lem: measuability}
                 Let $d,\fd\in \N$ and let $\cO=\{O\in\R^{d\times\fd}\colon ((OO^{\top}=I_{d})\vee (O^{\top}O=I_{\fd})\}$. Then
                 \begin{enumerate}[label=(\roman*)]
                   \item \label{item 1: measuability} there exists a measurable $\Phi\colon \R^{d\times\fd}\to\cO$ such that for all $\theta\in \R^{d\times\fd}$ it holds that  $ \textstyle \|\Phi(\theta)-\theta\|=\inf_{O\in \cO}\|O-\theta\|$ and
                     \item \label{item 2: measuability} it holds for every full-rank matrix $\theta\in\R^{d\times\fd}$ that there exists a unique $A\in \cO$ such that $
    \textstyle \|A-\theta\|=\inf_{O\in \cO}\|O-\theta\|$.
                 \end{enumerate}
                 \end{athm}
                 \end{tcolorbox}
                 \begin{aproof}
                 Throughout this proof let $P\colon \R^{d\times\fd}\to\R$ satisfy for all $\theta\in\R^{d\times\fd}$ that
                 \begin{equation}\llabel{def: P}
                   \textstyle  P(\theta)=\inf_{O\in\cO} \|O-\theta\|
                 \end{equation}
                 and let $H\subseteq \R^{d\times\fd}\times \cO$ satisfy
                 \begin{equation}\llabel{def: H}
                     H=\{(\theta,O)\in\R^{d\times\fd}\times \cO\colon \|\theta-O\|=P(\theta)\}.
                 \end{equation}
                 \argument{\cref{lem: Pi property 2};}{that for all $A\in\cO$ it holds that 
        \begin{equation}\llabel{eq1}
            \|A\|^2=\min\{d,\fd\}\leq d.
        \end{equation}}
        \argument{\lref{eq1};}{that $\cO$ is \llabel{arg1} bounded\dott}
        \argument{\lref{arg1};the fact that $\cO$ is closed; the fact that $\cO\subseteq\R^{d\times\fd}$}{that \llabel{arg2} $\cO$ is compact\dott}
           In the following we prove that for all $\theta\in \R^{d\times\fd}$ there exists $A\in\cO$ such that
           \begin{equation}\llabel{need to prove}
               \|\theta-A\|=P(\theta)=\textstyle\inf_{O\in\cO}\|\theta-\cO\|.
           \end{equation}
           We prove \lref{need to prove} by contradiction.
           In the following we assume that there exists $\theta\in \R^{d\times\fd}$ which satisfies for all $A\in \cO$ that
            \begin{equation}\llabel{assume 1}
                 \textstyle \|A-\theta\|>\inf_{O\in \cO}\|O-\theta\|.
            \end{equation}
            \startnewargseq
            \argument{;}{ that there exists $A_n\in \cO$, $n\in\N$, which satisfy
            \begin{equation}\llabel{def: A_n}
                \textstyle \lim_{n\to\infty }\|A_n-\theta\|=P(\theta)=\inf_{O\in \cO}\|O-\theta\|.
            \end{equation}}
            \startnewargseq
        \argument{the fact that $\cO$ is compact;}{that there exist $(n_k)_{k\in\N}\subseteq \N$ and $B\in\cO$ which satisfy
        \begin{equation}\llabel{eq2}
            \textstyle \limsup_{k\to\infty}\|A_{n_k}-B\|=0.
        \end{equation}}
        \startnewargseq
        \argument{\lref{assume 1};\lref{def: A_n};\lref{eq2};}{that
        \begin{equation}
           \textstyle \|B-\theta\|>\inf_{O\in \cO}\|O-\theta\|= \lim_{n\to\infty }\|A_n-\theta\|=\lim_{k\to\infty }\|A_{n_k}-\theta\|=\|B-\theta\|.
        \end{equation}}
        This contradiction shows \lref{need to prove}\dott
        \startnewargseq
        \argument{\lref{def: P};\lref{need to prove};}{that for all $\theta,\vartheta\in \R^{d\times\fd}$ there exist $A,B\in \cO$ such that
        \begin{equation}\llabel{eqt1}
            P(\theta)\leq \|\theta-B\| \leq \|\theta-\vartheta\|+\|\vartheta-B\|=P(\vartheta)+\|\theta-\vartheta\|.
        \end{equation}}
         \argument{\lref{def: P};\lref{need to prove};}{that for all $\theta,\vartheta\in \R^{d\times\fd}$ there exist $A,B\in \cO$ such that
        \begin{equation}\llabel{eqt2}
            P(\vartheta)\leq \|\vartheta-A\| \leq \|\vartheta-\theta\|+\|\theta-A\|=P(\theta)+\|\theta-\vartheta\|.
        \end{equation}}
        \argument{\lref{eqt2};\lref{eqt1};}{for all $\theta,\vartheta\in \R^{d\times\fd}$ that
        \begin{equation}\llabel{eqt3}
            |P(\theta)-P(\vartheta)|\leq \|\theta-\vartheta\|.
        \end{equation}}
        \argument{\lref{eqt3};}{that $P$ is \llabel{argt1} continuous\dott}
        \argument{\lref{argt1};\lref{def: H};}{that $H$ is \llabel{argt2} closed\dott}
        \argument{\lref{argt2};}{that \llabel{argt3} $H$ is measurable\dott}
       \argument{\lref{argt3};the fact that for all $\theta\in\R^{d\times\fd}$ it holds that $\{A\in\cO\colon\|\theta-A\|=\inf_{O\in\cO} \|\theta-A\|\}$ is compact;the Arsenin-Kunugui theorem (see, \eg, \cite[Theorem 18.18]{MR1321597})}{that there exists a measurable $C\subseteq H$ which satisfies for all $\theta\in\R^{d\times\fd}$ that 
       \begin{equation}\llabel{itemm1}
           \exists\, A\in\cO\colon (\theta,A)\in H \qquad \text{if and only if}  \qquad\exists!\, A\in\cO\colon (\theta,A)\in C.
       \end{equation}}
       \startnewargseq
       \argument{\lref{def: H};\lref{need to prove};}{that \llabel{arg3} for all $\theta\in \R^{d\times\fd}$ there exists $A\in\cO$ such that $(\theta,A)\in H$\dott}
       \argument{\lref{arg3};\lref{itemm1}}{that for all $\theta\in \R^{d\times\fd}$ there exists an unique $K_\theta\in\cO$ which satisfies
       \begin{equation}\llabel{eq6}
           (\theta,K_\theta)\in C.
       \end{equation}}
       \startnewargseq
       In the following let $\Phi\colon \R^{d\times\fd}\to\cO$ satisfy for all $\theta\in \R^{d\times\fd}$ that
       \begin{equation}\llabel{def: Phi}
           \Phi(\theta)=K_\theta.
       \end{equation}
       \argument{\lref{eq6};\lref{def: Phi};the fact that $C\subseteq H$;}{that for all $\theta\in \R^{d\times\fd}$ it holds that
       \begin{equation}\llabel{eq7}
           (\theta,\Phi(\theta))=(\theta,K_\theta)\in C \subseteq H.
       \end{equation}}
       \argument{\lref{eq7};\lref{def: H}}{for all $\theta\in\R^{d\times\fd}$ that
       \begin{equation}\llabel{eq8}
           \|\theta-\Phi(\theta)\|=\textstyle\inf_{O\in \cO}\|\theta-O\|.
       \end{equation}}
       \argument{\lref{eq6};\lref{eq7};}{that
       \begin{equation}\llabel{eq9}
           C=\{(\theta,\Phi(\theta))\colon \theta\in \R^{d\times\fd}\}.
       \end{equation}}
       \argument{\lref{eq9};}{that for every $\theta\in \R^{d\times\fd}$ and every measurable $B\subseteq\cO$ it holds that 
       \begin{equation}\llabel{arg7}
           \{O\in\cO\colon (\theta,O)\in C\cap (\R^{d\times\fd}\times B)\}=\begin{cases}
               \{\Phi(\theta)\}&\colon \Phi(\theta)\in B\\
               \emptyset &\colon \Phi(\theta)\notin B.
           \end{cases} 
       \end{equation}}
       \argument{\lref{arg7};}{that for every $\theta\in \R^{d\times\fd}$ and every measurable $B\subseteq\cO$ it holds that \llabel{arg8} $\{O\in\cO\colon (\theta,O)\in C\cap (\R^{d\times\fd}\times B)\}$ is compact\dott}
       \argument{\lref{arg8}; the fact that for all measurable $B\subseteq \cO$ it holds that $C\cap (\R^{d\times\fd}\times B)$ is measurable; the Arsenin-Kunugui theorem (see, \eg, \cite[Theorem 18.18]{MR1321597})}{that for all mearsurable $B\subseteq \cO$ it holds that
       \begin{equation}\llabel{arg9}
          \{\theta\in \R^{d\times\fd}\colon (\exists\, O\in \cO\colon (\theta,O)\in C\cap(\R^{d\times\fd}\times B))\} =\{\theta\in \R^{d\times\fd}\colon \Phi(\theta)\in B\}
       \end{equation}
       is measurable\dott}
       \argument{\lref{arg9};}{that \llabel{arg4} $\Phi$ is measurable\dott}
       \argument{\lref{arg4};\lref{eq8}}{\cref{item 1: measuability}\dott}
       \startnewargseq
        In the following we prove \cref{item 2: measuability} by contradiction. In the following we assume that there exist $A_1,A_2\in \cO$ and a full rank matrix $\theta\in \R^{d\times\fd}$ which satisfy that $A_1\neq A_2$ and
        \begin{equation}\llabel{assume 2}
            \textstyle \|A_1-\theta\|=\|A_2-\theta\|=\inf_{O\in \cO}\|O-\theta\|.
        \end{equation}
        \startnewargseq
        \argument{\unskip, \eg, \cite[Theorem 7.3.5]{horn_johnson_2012};}{that there exist $U\in \R^{d\times d}$, $V\in \R^{\fd\times\fd}$, $\Sigma=(\sigma_{i,j})_{(i,j)\in \{1,2,\dots,d\}\times \{1,2,\dots,\fd\}}\in \R^{d\times \fd}$ which satisfy that
\begin{enumerate}[label=(\Roman*)]
    \item \llabel{item 1} it holds that $UU^{\top}=U^{\top}U=I_d$ and $VV^{\top}=V^{\top}V=I_\fd$,
    \item \llabel{item 2} it holds for all $i\in \{1,2,\dots,d\}$, $j\in \{1,2,\dots,\fd\}$ with $i\neq j$ that $\sigma_{i,j}=0$,
    \item \llabel{item 3} it holds for all $i\in \{1,2,\dots, \min\{d,\fd\}\}$ that $\sigma_{i,i}\geq 0$, and
    \item \llabel{item 4} it holds that $\theta=U\Sigma V^{\top}$.
\end{enumerate}}
\startnewargseq
In the following let $\Delta=(\Delta_{i,j})_{(i,j)\in \{1,2,\dots,d\}\times \{1,2,\dots,\fd\}}\in \R^{d\times \fd}$ satisfy for all $i,j\in \{1,2,\dots,d\}$, $j=\{1,2,\dots,\fd\}$ that
\begin{equation}\llabel{def: Delta}
    \delta_{i,j}=\mathbbm 1_{\{i=j\}}
\end{equation}
and for every $k\in\{1,2\}$ let $W^k=(w^k_{i,j})_{(i,j)\in \{1,2,\dots,d\}\times \{1,2,\dots,\fd\}}\in \R^{d\times \fd}$ satisfy
\begin{equation}\llabel{def: W}
    W^k=U^\top A_k V.
\end{equation}
\startnewargseq
\argument{\lref{assume 2};}{that for all $k\in\{1,2\}$, $O\in\cO$ it holds that
\begin{equation}\llabel{eqtt2}
    \|A_k-\theta\|^2\leq \|O-\theta\|^2.
\end{equation}}
\argument{\lref{eqtt2};\cref{lem: Pi property 2}; the fact that for all $k\in\{1,2\}$ it holds that $A_k\in \cO$}{that for all $k\in\{1,2\}$, $O\in\cO$ it holds that
\begin{equation}\llabel{eqtt3}
\begin{split}
   \min\{d,\fd\}-2\spro{A_k,\theta}+\|\theta\|^2&= \|A_k\|^2-2\spro{A_k,\theta}+\|\theta\|^2\leq \|O\|^2-2\spro{O,\theta}+\|\theta\|^2\\
   &=\min\{d,\fd\}-2\spro{O,\theta}+\|\theta\|^2.
   \end{split}
\end{equation}}
\argument{\lref{eqtt3};\cref{lem: Pi property 1}}{for all $k\in\{1,2\}$, $O\in\cO$ it holds that
\begin{equation}\llabel{eqtt4}
   \Tr(O^{\top}\theta)= \spro{O,\theta}\leq \spro{A_k,\theta}=\Tr((A_k)^{\top}\theta).
\end{equation}}
\argument{\lref{eqtt4};\lref{item 1};\lref{item 4};the fact that $U\Delta V^{\top}\in \cO$; the fact that for all $A\in \R^{\fd\times d}$, $B\in \R^{d\times\fd}$ it holds that $\Tr(AB)=\Tr(BA)$}{that for all $k\in\{1,2\}$ it holds that
\begin{equation}\llabel{eqtt5}
    \begin{split}
    \Tr(\Delta^\top\Sigma)&=\Tr(V^\top V\Delta^\top \Sigma)   =\Tr(V\Delta^\top \Sigma V^\top)= \Tr((U\Delta V^\top)^\top U\Sigma V^{\top})\leq \Tr((A_k)^\top U\Sigma V^\top)\\
    &=\Tr(V^\top(A_k)^\top U\Sigma).
    \end{split}
\end{equation}}
\argument{\lref{def: W};the fact that for all $k\in\{1,2\}$ it holds that $W^k\in\cO$;}{that for all $i\in\{1,2,\dots,\min\{d,\fd\}\}$, $k\in\{1,2\}$ it holds that
\begin{equation}\llabel{eqtt6.5}
    \textstyle(w_{i,i}^k)^2\leq \min\bigl\{\sum_{j=1}^d(w_{j,i}^k)^2,\sum_{j=1}^\fd(w_{i,j}^k)^2\bigr\}\leq 1.
\end{equation}}
\argument{\lref{eqtt6.5};\lref{def: Delta};\lref{def: W};\lref{eqtt5};\lref{item 2};\lref{item 3};}{for all $k\in\{1,2\}$ that
\begin{equation}\llabel{eqtt6}
\begin{split}
  \textstyle\sum_{i=1}^{\min\{d,\fd\}}\sigma_{i,i}&= \textstyle \sum_{i=1}^{\min\{d,\fd\}}\delta_{i,i}\sigma_{i,i}=\spro{\Delta,\Sigma}=\Tr(\Delta^\top\Sigma)\leq \Tr((U^\top A_k V)^\top\Sigma)\\
  &\textstyle=\spro{U^\top A_k V ,\Sigma}=\sum_{i=1}^{\min\{d,\fd\}}w_{i,i}^k\sigma_{i,i}\leq \sum_{i=1}^{\min\{d,\fd\}}\sigma_{i,i}.
   \end{split}
\end{equation}}
\argument{the fact that $\theta$ is full rank;}{that for all $i\in\{1,2,\dots,\min\{d,\fd\}\}$ it holds that \llabel{eqtt1}
    $\sigma_{i,i}>0$\dott}
\argument{\lref{eqtt1};\lref{eqtt6};}{for all $i\in\{1,2,\dots,\min\{d,\fd\}\}$, $k\in\{1,2\}$ that
\begin{equation}\llabel{eqtt8}
    w_{i,i}^k=1.
\end{equation}}
\argument{\lref{def: W};the fact that for all $k\in \{1,2\}$ it holds that $U^\top A_k V\in\cO$;\cref{lem: Pi property 2};}{that for all $k\in\{1,2\}$ it holds that
\begin{equation}\llabel{eqtt7}
    \textstyle\sum_{i=1}^{\min\{d,\fd\}}(w_{i,i}^k)^2\leq \sum_{i=1}^d\sum_{j=1}^\fd (w_{i,j}^k)^2=\|W^k\|=\|U^\top  A_k V\|^2=\min\{d,\fd\}.
\end{equation}}
\argument{\lref{eqtt7};\lref{eqtt8};}{for all $i\in\{1,2,\dots,d\}$, $j\in\{1,2,\dots,\fd\}$, $k\in\{1,2\}$ with $i\neq j$ that \llabel{eqtt9} $w_{i,j}^k=0$\dott}
\argument{\lref{eqtt9};\lref{def: Delta};\lref{def: W};\lref{eqtt8}}{for all $k\in\{1,2\}$ that
\begin{equation}\llabel{eqtt10}
    U^\top A_kV=W^k=\Delta.
\end{equation}}
\argument{\lref{eqtt10};the assumption that $A_1\neq A_2$}{that
\begin{equation}\llabel{eqtt12}
    A_1\neq A_2=U\Delta V^\top=A_1.
\end{equation}}
This contradiction proves \cref{item 2: measuability}\dott
                 \end{aproof}
\subsection{Error estimates for MUON under boundedness assumptions}\label{subsec: muon under boundness}
 \begin{tcolorbox}[colback=white!95!gray,
                  colframe=black,
                  boxrule=0.5pt,
                  sharp corners,
                  enhanced,
                  breakable,
                 ]
\begin{athm}{lemma}{NS1}
 Let $K\in \N_0$, let $\nscoe=(\nscoe_{i,j})_{(i,j)\in (\N_0)^2}\colon(\N_0)^2\to \R$ satisfy $\min\{\nscoe_{0,0},\nscoe_{0,1}\}>0$ and $\#(\nscoe^{ - 1 }( \R\backslash\{0\} )) <\infty$, let $P_i\colon\R\to (0,\infty)$, $i\in\N$, satisfy for all $i\in\N\cap[0,K]$, $x\in\R$ that $P_i(x)=\sum_{j=0}^\infty \nscoe_{i,j}x^{2j}$, and let $x\colon \N_0\to\R$, $y\colon\N_0\to\R$, and $z\colon\N_0\to\R$ satisfy for all $n\in \N$ that 
 \begin{equation}\llabel{def: x}
    0\leq x_0\leq (\nscoe_{0,0})^{-1}=y_0,\qquad  x_n=x_{n-1}P_{n}(x_{n-1}),
 \end{equation}
 \begin{equation}\llabel{def: y}
      y_n=\textstyle|y_{n-1}|\bigl[\sum_{i=0}^\infty |\nscoe_{n,i}| (y_{n-1})^{2i}\bigr],
 \end{equation}
  \begin{equation}\llabel{def: z}
    \textstyle z_0=1,\qquad z_n=|z_{n-1}|\bigl[\sum_{i=0}^\infty|\nscoe_{n,i}| (y_{n-1})^{2i}\bigr].
 \end{equation}
 Then it holds for all $n\in \{0,1,\dots,K\}$ that
 \begin{equation}\llabel{conclude}
      \textstyle  \bigl[\prod_{k=1}^n\min\{\inf_{x\in \R}P_k(x),1\}\bigr]x_0\leq x_n\leq  \min\{y_n,z_n x_0\}.
 \end{equation}
\end{athm}
\end{tcolorbox}
\begin{aproof}
  \argument{the fact that for all $k\in\N\cap[0,K]$, $x\in\R$ it holds that $P_k(x)> 0$ ;the fact that for all $k\in\N\cap[0,K]$, $x\in\R$ it holds that $P_k(x)$ is a polynomial}{that for all $n\in\N\cap[0,K]$ it holds that
  \begin{equation}\llabel{eq0}
      \textstyle P_n(x_{n-1})\geq \inf_{x\in\R}P_n(x)>0.
  \end{equation}}
  \argument{\lref{eq0};\lref{def: x};induction;}{that for all $n\in \{0,1,\dots,K\}$ it holds that
  \begin{equation}\llabel{eq0.2}
      x_n\geq 0.
  \end{equation}}
  \argument{\lref{def: x};\lref{eq0.2};\lref{eq0};}{that for all $n\in \N\cap[0,K]$ it holds that
  \begin{equation}\llabel{eq1}
     \textstyle x_n\geq x_{n-1} \bigl( \inf_{x\in\R}P_n(x)\bigr).
  \end{equation}}
  We prove \lref{conclude} by induction on $n\in \{0,1,\dots,K\}$. For the base case $n=0$ note that the fact that $x_0\leq y_0$ and the fact that $z_0=1$ implies that
  \begin{equation}\llabel{base: eq1}
    \textstyle  \bigl[\prod_{k=1}^0\min\{\inf_{x\in \R}P_k(x),1\}\bigr]x_0=x_0\leq \min\{(\nscoe_{0,0})^{-1},x_0\}= \min\{y_0,z_0x_0\}.
  \end{equation}
  This ensures \lref{conclude} in the base case $n=0$. For the induction step we assume that there exists $n\in\N\cap (0,K]$ which satisfies
  \begin{equation}\llabel{ntp: assume}
       \textstyle   \bigl[\prod_{k=1}^{n-1}\min\{\inf_{x\in \R}P_k(x),1\}\bigr]x_0\leq x_{n-1}\leq  \min\{y_{n-1},z_{n-1}x_0\}.
  \end{equation}
  \startnewargseq
  \argument{\lref{def: x};\lref{def: y};\lref{eq0.2};\lref{ntp: assume}}{that
  \begin{equation}\llabel{eq7}
      x_n=x_{n-1}P_n(x_{n-1})=x_{n-1}\bigl[\textstyle\sum_{i=0}^\infty\nscoe_{n,i}(x_{n-1})^{2i}\bigr]\leq|y_{n-1}|\bigl[\sum_{i=0}^\infty|\nscoe_{n,i}|(y_{n-1})^{2i}\bigr] = y_n.
  \end{equation}}
  \argument{\lref{def: x};\lref{def: z};\lref{ntp: assume};}{that
  \begin{equation}\llabel{eq7.5}
  \begin{split}
      x_n&=x_{n-1}P_n(x_{n-1})=x_{n-1}\textstyle\bigl(\sum_{i=0}^\infty \nscoe_{n,i} (x_{n-1})^{2i}\bigr)\\
      &\leq\textstyle x_0z_{n-1}\bigl(\sum_{i=0}^\infty |\nscoe_{n,i}| (y_{n-1})^{2i}\bigr)=z_nx_0.
      \end{split}
  \end{equation}}
  \argument{\lref{def: x};\lref{eq0.2};\lref{eq1};}{that
  \begin{equation}\llabel{eq8}
  \begin{split}
      \textstyle x_n&=\textstyle x_{n-1} P_n(x_{n-1})\geq \bigl[\inf_{x\in\R} P_n(x)\bigr]  \bigl[\prod_{k=1}^{n-1}\min\{\inf_{x\in \R}P_k(x),1\}\bigr]x_0\\
      &\geq  \textstyle   \bigl[\prod_{k=1}^n\min\{\inf_{x\in \R}P_k(x),1\}\bigr]x_0.
      \end{split}
  \end{equation}}
  \argument{\lref{eq8};\lref{eq7};\lref{eq7.5}}{that
  \begin{equation}\llabel{eq9}
       \textstyle   \bigl[\prod_{k=1}^n\min\{\inf_{x\in \R}P_k(x),1\}\bigr]x_0\leq x_n\leq  \min\{y_n,z_nx_0\}.
  \end{equation}}
  \argument{\lref{eq9};\lref{ntp: assume};induction}{\lref{conclude}\dott}
    \end{aproof}
 \begin{tcolorbox}[colback=white!95!gray,
                  colframe=black,
                  boxrule=0.5pt,
                  sharp corners,
                  enhanced,
                  breakable,
                 ]
\begin{athm}{lemma}{NS error}[\textcolor{red}{Boundedness of the \NS\ functions}]
      Let $d,\fd\in \N$, $K\in \N_0$, $A \in \R^{d \times \fd}$, let $\nscoe=(\nscoe_{i,j})_{(i,j)\in (\N_0)^2}\colon(\N_0)^2\to \R$ satisfy $\min\{\nscoe_{0,0},\nscoe_{0,1}\}>0$ and $\#(\nscoe^{ - 1 }( \R\backslash\{0\} )) <\infty$, assume for all $i\in\N\cap[0,K]$, $x\in\R$ that $\sum_{j=0}^\infty \nscoe_{i,j}x^{2j}>0$, and let $y\colon\N_0\to\R$ and $z\colon\N_0\to\R$ satisfy for all $n\in \N$ that 
 \begin{equation}\llabel{def: y}
  y_0= (\nscoe_{0,0})^{-1},\qquad   y_n=\textstyle|y_{n-1}|\bigl[\sum_{i=0}^\infty |\nscoe_{n,i}| (y_{n-1})^{2i}\bigr],
 \end{equation}
  \begin{equation}\llabel{def: z}
    \textstyle z_0=1,\qqandqq z_n=|z_{n-1}|\bigl[\sum_{i=0}^\infty|\nscoe_{n,i}| (y_{n-1})^{2i}\bigr].
 \end{equation}
     Then it holds for all $n\in \{0,1,\dots,K\}$ that
     \begin{equation}\llabel{conclude}
         \lVert\TNewtonSchulzAlgorithm{\nscoe}{n}{A}\rVert\leq \min\{y_n[\min\{d,\fd\}]^{1/2},z_n(\nscoe_{0,1})^{-1}\|A\|\}
         \end{equation}
         \begin{equation}\llabel{conclude 2} \text{and}\textstyle\qquad\spro{\TNewtonSchulzAlgorithm{\nscoe}{n}{A},A}\geq   \bigl[\prod_{k=1}^n\min\{\inf_{x\in \R}(\sum_{j=0}^\infty \nscoe_{k,j}x^{2j}),1\}\bigr]\|A\|^2(\nscoe_{0,0}\|A\|+\nscoe_{0,1})^{-1}
     \end{equation}
     (cf.\ \cref{definition: NS}).
\end{athm}
\end{tcolorbox}
\begin{aproof}
    \argument{\unskip, \eg, \cite[Theorem 7.3.5]{horn_johnson_2012};}{that there exist $U\in \R^{d\times d}$, $V\in \R^{\fd\times\fd}$, $\Sigma=(\sigma_{i,j})_{(i,j)\in \{1,2,\dots,d\}\times\{1,2,\dots,\fd\}}\in \R^{d\times \fd}$ which satisfy that
\begin{enumerate}[label=(\Roman*)]
    \item \llabel{item 1} it holds that $UU^{\top}=U^{\top}U=I_d$ and $VV^{\top}=V^{\top}V=I_\fd$,
    \item \llabel{item 2} it holds for all $i\in \{1,2,\dots,d\}$, $j\in \{1,2,\dots,\fd\}$ with $i\neq j$ that $\sigma_{i,j}=0$,
    \item \llabel{item 3} it holds for all $i\in \{1,2,\dots, \min\{d,\fd\}\}$ that $\sigma_{i,i}\geq 0$, and
    \item \llabel{item 4} it holds that $A=U\Sigma V^{\top}$.
\end{enumerate}}
\startnewargseq
\argument{\lref{item 1};\lref{item 4};\cref{lem: Pi property 1}; the fact that for all $X,Y\in\R^\fd$ it holds that $\Tr(XY)=\Tr(YX)$}{that
\begin{equation}\llabel{eq1}
\begin{split}
    \|A\|^2&=\Tr(A^{\top}A)=\Tr((U\Sigma V^{\top})^{\top}U\Sigma V^{\top})=\Tr(V\Sigma^{\top}U^{\top}U\Sigma V^{\top})=\Tr(V\Sigma^{\top}\Sigma V^{\top})\\
    &=\Tr(V^{\top}V\Sigma^{\top}\Sigma )=\Tr(\Sigma^{\top}\Sigma).
    \end{split}
\end{equation}}
\argument{\lref{eq1};\lref{item 2};}{that
\begin{equation}\llabel{eq2}
    \|A\|^2=\sum_{i=1}^{\min\{d,\fd\}}(\sigma_{i,i})^2.
\end{equation}}
\startnewargseq
In the following let $\Delta_n=(\delta_{i,j}^n)_{(i,j)\in \{1,2,\dots,d\}\times \{1,2,\dots,\fd\}}\in \R^{d\times\fd}$, $n\in \N_0$, satisfy for all $n\in \N$, $i\in \{1,2,\dots,d\}$, $j\in \{1,2,\dots,\fd\}$ that
\begin{equation}\llabel{def: Delta}
    \textstyle \delta_{i,j}^0=\sigma_{i,j}(\nscoe_{0,0}\|A\|+\nscoe_{0,1})^{-1}\mathbbm 1_{\{i=j\}}\qqandqq \delta_{i,j}^n= \bigl[\sum_{k=0}^\infty\nscoe_{n,k}(\delta_{i,j}^{n-1})^{2k+1}\bigr]\mathbbm 1_{\{i=j\}}.
\end{equation}
In the following we prove that for all $n\in \{0,1,\dots,K\}$ it holds that
\begin{equation}\llabel{need to prove}
    \TNewtonSchulzAlgorithm{\nscoe}{n}{A}=U\Delta_n V^{\top}.
\end{equation}
We prove \lref{need to prove} by induction on $n\in \{0,1,\dots,K\}$. For the base case $n=0$ note that \lref{def: Delta} and  \lref{item 2} show that 
\begin{equation}
    \Delta_0=(\nscoe_{0,0}\|A\|+\nscoe_{0,1})^{-1}\Sigma.
\end{equation}
Combining this, \cref{T_B_D}, and \lref{item 4} implies that
\begin{equation}
\begin{split}
     \TNewtonSchulzAlgorithm{\nscoe}{0}{A}&=(\nscoe_{0,0}\|A\|+\nscoe_{0,1})^{-1}A=(\nscoe_{0,0}\|A\|+\nscoe_{0,1})^{-1}U\Sigma V^{\top}\\
&=U\bigl((\nscoe_{0,0}\|A\|+\nscoe_{0,1})^{-1}\Sigma\bigr) V^{\top}=U\Delta_0 V^{\top}.
    \end{split}
\end{equation}
This ensures \lref{need to prove} in the base case $n=0$. For the induction step we assume that there exists $n\in \N\cap(0,K]$ which satisfies that
\begin{equation}\llabel{ntp: assume}
      \TNewtonSchulzAlgorithm{\nscoe}{n-1}{A}=U\Delta_{n-1} V^{\top}.
\end{equation}
\startnewargseq
\argument{\lref{def: Delta};}{that 
\begin{equation}\llabel{ntp: eq0}
  \textstyle \sum_{i=0}^\infty \nscoe_{n,i}(\Delta_{n-1}(\Delta_{n-1})^{\top})^{i}\Delta_{n-1}=\Delta_n.
\end{equation}}
\argument{\lref{ntp: eq0};\cref{T_B_D};\lref{ntp: assume};\lref{item 1}}{that
\begin{equation}\llabel{ntp: eq1}
\begin{split}
     \TNewtonSchulzAlgorithm{\nscoe}{n}{A}\textstyle&=\textstyle\sum_{i=0}^\infty\nscoe_{n,i} (U\Delta_{n-1} V^{\top}(U\Delta_{n-1} V^{\top})^{\top})^iU\Delta_{n-1} V^{\top}\\
    &\textstyle=\sum_{i=0}^\infty\nscoe_{n,i} (U\Delta_{n-1} (\Delta_{n-1}) ^{\top}U^{\top})^iU\Delta_{n-1} V^{\top}\\
    &\textstyle =\sum_{i=0}^\infty \nscoe_{n,i} U(\Delta_{n-1}(\Delta_{n-1})^{\top})^{i} \Delta_{n-1}V^{\top}\\
    &=\textstyle U\bigl(\sum_{i=0}^\infty\nscoe_{n,i}(\Delta_{n-1}(\Delta_{n-1})^{\top})^{i} \Delta_{n-1}\bigr) V^{\top}\\
    &=U\Delta_n V^{\top}.
    \end{split}
\end{equation}}
\argument{\lref{ntp: eq1};\lref{ntp: assume};induction}{\lref{need to prove}\dott}
\startnewargseq
\argument{\lref{eq2};\lref{item 2}}{that for all $i\in \{1,2,\dots,d\}$, $j\in \{1,2,\dots,\fd\}$ with $i=j$ it holds that
\begin{equation}\llabel{eqt1}
    \sigma_{i,j}\leq \|A\|<(\nscoe_{0,0})^{-1}(\nscoe_{0,0}\|A\|+\nscoe_{0,1}).
\end{equation}}
\argument{\lref{def: Delta};\lref{item 2};\lref{item 3}}{that for all $i\in \{1,2,\dots,d\}$, $j\in \{1,2,\dots,\fd\}$ it holds that
\begin{equation}\llabel{eq3}
    0\leq \delta_{i,j}^0\leq (\nscoe_{0,0})^{-1}.
\end{equation}}
\argument{\lref{eq3};\lref{def: Delta};\cref{NS1};}{for all $i\in \{1,2,\dots,d\}$, $j\in \{1,2,\dots,\fd\}$, $n\in \{0,1,\dots,K\}$ that
\begin{equation}\llabel{eq4}
   \textstyle\bigl[\prod_{k=1}^n\min\{\inf_{x\in \R}(\sum_{j=0}^\infty \nscoe_{k,j}x^{2j}),1\}\bigr]\delta_{i,j}^0\leq \delta_{i,j}^n \leq \min\{y_n,z_n \delta_{i,j}^0\}.
\end{equation}}
\argument{\lref{def: Delta};\lref{eq4};}{that for all $n\in \{0,1,\dots,K\}$ it holds that
\begin{equation}\llabel{eq8'}
\begin{split}
  \spro{\Delta_n,\Delta_n}&=\sum_{i=1}^{\min\{d,\fd\}} (\delta_{i,i}^n)^2\leq \textstyle \min\biggl\{(y_n)^2\min\{d,\fd\},\sum\limits_{i=1}^{\min\{d,\fd\}}(z_n\delta^0_{i,i})^2\biggr\}\\
  &\leq   \textstyle\min\biggl\{(y_n)^2\min\{d,\fd\},(z_n)^2(\nscoe_{0,1})^{-2}\sum\limits_{i=1}^{\min\{d,\fd\}}(\sigma_{i,i})^2\biggr\}.
  \end{split}
\end{equation}}
\argument{\lref{eq8'};\lref{eq2}}{for all $n\in\{0,1,\dots,K\}$ that
\begin{equation}\llabel{eq8}
\begin{split}
    \spro{\Delta_n,\Delta_n}\leq \min\{(y_n)^2\min\{d,\fd\},(z_n)^2(\nscoe_{0,1})^{-2}\|A\|^2\}.
    \end{split}
\end{equation}}
\argument{\lref{eq8};\lref{item 1};\cref{lem: Pi property 1};the fact that for all $X,Y\in\R^\fd$ it holds that $\Tr(XY)=\Tr(YX)$}{that for all $n\in \{0,1,\dots,K\}$ it holds that
\begin{equation}\llabel{evidence 2}
\begin{split}
    \|U\Delta_n V^{\top}\|^2&=\Tr((U\Delta_n V^{\top})^{\top}U\Delta_nV^{\top})=\Tr(V(\Delta_n)^{\top}\Delta_n V^{\top})=\Tr(V^{\top}V(\Delta_n)^{\top}\Delta_n)\\
    &=\Tr((\Delta_n)^{\top}\Delta_n)=\spro{\Delta_n,\Delta_n}\leq \min\{(y_n)^2\min\{d,\fd\},(z_n)^2(\nscoe_{0,1})^{-2}\|A\|^2\}.
    \end{split}
\end{equation}}
\argument{\lref{eq2};\lref{def: Delta};\lref{item 2};\lref{item 3};\lref{eq4};}{that for all $n\in \{0,1,\dots,K\}$ it holds that
\begin{equation}\llabel{eq9}
\begin{split}
    \textstyle \spro{\Sigma,\Delta_n}\textstyle&=\textstyle\sum\limits_{i=1}^{\min \{d,\fd\}}\sigma_{i,i}\delta_{i,i}^n\geq  \sum\limits_{i=1}^{\min \{d,\fd\}}\sigma_{i,i}\bigl[\prod_{k=1}^n\min\{\inf_{x\in \R}(\sum_{j=0}^\infty \nscoe_{k,j}x^{2j}),1\}\bigr]\delta_{i,i}^0\\
    &\textstyle=\sum\limits_{i=1}^{\min \{d,\fd\}}\sigma_{i,i}\bigl[\prod_{k=1}^n\min\{\inf_{x\in \R}(\sum_{j=0}^\infty \nscoe_{k,j}x^{2j}),1\}\bigr]\sigma_{i,i}(\nscoe_{0,0}\|A\|+\nscoe_{0,1})^{-1}\\
    &\textstyle=\bigl[\prod_{k=1}^n\min\{\inf_{x\in \R}(\sum_{j=0}^\infty \nscoe_{k,j}x^{2j}),1\}\bigr]\|A\|^2(\nscoe_{0,0}\|A\|+\nscoe_{0,1})^{-1}.
    \end{split}
\end{equation}}
\argument{\lref{eq9};\lref{item 1};\lref{item 4};\cref{lem: Pi property 1};the fact that for all $X,Y\in\R^\fd$ it holds that $\Tr(XY)=\Tr(YX)$}{that for all $n\in \{0,1,\dots,K\}$ it holds that
\begin{equation}\llabel{evidence 3}
\begin{split}
    \spro{U\Delta_n V^{\top},A}&=\spro{U\Delta_n V^{\top},U\Sigma V^{\top}}=\Tr((U\Sigma V^{\top})^{\top}U\Delta_nV^{\top})=\Tr(V\Sigma^{\top}\Delta_n V^{\top})\\
    &=\Tr(V^{\top}V\Sigma^{\top}\Delta_n)
    =\textstyle \Tr(\Sigma^{\top}\Delta_n)=\spro{\Sigma,\Delta_n}\\
    &\textstyle\geq \bigl[\prod_{k=1}^n\min\{\inf_{x\in \R}(\sum_{j=0}^\infty \nscoe_{k,j}x^{2j}),1\}\bigr]\|A\|^2(\nscoe_{0,0}\|A\|+\nscoe_{0,1})^{-1}.
    \end{split}
\end{equation}}
\argument{\lref{need to prove};\lref{evidence 2};\lref{evidence 3}}{\lref{conclude,conclude 2}\dott}
\end{aproof}

   \begin{tcolorbox}[colback=white!95!gray,
                  colframe=black,
                  boxrule=0.5pt,
                  sharp corners,
                  enhanced,
                  breakable,
                 ]
                 \begin{definition}[\textcolor{red}{\MUON\ process}]\label{defifinition: MUON process}
                     Let $( \Omega, \cF,\P )$ be a probability space, let $M,\delta\in \N$, $d_1,d_2,\dots,d_\delta,\fd_1,\fd_2,\dots,\fd_\delta\in \N$, $\alpha,\chi\in \R$, let $\nscoe=(\nscoe_{i,j})_{(i,j)\in (\N_0)^2}\colon(\N_0)^2\to \R$ be a function, let $\cA=\times_{i=1}^\delta\R^{d_i\times\fd_i}$, let $(\setX,\cU)$ be a measurable space, let $\smalll=(\smalll(\theta,x))_{(\theta,x)\in \cA\times \setX} \allowbreak \colon \cA\times \setX\to\R$ be measurable, assume for all $x\in \setX$ that $(\cA\ni\theta\mapsto \smalll(\theta,x)\in \R)\in C^1(\cA,\R)$,
let $X_{ n,m } \colon \Omega \to \setX$, $(n,m) \in \N^2$, be random variables,
and let $\gamma\colon \N\to (0,\infty)$ be a function.
Then we say that $\Theta$ is the $\alpha$-$\gamma$-$\chi$-$\nscoe$-$M$-\MUON\ process for $\smalll$ on $( \Omega, \cF, \P )$ with respect to $(\setX,\cU)$ with data $( X_{ n,m } )_{ (n,m)\in \N^2 }$ (we say that $\Theta$ is the $\alpha$-$\gamma$-$\chi$-$\nscoe$-$M$-\MUON\ process for $\smalll$) if and only if there exist $K\in\N_0$ and $\bfm=(\bfm^{1},\dots,\bfm^{\delta})\colon \N_0\times\Omega\to\cA$ such that
\begin{enumerate}[label=(\roman*)]
\item it holds that $\Theta=(\Theta^1,\dots,\Theta^\delta)\colon \N_0\times\Omega\to\cA$ is a stochastic process,
    \item  it holds that $\Theta_{0}$ and $(X_{n,m})_{(n,m)\in \N^2}$ are independent, 
    \item \label{item 4: MUON process} it holds that $\gamma$ is non-increasing with $\limsup_{n\to\infty} (\gamma_n+(\gamma_n)^{-2}(\gamma_{n}-\gamma_{n+1}))=0$, 
    \item it holds that $\#(\nscoe^{ - 1 }( \R\backslash\{0\} )) <\infty$,
    \item it holds for all $i\in \N\cap [0,K]$, $x\in \R$ that $\sum_{j=0}^\infty \nscoe_{i,j}x^{2j}>0$,
    \item it holds that $\P( \| \Theta_0 \| \leq \chi ) =1+|\nscoe_{K+1,0}|<1+\min\{\nscoe_{0,0},\nscoe_{0,1}\}$, and
    \item it holds for all $n\in \N$, $i\in \{1,2,\dots,\delta\}$ that  
         \begin{equation}
            \textcolor{magenta}{\bfm_0=0},\qquad \textcolor{magenta}{\bfm_n= \alpha \bfm_{n-1}+(1-\alpha)\bigl[\textstyle \frac 1M \sum_{m=1}^M(\nabla_\theta\smalll)(\Theta_{n-1},X_{n,m})\bigr]},
            \end{equation}
            \begin{equation}\label{def: Theta: defifinition: MUON process}
            \text{and}\qquad  \textcolor{magenta}{\Theta_n^{i}=\Theta_{n-1}^{i}-\gamma_n\TNewtonSchulzAlgorithm{\nscoe}{K}{\bfm_n^{i}}}
         \end{equation}
\end{enumerate}
(cf.\ \cref{definition: NS}).
                 \end{definition}
                 \end{tcolorbox}

                 Note that the assumption that $\smalll$ is measurable, the assumption that for all $x\in \setX$ it holds that $(\cA\ni\theta\mapsto \smalll(\theta,x)\in \R)\in C^1(\cA,\R)$, and, \eg, \cite[Lemma 7.2.10]{ArBePhi2024} imply that $\cA\times\setX \ni (\theta,x) \mapsto ( \nabla_{ \theta } \smalll )( \theta, x ) \in \cA$ is measurable.

                 \begin{samepage}
                  \begin{tcolorbox}[colback=white!95!gray,
                  colframe=black,
                  boxrule=0.5pt,
                  sharp corners,
                  enhanced,
                  breakable,
                 ]
\begin{athm}{theorem}{muon bound convergence stochastic}[\textcolor{red}{\MUON\ error analysis under the assumption that \MUON\ stays bounded}]
 Let $(\Omega,\cF,\P)$ be a probability space, let $\delta\in \N$, $d_1,d_2,\dots,d_\delta,\fd_1,\fd_2,\dots,\fd_\delta\in \N$, $\alpha\in (0,1)$, let $\nscoe\colon(\N_0)^2\to \R$ be a function, let $(\setX,\cU)$ be a measurable space, let $X_{n,m}\colon \Omega\to\setX$, $(n,m)\in \N^2$, be \iid\ random variables, let $\cA=\times_{i=1}^\delta\R^{d_i\times\fd_i}$, let $\cK,p,\chi\in (0,\infty)$, let $\smalll=(\smalll(\theta,x))_{(\theta,x)\in \cA\times \setX}\colon \cA\times \setX\to\R$ be measurable, assume for all $x\in \setX$ that $(\cA\ni\theta\mapsto \smalll(\theta,x)\in \R)\in C^1(\cA,\R)$, assume for all $\theta,\vartheta\in \cA$, $x\in \setX$ with $\max\{\|\theta\|,\|\vartheta\|\}< \chi$ that
      \begin{equation}\llabel{assume: locally Lipschitz}
          \textcolor{magenta}{\|(\nabla_\theta\smalll)(\theta,x)-(\nabla_\theta\smalll)(\vartheta,x)\|\leq \cK\|\theta-\vartheta\|}\qqandqq \textcolor{magenta}{\|(\nabla_\theta\smalll)(\theta,x)\|\leq  \cK},
          \end{equation}
      let $\gamma\colon \N\to(0,\infty)$ and $\cL\colon \cA\to\R$ satisfy for all $\theta\in \cA$ that $\cL(\theta)=\E[\smalll(\theta,X_{1,1})]$, assume that $\cL$ is strongly convex, and let $\vartheta\in \cA$ satisfy $\|\vartheta\|< \chi$ and $\cL(\vartheta)=\inf_{\theta\in \cA}\cL(\theta)$.
     Then there exists $\fC\in \R$ such that for every $M,n\in \N$ and every \defmuon{\Theta}{\chi}{M}\ with $\inf_{m\in \N_0}\P(\|\Theta_m\|< \chi)=1$ it holds that
\begin{equation} \llabel{conclude}
     \textcolor{magenta}{\bigl(\E\bigl[\|\Theta_{n-1}-\vartheta\|^p\bigr]\bigr)^{1/p}\leq \fC (M^{-1/2}+\sqrt{\gamma_{n}})}\ifnocf.
\end{equation}
\cfout[.]
\end{athm}
\end{tcolorbox}
\end{samepage}
\begin{aproof}
Throughout this proof assume without loss of generality that there exist $N\in\N$ and an \defmuon{\Theta}{\chi}{N}.
\argument{the fact that there exist $N\in\N$ and an \defmuon{\Theta}{\chi}{N}}{that there exists an unique $K\in\N_0$ which satisfies that
\begin{enumerate}[label=(\roman*)]
\item \llabel{item 1} it holds that $\nscoe_{K+1,0}=0<\min\{\nscoe_{0,0},\nscoe_{0,1}\}$, 
\item \llabel{item 2} it holds that $\#(\nscoe^{ - 1 }( \R\backslash\{0\} )) <\infty$, and
    \item \llabel{item 3} it holds for all $n\in\N\cap[0,K]$, $x\in\R$ that
    $\sum_{i=0}^\infty \nscoe_{n,i}x^{2i}>0$
\end{enumerate}}
\startnewargseq
    \argument{\lref{item 1};\lref{item 2};\lref{item 3};\cref{NS error}}{that there exist $\scrc,\fC\in (0,\infty)$ such that for all $i\in \{1,2,\dots,\delta\}$, $A\in \R^{d_i\times\fd_i}$ it holds that
    \begin{equation}\llabel{eq1}
   \lVert\TNewtonSchulzAlgorithm{\nscoe}{K}{A}\rVert\leq \scrc\|A\| \qqandqq \spro{\TNewtonSchulzAlgorithm{\nscoe}{K}{A},A}\geq \frac{\fC\|A\|^2}{\|A\|+\nscoe_{0,1}(\nscoe_{0,0})^{-1}}.
    \end{equation}}
    \argument{\cref{T_B_D};induction}{that for all $i\in\{1,2,\dots,\delta\}$ it holds that \llabel{argt1} $\R^{d_i\times\fd_i}\ni A\mapsto\TNewtonSchulzAlgorithm{\nscoe}{K}{A}\in\R^{d_i\times\fd_i}$ is measurable\dott}
    \argument{\lref{argt1};\lref{eq1};\cref{theo: stochastic convergence}}{that there exists $\fC\in (0,\infty)$ such that for every $M\in \N$, every stochastic process $\bfm=(\bfm^{1},\dots,\bfm^{\delta})\colon \N_0\times\Omega\to\cA$, and every stochastic process $\Theta=(\Theta^1,\dots,\Theta^\delta)\colon \N_0\times\Omega\to\cA$ with the property that
    for all $n\in \N$, $i\in \{1,2,\dots,\delta\}$ it holds that
         \begin{equation}
           \bfm_0=0,\qquad \bfm_n= \alpha \bfm_{n-1}+(1-\alpha)\bigl[\textstyle \frac 1M \sum_{m=1}^M(\nabla_\theta\smalll)(\Theta_{n-1},X_{n,m})\bigr],
            \end{equation}
            \begin{equation}
           \P(\|\Theta_{n-1}\|< \chi)=1,\qquad \text{and}\qquad  \Theta_n^{i}=\Theta_{n-1}^{i}-\gamma_n\TNewtonSchulzAlgorithm{\nscoe}{K}{\bfm_n^{i}}
            \end{equation}
           and with the property that $\Theta_0$ and $(X_{n,m})_{(n,m)\in \N^2}$ are independent we have for all $n \in \N_0$ that
\begin{equation} \llabel{eq2}
     \bigl(\E\bigl[\|\Theta_n-\vartheta\|^p\bigr]\bigr)^{1/p}\leq \fC M^{-1/2}+\fC\sqrt{\gamma_{n+1}}.
\end{equation}}
   \argument{\lref{eq2};\lref{item 1};\lref{item 2};\lref{item 3}}{that there exists $\fC\in (0,\infty)$ such that for every $M\in \N$, $n\in \N_0$ and every \defmuon{\Theta}{\chi}{M}\ with $\inf_{m\in \N_0}\P(\|\Theta_m\|< \chi)=1$ we have that
   \begin{equation} \llabel{eq3}
     \bigl(\E\bigl[\|\Theta_n-\vartheta\|^p\bigr]\bigr)^{1/p}\leq \fC M^{-1/2}+\fC\sqrt{\gamma_{n+1}}.
\end{equation}}
  \argument{\lref{eq3}}{\lref{conclude}\dott}
\end{aproof}

In \cref{item 4: MUON process} in \cref{defifinition: MUON process} we employ the condition that the non-increasing learning rates $ \gamma_n \in (0,\infty) $, $ n \in \N $, satsfiy $\limsup_{n\to\infty} (\gamma_n+(\gamma_n)^{-2}(\gamma_{n}-\gamma_{n+1}))=0$. Loosely speaking, this condition ensures that the learning rates converge to zero but do not converge too quickly to zero (cf., \eg, \cite[Lemma 2.13 and Lemma 2.14]{DeDoArPhi2025}). In the following elementary fact, \cref{lem: equivalent learning rate}, we recall that this condition is equivalent to the condition that $\limsup_{ n \to \infty } ( \gamma_n + ( \gamma_n )^{ - 2 } ( \gamma_{ n - 1 } - \gamma_n ) ) = 0$. Only for completeness we include here a detailed proof for \cref{lem: equivalent learning rate}.
 \begin{tcolorbox}[colback=white!95!gray,
                  colframe=black,
                  boxrule=0.5pt,
                  sharp corners,
                  enhanced,
                  breakable,
                 ]
                 \begin{athm}{lemma}{lem: equivalent learning rate}
                     Let $(\gamma_n)_{n\in\N}\subseteq(0,\infty)$ be non-increasing, assume $\lim_{ n \to\infty } \gamma_n = 0$, and let $a \in \R$. Then the following two statements are equivalent:
                     \begin{enumerate}[label=(\roman*)]
                         \item \label{item 1: equivalent learning rate} It holds that $\limsup_{n\to\infty}|a-(\gamma_n)^{-2}(\gamma_{n-1}-\gamma_n)|=0$.
                           \item \label{item 2: equivalent learning rate} It holds that $\limsup_{n\to\infty} |a-(\gamma_n)^{-2}(\gamma_{n}-\gamma_{n+1})|=0$.
                     \end{enumerate}
                 \end{athm}
                 \end{tcolorbox}
                 \begin{aproof}
                     In the following we prove that (\ref{item 1: equivalent learning rate}$\rightarrow$\ref{item 2: equivalent learning rate}). For this we thus assume that
                     \begin{equation}\llabel{assume 1}
                         \textstyle\limsup_{n\to\infty}|a-(\gamma_n)^{-2}(\gamma_{n-1}-\gamma_n)|=0.
                     \end{equation}
                     \argument{\lref{assume 1}; }{that
                     \begin{equation}\llabel{as1: eq1}
                         \textstyle\lim_{n\to\infty}((\gamma_n)^{-2}(\gamma_{n-1}-\gamma_n))=a.
                         \end{equation}}
                      \argument{\lref{as1: eq1};the assumption that $\limsup_{n\to\infty}\gamma_n$; the fact that for all $n\in \N\backslash\{1\}$ it holds that $\gamma_{n-1}\geq \gamma_n$}{that
    \begin{equation}\llabel{as1: eq2}
        0\leq \textstyle\limsup_{n\to\infty} \bigl(\frac{\gamma_{n-1}-\gamma_n}{\gamma_{n}}\bigr)=0.
    \end{equation}}
    \argument{\lref{as1: eq2};}{that
    \begin{equation}\llabel{as1: eq3}
       \textstyle \lim_{n\to\infty} \bigl(\frac{\gamma_{n-1}-\gamma_n}{\gamma_{n}}\bigr)=0.
    \end{equation}}
    \argument{\lref{as1: eq3};}{that
    \begin{equation}\llabel{as1: eq4}
        \textstyle\lim_{n\to\infty} \bigl(\frac{\gamma_{n+1}}{\gamma_{n}}\bigr)=1.
    \end{equation}}
    \argument{\lref{as1: eq4};\lref{as1: eq1};}{that
    \begin{equation}\llabel{as1: eq6}
        \textstyle\lim_{n\to\infty} ((\gamma_{n})^{-2}(\gamma_{n}-\gamma_{n+1}))=\bigl[\lim_{n\to\infty} ((\gamma_{n})^{-2}(\gamma_{n+1})^2\bigr]\bigl[\lim_{n\to\infty} ((\gamma_{n+1})^{-2}(\gamma_{n}-\gamma_{n+1}))\bigr]=a.
    \end{equation}}
    \argument{\lref{as1: eq6};}{(\ref{item 1: equivalent learning rate}$\rightarrow$\ref{item 2: equivalent learning rate})\dott}
    \startnewargseq
       In the following we prove that (\ref{item 2: equivalent learning rate}$\rightarrow$\ref{item 1: equivalent learning rate}). For this we thus assume that
                     \begin{equation}\llabel{assume 2}
                        \textstyle\limsup_{n\to\infty} |a-(\gamma_n)^{-2}(\gamma_{n}-\gamma_{n+1})|=0.
                     \end{equation}
                     \argument{\lref{assume 2}; the fact that for all $n\in \N$ it holds that $(\gamma_n)^{-2}(\gamma_{n}-\gamma_{n+1})\geq 0$}{that
                     \begin{equation}\llabel{as2: eq1}
                         \lim_{n\to\infty} ((\gamma_n)^{-2}(\gamma_{n}-\gamma_{n+1}))=a.
                     \end{equation}}
                      \argument{\lref{as2: eq1};the assumption that $\limsup_{n\to\infty}\gamma_n=0$; the fact that for all $n\in \N$ it holds that $\gamma_{n}\geq \gamma_{n+1}$}{that
    \begin{equation}\llabel{as2: eq2}
        0\leq \textstyle\limsup_{n\to\infty} \bigl(\frac{\gamma_{n}-\gamma_{n+1}}{\gamma_{n}}\bigr)=0.
    \end{equation}}
    \argument{\lref{as2: eq2};}{that
    \begin{equation}\llabel{as2: eq3}
       \textstyle \lim_{n\to\infty} \bigl(\frac{\gamma_{n}-\gamma_{n+1}}{\gamma_{n}}\bigr)=0.
    \end{equation}}
    \argument{\lref{as2: eq3};}{that
    \begin{equation}\llabel{as2: eq4}
        \textstyle\lim_{n\to\infty} \bigl(\frac{\gamma_{n}}{\gamma_{n-1}}\bigr)=1.
    \end{equation}}
    \argument{\lref{as2: eq4};\lref{as2: eq1};}{that
    \begin{equation}\llabel{as2: eq6}
        \textstyle\lim_{n\to\infty} ((\gamma_{n})^{-2}(\gamma_{n-1}-\gamma_{n}))=\bigl[\lim_{n\to\infty} ((\gamma_{n})^{-2}(\gamma_{n-1})^2\bigr]\bigl[\lim_{n\to\infty} ((\gamma_{n-1})^{-2}(\gamma_{n-1}-\gamma_{n}))\bigr]=a.
    \end{equation}}
    \argument{\lref{as2: eq6};}{(\ref{item 2: equivalent learning rate}$\rightarrow$\ref{item 1: equivalent learning rate})\dott}
                 \end{aproof}
\subsection{Error estimates for MUON without boundedness assumptions}\label{subsec: without stays bound}
\begin{samepage}
 \begin{tcolorbox}[colback=white!95!gray,
                  colframe=black,
                  boxrule=0.5pt,
                  sharp corners,
                  enhanced,
                  breakable,
                 ]
\begin{athm}{theorem}{muon stochastic convergence}[\textcolor{red}{\MUON\ error analysis without the assumption that \MUON\ stays bounded}]
      Let $(\Omega,\cF,\P)$ be a probability space, let $\delta\in \N$, $d_1,d_2,\dots,d_\delta,\fd_1,\fd_2,\dots,\fd_\delta\in \N$, $\alpha\in (0,1)$,  $p,\chi\in (0,\infty)$, let $\nscoe\colon(\N_0)^2\to \R$ be a function,
  let $(\setX,\cU)$ be a measurable space, let $X_{n,m}\colon \Omega\to\setX$, $(n,m)\in \N^2$, be \iid\ random variables,  let $\cA=\times_{i=1}^\delta\R^{d_i\times\fd_i}$, let $\smalll=(\smalll(\theta,x))_{(\theta,x)\in \cA\times \setX}\colon \cA\times \setX\to\R$ be measurable, assume for all $x\in \setX$ that $(\cA\ni\theta\mapsto \smalll(\theta,x)\in \R)\in C^1(\cA,\R)$, assume for all  $\theta,\vartheta\in \cA$, $x,y\in \setX$ that $\textcolor{magenta}{ \|(\nabla_\theta\smalll)(\theta,x)-(\nabla_\theta\smalll)(\theta,y)\|+ |\smalll(0,x)|\leq  \chi}$ and
  \begin{equation}\llabel{assume}
     \textcolor{magenta}{ \|(\nabla_\theta\smalll)(\theta,x)-(\nabla_\theta\smalll)(\vartheta,x)\|^2\leq \chi\|\theta-\vartheta\|^2  \leq \chi^2\spro{(\nabla_{\theta}\smalll)(\theta,x)-(\nabla_{\theta}\smalll)(\vartheta,x),\theta-\vartheta}},
  \end{equation} 
      and let $\gamma\colon \N\to(0,\infty)$ be a function.
       Then 
    \begin{enumerate}[label=(\roman*)]
        \item \label{item 1: muon stochastic convergence}  there exists a unique $\vartheta\in \cA$ which satisfies 
         $\textcolor{magenta}{ \textstyle\E[\smalll(\vartheta,X_{1,1})]=\inf_{\theta\in \cA}\E[\smalll(\theta,X_{1,1})]}$ 
         and 
         \item \label{item 2: muon stochastic convergence} there exists $\fC\in \R$ such that for every $M,n\in \N$ and every \defmuon{\Theta}{\chi}{M}\ 
       it holds that
     \begin{equation} \llabel{conclude}
     \textcolor{magenta}{\bigl(\E\bigl[\|\Theta_{n-1}-\vartheta\|^p\bigr]\bigr)^{2/p}\leq \fC (M^{-1}+\gamma_{n})}\ifnocf.
\end{equation}
    \end{enumerate}
    \cfout[.]
\end{athm}
\end{tcolorbox}
\end{samepage}
\begin{aproof}
Throughout this proof let $\bfx\in \setX$.
\startnewargseq
\argument{\lref{assume};}{that there exists $\fC\in (0,\infty)$ such that for all $\theta\in\cA$, $x\in \setX$ it holds that
\begin{equation}\llabel{eqp1}
\begin{split}
    \|(\nabla_\theta\smalll)(\theta,x)\|&\leq \|(\nabla_\theta\smalll)(0,\bfx)\|+\|(\nabla_\theta\smalll)(0,x)-(\nabla_\theta\smalll)(0,\bfx)\|+\|(\nabla_\theta\smalll)(\theta,x)-(\nabla_\theta\smalll)(0,x)\|\\
    &\leq \|(\nabla_\theta\smalll)(0,\bfx)\|+\chi+\sqrt{\chi}\|\theta\|\leq \fC\|\theta\|+\fC.
    \end{split}
\end{equation}}
\argument{\lref{eqp1};\lref{assume};the fundamental theorem of calculus;}{that there exist $\fC_1,\fC_2\in (0,\infty)$ such that for all $\theta\in \cA$, $x\in\setX$ it holds that
\begin{equation}\llabel{eqp2}
\begin{split}
   | \smalll(\theta,x)|&\leq |\smalll(\theta,x)-\smalll(0,x)|+|\smalll(0,x)|\leq \biggl|\int_{0}^1 \spro{(\nabla_\theta\smalll)(t\theta,x),\theta}\,\d t\biggr|+\chi\\
   &\leq  \biggl|\int_{0}^1 \|(\nabla_\theta\smalll)(t\theta,x)\|\|\theta\|\,\d t\biggr|+\chi\leq  \biggl|\int_{0}^1(\fC_1\|\theta\|+\fC_1)\|\theta\|\,\d t\biggr|+\chi\\
   &\leq \fC_2\|\theta\|^2+\fC_2\|\theta\|+\fC_2.
   \end{split}
\end{equation}}
\argument{\lref{assume};\lref{eqp1};\lref{eqp2};the fact that $X_{1,1}\in \setX$;the dominated convergence theorem}{that for all $\theta\in \cA$ it holds that 
\begin{equation}\llabel{eqt1}
\nabla_\theta(\E[\smalll(\theta,X_{1,1})])=\E[(\nabla_\theta\smalll)(\theta,X_{1,1})].
\end{equation}}
    \argument{\lref{assume};\lref{eqt1};the fact that $X_{1,1}\in \setX$}{that for all $\theta,\vartheta\in \cA$ it holds that
    \begin{equation}\llabel{eqt2}
    \begin{split}
        &\spro{\nabla_\theta(\E[\smalll(\theta,X_{1,1})])-\nabla_\theta(\E[\smalll(\vartheta,X_{1,1})]),\theta-\vartheta}\\
        &=\spro{\E[(\nabla_\theta\smalll)(\theta,X_{1,1})]-\E[(\nabla_\theta\smalll)(\vartheta,X_{1,1})],\theta-\vartheta}\\&=\E[\spro{(\nabla_\theta\smalll)(\theta,X_{1,1})-(\nabla_\theta\smalll)(\vartheta,X_{1,1}),\theta-\vartheta}]\geq \chi^{-1}\|\theta-\vartheta\|^2.
        \end{split}
    \end{equation}}
\argument{\lref{eqt2};\unskip, \eg, \cite[Proposition 5.7.23]{ArBePhi2024}}{that \llabel{argt1} $\cA\ni\theta\mapsto \E[\smalll(\theta,X_{1,1})]\in \R$ is strongly convex\dott}
\argument{\lref{argt1};\unskip, \eg, \cite[Corollary 5.7.22]{ArBePhi2024}}{\cref{item 1: muon stochastic convergence}\dott}
In the following assume without loss of generality that there exist $N\in\N$ and an \defmuon{\Theta}{\chi}{N}.
\startnewargseq
\argument{the fact that there exists $N\in\N$ and an \defmuon{\Theta}{\chi}{N}}{that there exists an unique $K\in\N_0$ which satisfies that
\begin{enumerate}[label=(\Roman*)]
\item \llabel{item 1} it holds that $\nscoe_{K+1,0}=0<\min\{\nscoe_{0,0},\nscoe_{0,1}\}$,
\item \llabel{item 2} it holds that $\#(\nscoe^{ - 1 }( \R\backslash\{0\} )) <\infty$, and
    \item \llabel{item 3} it holds for all $n\in\N\cap[0,K]$, $x\in\R$ that
    $\sum_{i=0}^\infty \nscoe_{n,i}x^{2i}>0$
\end{enumerate}}
\startnewargseq
 \argument{\cref{T_B_D};\cref{NS error}}{that there exist $\scrc,\fC\in (0,\infty)$ such that for all $i\in \{1,2,\dots,\delta\}$, $A\in \R^{d_i\times\fd_i}$ it holds that
    \begin{equation}\llabel{eq1}
   \lVert\TNewtonSchulzAlgorithm{\nscoe}{K}{A}\rVert\leq \scrc \qqandqq \spro{\TNewtonSchulzAlgorithm{\nscoe}{K}{A},A}\geq \frac{\fC\|A\|^2}{\|A\|+\nscoe_{0,1}(\nscoe_{0,0})^{-1}}.
    \end{equation}}
    \argument{the fact that $\limsup_{n\to\infty}\gamma_n=0$; the fact that $\limsup_{n\to\infty} \bigl(\frac{\gamma_{n}-\gamma_{n+1}}{(\gamma_{n})^2}\bigr)<\infty$;}{that
    \begin{equation}\llabel{eq3'}
        \textstyle\limsup_{n\to\infty} \bigl(\frac{\gamma_{n}-\gamma_{n+1}}{\gamma_{n}}\bigr)=0.
    \end{equation}}
    \argument{\lref{eq3'};}{that
    \begin{equation}\llabel{eq3'.1}
        \textstyle\lim_{n\to\infty} \bigl(\frac{\gamma_{n}}{\gamma_{n+1}}\bigr)=1.
    \end{equation}}
    \argument{\lref{eq3'.1};the fact that $0<\alpha<1$}{that
    \begin{equation}\llabel{eq3''}
        \limsup_{n\to\infty}(\allowbreak(\gamma_{n+1})^{-1}\alpha\gamma_n)= \alpha<1.
    \end{equation}}
    \argument{\lref{eq1};\lref{eq3''};\cref{theo: MUON bound}; the fact that for all $n,m\in \N$ it holds that $X_{n,m}\in \setX$}{that there exists $\fC\in (0,\infty)$ such that for every $M\in \N$, every stochastic process $\bfm=(\bfm^{1},\dots,\bfm^{\delta})\colon \N_0\times\Omega\to\cA$, and every stochastic process $\Theta=(\Theta^{1},\dots,\Theta^{\delta})\colon \N_0\times\Omega\to\cA$ with the property that $\P(\|\Theta_0\|\leq \chi)=1$ and 
     with the property that for all $n\in \N$, $i\in \{1,2,\dots,\delta\}$ it holds that 
         \begin{equation}
            \bfm_0=0,\qquad \bfm_n= \alpha \bfm_{n-1}+(1-\alpha)\bigl[\textstyle \frac 1M\sum_{m=1}^M(\nabla_\theta\smalll)(\Theta_{n-1},X_{n,m})\bigr],
            \end{equation}
            \begin{equation}
            \text{and}\qquad  \Theta_n^{i}=\Theta_{n-1}^{i}-\gamma_n\TNewtonSchulzAlgorithm{\nscoe}{K}{\bfm_n^i}
         \end{equation}
         we have that
    \begin{equation}\llabel{eq2'}
        \textstyle\inf_{n\in\N_0}\P(\|\Theta_n\|< \fC)=1.
    \end{equation}}
    \argument{\lref{eq2'};\lref{item 1};\lref{item 2};\lref{item 3}}{that there exists $\fC\in (\max\{\|\vartheta\|,\chi\},\infty)$ which satisfies that for every $M\in \N$ and every \defmuon{\Theta}{\chi}{M} it holds that
         \begin{equation}\llabel{eq2}
        \textstyle\inf_{n\in\N_0}\P(\|\Theta_n\|< \fC)=1\ifnocf.
    \end{equation}}
    \startnewargseq
     \argument{\lref{eqp1};}{that there exists $\scrC\in (0,\infty)$ such that for all $\theta\in \cA$, $x\in \setX$ with $\|\theta\|<\fC$ it holds that
     \begin{equation}\llabel{eq4}
     \begin{split}
         \|(\nabla_{\theta}\smalll)(\theta,x)\|\leq \scrC.
         \end{split}
     \end{equation}}
    \argument{\lref{eq2};\lref{eq4};\cref{item 1: muon stochastic convergence};the fact that $\cA\ni\theta\mapsto \E[\smalll(\theta,X_{1,1})]\in \R$ is strongly convex;\cref{muon bound convergence stochastic}}{\lref{conclude}\dott}
\end{aproof}
\section{MUON applied to concrete examples SOPs}\label{sec: example}

In this section we illustrate \cref{muon bound convergence stochastic} and \cref{muon stochastic convergence} by means of several concrete example \SOPs. Specifically, in \cref{subsec: logistic} we apply \cref{muon stochastic convergence} from \cref{sec: convergence of stochastic muon} above to the situation of $ \ell_2 $-regularized logistic regression for binary classification (cf., \eg, \cite[Section 1]{ar1802.03081}), in \cref{subsec: quadradic muon} we apply \cref{muon stochastic convergence} to a class of quadratic \SOPs, and in \cref{subsec: deterministic muon} we apply \cref{muon bound convergence stochastic} from \cref{sec: convergence of stochastic muon} to the situation of deterministic optimization problems.
\subsection{MUON for logistic regression for binary classification}\label{subsec: logistic}

In this subsection we show in the elementary fact in \cref{lem: loss verify 2} that a class of $\ell_2$-regularized logistic regression problems satisfies the assumptions of the \SOP\ in \cref{muon stochastic convergence} and in \cref{main theorem 2 quadratic 4} we specialize \cref{muon stochastic convergence} to the situation where the considered \SOP\ belongs to this simple class of $\ell_2$-regularized logistic regression problems. \cref{main theorem 2 quadratic 4} follows directly from \cref{lem: loss verify 2} and \cref{muon stochastic convergence}.

\begin{tcolorbox}[colback=white!95!gray,
                  colframe=black,
                  boxrule=0.5pt,
                  sharp corners,
                  enhanced,
                  breakable,
                 ]
\begin{athm}{lemma}{lem: loss verify 2}
  Let  $d,\fd,v\in \N$, $\lambda_1\in [0,\infty)$, $\lambda_2\in (0,\infty)$, for every $i\in \{1,2,\dots,v\}$ let $M_i\colon\R^{d\times\fd}\to\R$ be affine functions, and let $\smalll=(\smalll(\theta,x))_{(\theta,x)\in \R^{d\times\fd}\times \{1,2,\dots,v\}}\colon\R^{d\times \fd}\times \{1,2,\dots,v\}\to\R$ satisfy for all $\theta\in \R^{d\times \fd}$, $i\in\{1,2,\dots,v\}$ that
    \begin{equation}\llabel{def: L}
   \textstyle\smalll(\theta,i) = \lambda_1 \textstyle\ln\bigl( 1 + \exp( M_i(\theta))\bigr)  + \lambda_2\| \theta \|^2 ,
    \end{equation}
    Then
    \begin{enumerate}[label=(\roman*)]
    \item \label{item 0: loss veify 2} it holds for all $i\in \{1,2,\dots,v\}$ that $(\R^{d\times\fd}\ni\theta\mapsto \smalll(\theta,i)\in \R)\in C^1(\R^{d\times\fd},\R)$,
        \item \label{item 1: loss verify 2} there exists $\cK\in (0,\infty)$ such that for all $\theta,\vartheta\in \R^{d\times\fd}$, $i,j\in \{1,2,\dots,v\}$ that $|\smalll(0,i)|\leq \cK$, $ \|(\nabla_\theta\smalll)(\theta,i)-(\nabla_\theta\smalll)(\vartheta,i)\|\leq \cK\|\theta-\vartheta\| $, and $ \|(\nabla_\theta\smalll)(\theta,i)-(\nabla_\theta\smalll)(\theta,j)\|\leq  \cK$, and
        \item \label{item 2: loss verify 2} there exists $\kappa\in (0,\infty)$ such that for all $\theta,\vartheta\in\R^{d\times\fd}$, $i\in \{1,2,\dots,v\}$ it holds that $\spro{(\nabla_{\theta}\smalll)(\theta,i)-(\nabla_{\theta}\smalll)(\vartheta,i),\theta-\vartheta}\geq \kappa\|\theta-\vartheta\|^2 $.
    \end{enumerate}
\end{athm}

\end{tcolorbox}
\begin{aproof}
    \argument{the fact that for all $i\in \{1,2,\dots,v\}$ it holds that $M_i$ is affine}{that for all $i\in \{1,2,\dots,v\}$ it holds that \llabel{arg1} $\bigl(\R^{d\times\fd}\ni\theta\mapsto \textstyle\lambda_1\ln\bigl( 1 + \exp( M_i(\theta))\bigr)  + \lambda_2\| \theta \|^2 \in \R\bigr)\in C^{1}(\R^{d\times\fd},\R)$\dott}
\argument{\lref{arg1};\lref{def: L}}{\cref{item 0: loss veify 2}\dott}
\startnewargseq
\argument{\lref{def: L};}{that there exists $\fC\in (0,\infty)$ such that for all $i\in\{1,2,\dots,v\}$ it holds that
\begin{equation}\llabel{eqt1}
 \textstyle   |\smalll(0,i)|\leq \lambda_1\ln\bigl(1+\exp\bigl(\sup_{j\in\{1,2,\dots,v\}}M_j(0)\bigr)\bigr)\leq\fC.
\end{equation}}
    \argument{\lref{def: L};}{that for all $\theta\in \R^{d\times\fd}$, $i\in \{1,2,\dots,v\}$ it holds that
    \begin{equation}\llabel{eq1}
        (\nabla_{\theta}\smalll)(\theta,i)=\frac{\lambda_1(\nabla M_i)(\theta)\exp(M_i(\theta))}{1+\exp(M_i(\theta))}+2\lambda_2\theta.
    \end{equation}}
    \argument{the fact that for all $i\in \{1,2,\dots,v\}$ it holds that $M_i$ is affine;}{that there exists $\cK\in (0,\infty)$ such that for all  $\theta,\vartheta\in \R^{d\times\fd}$ it holds that
    \begin{equation}\llabel{eq1'}
        (\nabla M_i)(\theta)=(\nabla M_i)(\vartheta),\qquad M_i(\theta)-M_i(\vartheta)=\spro{(\nabla M_i)(\theta),\theta-\vartheta},\qqandqq \| (\nabla M_i)(\theta)\|\leq \cK.
    \end{equation}}
    \argument{\lref{eq1};\lref{eq1'};the fact that $\R\ni x\mapsto \exp(x)(1+\exp(x))^{-1}\in \R$ is Lipschitz continuous}{that there exist $\cK_1,\cK_2,\cK_3\in (0,\infty)$ such that for all $\theta,\vartheta\in \R^{d\times\fd}$, $i\in \{1,2,\dots,v\}$ it holds that
    \begin{equation}\llabel{eq2}
    \begin{split}
        &\textstyle\| (\nabla_{\theta}\smalll)(\theta,i)- (\nabla_{\theta}\smalll)(\vartheta,i)\|\leq \textstyle\bigl\|\frac{\lambda_1(\nabla M_i)(\theta)\exp(M_i(\theta))}{1+\exp(M_i(\theta))}-\frac{\lambda_1(\nabla M_i)(\vartheta)\exp(M_i(\vartheta))}{1+\exp(M_i(\vartheta))}\bigr\|+\|2\lambda_2\theta-2\lambda_2\vartheta\|\\
        &=\lambda_1\textstyle\|(\nabla M_i)(\theta)\|\bigl|\frac{\exp(M_i(\theta))}{1+\exp(M_i(\theta))}-\frac{\exp(M_i(\vartheta))}{1+\exp(M_i(\vartheta))}\bigr|+2\lambda_2\|\theta-\vartheta\|\\
      &\leq \cK_1\cK_2\|M_i(\theta)-M_i(\vartheta)\|+2\lambda_2\|\theta-\vartheta\|\leq\cK_1\cK_2\|(\nabla M_i)(\theta)\|\|\theta-\vartheta\|+2\lambda_2\|\theta-\vartheta\|\leq \cK_3\|\theta-\vartheta\|. 
        \end{split}
    \end{equation}}
    \argument{\lref{eq1};\lref{eq1'}}{that there exists $\cK\in (0,\infty)$ such that for all $\theta\in \R^{d\times\fd}$, $i,j\in \{1,2,\dots,v\}$ it holds that
    \begin{equation}\llabel{eq3}
    \begin{split}
        \textstyle\| (\nabla_{\theta}\smalll)(\theta,i)- (\nabla_{\theta}\smalll)(\theta,j)\|&=\textstyle \bigl\|\frac{\lambda_1(\nabla M_i)(\theta)\exp(M_i(\theta))}{1+\exp(M_i(\theta))}-\frac{\lambda_1(\nabla M_j)(\theta)\exp(M_j(\theta))}{1+\exp(M_j(\theta))}\bigr\|\\
        &\textstyle\leq\lambda_1\|(\nabla M_i)(\theta)\|+\lambda_1\|(\nabla M_j)(\theta)\|\\
        &=\lambda_1\|(\nabla M_i)(0)\|+\lambda_1\|(\nabla M_j)(0)\|\leq \cK.
        \end{split}
    \end{equation}}
    \argument{\lref{eqt1};\lref{eq2};\lref{eq3}}{\cref{item 1: loss verify 2}\dott}
    \startnewargseq
    \argument{the fact that $\R\ni x\mapsto \ln(1+\exp(x))\in \R$ is convex;\unskip, \eg, \cite[Lemma 5.7.3]{ArBePhi2024}}{that for all $x,y\in \R$ it holds that
    \begin{equation}\llabel{eq5}
        \textstyle\bigl(\frac{\exp(x)}{1+\exp(x)}-\frac{\exp(y)}{1+\exp(y)}\bigr)(x-y)\geq 0.
    \end{equation}}
    \argument{\lref{eq1};\lref{eq1'};\lref{eq5}}{that there exists $\kappa\in (0,\infty)$ such that for all $\theta,\vartheta\in \R^{d\times\fd}$, $i\in\{1,2,\dots,v\}$ it holds that
    \begin{equation}\llabel{eq6}
    \begin{split}
        &\spro{(\nabla_{\theta}\smalll)(\theta,i)-(\nabla_{\theta}\smalll)(\vartheta,i),\theta-\vartheta}\\
        &=\textstyle \bigl\langle\frac{\lambda_1(\nabla M_i)(\theta)\exp(M_i(\theta))}{1+\exp(M_i(\theta))}-\frac{\lambda_1(\nabla M_i)(\vartheta)\exp(M_i(\vartheta))}{1+\exp(M_i(\vartheta))},\theta-\vartheta\bigr\rangle+2\lambda_2\spro{\theta-\vartheta,\theta-\vartheta}\\
        &=\textstyle\lambda_1\bigl(\frac{\exp(M_i(\theta))}{1+\exp(M_i(\theta))}-\frac{\exp(M_i(\vartheta))}{1+\exp(M_i(\vartheta))}\bigr)\spro{(\nabla M_i)(\theta),\theta-\vartheta}+2\lambda_2\|\theta-\vartheta\|^2\\
        &=\lambda_1\textstyle\bigl(\frac{\exp(M_i(\theta))}{1+\exp(M_i(\theta))}-\frac{\exp(M_i(\vartheta))}{1+\exp(M_i(\vartheta))}\bigr)(M_i(\theta)-M_i(\vartheta))+2\lambda_2\|\theta-\vartheta\|^2\\
        &\geq 2\lambda_2\|\theta-\vartheta\|^2.
        \end{split}
    \end{equation}}
    \argument{\lref{eq6};}{\cref{item 2: loss verify 2}\dott}
\end{aproof}
We now specialize \cref{muon stochastic convergence} to the situation where the loss function belongs to the class of quadratic loss functions studied in \cref{lem: loss verify}.

\begin{tcolorbox}[colback=white!95!gray,
                  colframe=black,
                  boxrule=0.5pt,
                  sharp corners,
                  enhanced,
                  breakable,
                 ]
\begin{athm}{cor}{main theorem 2 quadratic 4}[\textcolor{red}{\MUON\ for $\ell_2$-regularized logistic regression for binary classification}]
      Let $(\Omega,\cF,\P)$ be a probability space, let $d,\fd,v\in \N$, $x_1,x_2,\dots,x_v\in \R^{d\times\fd}$, $y_1,y_2,\dots,y_v\in \R$, $\alpha\in (0,1)$, $\lambda_1\in [0,\infty)$, $p,\chi,\lambda_2\in (0,\infty)$, let $\nscoe\colon(\N_0)^2\to \R$ be a function,
 let $X_{n,m}\colon \Omega\to\{1,2,\dots,v\}$, $(n,m)\in \N^2$, be \iid\ random variables, and let $\smalll=(\smalll(\theta,x))_{(\theta,x)\in \R^{d\times\fd}\times \{1,2,\dots,v\}}\colon\R^{d\times \fd}\times  \{1,2,\dots,v\}\to\R$ satisfy for all $\theta\in \R^{d\times \fd}$, $i\in\{1,2,\dots,v\}$ that
  \begin{equation}
   \textcolor{magenta}{ \smalll(\theta,i) = \lambda_1 \textstyle\ln\bigl( 1 + \exp( - y_i \spro{ x_i, \theta } )\bigr)  + \lambda_2 \| \theta \|^2 },
    \end{equation}
     and let $\gamma\colon \N\to(0,\infty)$ be a function.
    Then 
    \begin{enumerate}[label=(\roman*)]
        \item \label{item 1: main theorem 2 quadratic 4} there exists a unique $\vartheta\in \R^{d\times\fd}$ which satisfies $\textcolor{magenta}{ \E[\smalll(\vartheta,X_{1,1})]=\inf_{\theta\in \R^{d\times\fd}}\E[\smalll(\theta,X_{1,1})]}$ and
        \item \label{item 2: main theorem 2 quadratic 4} there exists $\fC\in \R$ such that for every $M,n\in \N$ and every \defmuon{\Theta}{\chi}{M}\ it holds that
\begin{equation} \llabel{conclude}
     \textcolor{magenta}{\bigl(\E\bigl[\|\Theta_{n-1}-\vartheta\|^p\bigr]\bigr)^{2/p}\leq \fC (M^{-1}+\gamma_{n})}\ifnocf.
\end{equation}
    \end{enumerate}
\cfout[.]
\end{athm}
\end{tcolorbox}
\begin{aproof}
    \argument{\cref{lem: loss verify 2};}{ that 
    \begin{enumerate}[label=(\Roman*)]
    \item \llabel{item 1} it holds for all $i\in \{1,2,\dots,v\}$ that $(\R^{d\times\fd}\ni\theta\mapsto \smalll(\theta,i)\in \R)\in C^1(\R^{d\times\fd},\R)$
        \item \llabel{item 2} there exists $\cK\in (0,\infty)$ such that for all $\theta,\vartheta\in \R^{d\times\fd}$, $i,j\in \{1,2,\dots,v\}$ that $|\smalll(0,i)|\leq \cK$, $ \|(\nabla_\theta\smalll)(\theta,i)-(\nabla_\theta\smalll)(\vartheta,i)\|\leq \cK\|\theta-\vartheta\| $, and $ \|(\nabla_\theta\smalll)(\theta,i)-(\nabla_\theta\smalll)(\theta,j)\|\leq  \cK$, and
        \item \llabel{item 3} there exists $\kappa\in (0,\infty)$ such that for all $i\in \{1,2,\dots,v\}$ it holds that $\spro{(\nabla_{\theta}\smalll)(\theta,i)-(\nabla_{\theta}\smalll)(\vartheta,i),\theta-\vartheta}\geq \kappa\|\theta-\vartheta\|^2 $.
    \end{enumerate}}
    \argument{\lref{item 3};\cref{muon stochastic convergence}}{\cref{item 1: main theorem 2 quadratic 4,item 2: main theorem 2 quadratic 4}\dott}
\end{aproof}
\subsection{MUON for quadratic SOPs} \label{subsec: quadradic muon}

In this subsection we show in the elementary fact in \cref{lem: loss verify} that a class of quadratic \SOPs\ satisfies the assumptions of the \SOP\ in \cref{muon stochastic convergence} and in \cref{main theorem 2 quadratic} we specialize \cref{muon stochastic convergence} to the situation where the considered \SOP\ belongs to this simple class of quadratic \SOPs. \cref{main theorem 2 quadratic} follows directly from \cref{lem: loss verify} and \cref{muon stochastic convergence}.
\begin{tcolorbox}[colback=white!95!gray,
                  colframe=black,
                  boxrule=0.5pt,
                  sharp corners,
                  enhanced,
                  breakable,
                 ]
\begin{athm}{lemma}{lem: loss verify}
    Let $\delta,\dimX,v\in \N$, $d_1,d_2,\dots,d_\delta,\fd_1,\fd_2,\dots,\fd_\delta\in \N$, $(\bfd_{i,j})_{(
    i,j)\in \N^2}\subseteq\N$, for every $i\in \{1,2,\dots,\delta\}$, $j\in\{1,2,\dots,v\}$ let  $q_{i,j}\in C( \R^\dimX,\R^{d_i\times \bfd_{i,j}})$ and let $M_{i,j}\colon\R^{d_i\times \fd_i}\to \R^{d_i\times \bfd_{i,j}}$ be an affine function, let $k\in \{1,2,\dots,v\}$ satisfy for all $i\in\{1,2,\dots,\delta\}$ that $M_{i,k}$ is invertible, let $\cA=\times_{i=1}^\delta\R^{d_i\times \fd_i}$, let $\smalll=(\smalll(\theta,x))_{(\theta,x)\in \cA\times \R^{\dimX}}\colon\cA\times \R^{\dimX}\to\R$ satisfy for all $\theta=(\theta_1,\dots,\theta_\delta)\in \cA$, $x\in \R^\dimX$ that
    \begin{equation}\llabel{def: L}
    \smalll(\theta,x)= \sum_{ i = 1 }^{ \delta }\sum_{j=1}^v \| M_{i,j} (\theta_i)- q_{i,j}( x ) \|^2,
    \end{equation}
    and let $\setX\subseteq\R^\dimX$ be bounded and measurable.
    Then
    \begin{enumerate}[label=(\roman*)]
    \item \label{item 0: loss veify} it holds for all $x\in \setX$ that $(\cA\ni\theta\mapsto \smalll(\theta,x)\in \R)\in C^1(\cA,\R)$,
        \item \label{item 1: loss verify} there exists $\cK\in(0,\infty)$ such that for all $\theta,\vartheta\in \cA$, $x,y\in \setX$ it holds that $|\smalll(0,x)|\leq \cK$, $\|(\nabla_\theta\smalll)(\theta,x)-(\nabla_\theta\smalll)(\vartheta,x)\|\leq \cK\|\theta-\vartheta\|$ and $ \|(\nabla_\theta\smalll)(\theta,x)-(\nabla_\theta\smalll)(\theta,y)\|\leq  \cK$.
        \item \label{item 2: loss verify} there exists $\kappa\in (0,\infty)$ such that for all $x\in \R^\dimX$ it holds that $\cA\ni\theta\mapsto \smalll(\theta,x)-\kappa\|\theta\|^2\in \R$ is convex.
    \end{enumerate}
\end{athm}
\end{tcolorbox}
\begin{aproof}
\argument{the fact that for all $i\in \{1,2,\dots,\delta\}$, $j\in\{1,2,\dots,v\}$ it holds that $M_{i,j}$ is affine}{that for all $i\in \{1,2,\dots,\delta\}$, $j\in\{1,2,\dots,v\}$, $x\in \R^\dimX$ it holds that \llabel{arg1} $(\cA\ni\theta=(\theta_1,\dots,\theta_\delta)\mapsto \|M_{i,j}(\theta_i)-q_{i,j}(x)\|^2\in \R)\in C^{1}(\cA,\R)$\dott}
\argument{\lref{arg1};\lref{def: L}}{\cref{item 0: loss veify}\dott}
\startnewargseq
\argument{\lref{def: L};the fact that for all $i\in\{1,2,\dots,\delta\}$, $j\in\{1,2,\dots,v\}$ it holds that $q_{i,j}$ is continuous; the fact that $\setX$ is bounded;}{that there exists $\fC\in (0,\infty)$ such that for all $x\in\setX$ it holds that
\begin{equation}\llabel{eqt1}
    |\smalll(0,x)|=\sum_{i=1}^\delta\sum_{j=1}^v\|M_{i,j}(0)-q_{i,j}(x)\|^2\leq \sum_{i=1}^\delta\sum_{j=1}^v\fC =\fC\delta v.
\end{equation}}
    \argument{\lref{def: L};}{that for all $\theta=(\theta_1,\dots,\theta_\delta)\in \cA$, $x\in \R^\dimX$, $i\in\{1,2,\dots,\delta\}$ it holds that
    \begin{equation}\llabel{eq1}
        (\nabla_{\theta_i}\smalll)(\theta,x)=\sum_{j=1}^v\bigl(2[(M_{i,j})'(\theta_i)]^{\top}(M_{i,j}(\theta_i)-q_{i,j}(x))\bigr).
    \end{equation}}
    \argument{the fact that for all $i\in \{1,2,\dots,\delta\}$, $j\in\{1,2,\dots,v\}$ it holds that $M_{i,j}$ is affine;}{that there exists $\cK\in (0,\infty)$ such that for all $i\in \{1,2,\dots,\delta\}$, $j\in \{1,2,\dots,v\}$, $\theta,\vartheta\in \R^{d_i\times\fd_i}$ it holds that
    \begin{equation}\llabel{eq1'}
        (M_{i,j})'(\theta)=(M_{i,j})'(\vartheta),\quad M_{i,j}(\theta)-M_{i,j}(\vartheta)=[(M_{i,j})'(\vartheta)](\theta-\vartheta),\qandq \| (M_{i,j})'(\theta)\|\leq \cK.
    \end{equation}}
    \argument{\lref{eq1};\lref{eq1'};the Cauchy-Schwarz inequality}{that there exist $\cK_1,\cK_2\in (0,\infty)$ such that for all $\theta=(\theta_1,\dots,\theta_\delta)$, $\vartheta=(\vartheta_1,\dots,\vartheta_\delta)\in \cA$, $x\in \R^\dimX$, $i\in\{1,2,\dots,\delta\}$ it holds that
    \begin{equation}\llabel{eq2}
    \begin{split}
        &\| (\nabla_{\theta_i}\smalll)(\theta,x)- (\nabla_{\theta_i}\smalll)(\vartheta,x)\|\\
        &\leq\sum_{j=1}^v\|2[(M_{i,j})'(\theta_i)]^{\top}(M_{i,j}(\theta_i)-q_{i,j}(x))-2[(M_{i,j})'(\vartheta_i)]^{\top}(M_{i,j}(\vartheta_i)-q_{i,j}(x))\|\\
        &=\sum_{j=1}^v\|2[(M_{i,j})'(\theta_i)]^{\top}(M_{i,j}(\theta_i)-q_{i,j}(x))-2[(M_{i,j})'(\theta_i)]^{\top}(M_{i,j}(\vartheta_i)-q_{i,j}(x))\|\\
        &\leq \cK_1\sum_{j=1}^v\|M_{i,j}(\theta_i)-q_{i,j}(x)-(M_{i,j}(\vartheta_i)-q_{i,j}(x))\|\\
        &= \cK_1\sum_{j=1}^v\|M_{i,j}(\theta_i)-M_{i,j}(\vartheta_i)\|\leq \cK_2\|\theta_i-\vartheta_i\|\leq \cK_2\|\theta-\vartheta\|.
        \end{split}
    \end{equation}}
    \argument{\lref{eq1};\lref{eq1'};the fact that $\sup_{i\in \{1,2,\dots,\delta\}}\sup_{j\in \{1,2,\dots,v\}}\sup_{x\in \setX}\|q_{i,j}(x)\|<\infty$;}{that there exists $\cK\in (0,\infty)$ such that for all $\theta=(\theta_1,\dots,\theta_\delta)\in \cA$, $x,y\in \setX$, $i\in\{1,2,\dots,\delta\}$ it holds that
    \begin{equation}\llabel{eq3}
    \begin{split}
       & \| (\nabla_{\theta_i}\smalll)(\theta,x)- (\nabla_{\theta_i}\smalll)(\theta,y)\|\\
       &\leq \sum_{j=1}^v\|2[(M_{i,j})'(\theta_i)]^{\top}(M_{i,j}(\theta_i)-q_{i,j}(x))-2[(M_{i,j})'(\theta_i)]^{\top}(M_{i,j}(\theta_i)-q_{i,j}(y))\|\\
        &= \sum_{j=1}^v\|2[(M_{i,j})'(\theta_i)]^{\top}(q_{i,j}(x)-q_{i,j}(y))\|\leq 2\sum_{j=1}^v\bigl(\|((M_{i,j})'(\theta_i))^{\top}\|\|q_{i,j}(x)-q_{i,j}(y)\|\bigr)\\
        &\leq 2\sum_{j=1}^v\cK=2v\cK.
        \end{split}
    \end{equation}}
    \argument{\lref{eqt1};\lref{eq3};\lref{eq2}}{\cref{item 1: loss verify}\dott}
    \startnewargseq
    \argument{the fact that for all $i\in \{1,2,\dots,\delta\}$, $j\in\{1,2,\dots,v\}$ it holds that $M_{i,j}$ is affine;the fact that for all $i\in \{1,2,\dots,\delta\}$ it holds that $M_{i,k}$ is invertible}{that for all $i\in \{1,2,\dots,\delta\}$, $\theta\in \R^{d_i\times\fd_i}$ it holds that
    \begin{equation}\llabel{eq5.1}
        \{\vartheta\in \R^{d_i\times\fd_i}\colon [(M_{i,k})'(\theta)]^{\top}\vartheta =0\}=\{0\}.
    \end{equation}}
    \argument{\lref{eq5.1};}{that there exists $\kappa\in (0,\infty)$ such that for all $i\in \{1,2,\dots,\delta\}$, $\theta\in \R^{d_i\times \fd_i}$ it holds that
    \begin{equation}\llabel{eq5}
        \theta^{\top}[(M_{i,k})'(\theta)]^{\top}[(M_{i,k})'(\theta)]\theta\geq \kappa\|\theta\|^2.
    \end{equation}}
    \argument{\lref{eq1};\lref{eq1'};\lref{eq5};}{that there exists $\kappa\in (0,\infty)$ such that for all $\theta=(\theta_1,\dots,\theta_\delta)$, $\vartheta=(\vartheta_1,\dots,\vartheta_\delta)\in \cA$, $x\in \R^\dimX$ it holds that
    \begin{equation}\llabel{eq6}
    \begin{split}
        &\spro{(\nabla_{\theta}\smalll)(\theta,x)-(\nabla_{\theta}\smalll)(\vartheta,x),\theta-\vartheta}=\sum_{i=1}^\delta \spro{(\nabla_{\theta_i}\smalll)(\theta,x)-(\nabla_{\theta_i}\smalll)(\vartheta,x),\theta_i-\vartheta_i}\\
        &=\sum_{i=1}^\delta \sum_{j=1}^v\spro{2[(M_{i,j})'(\theta_i)]^{\top}(M_{i,j}(\theta_i)-q_{i,j}(x))-2[(M_{i,j})'(\theta_i)]^{\top}(M_{i,j}(\vartheta_i)-q_{i,j}(x)),\theta_i-\vartheta_i}\\
        &=\sum_{i=1}^\delta \sum_{j=1}^v\spro{2[(M_{i,j})'(\theta_i)]^{\top}(M_{i,j}(\theta_i)-M_{i,j}(\vartheta_i)),\theta_i-\vartheta_i}\\
        &=\sum_{i=1}^\delta\sum_{j=1}^v 2(\theta_i-\vartheta_i)^{\top}[(M_{i,j})'(\theta_i)]^{\top}[(M_{i,j})'(\theta_i)](\theta_i-\vartheta_i)\\
        &\geq\sum_{i=1}^\delta 2(\theta_i-\vartheta_i)^{\top}[(M_{i,k})'(\theta_i)]^{\top}[(M_{i,k})'(\theta_i)](\theta_i-\vartheta_i)\geq \sum_{i=1}^\delta 2\kappa\|\theta_i-\vartheta_i\|^2=2\kappa\|\theta-\vartheta\|^2.
        \end{split}
    \end{equation}}
    \argument{\lref{eq6};\cite[Proposition 5.7.23]{ArBePhi2024}}{\cref{item 2: loss verify}\dott}
\end{aproof}
\begin{samepage}
 \begin{tcolorbox}[colback=white!95!gray,
                  colframe=black,
                  boxrule=0.5pt,
                  sharp corners,
                  enhanced,
                  breakable,
                 ]
\begin{athm}{cor}{main theorem 2 quadratic}[\textcolor{red}{\MUON\ for quadratic \SOPs}]
      Let $(\Omega,\cF,\P)$ be a probability space, let $\delta,\dimX,v\in \N$, $d_1,d_2,\dots,d_\delta,\fd_1,\fd_2,\dots,\fd_\delta\in \N$, $(\bfd_{i,j})_{(
    i,j)\in \N^2}\subseteq\N$, $\alpha\in (0,1)$,  $p,\chi\in (0,\infty)$, let $\nscoe\colon(\N_0)^2\to \R$ be a function, let $\setX\subseteq \R^\dimX$ be bounded and measurable,
 let $X_{n,m}\colon \Omega\to\setX$, $(n,m)\in \N^2$, be \iid\ random variables,  let $\cA=\times_{i=1}^\delta\R^{d_i\times \fd_i}$, for every $i\in \{1,2,\dots,\delta\}$, $j\in\{1,2,\dots,v\}$ let $q_{i,j}\in C( \R^\dimX,\R^{d_i\times \bfd_{i,j}})$ and let $M_{i,j}\colon\R^{d_i\times \fd_i}\to \R^{d_i\times \bfd_{i,j}}$ be an affine function, let $k\in \{1,2,\dots,v\}$ satisfy for all $i\in\{1,2,\dots,\delta\}$ that $M_{i,k}$ is invertible, let $\smalll=(\smalll(\theta,x))_{(\theta,x)\in \cA\times \R^{\dimX}}\colon\cA\times \R^{\dimX}\to\R$ satisfy for all $\theta=(\theta_1,\dots,\theta_\delta)\in \cA$, $x\in \R^\dimX$ that $
    \smalll(\theta,x)= \sum_{ i = 1 }^{ \delta } \sum_{j=1}^v\| M_{i,j} (\theta_i) - q_{i,j}( x ) \|^2$,
      and let $\gamma\colon\N\to(0,\infty)$ be a function.
       Then 
       \begin{enumerate}[label=(\roman*)]
           \item \label{item 1: main theorem 2 quadratic} there exists a unique  $\vartheta\in \cA$ which satisfies 
     $
       \textstyle  \textcolor{magenta}{ \E[\smalll(\vartheta,X_{1,1})]=\inf_{\theta\in \cA}\E[\smalll(\theta,X_{1,1})]}
   $ and \item \label{item 2: main theorem 2 quadratic} there exists $\fC\in \R$ such that for every $M,n\in \N$ and every \defmuon{\Theta}{\chi}{M}\ it holds that
\begin{equation} \llabel{conclude}
     \textcolor{magenta}{\bigl(\E\bigl[\|\Theta_{n-1}-\vartheta\|^p\bigr]\bigr)^{2/p}\leq \fC (M^{-1}+\gamma_n)}\ifnocf.
\end{equation}
       \end{enumerate}
\cfout[.]
\end{athm}
\end{tcolorbox}
\end{samepage}
\begin{aproof}
\argument{the assumption that for all  $\theta=(\theta_1,\dots,\theta_\delta)\in \cA$, $x\in \R^\dimX$ it holds that $
    \smalll(\theta,x)= \sum_{ i = 1 }^{ \delta } \sum_{j=1}^v\| M_{i,j} (\theta_i) - q_{i,j}( x ) \|^2$;\cref{lem: loss verify}}{that
\begin{enumerate}[label=(\Roman*)]
       \item \llabel{item 0} it holds for all $x\in \setX$ that $(\cA\ni\theta\mapsto \smalll(\theta,x)\in \R)\in C^1(\cA,\R)$
        \item \llabel{item 1} there exists $\cK\in(0,\infty)$ such that for all $\theta,\vartheta\in \cA$, $x,y\in \setX$ it holds that $|\smalll(0,x)|\leq \cK$, $\|(\nabla_\theta\smalll)(\theta,x)-(\nabla_\theta\smalll)(\vartheta,x)\|\leq \cK\|\theta-\vartheta\|$ and $ \|(\nabla_\theta\smalll)(\theta,x)-(\nabla_\theta\smalll)(\theta,y)\|\leq  \cK$.
        \item \llabel{item 2} there exists $\kappa\in (0,\infty)$ such that for all $x\in \R^\dimX$ it holds that $\cA\ni\theta\mapsto \smalll(\theta,x)-\kappa\|\theta\|^2\in \R$ is convex.
    \end{enumerate}}
    \argument{\lref{item 2};\cref{muon stochastic convergence}}{\cref{item 1: main theorem 2 quadratic,item 2: main theorem 2 quadratic}\dott}
\end{aproof}

\begin{samepage}
 \begin{tcolorbox}[colback=white!95!gray,
                  colframe=black,
                  boxrule=0.5pt,
                  sharp corners,
                  enhanced,
                  breakable,
                 ]
\begin{athm}{cor}{main theorem 2 quadratic 2}[\textcolor{red}{\MUON\ for quadratic \SOPs}]
      Let $(\Omega,\cF,\P)$ be a probability space, let $\delta,\dimX\in \N$, $d_1,d_2,\dots,d_\delta,\fd_1,\fd_2,\dots,\fd_\delta\in \N$, $\alpha\in (0,1)$,  $p,\chi\in (0,\infty)$, let $\nscoe\colon(\N_0)^2\to \R$ be a function, let $\setX\subseteq \R^\dimX$ be bounded and measurable,
 let $X_{n,m}\colon \Omega\to\setX$, $(n,m)\in \N^2$, be \iid\ random variables,  let $\cA=\times_{i=1}^\delta\R^{d_i\times \fd_i}$, for every $i\in \{1,2,\dots,\delta\}$ let $q_{i}\in C( \R^\dimX,\R^{d_i\times \fd_i})$ and let $M_{i}\colon\R^{d_i\times \fd_i}\to \R^{d_i\times \fd_{i}}$ be an invertible affine function, let $\smalll=(\smalll(\theta,x))_{(\theta,x)\in \cA\times \R^{\dimX}}\colon\cA\times \R^{\dimX}\to\R$ satisfy for all $\theta=(\theta_1,\dots,\theta_\delta)\in \cA$, $x\in \R^\dimX$ that $
    \smalll(\theta,x)= \sum_{ i = 1 }^{ \delta } \|M_i (\theta_i) - q_i( x ) \|^2 $,
      and let $\gamma\colon\N\to(0,\infty)$ be a function.
       Then 
       \begin{enumerate}[label=(\roman*)]
           \item \label{item 1: main theorem 2 quadratic 2} there exists a unique  $\vartheta\in \cA$ which satisfies 
     $
       \textstyle  \textcolor{magenta}{ \E[\smalll(\vartheta,X_{1,1})]=\inf_{\theta\in \cA}\E[\smalll(\theta,X_{1,1})]}
   $ and \item \label{item 2: main theorem 2 quadratic 2} there exists $\fC\in \R$ such that for every $M,n\in \N$ and every \defmuon{\Theta}{\chi}{M}\ it holds that
\begin{equation} \llabel{conclude}
     \textcolor{magenta}{\bigl(\E\bigl[\|\Theta_{n-1}-\vartheta\|^p\bigr]\bigr)^{2/p}\leq \fC (M^{-1}+\gamma_n)}\ifnocf.
\end{equation}
       \end{enumerate}
\cfout[.]
\end{athm}
\end{tcolorbox}
\end{samepage}
\begin{aproof}
    \argument{\cref{main theorem 2 quadratic} (applied with $v\curvearrowleft 1$,  $(\bfd_{i,1})_{i\in \{1,2,\dots,\delta\}}\curvearrowleft (\fd_i)_{i\in\{1,2,\dots,\delta\}}$, $(M_{i,1})_{i\in \{1,2,\dots,\delta\}}\curvearrowleft (M_i)_{i\in\{1,2,\dots,\delta\}}$, $(q_{i,1})_{i\in \{1,2,\dots,\delta\}}\curvearrowleft (q_i)_{i\in\{1,2,\dots,\delta\}}$, $k\curvearrowleft 1$ in the notation of \cref{main theorem 2 quadratic})}{\cref{item 1: main theorem 2 quadratic 2,item 2: main theorem 2 quadratic 2}\dott}
\end{aproof}

\subsection{MUON for deterministic optimization problems}
\label{subsec: deterministic muon}

In this subsection we apply \cref{muon bound convergence stochastic} above to the special situation of deterministic optimization problems. In \cref{deterministic muon bound convergence} we assume the boundedness of the considered \MUON\ optimization process and in \cref{deterministic muon convergence} we restrict ourselves to the situation where the gradient of the objective function is Lipschitz continuous but then we do not assume anymore that the \MUON\ optimization process is bounded.
\begin{samepage}
 \begin{tcolorbox}[colback=white!95!gray,
                  colframe=black,
                  boxrule=0.5pt,
                  sharp corners,
                  enhanced,
                  breakable,
                 ]
\begin{athm}{cor}{deterministic muon bound convergence}[\textcolor{red}{\MUON\ for deterministic optimization problems}]
      Let $\delta\in \N$, $K\in\N_0$, $d_1,d_2,\dots,d_\delta,\fd_1,\fd_2,\dots,\fd_\delta\in \N$, $\alpha\in (0,1)$, $\kappa\in (0,\infty)$, let $\nscoe=(\nscoe_{i,j})_{(i,j)\in (\N_0)^2}\colon(\N_0)^2\to \R$ satisfy $\textcolor{magenta}{\min\{\nscoe_{0,0},\nscoe_{0,1}\}>0}$ and $\textcolor{magenta}{\#(\nscoe^{ - 1 }( \R\backslash\{0\} )) <\infty}$, assume for all $i\in \N\cap[0,K]$, $x\in \R$ that $\textcolor{magenta}{\sum_{j=0}^\infty \nscoe_{i,j}x^{2j}>0}$,
  let $\cA=\times_{i=1}^\delta\R^{d_i\times\fd_i}$, $\smalll\in \allowbreak C^{1}(\cA,\R)$, assume that $\smalll$ is strongly convex, assume that $\nabla\smalll$ is locally Lipschitz continuous, let $\vartheta\in \cA$ satisfy $\smalll(\vartheta)=\inf_{\theta\in \cA}\smalll(\theta)$, let $(\gamma_n)_{n\in \N}\subseteq (0,\infty)$ be non-increasing, and assume $\limsup_{n\to\infty} (\gamma_n+(\gamma_n)^{-2}(\gamma_{n}-\gamma_{n+1}))=0$.
     Then there exists $\fC\in \R$ such that for every $\bfm=(\bfm^{1},\dots,\bfm^{\delta})\colon \N_0\to\cA$ and every $\Theta=(\Theta^{1},\dots,\Theta^{\delta})\colon \N_0\to\cA$ with the property that for all $n\in \N$, $i\in \{1,2,\dots,\delta\}$ it holds that 
         \begin{equation}\llabel{def: bfm}
            \textcolor{magenta}{\bfm_0=0},\qquad \textcolor{magenta}{\bfm_n= \alpha \bfm_{n-1}+(1-\alpha)(\nabla\smalll)(\Theta_{n-1})},
            \end{equation}
            \begin{equation}\llabel{def: Theta}
           \textcolor{magenta}{\|\Theta_{n-1}\|\leq \kappa},\qquad\text{and}\qquad  \textcolor{magenta}{\Theta_n^{i}=\Theta_{n-1}^{i}-\gamma_n\TNewtonSchulzAlgorithm{\nscoe}{K}{\bfm_n^i}}
         \end{equation}
     we have for all $n \in \N_0$ that
\begin{equation} \llabel{conclude}
     \textcolor{magenta}{\|\Theta_n-\vartheta\|\leq \fC \sqrt{\gamma_{n+1}}}
\end{equation}
(cf.\ \cref{definition: NS}).
\end{athm}
\end{tcolorbox}
\end{samepage}
\begin{aproof}
    Throughout this proof let $(\Omega,\cF,\P)$ be a probability space, for every $n,m$ let $X_{n,m}\colon \Omega\to\R$ satisfy for all $\omega\in \Omega$ that
    \begin{equation}\llabel{def: XY}
        X_{n,m}(\omega)=0,
    \end{equation}
    let $\fL=(\fL(\theta,x))_{(\theta,x)\in \cA\times\R}\colon \cA\times\R\to\R$ satisfy for all $\theta\in \cA$, $x\in \R$ that
    \begin{equation}\llabel{def: fL}
        \fL(\theta,x)=\smalll(\theta),
    \end{equation}
    let $\cL\colon \cA\to\R$ satisfy for all $\theta\in \cA$ that \begin{equation}\llabel{def: cL}
        \cL(\theta)=\E[\fL(\theta,X_{1,1})], 
    \end{equation}
    and  let $\bfa=(\bfa_{i,j})_{(i,j)\in (\N_0)^2}\colon(\N_0)^2\to \R$ satisfy for all $i,j\in\N_0$ that $\bfa_{i,j}=\nscoe_{i,j}\mathbbm 1_{\{0,1,\dots,K\}}(i)$.
    \argument{ \lref{def: XY};\lref{def: fL};\lref{def: cL}}{that for all $\theta\in \cA$, $x\in\R$ it holds that
    \begin{equation}\llabel{eq1}
       \cL(\theta)=\fL(\theta,0)=\smalll(\theta)=\fL(\theta,x). 
    \end{equation}}
    \argument{\lref{eq1};the assumption that $\smalll\in C^1(\cA)$;the assumption that $\smalll$ is strongly convex;the assumption that $\nabla \smalll$ is locally Lipschitz; the fact that $\smalll(\vartheta)=\inf_{\theta\in \cA}\smalll(\theta)$}{that
    \begin{enumerate}[label=(\roman*)]
\item \llabel{item 0} it holds that $\fL\in C^{1,0}(\cA\times\R,\R)$,
        \item \llabel{item 1} it holds that $\cL$ is strongly covex,
        \item \llabel{item 2} there exists $\cK\in (0,\infty)$ which satisfies for all $\theta_1,\theta_2\in \cA$, $x\in \R$ with $\max\{\|\theta_1\|,\|\theta_2\|\}\leq \kappa$ that $ \|(\nabla_\theta\fL)(\theta_1,x)-(\nabla_\theta\fL)(\theta_2,x)\|\leq \cK\|\theta_1-\theta_2\|$,
        and
        \item \llabel{item 3} it holds that $\cL(\vartheta)=\inf_{\theta\in \cA}\cL(\theta)$.
    \end{enumerate}}
    \argument{\lref{item 3};\cref{muon bound convergence stochastic} (applied with  $\cV\curvearrowleft \R$, $U\curvearrowleft \{0\}$, $p\curvearrowleft 1$, $\smalll\curvearrowleft \fL$, $\cL\curvearrowleft\cL$ in the notation of \cref{muon bound convergence stochastic};}{that there exists $\fC\in \R$ such that for every $M\in \N$, $n\in \N_0$ and every \defmuonfull{\Theta}{(\R,|\cdot|)}{(X_{v,w})_{(v,w)\in\N^2}}{2\kappa+2\|\vartheta\|}\ we have that
   \begin{equation}\llabel{eq3'}
\E\bigl[\|\Theta_n-\vartheta\|\bigr]\leq \fC (M^{-1/2}+\sqrt{\gamma_{n+1}})\ifnocf.
\end{equation} }
   \argument{\lref{eq3'};the assumption that $\min\{\nscoe_{0,0},\nscoe_{0,1}\}>0=\nscoe_{K+1,0}$;the assumption that $\#(\nscoe^{ - 1 }( \R\backslash\{0\} )) <\infty$} {that there exists $\fC\in \R$ such that for every $M\in \N$, every stochastic process $\bfm=(\bfm^{1},\dots,\bfm^{\delta})\colon \N_0\times\Omega\to\cA$, and every stochastic process $\Theta=(\Theta^{1},\dots,\Theta^{\delta})\colon \N_0\times\Omega\to\cA$ with the property that $\sup_{n\in \N_0}\|\Theta_n\|\leq \kappa$, with the property that for all $n\in \N$, $i\in \{1,2,\dots,\delta\}$ it holds that 
         \begin{equation}
            \bfm_0=0,\qquad \bfm_n= \alpha \bfm_{n-1}+(1-\alpha)\bigl[\textstyle \frac 1M \sum_{m=1}^M(\nabla_\theta\fL)(\Theta_{n-1},X_{n,m})\bigr],
            \end{equation}
            \begin{equation}
           \text{and}\qquad  \Theta_n^{i}=\Theta_{n-1}^{i}-\gamma_n\TNewtonSchulzAlgorithm{\bfa}{K}{\bfm_n^i},
         \end{equation}
         and with the property that  $\Theta_0$ and $(X_{n,m})_{(n,m)\in \N^2}$ are independent
     we have for all $n \in \N_0$ that
\begin{equation}\llabel{eq3''}
\E\bigl[\|\Theta_n-\vartheta\|\bigr]\leq \fC (M^{-1/2}+\sqrt{\gamma_{n+1}}).
\end{equation}}
 \argument{\lref{eq3''};} {that there exists $\fC\in \R$ such that for every $M\in \N$, every stochastic process $\bfm=(\bfm^{1},\dots,\bfm^{\delta})\colon \N_0\times\Omega\to\cA$, and every stochastic process $\Theta=(\Theta^{1},\dots,\Theta^{\delta})\colon \N_0\times\Omega\to\cA$ with the property that $\sup_{n\in \N_0}\|\Theta_n\|\leq 2\kappa+2\|\vartheta\|$, with the property that for all $n\in \N$, $i\in \{1,2,\dots,\delta\}$ it holds that 
         \begin{equation}
            \bfm_0=0,\qquad \bfm_n= \alpha \bfm_{n-1}+(1-\alpha)\bigl[\textstyle \frac 1M \sum_{m=1}^M(\nabla_\theta\fL)(\Theta_{n-1},X_{n,m})\bigr],
            \end{equation}
            \begin{equation}
           \text{and}\qquad  \Theta_n^{i}=\Theta_{n-1}^{i}-\gamma_n\TNewtonSchulzAlgorithm{\nscoe}{K}{\bfm_n^i},
         \end{equation}
         and with the property that  $\Theta_0$ and $(X_{n,m})_{(n,m)\in \N^2}$ are independent
     we have for all $n \in \N_0$ that
\begin{equation}\llabel{eq3}
\E\bigl[\|\Theta_n-\vartheta\|\bigr]\leq \fC (M^{-1/2}+\sqrt{\gamma_{n+1}}).
\end{equation}}
 \argument{\lref{def: XY};\lref{def: fL};}{that for all $\theta\in \cA$, $n,m\in \N$ it holds that
    \begin{equation}\llabel{evd1}
      \textstyle\frac 1M \sum_{m=1}^M(\nabla_\theta\fL)(\theta,X_{n,m})  = \textstyle\frac 1M \sum_{m=1}^M(\nabla_\theta\fL)(\theta,0)=(\nabla_\theta\fL)(\theta,0)=(\nabla\smalll)(\theta).
    \end{equation}}
    \argument{\lref{evd1};\lref{def: XY};\lref{eq3}}{that there exists $\fC\in \R$ such that for every $M\in \N$, every stochastic process $\bfm=(\bfm^{1},\dots,\bfm^{\delta})\colon \N_0\times\Omega\to\cA$, and every stochastic process $\Theta=(\Theta^{1},\dots,\Theta^{\delta})\colon \N_0\times\Omega\to\cA$ with the property that $\sup_{n\in \N_0}\|\Theta_n\|\leq 2\kappa+2\|\vartheta\|$ and with the property that for all $n\in \N$, $i\in \{1,2,\dots,\delta\}$ it holds that 
         \begin{equation}
            \bfm_0=0,\qquad \bfm_n= \alpha \bfm_{n-1}+(1-\alpha)(\nabla_\theta\smalll)(\Theta_{n-1}),
            \end{equation}
            \begin{equation}
           \text{and}\qquad  \Theta_n^{i}=\Theta_{n-1}^{i}-\gamma_n\TNewtonSchulzAlgorithm{\nscoe}{K}{\bfm_n^i}
         \end{equation}
     we have for all $n \in \N_0$ that
\begin{equation}\llabel{eq4}
\E\bigl[\|\Theta_n-\vartheta\|\bigr]\leq \fC (M^{-1/2}+\sqrt{\gamma_{n+1}}).
\end{equation}}
\argument{\lref{eq4};}{that there exists $\fC\in \R$ such that for every stochastic process $\bfm=(\bfm^{1},\dots,\bfm^{\delta})\colon \N_0\times\Omega\to\cA$ and every stochastic process $\Theta=(\Theta^{1},\dots,\Theta^{\delta})\colon \N_0\times\Omega\to\cA$ with the property that for all $n\in \N$, $i\in \{1,2,\dots,\delta\}$ it holds that 
         \begin{equation}
            \bfm_0=0,\qquad \bfm_n= \alpha \bfm_{n-1}+(1-\alpha)(\nabla_\theta\smalll)(\Theta_{n-1}),
            \end{equation}
            \begin{equation}
           \|\Theta_{n-1}\|\allowbreak\leq \kappa,\qquad \text{and}\qquad  \Theta_n^{i}=\Theta_{n-1}^{i}-\gamma_n\TNewtonSchulzAlgorithm{\nscoe}{K}{\bfm_n^i}
         \end{equation}
     we have for all $n \in \N_0$ that
\begin{equation}\llabel{eq5}
\E\bigl[\|\Theta_n-\vartheta\|\bigr]\leq \fC \sqrt{\gamma_{n+1}}.
\end{equation}}
\argument{\lref{eq5};}{\lref{conclude}\dott}
\end{aproof}
\begin{samepage}
\begin{tcolorbox}[colback=white!95!gray,
                  colframe=black,
                  boxrule=0.5pt,
                  sharp corners,
                  enhanced,
                  breakable,
                 ]
\begin{athm}{cor}{deterministic muon convergence}[\textcolor{red}{\MUON\ for deterministic optimization problems}]
      Let $\delta\in \N$, $K\in\N_0$, $d_1,d_2,\dots,d_\delta,\fd_1,\fd_2,\dots,\fd_\delta\in \N$, $\alpha\in (0,1)$, $\kappa\in (0,\infty)$, let $\nscoe=(\nscoe_{i,j})_{(i,j)\in (\N_0)^2}\colon(\N_0)^2\to \R$ satisfy $\textcolor{magenta}{\min\{\nscoe_{0,0},\nscoe_{0,1}\}>0}$ and $\textcolor{magenta}{\#(\nscoe^{ - 1 }( \R\backslash\{0\} )) <\infty}$,  assume for all $i\in \N\cap[0,K]$, $x\in \R$ that $\textcolor{magenta}{\sum_{j=0}^\infty\nscoe_{i,j}x^{2j}>0}$,
  let $\cA=\times_{i=1}^\delta\R^{d_i\times\fd_i}$, let $\smalll\in \allowbreak C^{1}(\cA,\R)$ be strongly convex, assume that $\nabla\smalll$ is Lipschitz continuous, let $\vartheta\in \cA$ satisfy $\smalll(\vartheta)=\inf_{\theta\in \cA}\smalll(\theta)$, let $(\gamma_n)_{n\in \N}\subseteq (0,\infty)$ be non-increasing, and assume $\limsup_{n\to\infty} (\gamma_n+(\gamma_n)^{-2}(\gamma_{n}-\gamma_{n+1}))=0$.
     Then there exists $\fC\in \R$ such that for every $\bfm=(\bfm^{1},\dots,\bfm^{\delta})\colon \N_0\to\cA$ and every $\Theta=(\Theta^{1},\dots,\Theta^{\delta})\colon \N_0\to\cA$ with the property that for all $n\in \N$, $i\in \{1,2,\dots,\delta\}$ it holds that 
         \begin{equation}\llabel{def: bfm}
            \textcolor{magenta}{\bfm_0=0},\qquad \textcolor{magenta}{\bfm_n= \alpha \bfm_{n-1}+(1-\alpha)(\nabla\smalll)(\Theta_{n-1})},
            \end{equation}
            \begin{equation}\llabel{def: Theta}
         \textcolor{magenta}{\|\Theta_0\|\leq \kappa},\qquad \text{and}\qquad  \textcolor{magenta}{\Theta_n^{i}=\Theta_{n-1}^{i}-\gamma_n\TNewtonSchulzAlgorithm{\nscoe}{K}{\bfm_n^i}}
         \end{equation}
     we have for all $n \in \N_0$ that
\begin{equation} \llabel{conclude}
     \textcolor{magenta}{\|\Theta_n-\vartheta\|\leq \fC \sqrt{\gamma_{n+1}}}
\end{equation}
(cf.\ \cref{definition: NS}).
\end{athm}
\end{tcolorbox}
\end{samepage}
\begin{aproof}
\argument{\cref{T_B_D};\cref{NS error}}{that there exist $\scrc,\fC\in (0,\infty)$ such that for all $i\in \{1,2,\dots,\delta\}$, $A\in \R^{d_i\times\fd_i}$ it holds that
    \begin{equation}\llabel{eq1}
   \lVert\TNewtonSchulzAlgorithm{\nscoe}{K}{A}\rVert\leq \scrc \qqandqq \spro{\TNewtonSchulzAlgorithm{\nscoe}{K}{A},A}\geq \frac{\fC\|A\|^2}{\|A\|+\nscoe_{0,1}(\nscoe_{0,0})^{-1}}.
    \end{equation}}
    \argument{the fact that $\limsup_{n\to\infty}\gamma_n=0$; the fact that $\limsup_{n\to\infty} \bigl(\frac{\gamma_{n}-\gamma_{n+1}}{(\gamma_{n})^2}\bigr)<\infty$;}{that
    \begin{equation}\llabel{eq3'}
        \textstyle\limsup_{n\to\infty} \bigl(\frac{\gamma_{n}-\gamma_{n+1}}{\gamma_{n}}\bigr)=0.
    \end{equation}}
    \argument{\lref{eq3'};}{that
    \begin{equation}\llabel{eq3'.1}
        \textstyle\lim_{n\to\infty} \bigl(\frac{\gamma_{n}}{\gamma_{n+1}}\bigr)=1.
    \end{equation}}
    \argument{\lref{eq3'.1};the fact that $0<\alpha<1$}{that
    \begin{equation}\llabel{eq3''}
        \limsup_{n\to\infty}(\allowbreak(\gamma_{n+1})^{-1}\alpha\gamma_n)= \alpha<1.
    \end{equation}}
    \argument{\lref{eq3''};\lref{eq1};\cref{theo: MUON bound};the assumption that $\smalll$ is strongly convex; the assumption that $\nabla\smalll$ is Lipschitz continuous}{that there exists $\fC\in \R$ such that for every $\bfm=(\bfm^{1},\dots,\bfm^{\delta})\colon \N_0\to\cA$ and every $\Theta=(\Theta^{1},\dots,\Theta^{\delta})\colon \N_0\to\cA$ with the property that for all $n\in \N$, $i\in \{1,2,\dots,\delta\}$ it holds that 
         \begin{equation}
          \bfm_0=0,\qquad \bfm_n= \alpha \bfm_{n-1}+(1-\alpha)(\nabla\smalll)(\Theta_{n-1}),
            \end{equation}
            \begin{equation}
    \|\Theta_0\|\leq \kappa,\qquad \text{and}\qquad  \Theta_n^{i}=\Theta_{n-1}^{i}-\gamma_n\TNewtonSchulzAlgorithm{\nscoe}{K}{\bfm_n^i}
         \end{equation}
     we have for all $n \in \N_0$ that
\begin{equation} \llabel{eq5}
   \|\Theta_n\|\leq \fC.
\end{equation}}
\argument{\cref{deterministic muon bound convergence};}{there exist $\fC_1,\fC_2\in \R$ such that for every $\bfm=(\bfm^{1},\dots,\bfm^{\delta})\colon\allowbreak \N_0\to\cA$ and every $\Theta=(\Theta^{1},\dots,\Theta^{\delta})\colon \N_0\to\cA$ with the property that for all $n\in \N$, $i\in \{1,2,\dots,\delta\}$ it holds that 
         \begin{equation}
            \bfm_0=0,\qquad \bfm_n= \alpha \bfm_{n-1}+(1-\alpha)(\nabla\smalll)(\Theta_{n-1}),
            \end{equation}
            \begin{equation}
       \|\Theta_{n-1}\|\leq \fC_1,\qquad \text{and}\qquad \Theta_n^{i}=\Theta_{n-1}^{i}-\gamma_n\TNewtonSchulzAlgorithm{\nscoe}{K}{\bfm_n^i}
         \end{equation}
     we have for all $n \in \N_0$ that
\begin{equation} \llabel{eq6}
\|\Theta_n-\vartheta\|\leq \fC_2 \sqrt{\gamma_{n+1}}.
\end{equation}}
\argument{\lref{eq5};\lref{eq6};}{\lref{conclude}\dott}
\end{aproof}
 \section{Approximative orthogonalization methods within MUON}\label{sec: practical}

In this section we show that the framework in \cref{main theorem} and \cref{muon stochastic convergence} is general enough to cover \MUON\ with five \NS\ steps with the original \MUON\ polynomial \cite{jordan2024muon} as implemented in the {\sc Pytorch} library \cite{pytorch_muon} (see \cref{lem: abc verify} in \cref{subsec: jordan muon} below) and to cover \MUON\  combined with the Polar Express method \cite{ar2505.16932} as implemented in the {\sc modded-nanoGPT} code \cite{Jordanmuoncode} (see \cref{lem: abc verify 2}  in \cref{subsec: polar express} below).

In addition, inspired by the smallness condition that for all $k \in \N_0$ it holds that $\delta_k \leq \delta < 1$ in \cite[Assumption 1, Theorem 2, and Theorem 4]{ar2510.19933} we establish in \cref{subsec: counter example}, in the situation of simple example optimization problems, lower bounds for the size of the error of orthogonalization methods within \MUON\ (see \cref{lem: counter example 1} and \cref{lem: counter example 3} below).
\subsection{Newton-Schulz iterations with the polynomial introduced in Jordan et al.\ (2024)}\label{subsec: jordan muon}
     \begin{tcolorbox}[colback=white!95!gray,
                  colframe=black,
                  boxrule=0.5pt,
                  sharp corners,
                  enhanced,
                  breakable,
                 ]
                 \begin{athm}{lemma}{lem: abc verify}[\textcolor{red}{5 \NS\ steps with the original \MUON\ polynomial \cite{jordan2024muon}}]
                    Let $\nscoe=(\nscoe_{i,j})_{(i,j)\in (\N_0)^2}\colon(\N_0)^2\to\R$ satisfy for all $i, j \in \N_0$, $k \in \{0,1,\dots,5\}$ with $\max\{ i, 2 j \} \geq 6$ that $\nscoe_{ i, j } = 0$ and
                   \begin{equation}\llabel{def: a}
                   \begin{split}
                      (\nscoe_{k,0},\nscoe_{k,1},\nscoe_{k,2})=
                       \begin{cases}
                       (1, 10^{-7},0) &\colon k=0,\\
(3.4445,-4.7750,2.0315) &\colon k>0.
                       \end{cases}
                       \end{split}
                   \end{equation}
                    Then 
                    \begin{enumerate}[label=(\roman*)]
                        \item \llabel{item 1} it holds that $\min\{\nscoe_{0,0},\nscoe_{0,1}\}>0=\nscoe_{6,0}$ and $\#(\nscoe^{ - 1 }( \R\backslash\{0\} )) <\infty$ and
                     \item \llabel{item 2} it holds for all $n\in \{1,2,\dots,5\}$, $x\in \R$ that $\sum_{i=0}^\infty \nscoe_{n,i}x^{2i}>0$.
                    \end{enumerate}
                 \end{athm}
                 \end{tcolorbox}
\begin{aproof}
     \argument{\lref{def: a}}{\lref{item 1}\dott} 
                    \startnewargseq
                    \argument{\lref{def: a}}{that for all $n\in\{1,2,\dots,5\}$ it holds that \llabel{arg2}$(\nscoe_{n,1})^2<4\nscoe_{n,1}\nscoe_{n,2}$\dott}
                    \argument{\lref{arg2};}{\lref{item 2}\dott}
\end{aproof}

\subsection{The Polar Express method introduced in Amsel et al.\ (2025)}\label{subsec: polar express}
  \begin{tcolorbox}[colback=white!95!gray,
                  colframe=black,
                  boxrule=0.5pt,
                  sharp corners,
                  enhanced,
                  breakable,
                 ]
                 \begin{athm}{lemma}{lem: abc verify 2}[\textcolor{red}{Polar Express method \cite{ar2505.16932}}]
                   Let $\nscoe=(\nscoe_{i,j})_{(i,j)\in (\N_0)^2}\colon(\N_0)^2\to\R$ satisfy for all $i, j \in \N_0$, $k \in \{0,1,\dots,5\}$ with $\max\{ i, 2 j \} \geq 6$ that $\nscoe_{ i, j } = 0$ and
                   \begin{equation}\llabel{def: a}
                   \begin{split}
                     & (\nscoe_{k,0},\nscoe_{k,1},\nscoe_{k,2})\\
                     &=
                       \begin{cases}
                       (1.02, 10^{-6},0) &\colon k=0,\\
(8.156554524902461,-22.48329292557795,15.878769915207462) &\colon k=1,\\
(4.042929935166739,-2.808917465908714,0.5000178451051316) &\colon k=2,\\
(3.8916678022926607,-2.772484153217685,0.5060648178503393) &\colon k=3,\\
(3.2857533657755655,-2.3681294933425376,0.46449024233003106) &\colon k=4,\\
(2.3465413258596377,-1.7097828382687081,0.42323551169305323) &\colon k=5.
                       \end{cases}
                       \end{split}
                   \end{equation}
                   
                   Then  \begin{enumerate}[label=(\roman*)]
                        \item \llabel{item 1} it holds that $\min\{\nscoe_{0,0},\nscoe_{0,1}\}>0=\nscoe_{6,0}$ and $\#(\nscoe^{ - 1 }( \R\backslash\{0\} )) <\infty$ and
                     \item \llabel{item 2} it holds for all $n\in \{1,2,\dots,5\}$, $x\in \R$ that $\sum_{i=0}^\infty \nscoe_{n,i}x^{2i}>0$.
                    \end{enumerate}
                 \end{athm}
                 \end{tcolorbox}
                 \begin{aproof}
                     \argument{\lref{def: a};the fact that $\nscoe_{6,0}=0$}{\lref{item 1}\dott} 
                    \startnewargseq
                    \argument{\lref{def: a}}{that for all $n\in\{1,2,\dots,5\}$ it holds that \llabel{arg2}$(\nscoe_{n,1})^2<4\nscoe_{n,1}\nscoe_{n,2}$\dott}
                    \argument{\lref{arg2};}{\lref{item 2}\dott}
                 \end{aproof}
\subsection{Error of orthogonalization methods within MUON}\label{subsec: counter example}

In the following results, \cref{lem: counter example 1} and \cref{lem: counter example 3}, we establish -- in the situation of simple example optimization problems -- lower bounds for the size of the error of orthogonalization methods within \MUON. Our conclusions of \cref{lem: counter example 1} and \cref{lem: counter example 3} are inspired by the smallness condition that for all $k \in \N_0$ it holds that $\delta_k \leq \delta < 1$ in \cite[Assumption 1, Theorem 2, and Theorem 4]{ar2510.19933}. Note that \cref{def: Phi: counter example 1} in \cref{lem: counter example 1} (cf.\ also \cref{def: Phi: counter example 3} in \cref{lem: counter example 3}) ensures that for all $\theta \in \R^{d\times\fd}$ it holds that $\Phi( \theta ) \in\operatorname{argmin}_{\vartheta\in\R^{d\times\fd},\,\nnorm{\vartheta}\leq 1}\spro{\theta,\vartheta}$.
                
\begin{tcolorbox}[colback=white!95!gray,
                  colframe=black,
                  boxrule=0.5pt,
                  sharp corners,
                  enhanced,
                  breakable,
                 ]
                 \begin{athm}{lemma}{lem: counter example 1}
                 Let $(\Omega,\cF,\P)$ be a probability space, let $d,\fd\in\N$, $K\in\N_0$, $\lambda,\fc\in (0,\infty)$, $\alpha\in (0,1)$, let $\smalll\colon \R^{d\times\fd}\times\R^{d\times\fd}\to\R$ satisfy for all $\theta,x\in \R^{d\times\fd}$ that $\smalll(\theta,x)=\lambda\|\theta-x\|^2$, for every $\theta \in \R^{ d \times \fd }$ let $\nnorm {\theta}  \in \R$ satisfy $\nnorm{ \theta}  = \sup_{ v \in \R^{ \fd },\, \|v\|\leq 1 } \| \theta v \| $ (spectral norm), let $\Phi\colon\R^{d\times\fd}\to\R^{d\times\fd}$ be measurable, assume for all $\theta\in\R^{d\times\fd}$ that
            \begin{equation}\label{def: Phi: counter example 1}
                    \textstyle \nnorm{ \Phi( \theta )} = 1 \qqandqq \spro{\theta, \Phi( \theta ) } = \inf_{ \vartheta\in\R^{d\times\fd},\,\nnorm{\vartheta}\leq 1 } \spro{ \theta, \vartheta }, 
                     \end{equation}
                     let $\nscoe=(\nscoe_{i,j})_{(i,j)\in (\N_0)^2}\colon(\N_0)^2\to \R$ satisfy $\min\{\nscoe_{0,0},\nscoe_{0,1}\}>0$ and $\#(\nscoe^{ - 1 }( \R\backslash\{0\} )) <\infty$, assume for all $i\in\N\cap[0,K]$, $x\in\R$ that $\sum_{j=0}^\infty \nscoe_{i,j}x^{2j}>0$, let $X_{n,m}\colon \Omega\to\R^{d\times\fd}$, $(n,m)\in \N^2$, be bounded \iid\ random variables, let $(\gamma_n)_{n\in \N}\subseteq (0,\infty)$ be non-increasing, assume $\limsup_{n\to\infty} (\gamma_n+(\gamma_n)^{-2}(\gamma_{n}-\gamma_{n+1}))=0$, for every $M\in\N$ let $\Theta^M\colon \N_0\times\Omega\to\R^{d\times\fd}$ and $\bfm^M\colon \N_0\times\Omega\to\R^{d\times\fd}$ satisfy for all
     $n\in \N$ that
         \begin{equation}\llabel{def: bfm}
            \bfm_0^M=0,\qquad \bfm_n^M= \alpha \bfm_{n-1}^M+(1-\alpha)\bigl[\textstyle \frac 1M \sum_{m=1}^M(\nabla_\theta\smalll)(\Theta_{n-1}^M,X_{n,m})\bigr],
            \end{equation}
            \begin{equation}\llabel{def: Theta}
         \P(\|\Theta_{0}^M\|< \fc)=1,\qquad \text{and}\qquad  \Theta_n^{M}=\Theta_{n-1}^{M}-\gamma_n\TNewtonSchulzAlgorithm{\nscoe}{K}{\bfm_n^{M}},
                \end{equation}
               assume for all $M\in\N$ that $\Theta_0^M$ and $(X_{n,m})_{(n,m)\in \N^2}$ are independent, and let $\delta\in (0,1)$ (cf.\ \cref{definition: NS}).
                     Then \begin{enumerate}[label=(\roman*)]
                         \item \label{item 1: counter example 1} there exists $\bfM\in\N$ such that for all $M\in\N\cap[\bfM,\infty)$ it holds that
                     \begin{equation}\llabel{conclude}
                       \textstyle \liminf_{n\to\infty} \P(\nnorm{\Phi(\bfm_n^M)+\TNewtonSchulzAlgorithm{\nscoe}{K}{\bfm_n^M}}>\delta)>\delta
                        \end{equation}
                        and
                        \item \label{item 2: counter example 1} there exist $M, n \in \N$ such that
                       $\P(\nnorm{\Phi(\bfm_n^M)+\TNewtonSchulzAlgorithm{\nscoe}{K}{\bfm_n^M}}\leq\delta)<1$.
                     \end{enumerate}
                 \end{athm}
                  \end{tcolorbox}
                 \begin{aproof}
\argument{\cref{NS error};}{that there exists $\fC\in (0,\infty)$ such that for all $\theta\in \R^{d\times\fd}$ that
\begin{equation}\llabel{eq3}
    \|\TNewtonSchulzAlgorithm{\nscoe}{K}{\theta}
    \|\leq \min\{\fC,\fC\|\theta\|\}.
\end{equation}}
\argument{\lref{eq3};}{that
\begin{equation}\llabel{eq4}
    \textstyle\limsup_{r\to 0} \sup_{\theta\in \R^{d\times\fd},\,\|\theta\|\leq r}\nnorm{\TNewtonSchulzAlgorithm{\nscoe}{K}{\theta}}=0.
\end{equation}}
\argument{\lref{eq4};\cref{def: Phi: counter example 1}}{that
\begin{equation}\llabel{eq5}
\begin{split}
 &  \textstyle \textstyle\liminf_{r\to 0} \inf_{\theta\in \R^{d\times\fd},\,\|\theta\|\leq r} \nnorm{\Phi(\theta)+\TNewtonSchulzAlgorithm{\nscoe}{K}{\theta}}\\
 &\textstyle\geq \bigl[\liminf_{r\to 0} \inf_{\theta\in \R^{d\times\fd},\,\|\theta\|\leq r} \nnorm{\Phi(\theta)}\bigr]-\bigl[\textstyle\limsup_{r\to 0} \sup_{\theta\in \R^{d\times\fd},\,\|\theta\|\leq r}\nnorm{\TNewtonSchulzAlgorithm{\nscoe}{K}{\theta}}\bigr]\\
   &=1-0=1.
   \end{split}
\end{equation}}
\argument{\lref{eq5};the fact that $\delta<1$}{that there exists $\varrho\in(0,\infty)$ which satisfies for all $\theta\in \R^{d\times\fd}$ with $\|\theta\|\leq \varrho$ it holds that
\begin{equation}\llabel{arg2} \nnorm{\Phi(\theta)+\TNewtonSchulzAlgorithm{\nscoe}{K}{\theta}}>\delta\dott
\end{equation}}
\startnewargseq
 \argument{the fact that $\limsup_{n\to\infty}\gamma_n=0$; the fact that $\limsup_{n\to\infty} \bigl(\frac{\gamma_{n}-\gamma_{n+1}}{(\gamma_{n})^2}\bigr)<\infty$;}{that
    \begin{equation}\llabel{eq3'}
        \textstyle\limsup_{n\to\infty} \bigl(\frac{\gamma_{n}-\gamma_{n+1}}{\gamma_{n}}\bigr)=0.
    \end{equation}}
    \argument{\lref{eq3'};}{that
    \begin{equation}\llabel{eq3'.1}
        \textstyle\limsup_{n\to\infty} \bigl(\frac{\gamma_{n}}{\gamma_{n+1}}\bigr)=1.
    \end{equation}}
    \argument{\lref{eq3'.1};the fact that $0<\alpha<1$}{that
    \begin{equation}\llabel{eq3''}
       \textstyle \limsup_{n\to\infty}(\allowbreak(\gamma_{n+1})^{-1}\alpha\gamma_n)= \alpha<1.
    \end{equation}}
\argument{\cref{theo: MUON bound};\lref{def: Theta};}{that there exists $\fC\in (0,\infty)$ such that for all $M\in\N$, $n\in\N_0$ it holds that 
\begin{equation}\llabel{eq6}
    \P(\|\Theta_n^M\|\leq \fC)=1.
\end{equation}}
\argument{\cref{def: Phi: counter example 1};\lref{def: Theta};\lref{eq3}}{that there exists $\fC\in (0,\infty)$ such that for all $M,n\in \N$ it holds that
    \begin{equation}\llabel{eqq3}
        \|\Theta_n^{M}-\Theta_{n-1}^{M}\|=\gamma_n\|\TNewtonSchulzAlgorithm{\nscoe}{K}{\bfm_n^{M}}\|\leq \gamma_n\fC.
    \end{equation}}
        \argument{\lref{def: Theta};induction}{that for all $M,n\in \N$ it holds that
        \begin{equation}\llabel{eqq6}
            \sigma(\Theta_{n}^{M})\subseteq\sigma \bigl(\Theta_0^{M},(X_{k,m})_{(k,m)\in \{1,2,\dots,n\}\times\N}\bigr).
        \end{equation}}
        \argument{\lref{eqq6};the assumption that for all $M\in \N$ it holds that $\Theta_{0}^{M}$ and $(X_{n,m})_{(n,m)\in \N^2}$ are independent; the fact that $X_{n,m}$, $(n,m)\in\N^2$, are independent}{that \llabel{argg1} for all $M,n\in \N$ it holds that $\Theta_{n-1}^{M}$ and $(X_{n,m})_{m\in \N}$ are independent\dott}
\argument{\lref{def: bfm};\lref{eq3''};\lref{eq6};\lref{eqq3};\lref{argg1};the fac that for all $\theta\in \R^{d\times\fd}$ it holds that $\smalll(\theta,x)=\lambda\|\theta-x\|^2$;\cref{lem: stochastic bfm analysis};}{that there exist $\fC\in (0,\infty)$ such that for all $M,n\in \N$ it holds that
    \begin{equation}\llabel{eq9}
       \textstyle\E\bigl[\|\bfm_{n}^{M}-2\lambda (\Theta_{n-1}^M-\E[X_{1,1}])\|\bigr]\leq \fC M^{-1/2}+\fC\gamma_n.
    \end{equation}}
\argument{\cref{muon stochastic convergence};}{that there exists $\fC\in (0,\infty)$ such that for all $M\in\N$, $n\in\N_0$ it holds that
\begin{equation}\llabel{eq10}
   \E\bigl[\|\Theta_{n}^M-\E[X_{1,1}]\|\bigr]\leq \fC (M^{-1/2}+\sqrt{\gamma_{n}}).
\end{equation}}
\argument{\lref{eq10};\lref{eq9};}{that there exists $\fC\in (0,\infty)$ which satisfies for all $M,n\in\N$ it holds that
\begin{equation}\llabel{def: fC}
\begin{split}
   \E\bigl[\|\bfm_n^M\|\bigr]&\leq \E\bigl[\|\bfm_{n}^{M}-2\lambda (\Theta_{n-1}^M-\E[X_{1,1}])\|\bigr]+\E\bigl[\|2\lambda(\Theta_{n-1}^M-\E[X_{1,1}])\|\bigr]\\
   &\leq \fC M^{-1/2}+\fC\sqrt{\gamma_{n}}+\fC\gamma_n.
   \end{split}
\end{equation}}
\startnewargseq
\argument{\lref{def: fC}; the Markov inequality}{that for all $M,n\in\N$ it holds that 
\begin{equation}\llabel{eq12}
    \P(\|\bfm_n^M\|> \varrho) \leq \varrho^{-1}\E\bigl[\|\bfm_n^M\|\bigr]\leq\fC\varrho^{-1} M^{-1/2}+\fC\varrho^{-1}\sqrt{\gamma_{n}}+\fC\varrho^{-1}\gamma_n.
\end{equation}}
\argument{\lref{eq12};\lref{arg2};}{that for all $M,n\in\N$ it holds that
\begin{equation}\llabel{eq13}
\begin{split}
    \P(\nnorm{\Phi(\bfm_n^M)+\TNewtonSchulzAlgorithm{\nscoe}{K}{\bfm_n^M}}>\delta)&\geq  \P(\|\bfm_n^M\|\leq  \varrho)=1- \P(\|\bfm_n^M\|> \varrho)\\
    &\geq 1-\fC\varrho^{-1} M^{-1/2}-\fC\varrho^{-1}\sqrt{\gamma_{n}}-\fC\varrho^{-1}\gamma_n.
    \end{split}
\end{equation}}
\argument{\lref{eq13};the fact that $\limsup_{n\to\infty}\gamma_n=0$}{thar for all $M\in\N$ it holds that
\begin{equation}\llabel{eq14}
   \liminf_{n\to\infty} \P(\nnorm{\Phi(\bfm_n^M)+\TNewtonSchulzAlgorithm{\nscoe}{K}{\bfm_n^M}}>\delta)\geq 1-\fC\varrho^{-1} M^{-1/2}.
\end{equation}}
\argument{\lref{eq14};}{\cref{item 1: counter example 1}\dott}
\startnewargseq
\argument{\cref{item 1: counter example 1};}{that there exists $M\in \N$ such that
\begin{equation}\llabel{eq15}
     \textstyle \liminf_{n\to\infty} \P(\nnorm{\Phi(\bfm_n^M)+\TNewtonSchulzAlgorithm{\nscoe}{K}{\bfm_n^M}}>\delta)>\delta>0.
\end{equation}}
\argument{\lref{eq15};}{that there exist $M, n\in\N$ such that
\begin{equation}\llabel{eq16}
  \textstyle  \P(\nnorm{\Phi(\bfm_n^M)+\TNewtonSchulzAlgorithm{\nscoe}{K}{\bfm_n^M}}>\delta)>0.
\end{equation}}
\argument{\lref{eq16};}{\cref{item 2: counter example 1}\dott}
                 \end{aproof}
                  \begin{tcolorbox}[colback=white!95!gray,
                  colframe=black,
                  boxrule=0.5pt,
                  sharp corners,
                  enhanced,
                  breakable,
                 ]
                 \begin{athm}{lemma}{lem: counter example 3}
                 Let $d,\fd\in\N$, $K\in\N_0$, $\lambda\in (0,\infty)$, $\vartheta\in \R^{d\times \fd}$, $\alpha\in (0,1)$, let $\smalll\colon \R^{d\times\fd}\to\R$ satisfy for all $\theta\in \R^{d\times\fd}$ that $\smalll(\theta)=\lambda\|\theta-\vartheta\|^2$, for every $\theta \in \R^{ d \times \fd }$ let $\nnorm {\theta}  \in \R$ satisfy $\nnorm{ \theta}  = \sup_{ v \in \R^{ \fd },\, \|v\|\leq  1 } \| \theta v \| $ (spectral norm), let $\Phi\colon\R^{d\times\fd}\to\R^{d\times\fd}$ be measurable, assume for all $\theta\in\R^{d\times\fd}$ that
            \begin{equation}\label{def: Phi: counter example 3}
                    \textstyle \nnorm{ \Phi( \theta )} = 1 \qqandqq \spro{\theta, \Phi( \theta ) } = \inf_{ v\in\R^{d\times\fd},\,\nnorm{v}\leq 1 } \spro{ \theta, v}, 
                     \end{equation}
                     let $\nscoe=(\nscoe_{i,j})_{(i,j)\in (\N_0)^2}\colon(\N_0)^2\to \R$ satisfy $\min\{\nscoe_{0,0},\nscoe_{0,1}\}>0$ and $\#(\nscoe^{ - 1 }( \R\backslash\{0\} )) <\infty$, assume for all $i\in\N\cap[0,K]$, $x\in\R$ that $\sum_{j=0}^\infty \nscoe_{i,j}x^{2j}>0$, let $(\gamma_n)_{n\in \N}\subseteq (0,\infty)$ be non-increasing, assume $\limsup_{n\to\infty} (\gamma_n+(\gamma_n)^{-2}(\gamma_{n}-\gamma_{n+1}))=0$, and let $\Theta\colon \N_0\to\R^{d\times\fd}$ and $\bfm\colon \N_0\to\R^{d\times\fd}$ satisfy for all
     $n\in \N$ that
         \begin{equation}\llabel{def: bfm}
            \bfm_0=0,\qquad \bfm_n= \alpha \bfm_{n-1}+(1-\alpha)(\nabla\smalll)(\Theta_{n-1}),
            \end{equation}
            \begin{equation}\llabel{def: Theta}
          \text{and}\qquad  \Theta_n=\Theta_{n-1}-\gamma_n\TNewtonSchulzAlgorithm{\nscoe}{K}{\bfm_n}
                \end{equation}
               (cf.\ \cref{definition: NS}).
                     Then
                     \begin{equation}\llabel{conclude}
                    \textstyle   \liminf_{n\to\infty}\nnorm{\Phi(\bfm_n)+\TNewtonSchulzAlgorithm{\nscoe}{K}{\bfm_n}}\geq 1.
                     \end{equation}
     \end{athm} 
     \end{tcolorbox}
    \begin{aproof}
          \argument{\cref{NS error};}{that there exists $\fC\in (0,\infty)$ such that for all $\theta\in \R^{d\times\fd}$ that
\begin{equation}\llabel{eq3}
    \|\TNewtonSchulzAlgorithm{\nscoe}{K}{\theta}
    \|\leq \fC\|\theta\|.
\end{equation}}
\argument{\lref{eq3};}{that
\begin{equation}\llabel{eq4}
    \textstyle\limsup_{r\to 0} \sup_{\theta\in \R^{d\times\fd},\,\|\theta\|\leq r}\nnorm{\TNewtonSchulzAlgorithm{\nscoe}{K}{\theta}}=0.
\end{equation}}
\argument{\lref{eq4};\cref{def: Phi: counter example 3}}{that
\begin{equation}\llabel{eq5}
\begin{split}
   &\textstyle \textstyle\liminf_{r\to 0} \inf_{\theta\in \R^{d\times\fd},\,\|\theta\|\leq r} \nnorm{\Phi(\theta)+\TNewtonSchulzAlgorithm{\nscoe}{K}{\theta}}\\
   &\textstyle\geq \bigl[\liminf_{r\to 0} \inf_{\theta\in \R^{d\times\fd},\,\|\theta\|\leq r} \nnorm{\Phi(\theta)}\bigr]
   -\textstyle\bigl[\limsup_{r\to 0} \sup_{\theta\in \R^{d\times\fd},\,\|\theta\|\leq r}\nnorm{\TNewtonSchulzAlgorithm{\nscoe}{K}{\theta}}\bigr]\\
   &=1-0=1.
   \end{split}
\end{equation}}
\argument{\cref{deterministic muon convergence}}{that there exists $\fC\in(0,\infty)$ such that for all $n\in\N_0$ it holds that
\begin{equation}\llabel{eq6}
    \|\Theta_n-\vartheta\|\leq \fC\sqrt{\gamma_{n+1}}.
\end{equation}}
\argument{\lref{eq6}; the fact that $\limsup_{n\to\infty} \gamma_n=0$;}{that
\begin{equation}\llabel{eq7}
   \textstyle \limsup_{n\to\infty}\|\Theta_n-\vartheta\|=0.
\end{equation}}
\argument{\lref{eq7};}{that there exists $\fC\in (0,\infty)$ which satisfies for all $n\in\N_0$ that
\begin{equation}\llabel{eq8}
    \|\Theta_n-\vartheta\|\leq \fC.
\end{equation}}
In the following we prove that 
\begin{equation}\llabel{need to prove}
  \textstyle  \limsup\limits_{n\to\infty}\biggl[\sum\limits_{k=1}^n\alpha^{n-k}\textstyle\|\Theta_{k-1}-\vartheta\|\biggr]=0. 
\end{equation}
\startnewargseq 
\argument{the fact that $0<\alpha<1$}{that for all $\varepsilon\in (0,\infty)$ there exists $N_\varepsilon\in\N$ such that for all $n\in\N\cap[N_\varepsilon,\infty)$ it holds that
\begin{equation}\llabel{def: N}
   \fC(1-\alpha)^{-1} \alpha^{N_\varepsilon}<\frac{\varepsilon}{2}.
\end{equation}}
\startnewargseq
\argument{\lref{eq7};}{that for all $\varepsilon \in (0,\infty)$ there exists $M_\varepsilon\in \N\cap[N_\varepsilon,\infty)$ which satisfies for all $n\in \N\cap[M_\varepsilon,\infty)$ that
\begin{equation}\llabel{eq9}
    \|\Theta_n-\vartheta\|\leq \frac{\varepsilon(1-\alpha)}{2}.
\end{equation}}
\startnewargseq
\argument{\lref{eq8};\lref{def: N};\lref{eq9}}{that for all $\varepsilon\in (0,\infty)$, $n\in \N\cap[2M_\varepsilon+1,\infty)$ it holds that
\begin{equation}\llabel{eq10}
\begin{split}
   &\textstyle\sum\limits_{k=1}^n\alpha^{n-k}\textstyle\|\Theta_{k-1}-\vartheta\|=\textstyle\biggl[\sum\limits_{k=1}^{n-M_\varepsilon}\alpha^{n-k}\textstyle\|\Theta_{k-1}-\vartheta\|\biggr]+\biggl[\sum\limits_{k=n-M_\varepsilon+1}^n\alpha^{n-k}\textstyle\|\Theta_{k-1}-\vartheta\|\biggr]\\
   &\leq\textstyle \biggl[\sum\limits_{k=1}^{n-M_\varepsilon}\alpha^{n-k}\textstyle\fC\biggr]+\biggl[\sum\limits_{k=n-M_\varepsilon+1}^n\alpha^{n-k}\textstyle\frac{\varepsilon(1-\alpha)}{2}\biggr]=\textstyle \fC\biggl[\sum\limits_{k=M_\varepsilon}^{n-1}\alpha^{k}\textstyle\biggr]+\displaystyle\frac{\varepsilon(1-\alpha)}{2}\textstyle\biggl[\sum\limits_{k=0}^{M_\varepsilon-1}\alpha^{k}\textstyle\biggr]\\
   &\leq \fC (1-\alpha)^{-1}\alpha^{M_\varepsilon}+\frac{\varepsilon(1-\alpha)}{2(1-\alpha)}\leq \fC (1-\alpha)^{-1}\alpha^{N_\varepsilon}+\frac{\varepsilon(1-\alpha)}{2(1-\alpha)}\leq \frac{\varepsilon}{2}+ \frac{\varepsilon}{2}=\varepsilon.
   \end{split}
\end{equation}}
\argument{\lref{eq10};}{\lref{need to prove}\dott}
\startnewargseq
 \argument{\lref{def: bfm};the assumption that for all $\theta\in\R^{d\times\fd}$ it holds that $\smalll(\theta)=\lambda|\theta-\vartheta|^2$;\cref{momentum:representation}}{that for all $n\in \N_0$ it holds that
    \begin{equation}\llabel{eqt1}
    \begin{split}
        \bfm_{n}&=\textstyle(1-\alpha)\biggl[\sum\limits_{k=1}^n\alpha^{n-k}(\nabla_\theta\smalll)(\Theta_{k-1})\biggr]=\textstyle(1-\alpha)\biggl[\sum\limits_{k=1}^n\alpha^{n-k}\textstyle2\lambda(\Theta_{k-1}-\vartheta)\biggr]\\
        &=\textstyle2\lambda(1-\alpha)\biggl[\sum\limits_{k=1}^n\alpha^{n-k}\textstyle(\Theta_{k-1}-\vartheta)\biggr].
        \end{split}
        \end{equation}}
\argument{\lref{eqt1};}{for all $n\in\N_0$ that
\begin{equation}\llabel{eqt2}
    \|\bfm_n\|\leq \textstyle2\lambda(1-\alpha)\biggl[\sum\limits_{k=1}^n\alpha^{n-k}\textstyle\|\Theta_{k-1}-\vartheta\|\biggr].
\end{equation}}
\argument{\lref{eqt2};\lref{need to prove}}{that
\begin{equation}\llabel{eqt3}
    \textstyle\limsup_{n\to\infty}\|\bfm_n\|=0.
\end{equation}}
\argument{\lref{eqt3};\lref{eq5}}{\lref{conclude}\dott}
\end{aproof}
\section{Non-convergence of MUON}\label{sec: non convergence}

The main purpose of this section is to prove \cref{MUON non-convergence 1d big batch} in \cref{subsec: muon non convergence} below. \cref{MUON non-convergence 1d big batch} shows for \MUON\ optimization methods applied to a simple one-dimensional quadratic \SOP\ that for almost every mini-batch size $M \in \N$ we have that the optimizer fails to converge to the solution of the \SOP\ (the unique global minimizer of the \SOP) as the number of gradient steps $n$ goes to infinity. \cref{MUON non-convergence 1d big batch} covers \MUON\ with an arbitrary finite number of \NS\ steps with the original \MUON\ polynomial \cite{jordan2024muon} as implemented in the {\sc Pytorch} library \cite{pytorch_muon} (see \cref{lem: abc verify} in \cref{subsec: jordan muon} above) and \MUON\ combined with the Polar Express method \cite{ar2505.16932} as implemented in the {\sc modded-nanoGPT} code \cite{Jordanmuoncode} (see \cref{lem: abc verify 2} in \cref{subsec: polar express} above) as special cases. 

Our overall strategy of the proof of \cref{MUON non-convergence 1d big batch} is inspired by the arguments in our proof of the non-convergence part in the \emph{\Adam\ symmetry theorem} in our preliminary work \cite[Theorem 1.1]{DeDoArPhi2025}. Specifically, we first establish in \cref{vector field} in \cref{subsec: muon limit theorem} a special case of a certain kind of \emph{\MUON\ limit theorem}, which shows for a simple one-dimensional quadratic \SOP\ that if \MUON\ converges to a vector $\vartheta$, then this vector $\vartheta$ must be a zero of the \emph{\MUON\ vector field} introduced in \cref{def: f: vector field} in \cref{vector field}. 

The purpose of \cref{subsec: property of stationary momentum}, \cref{subsec: sharp convergence rate smrv}, and \cref{subsec: non vanishing} is then to prove \cref{non-zero vector field 3} in \cref{subsec: non vanishing}, which shows for almost every mini-batch size $M \in\N$ that the solution of the \SOP\ is not a zero of the \emph{\MUON\ vector field}. In our proof of \cref{MUON non-convergence 1d big batch} in \cref{subsec: muon non convergence} we then combine the \emph{\MUON\ limit theorem} from \cref{vector field} with the fact that for almost every mini-batch size $M$ the solution of the \SOP\ is not a zero of the \emph{\MUON\ vector field} according to \cref{non-zero vector field 3} to see that that for almost every mini-batch size we have that \MUON\ can not converge to the solution of the \SOP.
\subsection{MUON-type limit theorem}\label{subsec: muon limit theorem}

In the next auxiliary result, \cref{Casero}, we briefly recall a simple consequence of the Stolz--Cesàro theorem on
Cesàro summation limits with weights. Only for the sake of completeness we include here a short self-contained proof
for \cref{Casero}. We employ \cref{Casero} in our proof of the \emph{limit theorem for \MUON-type methods} in \cref{vector field} below.

\begin{tcolorbox}[colback=white!95!gray,
                  colframe=black,
                  boxrule=0.5pt,
                  sharp corners,
                  enhanced,
                  breakable,
                 ]
\begin{athm}{lemma}{Casero}[\textcolor{red}{Cesàro summation limits with weights}]
let $d\in \N$, $\bfe\in \R^d$ and let $\gamma\colon\N\to(0,\infty)$ and $e\colon \N_0\to\R^d$ satisfy $\sum_{n=1}^\infty\gamma_n=\infty$ and  $\limsup_{n\to\infty}\|e_n-\bfe\|=0$. Then\begin{equation}\llabel{conclude}
   \textstyle \limsup_{n\to\infty}\big[\bigl(\sum_{k=1}^n\gamma_k\bigr)^{-1}\bigl(\sum_{k=1}^n\gamma_k \|e_k-\bfe\|\bigr)\bigr]=0.
\end{equation}
\end{athm}
\end{tcolorbox}
\begin{aproof}
\argument{the assumption that $\limsup_{n\to\infty}\|e_n-\bfe\|=0$}{that there exists $N \colon (0,\infty) \to \N$ which satisfies for all $\varepsilon\in (0,\infty)$ that 
\begin{equation}\llabel{def: N}
   \textstyle\sup_{n\in\N\cap[N_\varepsilon,\infty)} \|e_n-\bfe\|\leq \varepsilon.
\end{equation}}
\startnewargseq
\argument{\lref{def: N};the assumption that $\sum_{k=1}^\infty\gamma_k=\infty$}{that for all $\varepsilon\in (0,\infty)$ it holds that
\begin{equation}\llabel{eq1}
    \begin{split}
        &\textstyle \textstyle \limsup_{n\to\infty}\big[\bigl(\sum_{k=1}^n\gamma_k\bigr)^{-1}\bigl(\sum_{k=1}^n\gamma_k \|e_k-\bfe\|\bigr)\bigr]\\
&=\textstyle\limsup_{n\to\infty}\big[\bigl(\sum_{k=1}^{n}\gamma_k\bigr)^{-1}\bigl(\sum_{k=1}^{N_\varepsilon}\gamma_k \|e_k-\bfe\|\bigr)\bigr]\\
&\textstyle+\limsup_{n\to\infty}\big[\bigl(\sum_{k=1}^n\gamma_k\bigr)^{-1}\bigl(\sum_{k=N_\varepsilon+1}^{n}\gamma_k \|e_k-\bfe\|\bigr)\bigr]\\
&\leq \textstyle\limsup_{n\to\infty}\big[\bigl(\sum_{k=1}^{n}\gamma_k\bigr)^{-1}N_\varepsilon\bigl(\sup_{k\in \{1,2,\dots,N_\varepsilon\}}\gamma_k \|e_k-\bfe\|\bigr)\bigr]\\
&\textstyle+\limsup_{n\to\infty}\big[\bigl(\sum_{k=1}^n\gamma_k\bigr)^{-1}\bigl(\sum_{k=N_\varepsilon+1}^{n}\gamma_k \varepsilon\bigr)\bigr]\leq 0+\varepsilon=\varepsilon.
    \end{split}
\end{equation}}
\argument{\lref{eq1};}{\lref{conclude}\dott}
\end{aproof}

The work \cite[Corollary 1.10]{DereichJentzenKassing2025} (cf.\ \cite[Theorem 3.2]{DeDoArPhi2025}) establishes a result, which we refer to as \emph{\Adam\ limit theorem} that shows in a very wide generality that if the \Adam\ optimizer converges in probability to any random variable $\xi \colon \Omega \to\R^d$, then this random variable $\xi$ must almost surely be a zero of the \emph{\Adam\ vector field} introduced in \cite[Definition 2.4]{DereichAdamconvergence2024}. This \Adam\ limit theorem is then used in \cite{DeDoArPhi2025} to prove \emph{non-convergence for \Adam} in the \emph{\Adam\ symmetry theorem} in \cite[Theorem 1.1]{DeDoArPhi2025}. 

Inspired by this \Adam\ limit theorem and its use to derive non-convergence properties for \Adam, we establish in the next result, \cref{vector field}, in a very special case such a \emph{limit theorem for \MUON} and related methods. 

The \SGD\ optimization methods in \cref{vector field} are a general class of \MUON-type methods and cover the \MUON\ optimizer as a special case with an arbitrary number of \NS\ steps and the \NS\ polynomial as originally proposed in \cite{jordan2024muon} and as used in practice (see \cref{def: Pi: vector field} below, \cref{NS error} above, and \cref{NSLipschitz} below). The optimization problem class in \cref{vector field} is, however, very restrictive and only applies to a very specific one-dimensional quadratic example class of optimization problems.
\begin{samepage}
 \begin{tcolorbox}[colback=white!95!gray,
                  colframe=black,
                  boxrule=0.5pt,
                  sharp corners,
                  enhanced,
                  breakable,
                 ]
\begin{athm}{prop}{vector field}[\textcolor{red}{\MUON-type limit theorem}]
Let $(\Omega,\cF,\P)$ be a probability space, $\lambda\in (0,\infty)$, let $\smalll=(\smalll(\theta,x))_{(\theta,x)\in \R\times \R}\colon\R\times\R$ satisfy for all $\theta,x\in \R$ that $\smalll(\theta,x)=\lambda|\theta-x|^2$, let $\alpha\in (0,1)$, $M\in \N$, $\cK\in\R$, let $X_{n,m}\colon \Omega\to\R$, $(n,m)\in\N^2$, be bounded \iid\ random variables, let $\sgn\colon \R\to\R$ satisfy for all $x,y\in \R$ that
\begin{equation}\label{def: Pi: vector field}
 \textcolor{magenta}{ |\sgn(x)|\leq \cK}\qqandqq \textcolor{magenta}{|\sgn(x)-\sgn(y)|\leq \cK|x-y|},
\end{equation}
let $(\gamma_n)_{n\in\N}\subseteq(0,\infty)$ satisfy $\sup_{n\in \N}\gamma_n<\infty=\sum_{n=1}^\infty\gamma_n$,
 let $\Theta \colon\N_0\times\Omega\to\R$ and $\bfm\colon \N_0\times\Omega\to\R$ be stochastic processes, assume for all $n\in \N$ that 
         \begin{equation}\llabel{def: bfm}
          \textstyle\textcolor{magenta}{\bfm_0=0},\qquad \textcolor{magenta}{\bfm_n= \alpha \bfm_{n-1}+(1-\alpha)\bigl[\frac 1M\sum_{m=1}^M(\nabla_\theta\smalll)(\Theta_{n-1},X_{n,m})\bigr]},
            \end{equation}
            \begin{equation}\llabel{def: Theta}
           \textcolor{magenta}{|\Theta_0|\leq \cK},\qquad \text{and}\qquad  \textcolor{magenta}{\Theta_n=\Theta_{n-1}-\gamma_n\sgn(\bfm_n)},
         \end{equation}
         let $f\colon \R\to\R$ satisfy for all $\theta\in \R$ that
         \begin{equation}\label{def: f: vector field}
             \textcolor{magenta}{\textstyle f(\theta)=\E\bigl[\sgn\bigl((1-\alpha)\sum_{k=0}^\infty\alpha^k\bigl[\frac 1M\sum_{m=1}^M(\nabla_\theta\smalll)(\theta,X_{k+1,m})\bigr]\bigr)\bigr]},
         \end{equation}
         and let $\vartheta\in \R$ satisfy
           $\limsup_{n\to\infty} \E\bigl[|\Theta_n-\vartheta|\bigr]=0$.
Then 
  $\textcolor{magenta}{f(\vartheta)=0}$.
\end{athm}
\end{tcolorbox}
\end{samepage}
\begin{aproof}
    Throughout this proof for every $n\in \N$ let $\bbM_n\colon \Omega\to\R$ satisfy
    \begin{equation}\llabel{def: bbM}
        \textstyle\bbM_n=2\lambda(1-\alpha)\biggl[\sum\limits_{k=1}^{n}\alpha^{n-k}\bigl[\frac 1M\sum_{m=1}^M(\vartheta-X_{k,m})\bigr]\biggr],
    \end{equation} 
    for every $n\in \N\cup\{\infty\}$ let $\bfM_n\colon\Omega\to\R$ satisfy
    \begin{equation}\llabel{def: bfM}
       \textstyle \bfM_n=2\lambda(1-\alpha)\biggl[\sum\limits_{k=0}^{n-1}\alpha^{k}\bigl[\frac 1M\sum_{m=1}^M(\vartheta-X_{k+1,m})\bigr]\biggr],
    \end{equation}
    and let $\Gamma\in (0,\infty)$ satisfy
    \begin{equation}\llabel{def: Gamma}
       \textstyle \Gamma=\sup_{n\in \N}\gamma_n.
    \end{equation}
    \argument{\lref{def: bfm};the assumption that for all $\theta,x\in \R$ it holds that $\smalll(\theta,x)=\lambda|\theta-x|^2$;\cref{momentum:representation}}{that for all $n\in \N_0$ it holds that
    \begin{equation}\llabel{eq1}
    \begin{split}
        \bfm_{n}&=\textstyle(1-\alpha)\biggl[\sum\limits_{k=1}^n\alpha^{n-k}\textstyle \frac 1M \bigl[\textstyle\sum_{m=1}^M(\nabla_\theta\smalll)(\Theta_{k-1},X_{k,m})\bigr]\biggr]\\
        &=\textstyle(1-\alpha)\biggl[\sum\limits_{k=1}^n\alpha^{n-k}\textstyle \frac 1M \bigl[\textstyle\sum_{m=1}^M[2\lambda(\Theta_{k-1}-X_{k,m})]\bigr]\biggr]\\
        &=\textstyle2\lambda(1-\alpha)\biggl[\sum\limits_{k=1}^n\alpha^{n-k}\textstyle \frac 1M \bigl[\textstyle\sum_{m=1}^M(\Theta_{k-1}-X_{k,m})\bigr]\biggr].
        \end{split}
        \end{equation}}
        \argument{\lref{eq1};\lref{def: bbM};}{for all $n\in\N $ that
        \begin{equation}\llabel{eq2}
            \bfm_n-\bbM_n=\textstyle2\lambda(1-\alpha)\biggl[\sum\limits_{k=1}^n\alpha^{n-k}(\Theta_{k-1}-\vartheta)\biggr].
        \end{equation}}
        \argument{\lref{eq2};}{for all $n\in \N$ that
        \begin{equation}\llabel{eq3}
            \E\bigl[|\bfm_n-\bbM_n|\bigr]\leq \textstyle2\lambda(1-\alpha)\biggl[\sum\limits_{k=1}^n\alpha^{n-k}\E\bigl[|\Theta_{k-1}-\vartheta|\bigr]\biggr].
        \end{equation}}
        \argument{\cref{def: Pi: vector field};\lref{def: Theta};\lref{def: Gamma}}{for all $n\in \N$ that
        \begin{equation}\llabel{eq4}
            |\Theta_n-\Theta_{n-1}|=\gamma_n|\sgn(\bfm_n)|\leq \cK\gamma_n\leq \cK\Gamma.
        \end{equation}}
        \argument{\lref{eq4};\lref{def: Theta};induction;}{for all $n\in \N_0$ that
        \begin{equation}\llabel{eq5}
            |\Theta_n|\leq \cK\Gamma n+|\Theta_0|\leq \cK\Gamma n+\cK.
        \end{equation}}
        \argument{\lref{eq5};}{for all $n\in \N_0$ that
        \begin{equation}\llabel{eq6}
            |\Theta_n-\vartheta|\leq |\Theta_n|+|\vartheta|\leq  \cK\Gamma n+\cK+|\vartheta|.
        \end{equation}}
         In the following we prove that 
        \begin{equation}\llabel{need to prove}
            \limsup_{n\to\infty}\E\bigl[|\bfm_n-\bbM_n|\bigr]=0.
        \end{equation}
        \startnewargseq
          \argument{the fact that $\limsup_{n\to\infty} \E\bigl[|\Theta_n-\vartheta|\bigr]=0$;}{that for all $\varepsilon\in (0,\infty)$ there exists $\fN_\varepsilon\in \N$ which satisfies for all $n\in \N\cap[\fN_\varepsilon,\infty)$ that
        \begin{equation}\llabel{def: fN}
            \E\bigl[|\Theta_n-\vartheta|\bigr]\leq \frac{\varepsilon}{4\lambda}.
        \end{equation}}
        \startnewargseq
        \argument{\lref{eq6};\lref{def: fN}}{that there exists $\bfK\in (0,\infty)$ which satisfies for all $n\in\N_0$ that
        \begin{equation}\llabel{def: bfK}
            \E[|\Theta_n-\vartheta|]\leq \textstyle\max \bigl\{\frac{1}{4\lambda},\cK\Gamma \fN_1+\cK+|\vartheta|\bigr\}=\bfK.
        \end{equation}}
        \startnewargseq
   \argument{the fact that $0<\alpha<1$;}{that for all $\varepsilon\in (0,\infty)$ there exists $\bfN_\varepsilon\in \N$ which satisfies 
        \begin{equation}\llabel{def: bfN}
            4\lambda\bfK \alpha^{\bfN_\varepsilon}\leq \varepsilon.
        \end{equation}}
        \startnewargseq
        \argument{\lref{def: bfK};\lref{def: bfN}}{that for all $\varepsilon\in (0,\infty)$, $n\in \N\cap[\bfN_\varepsilon,\infty)$ it holds that
        \begin{equation}
            \begin{split}\llabel{ntp: eq1}
               & \textstyle2\lambda(1-\alpha)\biggl[\sum\limits_{k=1}^
                {n-\bfN_\varepsilon}\alpha^{n-k}\E\bigl[|\Theta_{k-1}-\vartheta|\bigr]\biggr]
                \leq \textstyle2\lambda(1-\alpha)\biggl[\sum\limits_{k=1}^
                {n-\bfN_\varepsilon}\alpha^{n-k}\bfK\biggr]\\
                &=\textstyle2\bfK\lambda(1-\alpha)\biggl[\sum\limits_{k=\bfN_\varepsilon}^
                {n-1}\alpha^{k}\biggr]\leq 2\lambda\bfK \alpha^{\bfN_\varepsilon}\leq \frac{\varepsilon}{2}.
            \end{split}
        \end{equation}}
        \argument{\lref{def: fN};}{that for all $\varepsilon\in (0,\infty)$, $n\in \N\cap[\bfN_\varepsilon+\fN_\varepsilon,\infty)$ it holds that
        \begin{equation}\llabel{ntp: eq2}
        \begin{split}
            \textstyle2\lambda(1-\alpha)\biggl[\sum\limits_{k=n-\bfN_\varepsilon+1}^
                {n}\alpha^{n-k}\E\bigl[|\Theta_{k-1}-\vartheta|\bigr]\biggr]\leq \textstyle2\lambda(1-\alpha)\biggl[\sum\limits_{k=n-\bfN_\varepsilon+1}^
                {n}\alpha^{n-k}\bigl(\frac{\varepsilon}{4\lambda}\bigr)\biggr]\leq \frac{\varepsilon}{2}.
                \end{split}
        \end{equation}}
        \argument{\lref{ntp: eq2};\lref{ntp: eq1}}{that for all $\varepsilon\in (0,\infty)$, $n\in \N\cap[\bfN_\varepsilon+\fN_\varepsilon,\infty)$ it holds that
        \begin{equation}\llabel{ntp: eq3}
        \begin{split}
              \textstyle2\lambda(1-\alpha)\biggl[\sum\limits_{k=1}^
                {n}\alpha^{n-k}\E\bigl[|\Theta_{k-1}-\vartheta|\bigr]\biggr]
                &\textstyle\leq   \textstyle2\lambda(1-\alpha)\biggl[\sum\limits_{k=1}^
                {n-\bfN_\varepsilon}\alpha^{n-k}\E\bigl[|\Theta_{k-1}-\vartheta|\bigr]\biggr]
                \\& +\textstyle2\lambda(1-\alpha)\biggl[\sum\limits_{k=n-\bfN_\varepsilon+1}^
                {n}\alpha^{n-k}\E\bigl[|\Theta_{k-1}-\vartheta|\bigr]\biggr]\\
                &\leq \frac{\varepsilon}{2}+\frac{\varepsilon}{2}=\varepsilon.
                \end{split}
        \end{equation}}
        \argument{\lref{ntp: eq3};\lref{eq3};}{for all $\varepsilon\in (0,\infty)$, $n\in \N\cap[\bfN_\varepsilon+\fN_\varepsilon,\infty)$ that
        \begin{equation}\llabel{ntp: eq4}
            \E\bigl[|\bfm_n-\bbM_n|\bigr]\leq \textstyle2\lambda(1-\alpha)\biggl[\sum\limits_{k=1}^n\alpha^{n-k}\E\bigl[|\Theta_{k-1}-\vartheta|\bigr]\biggr]\leq \varepsilon.
        \end{equation}}
        \argument{\lref{ntp: eq4};}{\lref{need to prove}\dott}
        \startnewargseq
         \argument{\cref{def: Pi: vector field};\lref{need to prove};}{that
        \begin{equation}\llabel{eq9}
             \limsup_{n\to\infty}\E\bigl[|\sgn(\bfm_n)-\sgn(\bbM_n)|\bigr]\leq  \limsup_{n\to\infty}\E\bigl[\cK|\bfm_n-\bbM_n|\bigr]=0. 
        \end{equation}}
        \argument{\lref{def: bfM};the assumption that $X_{n,m}$, $n,m\in \N^2$, are bounded \iid\@}{that
        \begin{equation}\llabel{eq7'}
            \limsup_{n\to\infty}\E\bigl[|\bfM_n-\bfM_\infty|\bigr]=0.
        \end{equation}}
        \argument{\lref{eq7'};\cref{def: Pi: vector field};}{that
        \begin{equation}\llabel{eq7}
            \limsup_{n\to\infty}\E\bigl[|\sgn(\bfM_n)-\sgn(\bfM_\infty)|\bigr]\leq   \limsup_{n\to\infty}\E\bigl[\cK|\bfM_n-\bfM_\infty|\bigr]=0.
        \end{equation}}
        \argument{\lref{def: bbM};\lref{def: bfM};the fact that $X_{n,m}$, $n,m\in \N^2$, are \iid\@}{that \llabel{arg1} for all $n\in \N$ it holds that $\bfM_n$ and $\bbM_n$ are identically distributed\dott}
        \argument{\lref{arg1};}{that for all $n\in \N$ it holds that
        \begin{equation}\llabel{eq7.1}
            \E[\sgn(\bbM_n)]=\E[\sgn(\bfM_n)].
        \end{equation}}
        \argument{\lref{eq9};\lref{eq7};\lref{eq7.1};}{that
        \begin{equation}\llabel{eq8}
        \begin{split}
            &\textstyle \limsup_{n\to\infty}\bigl|\E[\sgn(\bfm_n)]-\E[\sgn(\bfM_\infty)]\bigr|\\&\textstyle\leq \bigl(\limsup_{n\to\infty}\bigl|\E[\sgn(\bfm_n)]-\E[\sgn(\bbM_n)]\bigr|\bigr)+\bigl(\limsup_{n\to\infty}\bigl|\E[\sgn(\bbM_n)]-\E[\sgn(\bfM_\infty)]\bigr|\bigr)\\
            &\textstyle \textstyle\leq  \bigl(\limsup_{n\to\infty}\E\bigl[|\sgn(\bfm_n)-\sgn(\bbM_n)|\bigr]\bigr)+\bigl(\limsup_{n\to\infty}\bigl|\E[\sgn(\bfM_n)]-\E[\sgn(\bfM_\infty)]\bigr|\bigr)\\
            & \textstyle \textstyle\leq  \bigl(\limsup_{n\to\infty}\E\bigl[|\sgn(\bfm_n)-\sgn(\bbM_n)|\bigr]\bigr)+\bigl(\limsup_{n\to\infty}\E\bigl[|\sgn(\bfM_n)-\sgn(\bfM_\infty)|\bigr]\bigr)\\
            &=0.
             \end{split}
        \end{equation}}
         \argument{\lref{eq8};\cref{Casero}; the fact that $\sum_{n=1}^\infty\gamma_n=\infty$}{that
        \begin{equation}\llabel{eq12}
           \textstyle \limsup_{n\to\infty}\bigl(\sum_{k=1}^n\gamma_k)^{-1}\bigl(\sum_{k=1}^n\gamma_k\bigl|\E[\sgn(\bfm_k)]-\E[\sgn(\bfM_\infty)]\bigr|\bigr)=0.
        \end{equation}}
        \argument{\lref{eq12};}{that
        \begin{equation}\llabel{eq13}
        \begin{split}
        \textstyle \lim_{n\to\infty}\bigl(\sum_{k=1}^n\gamma_k)^{-1}\bigl(\sum_{k=1}^n\gamma_k\E[\sgn(\bfm_k)]\bigr)&=\textstyle \lim_{n\to\infty}\bigl(\sum_{k=1}^n\gamma_k)^{-1}\bigl(\sum_{k=1}^n\gamma_k\E[\sgn(\bfM_\infty)]\bigr)\\
        &=\E[\sgn(\bfM_\infty)].
        \end{split}
        \end{equation}}
        \argument{\lref{def: Theta};induction}{that for all $n\in \N$
 it holds that
 \begin{equation}\llabel{eq10}
           \textstyle \Theta_n-\Theta_0=-\sum_{k=1}^n\gamma_k\sgn(\bfm_k).
        \end{equation}}
        \argument{\lref{def: Theta};the fact that $\limsup_{n\to\infty} \E\bigl[|\Theta_n-\vartheta|\bigr]=0$;}{that
        \begin{equation}\llabel{eq10.1}
        \begin{split}
           &\textstyle \limsup_{n\to\infty}\E\bigl[|\Theta_n-\Theta_0|\bigr]\leq  \bigl(\limsup_{n\to\infty}\E\bigl[|\Theta_n-\vartheta|\bigr]\bigr)+ \bigl(\limsup_{n\to\infty}\E\bigl[|\vartheta-\Theta_0|\bigr]\bigr)\\
           &\textstyle= \limsup_{n\to\infty}\E\bigl[|\vartheta-\Theta_0|\bigr]= \E\bigl[|\vartheta-\Theta_0|\bigr]\leq \E\bigl[|\vartheta|\bigr]+\E\bigl[|\Theta_0|\bigr] \leq \|\vartheta\|+\cK<\infty.
           \end{split}
        \end{equation}}
        \argument{\lref{eq10};\lref{eq10.1};the fact that $\sum_{n=1}^\infty\gamma_n=\infty$}{that
        \begin{equation}\llabel{eq11}
\begin{split}
&\textstyle\limsup_{n\to\infty}\bigl(\sum_{k=1}^n\gamma_k)^{-1}\bigl|\sum_{k=1}^n\gamma_k\E[\sgn(\bfm_k)]\bigr|\\
&=\textstyle\limsup_{n\to\infty}\bigl(\sum_{k=1}^n\gamma_k)^{-1}\Bigl|\E\bigl[\sum_{k=1}^n\gamma_k\sgn(\bfm_k)\bigr]\Bigr|=\textstyle\limsup_{n\to\infty}\bigl(\sum_{k=1}^n\gamma_k)^{-1}\bigl|\E[\Theta_n-\Theta_0]\bigr|\\
&\leq \textstyle\limsup_{n\to\infty}\bigl(\sum_{k=1}^n\gamma_k)^{-1}\E\bigl[|\Theta_n-\Theta_0|\bigr]=0.
\end{split}
        \end{equation}}
        \argument{\lref{eq11};}{that
        \begin{equation}\llabel{eq11.1}
    \textstyle  \lim_{n\to\infty}\bigl(\sum_{k=1}^n\gamma_k)^{-1}\bigl(\sum_{k=1}^n\gamma_k\E[\sgn(\bfm_k)]\bigr)=0.
        \end{equation}}
        \argument{\lref{eq13};\lref{eq11.1};}{that
        \begin{equation}\llabel{eq14}
            \E[\sgn(\bfM_\infty)]=0.
        \end{equation}}
        \argument{\lref{eq14};\cref{def: f: vector field};\lref{def: bfM};the assumption that for all $\theta,x\in \R$ it holds that $\smalll(\theta,x)=\lambda|\theta-x|^2$}{that
        \begin{equation}
            f(\vartheta)=\E[\sgn(\bfM_\infty)]=0.
        \end{equation}}
\end{aproof}
\subsection{Properties of stationary momentum random variables (SMRVs)}\label{subsec: property of stationary momentum}

\begin{tcolorbox}[colback=white!95!gray,
                  colframe=black,
                  boxrule=0.5pt,
                  sharp corners,
                  enhanced,
                  breakable,
                 ]
\begin{athm}{lemma}{lem: sum of 2}
    Let $(\Omega,\cF,\P)$ be a probability space, let $M\in \N$, $N\in \N_0$, let $X_n\colon \Omega\to \R$, $n\in \N$, be \iid\ random variables, assume $\E\bigl[|X_{1}|^{2N+3}\bigr]<\infty$, and assume for all $n\in \{0,1,\dots,N\}$ that $\E[(X_1)^{2n+1}]=0$. Then 
    \begin{enumerate}[label=(\roman*)]
        \item \label{item 1: sum of 2} it holds for all $p\in \{0,1,\dots,N\}$ that
$
      \textstyle \E\bigl[\bigl(\sum_{m=1}^MX_m\bigr)^{2p+1}\bigr]=0$
    and 
    \item \label{item 2: sum of 2} it holds that 
$
      \textstyle \E\bigl[\bigl(\sum_{m=1}^MX_m\bigr)^{2N+3}\bigr]=M\E[(X_1)^{2N+3}]$.
    \end{enumerate}
\end{athm}
\end{tcolorbox}
\begin{aproof}
    \argument{the assumption that $\E\bigl[|X_{1}|^{2N+3}\bigr]<\infty$;the fact that $X_n$, $n\in \N$, are \iid\@;the multinomial theorem}{that for all $p\in \{0,1,\dots,N+1\}$ it holds that
    \begin{equation}\llabel{eq2}
    \begin{split}
&\textstyle\E\biggl[\biggl(\sum\limits_{m=1}^MX_m\biggr)^{2p+1}\biggr]\\
&\textstyle=\sum\limits_{k_1,k_2,\dots,k_M\in \N_0,\, k_1+k_2+\dots+k_M=2p+1}\frac{(2p+1)!}{(k_1)!(k_2)!\dots (k_M)!}\E[(X_1)^{k_1}(X_2)^{k_2}\dots(X_M)^{k_M}]\\
&\textstyle=\sum\limits_{k_1,k_2,\dots,k_M\in \N_0,\, k_1+k_2+\dots+k_M=2p+1}\frac{(2p+1)!}{(k_1)!(k_2)!\dots (k_M)!}\E[(X_1)^{k_1}]\,\E[(X_2)^{k_2}]\dots\E[(X_M)^{k_M}].
\end{split}
    \end{equation}}
    \argument{the fact that for all $q\in \{0,1,\dots,N\}$, $i\in\{1,2,\dots,M\}$ it holds that $\E[(X_i)^{2q+1}]=0$;the fact that for all $p\in \{0,1,\dots, N+1\}$, $k_1,k_2,\dots,k_M\in\N_0$ with $k_1+k_2+\dots+k_M=2p+1$ and $\max\{k_1,k_2,\dots,k_M\}<2N+3$ there exist $n\in \{1,2,\dots,M\}$, $q\in \{0,1,\dots,N\}$ such that $k_n=2q+1$}{that for all $p\in \{0,1,\dots,N+1\}$, $k_1,k_2,\dots,k_M\in \N_0$ with $k_1+k_2+\dots+k_M=2p+1$ and $\max\{k_1,k_2,\dots,k_M\}<2N+3$ it holds that
    \begin{equation}\llabel{eq3}
        \E[(X_1)^{k_1}]\,\E[(X_2)^{k_2}]\dots\E[(X_M)^{k_M}]=0.
    \end{equation}}
      \argument{\lref{eq3};\lref{eq2}; the assumption that $X_n$, $n\in \N$, are \iid\@}{that for all $p\in \{0,1,\dots,N\}$ it holds that
    \begin{equation}\llabel{eq4'}
    \begin{split}
&\textstyle\E\bigl[\bigl(\sum_{m=1}^MX_m\bigr)^{2p+1}\bigr]\\
&=\textstyle\sum\limits_{\substack{k_1,k_2,\dots,k_M\in \N_0,\, k_1+k_2+\dots+k_M=2p+1,\\
\max\{k_1,k_2,\dots,k_M\}<2N+3}}\frac{(2p+1)!}{(k_1)!(k_2)!\dots (k_M)!}\E[(X_1)^{k_1}]\,\E[(X_2)^{k_2}]\dots\E[(X_M)^{k_M}]\\&=0.
\end{split}
    \end{equation}}
    \argument{\lref{eq4'};}{\cref{item 1: sum of 2}\dott}
    \startnewargseq
    \argument{\lref{eq2};\lref{eq3}; the assumption that $X_n$, $n\in \N$, are \iid\@}{that
    \begin{equation}\llabel{eq4}
    \begin{split}
&\textstyle\E\bigl[\bigl(\sum_{m=1}^MX_m\bigr)^{2N+3}\bigr]\\
&=\textstyle\sum\limits_{\substack{k_1,k_2,\dots,k_M\in \N_0,\, k_1+k_2+\dots+k_M=2N+3,\\
\max\{k_1,k_2,\dots,k_M\}<2N+3}}\frac{(2N+3)!}{(k_1)!(k_2)!\dots (k_M)!}\E[(X_1)^{k_1}]\,\E[(X_2)^{k_2}]\dots\E[(X_M)^{k_M}]\\
&+\textstyle\sum\limits_{\substack{k_1,k_2,\dots,k_M\in \N_0,\, k_1+k_2+\dots+k_M=2N+3,\\
\max\{k_1,k_2,\dots,k_M\}=2N+3}}\frac{(2N+3)!}{(k_1)!(k_2)!\dots (k_M)!}\E[(X_1)^{k_1}]\,\E[(X_2)^{k_2}]\dots\E[(X_M)^{k_M}]\\
&\textstyle=\sum\limits_{\substack{k_1,k_2,\dots,k_M\in \N_0,\, k_1+k_2+\dots+k_M=2N+3,\\
\max\{k_1,k_2,\dots,k_M\}=2N+3}}\frac{(2N+3)!}{(k_1)!(k_2)!\dots (k_M)!}\E[(X_1)^{k_1}]\,\E[(X_2)^{k_2}]\dots\E[(X_M)^{k_M}]\\
&\textstyle=\sum_{m=1}^M\E[(X_m)^{2N+3}]=M\,\E[(X_1)^{2N+3}].
\end{split}
    \end{equation}}
    \argument{\lref{eq4};}{\cref{item 2: sum of 2}\dott}
\end{aproof}
\begin{tcolorbox}[colback=white!95!gray,
                  colframe=black,
                  boxrule=0.5pt,
                  sharp corners,
                  enhanced,
                  breakable,
                 ]
\begin{athm}{cor}{lem: non-negative moment}
    Let $(\Omega,\cF,\P)$ be a probability space, let $M\in \N$, $N\in \N_0$, $\alpha\in [0,1)$, let $X_{n,m}\colon \Omega\to\R$, $(n,m)\in \N^2$, be \iid\ random variables, assume $\E\bigl[|X_{1,1}|^{2N+3}\bigr]<\infty$, and assume $\sum_{p=0}^N|\E[(X_{1,1})^{2p+1}]|=0$. Then 
    \begin{enumerate}[label=(\roman*)]
    \item \label{item 0: non-negative moment} it holds that 
    \begin{equation}
          \textstyle\sum_{p=0}^{N+1} \E\bigl[\textstyle \bigl(\sum_{k=0}^{\infty}\alpha^k|\frac{1}{M}\sum_{m=1}^MX_{k+1,m}|\bigr)^{2p+1}\bigr]<\infty,
    \end{equation}
        \item \label{item 1: non-negative moment} it holds that
      \begin{equation}
           \textstyle\sum_{p=0}^N\bigl| \E\bigl[\textstyle \bigl(\sum_{k=0}^{\infty}\alpha^k[\frac{1}{M}\sum_{m=1}^MX_{k+1,m}]\bigr)^{2p+1}\bigr]\bigr|=0,
      \end{equation}
    and 
    \item \label{item 2: non-negative moment} it holds that 
  \begin{equation}
  \begin{split}
            &\E\bigl[\textstyle \bigl(\sum_{k=0}^{\infty}\alpha^k[\frac{1}{M}\sum_{m=1}^MX_{k+1,m}]\bigr)^{2N+3}\bigr]=M^{-2N-2}(1-\alpha^{2N+3})^{-1}\E[(X_{1,1})^{2N+3}].
            \end{split}
      \end{equation}
    \end{enumerate}
\end{athm}
\end{tcolorbox}
\begin{aproof}
    Throughout this proof for every $k\in \N_0$ let $\bbX_k\colon \Omega\to\R$ satisfy
    \begin{equation}\llabel{def: bbX}
       \displaystyle \bbX_k=\displaystyle\frac{1}{M} \textstyle\sum\limits_{m=1}^MX_{k+1,m}.
    \end{equation}
    \argument{\lref{def: bbX};the assumption that $X_{n,m}$, $(n,m)\in\N^2$, are \iid\@; the assumption that $\E[|X_{1,1}|^{2N+3}]\allowbreak<\infty$; the triangle inequality}{that there exists $\fC\in (0,\infty)$ such that for all $p\in\{1,2,\dots,2N+3\}$, $k\in\N_0$ it holds that
    \begin{equation}\llabel{eqt1}
    \begin{split}
\textstyle(\E[|\bbX_k|^p])^{1/p}&\textstyle=\bigl(\E\bigl[\bigl|\frac{1}{M} \textstyle\sum_{m=1}^MX_{k+1,m}\bigr|^p\bigr]\bigr)^{1/p}\leq \sum_{m=1}^M(\E[|M^{-1}X_{k+1,m}|^p])^{1/p}\\
        &\textstyle=M(\E[|M^{-1}X_{k+1,1}|^p])^{1/p}\leq \fC.
        \end{split}
    \end{equation}}
    \argument{\lref{eqt1};the triangle inequality}{that there exist $\fC_1,\fC_2\in (0,\infty)$ such that for all $p\in\{1,2,\dots,2N+3\}$, $K\in\N$ it holds that
    \begin{equation}\llabel{eqt2}
\textstyle\bigl(\E\bigl[\bigl(\sum_{k=0}^K\alpha^k|\bbX_k|\bigr)^p\bigr]\bigr)^{1/p}\leq \sum_{k=0}^K(\E[|\alpha^k\bbX_k|^p])^{1/p}\leq \sum_{k=0}^K\alpha^k\fC_1\leq\fC_2.  
\end{equation}}
\argument{\lref{eqt2};the monotone convergence theorem;}{that there exists $\fC\in (0,\infty)$ such that for all $p\in\{1,2,\dots,2N+3\}$ it holds that
\begin{equation}\llabel{eqt3'}
\textstyle\E\bigl[\bigl(\sum_{k=0}^\infty\alpha^k|\bbX_k|\bigr)^{p}\bigr]\leq \fC.
\end{equation}}
\argument{\lref{eqt3'}; the fact that $\textstyle\E\bigl[\bigl(\sum_{k=0}^\infty\alpha^k|\bbX_k|\bigr)^{0}\bigr]=(1-\alpha)^{-1}$}{that there exists $\fC\in (0,\infty)$ such that for all $p\in\{0,1,\dots,2N+3\}$ it holds that
\begin{equation}\llabel{eqt3}
\textstyle\E\bigl[\bigl(\sum_{k=0}^\infty\alpha^k|\bbX_k|\bigr)^{p}\bigr]\leq \fC.
\end{equation}}
\argument{\lref{eqt3};\lref{def: bbX}}{\cref{item 0: non-negative moment}\dott}
\startnewargseq
    \argument{\cref{item 0: non-negative moment};the multinomial theorem;the fact that $\bbX_k$, $k\in \N_0$, are \iid;}{that for all $p\in\{0,1,\dots,N+1\}$ it holds that
    \begin{equation}\llabel{eq0}
    \begin{split}
        \E\biggl[\textstyle \biggl(\sum\limits_{k=0}^{\infty}\alpha^k\bbX_k\biggr)^{2p+1}\biggr]
        &=\E\biggl[\textstyle \sum\limits_{(k_i)_{i\in \N_0}\subseteq\N_0,\, \sum_{i=0}^\infty k_i=2p+1}\frac{(2p+1)!\alpha^{\sum_{i=0}^\infty ik_i}}{\prod_{i=0}^\infty(k_i)!}\biggl(\prod\limits_{i=0}^\infty(\bbX_i)^{k_i}\biggr)\biggr]\\
       &=\textstyle \sum\limits_{(k_i)_{i\in \N_0}\subseteq\N_0,\, \sum_{i=0}^\infty k_i=2p+1}\frac{(2p+1)!\alpha^{\sum_{i=0}^\infty ik_i}}{\prod_{i=0}^\infty(k_i)!}\E\biggl[\prod\limits_{i=0}^\infty(\bbX_i)^{k_i}\biggr].
    \end{split}
    \end{equation}}
    \argument{\lref{def: bbX}; the assumption that $\E\bigl[|X_{1,1}|^{2N+3}\bigr]<\infty$;the assumption that for all $n\in \{0,1,\dots,N\}$ it holds that $\E[(X_{1,1})^{2n+1}]=0$;\cref{item 1: sum of 2} in \cref{lem: sum of 2}}{that for all $p\in \{0,1,\dots,N\}$, $k\in \N_0$ it holds that
    \begin{equation}\llabel{eq1}
      \begin{split}
    \E\bigl[(\bbX_k)^{2p+1}\bigr]&=\textstyle\E\bigl[\bigl(\frac{1}{M} \textstyle\sum_{m=1}^MX_{k+1,m}\bigr)^{2p+1}\bigr]=M^{-2p-1}\E\bigl[\bigl( \textstyle\sum_{m=1}^MX_{k+1,m}\bigr)^{2p+1}\bigr]=0.
      \end{split}  
    \end{equation}}
    \argument{the fact that $\bbX_k$, $k\in \N_0$, are \iid\@;\lref{eq1}}{that for all $p\in \{0,1,\dots,N\}$, $(k_i)_{i\in\N_0} \subseteq\N_0$ with $\sum_{i=0}^\infty k_i=2p+1$ there exist $j\in \N_0$, $r\in\{0,1,\dots,N\}$ such that $k_j=2r+1$ and 
    \begin{equation}\llabel{eq1.2}
    \begin{split}
\textstyle\E\bigl[\prod_{i=0}^\infty(\bbX_i)^{k_i}\bigr]
&=\textstyle\E\bigl[\bigl(\prod_{i\in \N_0\backslash\{j\}}(\bbX_i)^{k_i}\bigr)(\bbX_j)^{k_j}\bigr]=\textstyle\E\bigl[\bigl(\prod_{i\in \N_0\backslash\{j\}}(\bbX_i)^{k_i}\bigr)(\bbX_j)^{2r+1}\bigr]\\&=\textstyle\E\bigl[\bigl(\prod_{i\in \N_0\backslash\{j\}}(\bbX_i)^{k_i}\bigr)\bigr]\E[(\bbX_j)^{2r+1}]=0.
\end{split}
    \end{equation}}
      \argument{\lref{eq0};\lref{eq1.2}}{that for all $p\in \{0,1,\dots,N\}$ it holds that
    \begin{equation}\llabel{eq3}
    \begin{split}
       & \E\biggl[\textstyle \biggl(\sum\limits_{k=0}^{\infty}\alpha^k\bbX_k\biggr)^{2p+1}\biggr]=\textstyle \sum\limits_{(k_i)_{i\in \N_0}\subseteq\N_0,\, \sum_{i=0}^\infty k_i=2p+1}\frac{(2p+1)!\alpha^{\sum_{i=0}^\infty ik_i} 0}{\prod_{i=0}^\infty(k_i)!}
        =0.
    \end{split}
    \end{equation}}
    \argument{\lref{eq3};\lref{def: bbX}}{\cref{item 1: non-negative moment}\dott}
    \startnewargseq
    \argument{\lref{def: bbX}; the assumption that $\E\bigl[|X_{1,1}|^{2N+3}\bigr]<\infty$;the assumption that for all $n\in \{0,1,\dots,N\}$ it holds that $\E[(X_{1,1})^{2n+1}]=0$;\cref{item 2: sum of 2} in \cref{lem: sum of 2}}{that for all $k\in \N_0$ it holds that
    \begin{equation}\llabel{eq2}
      \begin{split}
    \E\bigl[(\bbX_k)^{2N+3}\bigr]&=\textstyle\E\bigl[\bigl(\frac{1}{M} \textstyle\sum_{m=1}^MX_{k+1,m}\bigr)^{2N+3}\bigr]=M^{-2N-3}\E\bigl[\bigl( \textstyle\sum_{m=1}^MX_{k+1,m}\bigr)^{2N+3}\bigr]\\
          &=M^{-2N-3}M\E\bigl[(X_{1,1})^{2N+3}]=M^{-2N-2}\E\bigl[(X_{1,1})^{2N+3}].
      \end{split}  
    \end{equation}
    }
    \argument{\lref{eq1};the fact that $\bbX_k$, $k\in \N_0$, are \iid\@}{that for all $(k_i)_{i\in\N_0} \subseteq\{0,1,\dots,2N+2\}$ with $\sum_{i=0}^\infty k_i=2N+3$ there exist $j\in \N_0$, $r\in\{0,1,\dots,N\}$ such that $k_j=2r+1$ and 
    \begin{equation}\llabel{eq2.2}
    \begin{split}
\textstyle\E\bigl[\prod_{i=0}^\infty(\bbX_i)^{k_i}\bigr]&=\textstyle\E\bigl[\bigl(\prod_{i\in \N_0\backslash\{j\}}(\bbX_i)^{k_i}\bigr)(\bbX_j)^{k_j}\bigr]=\textstyle\E\bigl[\bigl(\prod_{i\in \N_0\backslash\{j\}}(\bbX_i)^{k_i}\bigr)(\bbX_j)^{2r+1}\bigr]\\&=\textstyle\E\bigl[\bigl(\prod_{i\in \N_0\backslash\{j\}}(\bbX_i)^{k_i}\bigr)\bigr]\E[(\bbX_j)^{2r+1}]=0.
\end{split}
    \end{equation}}
\argument{\lref{eq0};\lref{eq2.2};\lref{eq2}}{that
\begin{equation}\llabel{eq4}
\begin{split}
    &\E\biggl[\textstyle \biggl(\sum\limits_{k=0}^{\infty}\alpha^k\bbX_k\biggr)^{2N+3}\biggr]\\
       &= \textstyle \sum\limits_{(k_i)_{i\in \N_0}\subseteq\N_0,\, \sum_{i=0}^\infty k_i=2N+3,\,\max_{i\in\N_0}k_i<2N+3}\frac{(2N+3)!\alpha^{\sum_{i=0}^\infty ik_i}}{\prod_{i=0}^\infty(k_i)!}\E\biggl[\prod\limits_{i=0}^\infty(\bbX_i)^{k_i}\biggr]\\
       &\textstyle+\sum\limits_{(k_i)_{i\in \N_0}\subseteq\N_0,\, \sum_{i=0}^\infty k_i=2N+3,\,\max_{i\in\N_0}k_i=2N+3}\frac{(2N+3)!\alpha^{\sum_{i=0}^\infty ik_i}}{\prod_{i=0}^\infty(k_i)!}\E\biggl[\prod\limits_{i=0}^\infty(\bbX_i)^{k_i}\biggr]\\
       &\textstyle=\sum\limits_{(k_i)_{i\in \N_0}\subseteq\N_0,\, \sum_{i=0}^\infty k_i=2N+3,\,\max_{i\in\N_0}k_i=2N+3}\frac{(2N+3)!\alpha^{\sum_{i=0}^\infty ik_i}}{\prod_{i=0}^\infty(k_i)!}\E\biggl[\prod\limits_{i=0}^\infty(\bbX_i)^{k_i}\biggr]\\
       &\textstyle=\sum\limits_{m=0}^\infty \bigl(\alpha^{m(2N+3)}\E\bigl[(\bbX_m)^{2N+3}\bigr]\bigr)=\textstyle\sum\limits_{m=0}^\infty \bigl(\alpha^{m(2N+3)} M^{-2N-2}\E\bigl[(X_{1,1})^{2N+3}]\bigr)\\
       &\textstyle=(1-\alpha^{2N+3})^{-1}M^{-2N-2}\E\bigl[(X_{1,1})^{2N+3}]. 
       \end{split}
\end{equation}}
\argument{\lref{eq4};\lref{def: bbX}}{\cref{item 2: non-negative moment}\dott}
\end{aproof}
\begin{tcolorbox}[colback=white!95!gray,
                  colframe=black,
                  boxrule=0.5pt,
                  sharp corners,
                  enhanced,
                  breakable,
                 ]
\begin{athm}{lemma}{lem: upper order of Adam coefficient 1}
    Let $(\Omega,\cF,\P)$ be a probability space, let $p\in \N\backslash\{1\}$, $\alpha\in [0,1)$, $\varepsilon\in (0,\infty)$, $\delta\in \N_0$, let $X_{n,m}\colon \Omega\to\R$, $(n,m)\in \N^2$, be bounded centered \iid\ random variables, and let $\psi\colon \R\to\R$ satisfy for all $x,y\in \R$, $\kappa\in (0,\infty)$ that $|\psi(x)-\psi(y)|\leq |x-y|$, $|\psi(x)|\leq |x|^\delta$, and $\psi(\kappa x)=\kappa^\delta \psi(x)$. Then
    \begin{enumerate}[label=(\roman*)]
        \item \label{item 1: upper order of Adam coefficient 1} it holds that $\sup_{M\in\N}\E\bigl[\big((1-\alpha)\sum_{k=0}^\infty \alpha^k|M^{-1/2}\sum_{m=1}^MX_{k+1,m}|\bigr)^{p}\bigr]<\infty$ and
        \item \label{item 2: upper order of Adam coefficient 1} there exists $\bfK\in \N$ which satisfies for all $K\in \N\cap[\bfK,\infty)$, $M\in \N$ that
      \begin{equation}\llabel{conclude}
      \begin{split}
           &M^{\nicefrac{(p+\delta)}{2}}\bigl|\E\bigl[\textstyle \psi\bigl(\sum_{k=0}^{\infty}\alpha^k(\frac{1}{M}\sum_{m=1}^MX_{k+1,m})\bigr)|(1-\alpha)\sum_{k=0}^\infty \alpha^k[\frac{1}{M}\sum_{m=1}^MX_{k+1,m}]|^{p}\bigr]\\
           &-\E\bigl[\textstyle \psi\bigl(\sum_{k=0}^{K}\alpha^k(\frac{1}{M}\sum_{m=1}^MX_{k+1,m})\bigr)|(1-\alpha)\sum_{k=0}^K \alpha^k[\frac{1}{M}\sum_{m=1}^MX_{k+1,m}]|^{p}\bigr]\bigr|\leq \varepsilon.
           \end{split}
      \end{equation}
    \end{enumerate}
\end{athm}
\end{tcolorbox}
\begin{aproof}
      Throughout this proof assume without loss of generality that $\alpha>0$, let $c\in (0,\infty)$ satisfy
   \begin{equation}\llabel{def: c}
       \P(|X_{1,1}|\leq c)=1,
   \end{equation}
   for every $M\in\N$, $k\in \N_0$ let $\bbX_k^M\colon \Omega\to\R$ satisfy
    \begin{equation}\llabel{def: bbX}
       \displaystyle \bbX_k^M=M^{-1/2} \textstyle\sum_{m=1}^MX_{k+1,m},
    \end{equation} 
   and for every $K\in \N\cup\{\infty\}$, $M\in \N$ let $A_K^M\colon \Omega\to\R$  satisfy
   \begin{equation}\llabel{def: AB}
      \textstyle A_K^M=\sum_{k=0}^{K}\alpha^k\bbX_k^M.
   \end{equation}
    \argument{\lref{def: bbX};\lref{def: AB}; the fact that for all $\kappa\in [0,\infty)$, $x\in \R$ it holds that $\psi(\kappa x)=\kappa^\delta \psi(x)$}{that for all $M\in \N$, $K\in \N\cup\{\infty\}$ it holds that
    \begin{equation}\llabel{eq1}
    \begin{split}
          &M^{\nicefrac{(p+\delta)}{2}}\Bigl(\E\Bigl[\textstyle \psi\bigl(\sum_{k=0}^{K}\alpha^k(\frac{1}{M}\sum_{m=1}^MX_{k+1,m})\bigr)\bigl|(1-\alpha)\sum_{k=0}^K \alpha^k[\frac{1}{M}\sum_{m=1}^MX_{k+1,m}]\bigr|^{p}\Bigr]\Bigr)\\
           &=\E\Bigl[\textstyle \psi\bigl(\sum_{k=0}^{K}\alpha^k(M^{-1/2}\sum_{m=1}^MX_{k+1,m})\bigr)\bigl|(1-\alpha)\sum_{k=0}^K\alpha^k[M^{-1/2}\sum_{m=1}^MX_{k+1,m}]\bigr|^{p}\Bigr]\\
          &= \E\Bigl[\textstyle \psi\bigl(\sum_{k=0}^{K}\alpha^k\bbX_k^M\bigr)\bigl|(1-\alpha)\sum_{k=0}^K \alpha^k\bbX_k^M\bigr|^{p}\Bigr]= (1-\alpha)^{p}\E\bigl[\psi(A^M_K)|A^M_K|^{p}\bigr].
          \end{split}
    \end{equation}}
    \argument{\lref{def: AB};the assumption that for all $x,y\in \R$ it holds that $|\psi(x)-\psi(y)|\leq |x-y|$;the fact that for all $x$ it holds that $|\psi(x)|\leq |x|^\delta$;the fact that for all $x,y\in [0,\infty)$ it holds that $||x|^{p}-|y|^{p}|\leq p|x-y|(x^{p-1}+y^{p-1})$;the triangle inequality;the Cauchy-Schwarz inequality}{that for all $M,K\in\N $ it holds that
    \begin{equation}\llabel{eq2}
    \begin{split}
       &\Bigl|\E\bigl[\psi(A^M_\infty)|A^M_\infty|^{p}\bigr]-\E\bigl[\psi(A^M_K)
    |A^M_K|^{p}\bigr]\Bigr|\\
       &\leq \Bigl|\E\bigl[\psi(A^M_\infty)|A^M_\infty|^{p}\bigr]-\E\bigl[\psi(A^M_K)|A^M_\infty|^{p}\bigr]\Bigr|+\Bigl|\E\bigl[\psi(A^M_K)|A^M_\infty|^{p}\bigr]-\E\bigl[\psi(A^M_K)|A^M_K|^{p}\bigr]\Bigr|\\
       &\leq \E\bigl[|A^M_\infty-A_K^M||A^M_\infty|^{p}\bigr]+\E\bigl[|A^M_K|^{\delta}\bigl||A^M_\infty|^{p}-|A^M_K|^{p}\bigr|\bigr]\\
       &\leq \E\bigl[|A^M_\infty-A_K^M||A^M_\infty|^{p}\bigr]+\E\bigl[p|A^M_K|^\delta(|A_\infty^M|^{p-1}+|A_K^M|^{p-1})|A_\infty^M-A_K^M|\bigr]\\
       &\leq \bigl(\E[(A^M_\infty-A_K^M)^2]\bigr)^{1/2}\bigl(\E[(A^M_\infty)^{2p}]\bigr)^{1/2}\\
       &+p\bigl(\E[|A^M_K|^{4\delta}]\bigr)^{1/4}\bigl(\E[(|A_\infty^M|^{p-1}+|A_K^M|^{p-1})^4]\bigr)^{1/4}\bigl(\E\bigl[(A^M_\infty-A^M_K)^{2}\bigr]\bigr)^{1/2}\\
       &\leq\bigl(\E[(A^M_\infty-A_K^M)^2]\bigr)^{1/2}\bigl(\E[(A^M_\infty)^{2p}]\bigr)^{1/2}\\
       &+p\bigl(\E[|A^M_K|^{4\delta}]\bigr)^{1/4}\Bigl(\bigl(\E[(A_\infty^M)^{4p-4}]\bigr)^{1/4}+\bigl(\E[(A_K^M)^{4p-4}]\bigr)^{1/4}\Bigr)\bigl(\E\bigl[(A^M_\infty-A^M_K)^{2}\bigr]\bigr)^{1/2}.
    \end{split}
    \end{equation}}
    \argument{the assumption that $X_{n,m}$, $(n,m)\in \N^2$, are \iid;the fact that $\E[X_{1,1}]=0$;\lref{def: c};\lref{def: bbX};\lref{def: AB}; the multinomial theorem}{that for all $M,K\in\N$ it holds that
    \begin{equation}\llabel{evd1}
    \begin{split}
       & \E[(A^M_\infty-A_K^M)^2]=\E\textstyle\bigl[\bigl(\sum_{k=K+1}^\infty \alpha^k\bbX_k^M\bigr)^2\bigr]=\E\bigl[\sum_{k,i=K+1}^\infty\alpha^k\bbX_k^M\alpha^i\bbX_i^M\bigr]\\
&\textstyle =\sum_{k,i=K+1}^\infty \E\bigl[\alpha^k\bbX_k^M\alpha^i\bbX_i^M\bigr]=\sum_{k=K+1}^\infty \alpha^{2k}\E[(\bbX_k^M)^2] =\alpha^{2K+2}(1-\alpha^2)^{-1}\E\bigl[(\bbX_1^M)^2\bigr]\\
&=\textstyle\alpha^{2K+2}(1-\alpha^2)^{-1}M^{-1}\E\bigl[\bigl(\sum_{m=1}^M X_{1,m}\bigr)^2\bigr]=\alpha^{2K+2}(1-\alpha^2)^{-1}\E[(X_{1,1})^2]\\
&\leq \alpha^{2K+2}(1-\alpha^2)^{-1}c^2.
        \end{split}
    \end{equation}}
    \argument{the assumption that $X_{n,m}$, $(n,m)\in \N^2$, are \iid;the fact that $\E[X_{1,1}]=0$;\lref{def: c};\cite[Theorem 3]{MR271721}}{that there exists $\fC\colon [0,\infty)\to (0,\infty)$ which satisfies for all $k\in \N_0$, $M\in \N$, $q\in [2,\infty)$ that
    \begin{equation}\llabel{def: fC}
    \begin{split}
      &\textstyle  \E\Bigl[\bigl|\sum_{m=1}^MX_{k+1,m}\bigr|^q\Bigr]\textstyle\leq \fC(q)\max\bigl\{\sum_{m=1}^M\E\bigl[|X_{k+1,m}|^q\bigr],\bigl(\sum_{m=1}^M \E[(X_{k+1,m})^2]\bigr)^{q/2}\bigr\} \\
      &=\textstyle \fC(q)\max\bigl\{M\E\bigl[|X_{k+1,1}|^q\bigr],\bigl(M\E[(X_{k+1,1})^2]\bigr)^{q/2}\bigr\}
      \leq \fC(q)\max\bigl\{Mc^q,\bigl(Mc^2\bigr)^{q/2}\bigr\}\\
      &\leq \textstyle\fC(q)\max\bigl\{Mc^q,\bigl(Mc^2\bigr)^{q/2}\bigr\}\leq \fC(q) M^{q/2}c^q.
      \end{split}
    \end{equation}}
    \startnewargseq
    \argument{\lref{def: bbX};\lref{def: fC}}{that for all $k\in \N_0$, $M\in \N$, $q\in [2,\infty)$ it holds that
    \begin{equation}\llabel{eq3}
\E\bigl[|\bbX_k^M|^q\bigr]=\E\Bigl[\textstyle\bigl|M^{-1/2}\sum_{m=1}^MX_{k+1,m}\bigr|^{q}\Bigr]=M^{-q/2}\E\Bigl[\textstyle\bigl|\sum_{m=1}^MX_{k+1,m}\bigr|^{q}\Bigr]\leq \fC(q)c^q.
    \end{equation}}
    \argument{\lref{eq3};the multinomial theorem; the fact that for all $M\in \N$ it holds that $\bbX_k^M$, $k\in \N_0$, are \iid; the Holder inequality}{that for all $q\in \N\backslash\{1\}$, $M\in \N$, $K_1\in \N_0$, $K_2\in (\N\cap[K_1,\infty])\cup\{\infty\}$ it holds that
    \begin{equation}\llabel{eq4}
    \begin{split}
&\textstyle\E\biggl[\biggl(\sum\limits_{k=K_1}^{K_2}\alpha^k|\bbX_k^M|\biggr)^q\biggr]
=\textstyle\E\biggl[\sum\limits_{(k_i)_{i\in \N_0\cap[K_1,K_2]}\subseteq\N_0,\, \sum_{i=K_1}^{K_2} k_i=q}\frac{q!}{\prod_{i=K_1}^{K_2}(k_i)!}\biggl(\prod\limits_{i=K_1}^{K_2}\alpha^{ik_i}|\bbX_i^M|^{k_i}\biggr)\biggr]\\
&=\textstyle\E\biggl[\sum\limits_{(k_i)_{i\in \N_0\cap[K_1,K_2]}\subseteq\N_0,\, \sum_{i=K_1}^{K_2} k_i=q}\frac{q!\alpha^{\sum_{i=K_1}^{K_2}ik_i}}{\prod_{i=K_1}^{K_2}(k_i)!}\biggl(\prod\limits_{i=K_1}^{K_2}|\bbX_i^M|^{k_i}\biggr)\biggr]\\
&=\textstyle\sum\limits_{(k_i)_{i\in \N_0\cap[K_1,K_2]}\subseteq\N_0,\, \sum_{i=K_1}^{K_2} k_i=q}\frac{q!\alpha^{\sum_{i=K_1}^{K_2}ik_i}}{\prod_{i=K_1}^{K_2}(k_i)!}\biggl(\prod\limits_{i=K_1}^{K_2}\E\bigl[|\bbX_i^M|^{k_i}\bigr]\biggr)\\
&\leq\textstyle\sum\limits_{(k_i)_{i\in \N_0\cap[K_1,K_2]}\subseteq\N_0,\, \sum_{i=K_1}^{K_2} k_i=q}\frac{q!\alpha^{\sum_{i=K_1}^{K_2}ik_i}}{\prod_{i=K_1}^{K_2}(k_i)!}\biggl(\prod\limits_{i=K_1}^{K_2}\bigl(\E\bigl[|\bbX_1^M|^{q}\bigr]\bigr)^{k_i/q}\biggr)\\
&=\textstyle\sum\limits_{(k_i)_{i\in \N_0\cap[K_1,K_2]}\subseteq\N_0,\, \sum_{i=K_1}^{K_2} k_i=q}\frac{q!\alpha^{\sum_{i=K_1}^{K_2}ik_i}}{\prod_{i=K_1}^{K_2}(k_i)!}\E\bigl[|\bbX_1^M|^{q}\bigr]\\
&\leq\textstyle\sum\limits_{(k_i)_{i\in \N_0\cap[K_1,K_2]}\subseteq\N_0,\, \sum_{i=K_1}^{K_2} k_i=q}\frac{q!\alpha^{\sum_{i=K_1}^{K_2}ik_i}}{\prod_{i=K_1}^{K_2}(k_i)!}\fC(q)c^q= \fC(q)c^q\biggl(\sum\limits_{i=K_1}^{K_2}\alpha^i\biggr)^q\\
&\leq \fC(q)c^q\alpha^{qK_1}(1-\alpha)^{-q}.
\end{split}
    \end{equation}}
    \argument{\lref{eq4};\lref{def: bbX};}{\cref{item 1: upper order of Adam coefficient 1}\dott}
    \startnewargseq
    \argument{\lref{def: bbX};\lref{def: AB};\lref{eq4}}{that for all $K\in \N\cup\{\infty\}$, $M\in \N$, $q\in \N\backslash\{1\}$ it holds that
    \begin{equation}\llabel{eq5}
    \begin{split}
\E[|A_K^M|^q]&=\textstyle\E\bigl[|\sum_{k=0}^K\alpha^k\bbX_k^M|^{q}\bigr]\leq\E\bigl[\bigl(\sum_{k=0}^K\alpha^k|\bbX_k^M|\bigr)^{q}\bigr] \leq \fC(q)c^q(1-\alpha)^{-q}.
\end{split}
    \end{equation}}
    \argument{\lref{eq5};Holder inequality;the fact that $\delta\in \N_0$}{that for all $M\in \N$, $K\in \N$ it holds that
    \begin{equation}\llabel{evd3}
   \E[(A_K^M)^{4p-4}]\leq \fC(4p-4)c^{4p-4}(1-\alpha)^{4-4p},\qquad\E[(A_\infty^M)^{4p-4}]\leq \fC(4p-4)c^{4p-4}(1-\alpha)^{4-4p},
   \end{equation}
   \begin{equation}\llabel{evd4}
  \E[(A_\infty^M)^{2p}]\leq \fC(2p)c^{2p}(1-\alpha)^{-2p},\qquad\text{and}\qquad \E[(A_K^M)^{4\delta}]\leq 1+\fC(4\delta)c^{4\delta}(1-\alpha)^{-4\delta}.
\end{equation}}
\argument{\lref{eq2};\lref{evd1};\lref{evd4}}{that there exists $\bfC\in (0,\infty)$ such that for all $M,K\in \N$ it holds that
\begin{equation}\llabel{eq7}
    \begin{split}
        &\bigl|\E\bigl[\psi(A^M_\infty)|A^M_\infty|^{p}\bigr]-\E\bigl[\psi(A^M_K)|A^M_K|^{p}\bigr]\bigr|\leq\bigl(\E[(A^M_\infty-A_K^M)^2]\bigr)^{1/2}\bigl(\E[(A^M_\infty)^{2p}]\bigr)^{1/2}\\
       &+p\bigl(\E[(A^M_K)^{4\delta}]\bigr)^{1/4}\Bigl(\bigl(\E[(A_\infty^M)^{4p-4}]\bigr)^{1/4}+\bigl(\E[(A_K^M)^{4p-4}]\bigr)^{1/4}\Bigr)\bigl(\E\bigl[(A^M_\infty-A^M_K)^{2}\bigr]\bigr)^{1/2}\\
       &\leq \alpha^{K+1}(1-\alpha^2)^{-1/2}c[\fC(2p)]^{1/2}c^p(1-\alpha)^{-p}\\
       &+p[1+\fC(4\delta)c^{4\delta}(1-\alpha)^{-4\delta}]^{1/4}(1-\alpha)^{1-p}\bigl([\fC(4p-4)]^{1/4} c^{p-1}+[\fC(4p-4)]^{1/4} c^{p-1}\bigr)\\
       &\quad\cdot\alpha^{K+1}(1-\alpha^2)^{-1/2}c\\
       &=\bfC \alpha^{K+1}.
    \end{split}
\end{equation}}
\argument{\lref{eq7};\lref{eq1}}{that there exists $\bfC\in (\varepsilon,\infty)$ which satisfies for all $M,K\in \N$ that
\begin{equation}\llabel{eq8}
     \begin{split}
           &M^{\nicefrac{(p+\delta)}{2}}\Bigl|\E\bigl[\psi\textstyle \bigl(\sum_{k=0}^{\infty}\alpha^k(\frac{1}{M}\sum_{m=1}^MX_{k+1,m})\bigr)\bigl|(1-\alpha)\sum_{k=0}^\infty \alpha^k[\frac{1}{M}\sum_{m=1}^MX_{k+1,m}]\bigr|^{p}\bigr]\\
           &-\E\bigl[\textstyle \psi\bigl(\sum_{k=0}^{K}\alpha^k(\frac{1}{M}\sum_{m=1}^MX_{k+1,m})\bigr)\bigl|(1-\alpha)\sum_{k=0}^K \alpha^k[\frac{1}{M}\sum_{m=1}^MX_{k+1,m}]\bigr|^{p}\bigr]\Bigr|\\
           &= (1-\alpha)^p\Bigl|\E\bigl[\psi(A^M_\infty)|A^M_\infty|^{p}\bigr]-\E\bigl[\psi(A^M_K)|A^M_K|^{p}\bigr]\Bigr|\leq \bfC\alpha^K.
           \end{split}
\end{equation}}
\startnewargseq
\argument{\lref{eq8}; the fact that $0< \alpha<1$}{that for all $M\in \N$, $K\in \N\cap[ \log_{\alpha}(\varepsilon \bfC^{-1}),\infty)$ it holds that
\begin{equation}\llabel{eq9}
     \begin{split}
          &M^{\nicefrac{(p+\delta)}{2}}\Bigl|\E\bigl[\textstyle \psi\bigl(\sum_{k=0}^{\infty}\alpha^k(\frac{1}{M}\sum_{m=1}^MX_{k+1,m})\bigr)\bigl|(1-\alpha)\sum_{k=0}^\infty \alpha^k[\frac{1}{M}\sum_{m=1}^MX_{k+1,m}]\bigr|^{p}\bigr]\\
           &-\E\bigl[\textstyle \psi\bigl(\sum_{k=0}^{K}\alpha^k(\frac{1}{M}\sum_{m=1}^MX_{k+1,m})\bigr)\bigl|(1-\alpha)\sum_{k=0}^K \alpha^k[\frac{1}{M}\sum_{m=1}^MX_{k+1,m}]\bigr|^{p}\bigr]\Bigr|\\
           &\leq \bfC\alpha^{\log_{\alpha}(\varepsilon \bfC^{-1})}= \varepsilon.
           \end{split}
\end{equation}}
\argument{\lref{eq9};}{\lref{conclude}\dott}
\end{aproof}

In \cref{conclude 1: upper order of Adam coefficient} in \cref{cor: upper order of Adam coefficient} below we establish sharp convergence rates for the expectations of signed absolute powers of certain stationary momentum random variables. In our proof of \cref{cor: upper order of Adam coefficient} we employ the elementary fact on weakly convergent probability measures in \cref{lem: weak converge} below. Only for completeness we include here a detailed proof for \cref{lem: weak converge}.
\begin{tcolorbox}[colback=white!95!gray,
                  colframe=black,
                  boxrule=0.5pt,
                  sharp corners,
                  enhanced,
                  breakable,
                 ]
                  \begin{athm}{lemma}{lem: weak converge}
                      Let $d\in\N$, let $\mu_n\colon \cB(\R^d)\to[0,\infty]$, $n\in \N_0$, be probability measures, assume that $(\mu_n)_{n\in \N}$ converges weakly to $\mu_0$, assume for all $p\in \N$ that $\limsup_{n\to\infty}\int_{\R^d}\|x\|^p\,\mu_n(\d x)<\infty$, and let $\varphi\in C(\R^d,\R)$ be at most polynomial growing. Then
                      \begin{equation}\llabel{conclude}
                          \textstyle\limsup_{n\to\infty} \bigl|\int_{\R^d}\varphi(x)\,\mu_n(\d x)-\int_{\R^d}\varphi(x)\,\mu_0(\d x)\bigr|=0.
                      \end{equation}
                  \end{athm} 
                    \end{tcolorbox}
                  \begin{aproof}
                Throughout this proof for every $R\in (0,\infty)$ let $g_R\colon\R\to\R$ satisfy for all $x\in\R$ that
                \begin{equation}\llabel{def: g}
                    g_R(x)=\max\{-R,\min\{x,R\}\}
                \end{equation}
                and let $p\in\N$, $\scrc\in (0,\infty)$ satisfy for all $x\in\R^d$ that
                \begin{equation}\llabel{def: p}
                    |\varphi(x)|\leq \scrc+\scrc \|x\|^p.
                \end{equation}
                \argument{the assumption that for all $q\in \N$ it holds that $\limsup_{n\to\infty}\int_{\R^d}\|x\|^q\,\mu_n(\d x)<\infty$;}{that for all $q\in\N$ there exists $\bfM_q\in\N$ which satisfies
                \begin{equation}\llabel{def: bfM}
                   \textstyle \sup_{n\in \N\cap[\bfM_q,\infty)}\int_{\R^d}\|x\|^q\,\mu_n(\d x)<\infty.
                \end{equation}}
                \startnewargseq
                \argument{the assumption that $(\mu_n)_{n\in\N}$ converges weakly to $\mu_0$;}{that for all $R\in (0,\infty)$, $q\in \N$ it holds that
                \begin{equation}\llabel{eqt1}
                 \int_{\R^d}\min\{\|x\|^q,R\}\,\mu_0(\d x)=\lim_{n\to\infty}\int_{\R^d}\min\{\|x\|^q,R\}\,\mu_n(\d x)\leq \limsup_{n\to\infty}\int_{\R^d}\|x\|^q\,\mu_n(\d x).
                \end{equation}}
                \argument{\lref{eqt1};monotone convergence theorem}{that \llabel{arg1} for all $q\in \N$ it holds that $ \int_{\R^d}\|x\|^q\,\mu_0(\d x)<\infty$\dott}
                \argument{\lref{arg1};\lref{def: bfM};}{that for all $q\in\N$ it hold that
                \begin{equation}\llabel{eqt3}
                    \textstyle\sup_{n\in\{0\}\cup(\N\cap[\bfM_q,\infty))}\int_{\R^d}\|x\|^q\,\mu_n(\d x)<\infty.
                \end{equation}}
                \argument{\lref{def: g};\lref{def: p};\lref{def: bfM}}{that for all $n\in \N\cap[\bfM_{2p},\infty)$, $R\in (0,\infty)$ it holds that
                \begin{equation}\llabel{eq2}
                \begin{split}
               &  \textstyle   \bigl|\int_{\R^d}g_R(\varphi(x))\,\mu_n(\d x)-\int_{\R^d}\varphi(x)\mu_n(\d x)\bigr|\leq \int_{\R^d}|\varphi(x)|\mathbbm 1_{\{|\varphi(x)| \geq R\}}\,\mu_n(\d x)\\
               &\leq \int_{\R^d}|\varphi(x)|^2 R^{-1}\,\mu_n(\d x)\leq \int_{\R^d}R^{-1}(\scrc+\scrc\|x\|^p)^2\,\mu_n(\d x)\leq 2R^{-1}\int_{\R^d}(\scrc^2+\scrc^2\|x\|^{2p})\,\mu_n(\d x).
                    \end{split}
                \end{equation}}
                \argument{\lref{eq2};\lref{eqt3};the fact that for all $n\in\N_0$ it holds that $\mu_n(\R^d)\leq 1$}{that there exists $\fC\in (0,\infty)$ which satisfies for all $R\in (0,\infty)$ that
                \begin{equation}\llabel{def: fC}
                  \textstyle  \sup_{n\in\{0\}\cup(\N\cap[\bfM_{2p},\infty))}\bigl|\int_{\R^d}g_R(\varphi(x))\,\mu_n(\d x)-\int_{\R^d}\varphi(x)\mu_n(\d x)\bigr|\leq \fC R^{-1}.
                \end{equation}}
                \startnewargseq
                \argument{the assumption that $(\mu_n)_{n\in\N}$ converges weakly to $\mu_0$;the fact that for all $R\in (0,\infty)$ it holds that $g_R\circ \varphi$ is bounded;the fact that for all $R\in (0,\infty)$ it holds that $g_R\circ\varphi$ is continuous;}{that for all $R\in (0,\infty)$ it holds that
                \begin{equation}\llabel{eq1}
                   \textstyle \lim_{n\to\infty} \int_{\R^d}g_R(\varphi(x))\,\mu_n(\d x)=\int_{\R^d}g_R(\varphi(x))\,\mu_0(\d x).
                \end{equation}}
                \argument{\lref{eq1};}{that for $\varepsilon\in (0,\infty)$ there exists $\bfK_\varepsilon\in (0,\infty)$ which satisfies for all $n\in \N\cap[\bfK_\varepsilon,\infty)$ that
                \begin{equation}\llabel{def: K}
                   \textstyle \bigl|\int_{\R^d}g_{\fC\varepsilon^{-1}}(\varphi(x))\,\mu_0(\d x)-\int_{\R^d}g_{\fC\varepsilon^{-1}}(\varphi(x))\,\mu_n(\d x)\bigr|\leq\varepsilon.
                \end{equation}}
                \startnewargseq
                \argument{\lref{def: fC};\lref{def: K}}{that for all $\varepsilon\in(0,\infty)$, $n\in \N\cap[\max\{\bfM_{2p},\bfK_\varepsilon\},\infty)$ it holds that
                \begin{equation}\llabel{eq3}
                \begin{split}
                  \textstyle  \bigl|\int_{\R^d}\varphi(x)\,\mu_n(\d x)-\int_{\R^d}\varphi(x)\,\mu_0(\d x)\bigr|
                  &\leq \textstyle \bigl|\int_{\R^d}g_{\fC\varepsilon^{-1}}(\varphi(x))\,\mu_n(\d x)-\int_{\R^d}\varphi(x)\,\mu_n(\d x)\bigr|\\
                  &\textstyle+\bigl|\int_{\R^d}g_{\fC\varepsilon^{-1}}(\varphi(x))\,\mu_n(\d x)-\int_{\R^d}g_{\fC\varepsilon^{-1}}(\varphi(x))\,\mu_0(\d x)\bigr|\\
                &\textstyle+\bigl|\int_{\R^d}g_{\fC\varepsilon^{-1}}(\varphi(x))\,\mu_0(\d x)-\int_{\R^d}\varphi(x)\,\mu_0(\d x)\bigr|\\
                   &\leq \varepsilon+\varepsilon+\varepsilon=3\varepsilon.
                    \end{split}
                \end{equation}}
                \argument{\lref{eq3};}{\lref{conclude}\dott}
                  \end{aproof}
                 
\begin{tcolorbox}[colback=white!95!gray,
                  colframe=black,
                  boxrule=0.5pt,
                  sharp corners,
                  enhanced,
                  breakable,
                 ]
\begin{athm}{cor}{cor: upper order of Adam coefficient}
    Let $(\Omega,\cF,\P)$ be a probability space, let $\alpha\in [0,1)$, $p\in \N\backslash\{1\}$, $\delta\in \N_0$, let $X_{n,m}\colon \Omega\to\R$, $(n,m)\in \N^2$, be bounded centered \iid\ random variables, let $\psi\colon \R\to\R$ satisfy for all $x,y\in \R$, $\kappa\in (0,\infty)$ that $|\psi(x)-\psi(y)|\leq |x-y|$, $|\psi(x)|\leq |x|^\delta$, and $\psi(\kappa x)=\kappa^\delta \psi(x)$, and for every $k\in \N_0$, $M \in \N$ let $\bbX_{k,M} \colon \Omega \to \R$ satisfy $\bbX_{k,M} = \frac{ 1 }{ M } \sum_{ m=1 }^{ M } X_{ k+1, m }$. Then
      \begin{equation}\label{conclude 1: upper order of Adam coefficient}
            \textstyle \limsup_{M\to\infty}\bigl|M^{\nicefrac{(p+1)}{2}}\,\E\bigl[\textstyle \bigl(\sum_{k=0}^{\infty}\alpha^k\bbX_{k,M}\bigr)|(1-\alpha)\sum_{k=0}^\infty \alpha^k\bbX_{k,M}|^p\bigr]\bigr|=0
      \end{equation}
    and 
  \begin{equation}\label{conclude 2: upper order of Adam coefficient}
  \begin{split}
            \textstyle \limsup_{M\to\infty}\bigl(M^{\nicefrac{(p+\delta)}{2}}\,\E\bigl[\textstyle \psi\bigl(\sum_{k=0}^{\infty}\alpha^k\bbX_{k,M}\bigr)|(1-\alpha)\sum_{k=0}^\infty \alpha^k\bbX_{k,M}|^{p}\bigr]\bigr)<\infty.
            \end{split}
      \end{equation}
\end{athm}
\end{tcolorbox}
\begin{aproof}
Throughout this proof for every $k\in \N_0$, $M\in\N$ let $\bfX_k^M\colon \Omega\to\R$ satisfy
    \begin{equation}\llabel{def: bbX}
       \displaystyle \bfX_k^M=M^{-1/2} \textstyle\sum_{m=1}^MX_{k+1,m}.
    \end{equation}
    \argument{\lref{def: bbX};the fact that for all $\kappa\in (0,\infty)$, $x\in \R$ it holds that $\psi(\kappa x)=\kappa^\delta \psi(x)$}{that for all $M\in \N$, $K\in \N\cup\{\infty\}$ it holds that
    \begin{equation}\llabel{eq1}
    \begin{split}
    &M^{\nicefrac{(p+1)}{2}}\,\E\bigl[\textstyle \bigl(\sum_{k=0}^{K}\alpha^k\bbX_{k,M}\bigr)\bigl|(1-\alpha)\sum_{k=0}^K \alpha^k\bbX_{k,M}\bigr|^p\bigr]\\
          &=M^{\nicefrac{(p+1)}{2}}\Bigl(\E\bigl[\textstyle \bigl(\sum_{k=0}^{K}\alpha^k[\frac{1}{M}\sum_{m=1}^MX_{k+1,m}]\bigr)\bigl|(1-\alpha)\sum_{k=0}^K \alpha^k[\frac{1}{M}\sum_{m=1}^MX_{k+1,m}]\bigr|^{p}\bigr]\Bigr)\\
           &=\E\bigl[\textstyle \bigl(\sum_{k=0}^{K}\alpha^k[M^{-1/2}\sum_{m=1}^MX_{k+1,m}]\bigr)\bigl|(1-\alpha)\sum_{k=0}^K \alpha^k[M^{-1/2}\sum_{m=1}^MX_{k+1,m}]\bigr|^{p}\bigr]\\
          &= \E\bigl[\textstyle \bigl(\sum_{k=0}^{K}\alpha^k\bfX_k^M\bigr)\bigl|(1-\alpha)\sum_{k=0}^K \alpha^k\bfX_k^M\bigr|^{p}\bigr]
          \end{split}
    \end{equation}
    and 
    \begin{equation}\llabel{eq2}
    \begin{split}
    &M^{\nicefrac{(p+\delta)}{2}}\,\E\bigl[\textstyle \psi\bigl(\sum_{k=0}^{K}\alpha^k\bbX_{k,M}\bigr)\bigl|(1-\alpha)\sum_{k=0}^K \alpha^k\bbX_{k,M}\bigr|^{p}\bigr]\\
          &=M^{\nicefrac{(p+\delta)}{2}}\biggl(\E\bigl[\textstyle \psi\bigl(\sum_{k=0}^{K}\alpha^k\bigl(\frac{1}{M}\sum_{m=1}^MX_{k+1,m}\bigr)\bigr)\bigl|(1-\alpha)\sum_{k=0}^K \alpha^k[\frac{1}{M}\sum_{m=1}^MX_{k+1,m}]\bigr|^{p}\bigr]\biggr)\\
           &=\E\bigl[\textstyle \psi\bigl(\sum_{k=0}^{K}\alpha^k\bigl(M^{-1/2}\sum_{m=1}^MX_{k+1,m}\bigr)\bigr)\bigl|(1-\alpha)\sum_{k=0}^K \alpha^k[M^{-1/2}\sum_{m=1}^MX_{k+1,m}]\bigr|^{p}\bigr]\\
          &= \E\bigl[\textstyle \psi\bigl(\sum_{k=0}^{K}\alpha^k\bfX_k^M\bigr)\bigl|(1-\alpha)\sum_{k=0}^K \alpha^k\bfX_k^M\bigr|^{p}\bigr].
          \end{split}
    \end{equation}}
\argument{\cref{item 2: upper order of Adam coefficient 1} in \cref{lem: upper order of Adam coefficient 1};\lref{eq2}}{that for all $\varepsilon\in (0,\infty)$ there exists $\bfK_\varepsilon\in \N$ which satisfy for all $M\in \N$ that
\begin{equation}\llabel{def: K1}
           \bigl|\E\bigl[\textstyle \bigl(\sum_{k=0}^{\infty}\alpha^k\bfX_k^M\bigr)|(1-\alpha)\sum_{k=0}^\infty \alpha^k\bfX_k^M|^{p}\bigr]
           -\E\bigl[\textstyle \bigl(\sum_{k=0}^{\bfK_\varepsilon}\alpha^k\bfX_k^M\bigr)|(1-\alpha)\sum_{k=0}^{\bfK_\varepsilon} \alpha^k\bfX_k^M|^{p}\bigr]\bigr|\leq \frac{\varepsilon}{2}.
\end{equation}
and 
\begin{equation}\llabel{def: K2}
\begin{split}
          & \bigl|\E\bigl[\psi\textstyle \bigl(\sum_{k=0}^{\infty}\alpha^k\bfX_k^M\bigr)|(1-\alpha)\sum_{k=0}^\infty \alpha^k\bfX_k^M|^{p}\bigr]\\
         &\quad-\E\bigl[\textstyle \psi\bigl(\sum_{k=0}^{\bfK_\varepsilon}\alpha^k\bfX_k^M\bigr)|(1-\alpha)\sum_{k=0}^{\bfK_\varepsilon} \alpha^k\bfX_k^M|^{p}\bigr]\bigr|\leq \varepsilon.
         \end{split}
\end{equation}}
\startnewargseq
\argument{\lref{def: bbX};the assumption that $X_{n,m}$, $(n,m)\in \N^2$, are \iid; the central limit theorem (see, \eg, \cite[Theorem 15.58]{klenkeprobability})}{that there exist random variables $N_{k,j}\colon\Omega\to\R$, $(k,j)\in(\N_0)^2$, which satisfy that
\begin{enumerate}[label=(\roman*)]
    \item \llabel{item 1} it holds for all $K\in\N$, $j\in\{0,1,\dots,K\}$ that $N_{K,j}$ is a normal distribution,
    \item \llabel{item 2} it holds for all $K\in\N$, $i,j\in\{0,1,\dots,K\}$ that
       $ \E[N_{K,i}]=0$ and $\operatorname{cov}(N_{K,i},N_{K,j})=\E[(X_{1,1})^2]\mathbbm 1_{\{i=j\}}$,
    and
    \item \llabel{item 3} it holds for all $K\in \N$ that
$\textstyle (\bfX^M_{0},\bfX^M_1,\dots,\bfX^M_K)$ converges in distribution to $(N_{K,0},N_{K,1},\dots,N_{K,K})$ as $M$ goes to infinity\dott
\end{enumerate}}
\argument{\lref{item 3};}{that for all $K\in \N$ it holds that
\begin{equation}\llabel{eq3}
   \textstyle \sum_{k=0}^K\alpha^k\bfX^M_k \rightarrow \sum_{k=0}^K\alpha^kN_{K,k}
\end{equation}
in distribution as $M\to\infty$\dott}
\argument{\cref{item 1: upper order of Adam coefficient 1} in \cref{lem: upper order of Adam coefficient 1};\lref{def: bbX}}{that for all $K,q\in \N$ it holds that
\begin{equation}\llabel{eqt1}
  \textstyle  \sup_{M\in\N} \E\bigl[\textstyle \bigl|\sum_{k=0}^{K} \alpha^k\bfX_k^M\bigr|^{q}\bigr]\leq \textstyle  \sup_{M\in\N} \E\bigl[\textstyle \bigl(\sum_{k=0}^{\infty} \alpha^k|\bfX_k^M|\bigr)^{q}\bigr]<\infty.
\end{equation}}
\argument{\lref{eqt1};\lref{eq3};
\cref{lem: weak converge}}{that for all $K\in\N$ it holds that
\begin{equation}\llabel{eqt2}
    \begin{split}
        &\lim_{M\to\infty}\E\bigl[\textstyle \bigl(\sum_{k=0}^{K}\alpha^k\bfX_k^M\bigr)\bigl|(1-\alpha)\sum_{k=0}^{K} \alpha^k\bfX_k^M\bigr|^{p}\bigr]=\E\bigl[\textstyle \bigl(\sum_{k=0}^{K}\alpha^kN_{K,k}\bigr)\bigl|(1-\alpha)\sum_{k=0}^{K} \alpha^kN_{K,k}\bigr|^{p}\bigr].
    \end{split}
\end{equation}}
\argument{the fact that for all $K\in \N$ it holds that $(N_{K,0},N_{K,1}\allowbreak,\allowbreak\dots,N_{K,K})$ and $(-N_{K,0},-N_{K,1},\allowbreak\dots,\allowbreak-N_{K,K})$ are identically distributed;}{that for all $K\in\N$ and all measurable $\varphi \colon \R^{K+1} \to \R$ it holds that \llabel{argg1} $\varphi( N_{K,0},N_{K,1},\dots,N_{K,K} )$ and $\varphi( -N_{K,0},-N_{K,1},\dots,\allowbreak-N_{K,K})$ are identically distributed\dott}
\argument{\lref{eqt2};\lref{argg1}}{that for all $\varepsilon\in (0,\infty)$ it holds that
\begin{equation}\llabel{eq4}
\begin{split}
    &\lim_{M\to\infty}\E\bigl[\textstyle \bigl(\sum_{k=0}^{\bfK_\varepsilon}\alpha^k\bfX_k^M\bigr)\bigl|(1-\alpha)\sum_{k=0}^{\bfK_\varepsilon} \alpha^k\bfX_k^M\bigr|^{p}\bigr]\\
    &=\E\bigl[\textstyle \bigl(\sum_{k=0}^{\bfK_\varepsilon}\alpha^kN_{\bfK_\varepsilon,k}\bigr)\bigl|(1-\alpha)\sum_{k=0}^{\bfK_\varepsilon} \alpha^kN_{\bfK_\varepsilon,k}\bigr|^{p}\bigr]\\
    &=\frac{1}{2}\Bigl(\E\bigl[\textstyle \bigl(\sum_{k=0}^{\bfK_\varepsilon}\alpha^kN_{\bfK_\varepsilon,k}\bigr)\bigl|(1-\alpha)\sum_{k=0}^{\bfK_\varepsilon} \alpha^kN_{\bfK_\varepsilon,k}\bigr|^{p}\bigr]\\
    &\quad+\E\bigl[\textstyle \bigl(\sum_{k=0}^{\bfK_\varepsilon}\alpha^k(-N_{\bfK_\varepsilon,k})\bigr)\bigl|(1-\alpha)\sum_{k=0}^{\bfK_\varepsilon} \alpha^k(-N_{\bfK_\varepsilon,k})\bigr|^{p}\bigr]\Bigr)\\
    &=0.
    \end{split}
\end{equation}}
\argument{\lref{eq4};}{that for all $\varepsilon\in (0,\infty)$ there exists $\bfM_\varepsilon\in \N$ which satisfies for all $M\in \N\cap[\bfM_\varepsilon,\infty)$ that
\begin{equation}\llabel{def: bfM}
    \Bigl|\E\bigl[\textstyle \bigl(\sum_{k=0}^{\bfK_\varepsilon}\alpha^k\bfX_k^M\bigr)\bigl|(1-\alpha)\sum_{k=0}^{\bfK_\varepsilon} \alpha^k\bfX_k^M\bigr|^{p}\bigr]\Bigr|<\frac{\varepsilon}{2}.
\end{equation}}
\startnewargseq
\argument{\lref{def: K1};\lref{def: bfM}}{that for all $\varepsilon\in (0,\infty)$, $M\in \N\cap[\bfM_\varepsilon,\infty)$ it holds that
\begin{equation}\llabel{eq7}
    \begin{split}
           & \Bigl|\E\bigl[\textstyle \bigl(\sum_{k=0}^{\infty}\alpha^k\bfX_k^M\bigr)\bigl|(1-\alpha)\sum_{k=0}^{\infty} \alpha^k\bfX_k^M\bigr|^{p}\bigr]\Bigr|\\
           & \leq \Bigl|\E\bigl[\textstyle \bigl(\sum_{k=0}^{\bfK_\varepsilon}\alpha^k\bfX_k^M\bigr)\bigl|(1-\alpha)\sum_{k=0}^{\bfK_\varepsilon} \alpha^k\bfX_k^M\bigr
           |^{p}\bigr]\Bigr|\\
          & +  \Bigl|\E\bigl[\textstyle \bigl(\sum_{k=0}^{\infty}\alpha^k\bfX_k^M\bigr)\bigl|(1-\alpha)\sum_{k=0}^\infty \alpha^k\bfX_k^M\bigr|^{p}\bigr] -\E\bigl[\textstyle \bigl(\sum_{k=0}^{\bfK_\varepsilon}\alpha^k\bfX_k^M\bigr)\bigl|(1-\alpha)\sum_{k=0}^{\bfK_\varepsilon} \alpha^k\bfX_k^M\bigr|^{p}\bigr]\Bigr|\\
        &\leq \frac{\varepsilon}{2}+\frac{\varepsilon}{2}=\varepsilon.
    \end{split}
\end{equation}}
\argument{\lref{eq7};}{that
\begin{equation}\llabel{eq8}
    \limsup_{M\to\infty}\Bigl|\E\bigl[\textstyle \bigl(\sum_{k=0}^{\infty}\alpha^k\bfX_k^M\bigr)\bigl|(1-\alpha)\sum_{k=0}^{\infty} \alpha^k\bfX_k^M\bigr|^{p}\bigr]\Bigr|=0.
\end{equation}}
\argument{\lref{eq8};\lref{eq1}}{\cref{conclude 1: upper order of Adam coefficient}\dott}
\startnewargseq
\argument{\lref{eq3};\lref{eqt1};the assumption that for all $x\in \R$ it holds that $|\psi(x)|\leq |x|^\delta$;
\cref{lem: weak converge}}{that for all $K\in\N$ it holds that
\begin{equation}\llabel{eqt3}
    \begin{split}
        &\lim_{M\to\infty}\E\bigl[\textstyle \psi\bigl(\sum_{k=0}^{K}\alpha^k\bfX_k^M\bigr)\bigl|(1-\alpha)\sum_{k=0}^{K} \alpha^k\bfX_k^M\bigr|^{p}\bigr]\\
        &=\E\bigl[\psi\textstyle \bigl(\sum_{k=0}^{K}\alpha^kN_{K,k}\bigr)\bigl|(1-\alpha)\sum_{k=0}^{K} \alpha^kN_{K,k}\bigr|^{p}\bigr].
    \end{split}
\end{equation}}
\argument{\lref{eqt3};the fact that for all $k\in\{0,1,\dots,\bfK_1\}$ it holds that $N_{\bfK_1,k}$ is a normal distribution;the fact that $0\leq\alpha<1$;the Cauchy-Schwarz inequality; the assumption that for all $x\in \R$ it holds that $|\psi(x)|\leq |x|^\delta$}{that 
\begin{equation}\llabel{eq4'}
\begin{split}
    &\lim_{M\to\infty}\E\bigl[\textstyle\psi \bigl(\sum_{k=0}^{\bfK_1}\alpha^k\bfX_k^M\bigr)\bigl|(1-\alpha)\sum_{k=0}^{\bfK_1} \alpha^k\bfX_k^M\bigr|^{p}\bigr]\\
    &=\E\bigl[\psi\textstyle \bigl(\sum_{k=0}^{\bfK_1}\alpha^kN_{\bfK_1,k}\bigr)\bigl|(1-\alpha)\sum_{k=0}^{\bfK_1} \alpha^kN_{\bfK_1,k}\bigr|^{p}\bigr]\\
    &\textstyle\leq \bigl(\E\bigl[\bigl|\psi\bigl(\sum_{k=0}^{\bfK_1}\alpha^kN_{\bfK_1,k}\bigr)\bigr|^2\bigr]\bigr)^{1/2}\bigl(\E\bigl[\bigl|(1-\alpha)\sum_{k=0}^{\bfK_1} \alpha^kN_{\bfK_1,k}\bigr|^{2p}\bigr]\bigr)^{1/2}\\
    & \textstyle\leq \bigl(\E\bigl[\bigl|\sum_{k=0}^{\bfK_1}\alpha^kN_{\bfK_1,k}\bigr|^{2\delta}\bigr]\bigr)^{1/2}\bigl(\E\bigl[\bigl|(1-\alpha)\sum_{k=0}^{\bfK_1} \alpha^kN_{\bfK_1,k}\bigr|^{2p}\bigr]\bigr)^{1/2}\\
    &<\infty.
    \end{split}
\end{equation}}
\argument{\lref{eq4'};}{that there exist $\fC,\fM\in (0,\infty)$ which satisfy for all $M\in \N\cap[\fM,\infty)$ that
\begin{equation}\llabel{def: bfM'}
    \Bigl|\E\bigl[\textstyle \psi\bigl(\sum_{k=0}^{\bfK_1}\alpha^k\bfX_k^M\bigr)\bigl|(1-\alpha)\sum_{k=0}^{\bfK_1} \alpha^k\bfX_k^M\bigr|^{p}\bigr]\Bigr|\leq \fC.
\end{equation}}
\startnewargseq
\argument{\lref{def: K2};\lref{def: bfM'}}{that for all $M\in \N\cap[\fM,\infty)$ it holds that
\begin{equation}\llabel{eq7'}
    \begin{split}
           & \Bigl|\E\bigl[\textstyle \psi\bigl(\sum_{k=0}^{\infty}\alpha^k\bfX_k^M\bigr)\bigl|(1-\alpha)\sum_{k=0}^{\infty} \alpha^k\bfX_k^M\bigr|^{p}\bigr]\Bigr|\\
           & \leq \Bigl|\E\bigl[\textstyle \psi\bigl(\sum_{k=0}^{\bfK_1}\alpha^k\bfX_k^M\bigr)\bigl|(1-\alpha)\sum_{k=0}^{\bfK_1} \alpha^k\bfX_k^M\bigr|^{p}\bigr]\Bigr|\\
          & +  \Bigl|\E\bigl[\textstyle \psi\bigl(\sum_{k=0}^{\infty}\alpha^k\bfX_k^M\bigr)\bigl|(1-\alpha)\sum_{k=0}^{\infty} \alpha^k\bfX_k^M\bigr|^{p}\bigr] -\E\bigl[\textstyle \psi\bigl(\sum_{k=0}^{\bfK_1}\alpha^k\bfX_k^M\bigr)\bigl|(1-\alpha)\sum_{k=0}^{\bfK_1} \alpha^k\bfX_k^M\bigr|^{p}\bigr]\Bigr|\\
        &\leq 1+\fC.
    \end{split}
\end{equation}}
\argument{\lref{eq7'};}{that
\begin{equation}\llabel{eq8'}
    \limsup_{M\to\infty}\E\bigl[\textstyle \psi\bigl(\sum_{k=0}^{\infty}\alpha^k\bfX_k^M\bigr)\bigl|(1-\alpha)\sum_{k=0}^{\infty} \alpha^k\bfX_k^M\bigr|^{p}\bigr]<\infty.
\end{equation}}
\argument{\lref{eq8'};\lref{eq2}}{\cref{conclude 2: upper order of Adam coefficient}\dott}
\end{aproof}
\subsection{Sharp convergence rates for the expectation of the signed square of SMRVs}\label{subsec: sharp convergence rate smrv}
\begin{tcolorbox}[colback=white!95!gray,
                  colframe=black,
                  boxrule=0.5pt,
                  sharp corners,
                  enhanced,
                  breakable,
                 ]
                 \begin{athm}{lemma}{lem: tg1}
                     Let $\chi\in \R\backslash\{0\}$, $B,t\in (0,\infty)$ and let $r\colon\N\to \bbC$ satisfy $\limsup_{M\to\infty} (M^{3/2}|r_M|)=0$. Then\footnotemark
                     \begin{equation}\llabel{conclude}
                      \textstyle   \limsup_{M\to\infty}\bigl|\sqrt{M}\imagine\bigl(\bigl(1-\frac{Bt^2 }{2M}-\frac{\chi t^3\imag}{6M\sqrt{M}}+r_M\bigr)^M\bigr)+\frac{\chi t^3 \exp(-Bt^2/2)}{6}\bigr|=0.
                     \end{equation}
                 \end{athm}
                 \end{tcolorbox}
                 \footnotetext{Note that for all $a,b\in \R$ it holds that $\imagine(a+ b\imag)=b$ (imagine part of complex number).}
                 \begin{aproof}
                     Throughout this proof for every $M\in \N$ let $u_M,v_M\in \R$ satisfy 
                     \begin{equation}\llabel{def: uv}
                         r_M=u_M+v_M\imag.
                     \end{equation}
                     \argument{\lref{def: uv};the assumption that $\limsup_{M\to\infty} (M^{3/2}|r_M|)=0$}{that
                     \begin{equation}\llabel{eq1}
                         \textstyle\limsup_{M\to\infty} (M^{3/2}|u_M|)=\limsup_{M\to\infty} (M^{3/2}|v_M|)=0.
                     \end{equation}}
                     \argument{\lref{eq1};}{that
                     \begin{equation}\llabel{eq2}
                     \begin{split}
                        & \textstyle \lim_{M\to\infty}\bigl[M^{3/2}\bigl(\frac{\chi t^3}{6M\sqrt{M}}-v_M\bigr)\bigl(1-\frac{Bt^2}{2M}+u_M\bigr)^{-1}\bigr]\\
                        &=\textstyle \textstyle \lim_{M\to\infty}\bigl[\bigl(\frac{\chi t^3}{6}-M^{3/2}v_M\bigr)\bigl(1-\frac{Bt^2}{2M}+u_M\bigr)^{-1}\bigr]=\textstyle\bigl(\frac{\chi t^3}{6}-0\bigr) (1-0-0)^{-1}=\frac{\chi t^3}{6}.
                         \end{split}
                     \end{equation}}
                     \argument{\lref{eq2};}{that there exist $\fM\in\N$ and $s_M\in \R$, $M\in \N$, which satisfy for all $M\in\N\cap[\fM,\infty)$ that
                     \begin{equation}\llabel{eq3}
                         \limsup_{N\to\infty} N^{3/2}|s_N|=0\qqandqq \textstyle \bigl(\frac{\chi t^3}{6M\sqrt{M}}-v_M\bigr)\bigl(1-\frac{Bt^2}{2M}+u_M\bigr)^{-1}=\frac{\chi t^3}{6M\sqrt{M}}+s_M.
                     \end{equation}}
                     \startnewargseq
                     \argument{\lref{eq3};}{that for all $M\in\N\cap[\fM,\infty)$ it holds that
                     \begin{equation}\llabel{eq4}
                      \textstyle   \frac{\chi t^3}{6M\sqrt{M}}-v_M=\bigl(\frac{\chi t^3}{6M\sqrt{M}}+s_M\bigr)\bigl(1-\frac{Bt^2}{2M}+u_M\bigr).
                     \end{equation}}
                     \argument{\lref{eq4};\lref{def: uv};the multinomial theorem}{that for all $M\in\N\cap[\fM,\infty)$ it holds that
                     \begin{equation}\llabel{eq5}
                     \begin{split}
                        &\textstyle \bigl(1-\frac{Bt^2}{2M}-\frac{\chi t^3\imag}{6M\sqrt{M}}+r_M\bigr)^M=\bigl(1-\frac{Bt^2}{2M}+u_M-\bfi\bigl(\frac{\chi t^3}{6M\sqrt{M}}-v_M\bigr)\bigr)^M\\
                        &\textstyle=\bigl(1-\frac{Bt^2}{2M}+u_M-\bfi\bigl(\frac{\chi t^3}{6M\sqrt{M}}+s_M\bigr)\bigl(1-\frac{Bt^2}{2M}+u_M\bigr)\bigr)^M\\
                        &=\textstyle\bigl(1-\frac{Bt^2}{2M}+u_M\bigr)^M\bigl(1-\bfi\bigl(\frac{\chi t^3}{6M\sqrt{M}}+s_M\bigr)\bigr)^M\\
                        &\textstyle
                       =\bigl(1-\frac{Bt^2}{2M}+u_M\bigr)^M \biggl[\sum\limits_{k=0}^M\frac{M!}{k!(M-k)!}(-\imag)^k\bigl(\frac{\chi t^3}{6M\sqrt{M}}+s_M\bigr)^k\biggr]
                        \end{split}
                     \end{equation}}
                     \argument{\lref{eq5};}{for all $M\in\N\cap[\fM,\infty)$ that
                     \begin{equation}\llabel{eq6}
                     \begin{split}
                        & \textstyle\imagine\bigl(\bigl(1-\frac{Bt^2}{2M}-\frac{\chi t^3\imag}{6M\sqrt{M}}+r_M\bigr)^M\bigr)\\&=\textstyle\bigl(1-\frac{Bt^2}{2M}+u_M\bigr)^M \biggl[\sum\limits_{k=0}^{\lfloor (M-1)/2\rfloor}\frac{M!}{(2k+1)!(M-2k-1)!}(-\imag)^{2k}\bigl(\frac{\chi t^3}{6M\sqrt{M}}+s_M\bigr)^{2k+1}\biggr]\\
                         &=\textstyle\bigl(1-\frac{Bt^2}{2M}+u_M\bigr)^M \biggl[\sum\limits_{k=0}^{\lfloor (M-1)/2\rfloor}\frac{M!}{(2k+1)!(M-2k-1)!}(-1)^{k}\bigl(\frac{\chi t^3}{6M\sqrt{M}}+s_M\bigr)^{2k+1}\biggr]\\
                         &=\textstyle -\bigl(1-\frac{Bt^2}{2M}+u_M\bigr)^M \bigl(\frac{\chi t^3}{6\sqrt{M}}+Ms_M\bigr)\\
                         &-\textstyle\bigl(1-\frac{Bt^2}{2M}+u_M\bigr)^M \biggl[\sum\limits_{k=1}^{\lfloor (M-1)/2\rfloor}\frac{M!}{(2k+1)!(M-2k-1)!}(-1)^{k}\bigl(\frac{\chi t^3}{6M\sqrt{M}}+s_M\bigr)^{2k+1}\biggr].
                         \end{split}
                     \end{equation}}
                     \argument{\lref{eq6}}{for all $M\in \N\cap[\fM,\infty)$ that
                     \begin{equation}\llabel{eq7}
                     \begin{split}
                         &\textstyle \sqrt{M}\imagine\bigl(\bigl(1-\frac{Bt^2}{2M}-\frac{\chi t^3\imag}{6M\sqrt{M}}+r_M\bigr)^M\bigr)=-\textstyle \bigl(1-\frac{Bt^2}{2M}+u_M\bigr)^M \bigl(\frac{\chi t^3}{6}+M\sqrt{M}s_M\bigr)\\
                         &-\textstyle\bigl(1-\frac{Bt^2}{2M}+u_M\bigr)^M \biggl[\sum\limits_{k=1}^{\lfloor (M-1)/2\rfloor}\frac{M!}{(2k+1)!(M-2k-1)!}(-1)^{k}\bigl(\frac{\chi t^3}{6M\sqrt{M}}+s_M\bigr)^{2k+1}\sqrt{M}\biggr].
                         \end{split}
                     \end{equation}}
                     \argument{\lref{eq1};}{that  \begin{equation}\llabel{arg1}
                        \textstyle \limsup_{M\to\infty}(M(-\frac{Bt^2}{2M}+u_M))=\frac{-Bt^2}{2}<0=\limsup_{M\to\infty}(\frac{-Bt^2}{2M}+u_M).
                     \end{equation}}
                     \argument{\lref{arg1};the fact that $\lim_{x\to 0} (x^{-1}\log(1+x))=1$}{that
                     \begin{equation}\llabel{eq8}
                     \begin{split}
                         &\textstyle \lim_{M\to\infty} \bigl(M\log\bigl(1-\frac{Bt^2}{2M}+u_M\bigr)\bigr)\\
                         &\textstyle=\bigl[\lim_{M\to\infty} \bigl(-\frac{Bt^2}{2M}+u_M\bigr)^{-1}\log\bigl(1-\frac{Bt^2}{2M}+u_M\bigr)\bigr]\bigl[\lim_{M\to\infty} M\bigl(-\frac{Bt^2}{2M}+u_M\bigr)\bigr]\\
                         &=\textstyle 1\bigl(-\frac{Bt^2}{2}\bigr)=\frac{-Bt^2}{2}.
                         \end{split}
                     \end{equation}}
                     \argument{\lref{eq8};}{that
                     \begin{equation}\llabel{eq9}
                         \lim_{M\to\infty}\textstyle\bigl(1-\frac{Bt^2}{2M}+u_M\bigr)^M=\exp(-Bt^2/2).
                     \end{equation}}
                     \argument{\lref{eq9};\lref{eq3};}{that
                     \begin{equation}\llabel{evd1}
                     \begin{split}
                        &\textstyle \lim_{M\to\infty}\bigl[\bigl(1-\frac{Bt^2}{2M}+u_M\bigr)^M \bigl(\frac{\chi t^3}{6}+M\sqrt{M}s_M\bigr)\bigr]\\
                        &\textstyle=\bigl[\lim_{M\to\infty}\bigl(1-\frac{Bt^2}{2M}+u_M\bigr)^M\bigr]\lim_{M\to\infty}\bigl[ \frac{\chi t^3}{6}+M\sqrt{M}s_M\bigr]=\displaystyle\frac{\exp(-Bt^2/2)\chi t^3}{6}.
                        \end{split}
                     \end{equation}}
                     \argument{\lref{eq3}}{that there exists $\bfM\in \N$ such that for all $M\in \N\cap[\bfM,\infty)$ it holds that \llabel{arg2} $|s_M|\leq M^{-3/2}$\dott}
\argument{\lref{arg2};the fact that for all $y\in (0,\infty)$, $M\in \N$ it holds that $(1+\frac yM)^M\leq \exp(y)$}{that there exists $\bfM\in \N$ such that for all $M\in \N\cap[\bfM,\infty)$ it holds that
\begin{equation}\llabel{eq10}
\begin{split}
  &\textstyle  \biggl|\sum\limits_{k=1}^{\lfloor (M-1)/2\rfloor}\bigl[\frac{M!}{(2k+1)!(M-2k-1)!}(-1)^{k}\bigl(\frac{\chi t^3}{6M\sqrt{M}}+s_M\bigr)^{2k+1}\sqrt{M}\bigr]\biggr|\\
  &\textstyle \leq \sum\limits_{k=1}^{\lfloor (M-1)/2\rfloor}\bigl[\frac{M!}{(2k+1)!(M-2k-1)!}\bigl(\bigl|\frac{\chi t^3}{6M\sqrt{M}}\bigr|+|s_M|\bigr)^{2k+1}\sqrt{M}\bigr]\\
  &\textstyle  \leq \sum\limits_{k=1}^{\lfloor (M-1)/2\rfloor}\bigl[\frac{M!}{(2k+1)!(M-2k-1)!}\bigl(\bigl|\frac{\chi t^3}{6M\sqrt{M}}\bigr|+|M^{-3/2}|\bigr)^{2k+1}\sqrt{M}\bigr]\\
  &\leq\textstyle\sum\limits_{k=1}^{\lfloor (M-1)/2\rfloor}\bigl[\frac{M!}{(2k+1)!(M-2k-1)!}\bigl(\frac{|\chi t^3|+1}{M\sqrt{M}}\bigr)^{2k+1}\sqrt{M}\bigr]\leq \textstyle\sum\limits_{k=3}^{M}\bigl[\frac{M!}{(2k+1)!(M-2k-1)!}\bigl(\frac{|\chi t^3|+1}{M\sqrt{M}}\bigr)^{k}\sqrt{M}\bigr]\\
  &\leq \textstyle\sum\limits_{k=3}^{M}\bigl[\frac{M!}{(2k+1)!(M-2k-1)!}\bigl(\frac{|\chi t^3|+1}{M}\bigr)^{k}M^{-1}\bigr]\leq \textstyle\sum\limits_{k=1}^{M}\bigl[\frac{M!}{(2k+1)!(M-2k-1)!}\bigl(\frac{|\chi t^3|+1}{M}\bigr)^{k}M^{-1}\bigr]\\
  &=\textstyle M^{-1}\biggl[\sum\limits_{k=1}^{M}\frac{M!}{(2k+1)!(M-2k-1)!}\bigl(\frac{|\chi t^3|+1}{M}\bigr)^{k}\biggr]\leq M^{-1}\bigl(1+\frac{|\chi t^3|+1}{M}\bigr)^M\leq M^{-1}\exp(|\chi t^3|+1).
  \end{split}
\end{equation}}
\argument{\lref{eq10};}{that
\begin{equation}\llabel{eq11}
\textstyle\lim\limits_{M\to\infty}\biggl[\sum\limits_{k=1}^{\lfloor (M-1)/2\rfloor}\frac{M!}{(2k+1)!(M-2k-1)!}(-1)^{k}\bigl(\frac{\chi t^3}{6M\sqrt{M}}+s_M\bigr)^{2k+1}\sqrt{M}\biggr]=0.
\end{equation}}
\argument{\lref{eq11};\lref{eq9};}{that 
\begin{equation}\llabel{evd2}
   \textstyle \lim\limits_{M\to\infty}\biggl(\bigl(1-\frac{Bt^2}{2M}+u_M\bigr)^M \biggl[\sum\limits_{k=1}^{\lfloor (M-1)/2\rfloor}\frac{M!}{(2k+1)!(M-2k-1)!}(-1)^{k}\bigl(\frac{\chi t^3}{6M\sqrt{M}}+s_M\bigr)^{2k+1}\sqrt{M}\biggr]\biggr)=0.
\end{equation}}
\argument{\lref{evd2};\lref{eq7};\lref{evd1}}{that
\begin{equation}\llabel{eq12}
  \textstyle  \lim_{M\to\infty}\bigl[\sqrt{M}\imagine\bigl(\bigl(1-\frac{Bt^2}{2M}-\frac{\chi t^3\imag}{6M\sqrt{M}}+r_M\bigr)^M\bigr)\bigr]=\displaystyle\frac{-\exp(-Bt^2/2)\chi t^3}{6}.
\end{equation}}
\argument{\lref{eq12};}{\lref{conclude}\dott}
                 \end{aproof}

                 In the next elementary result, \cref{lem: integral}, we provide an integral representation for the signed square function $\R \ni x \mapsto x |x| \in \R$. We employ \cref{lem: integral} in our proof of \cref{lem: tg2} below, which establishes sharp convergence rates for the expectations of the signed square function evaluated at stationary momentum random variables.
                 \begin{tcolorbox}[colback=white!95!gray,
                  colframe=black,
                  boxrule=0.5pt,
                  sharp corners,
                  enhanced,
                  breakable,
                 ]
                 \begin{athm}{lemma}{lem: integral}[\textcolor{red}{An integral representation of the signed square function}]
                 Let $x \in \R$ and let $f \colon (0,\infty) \to \R$  satisfy for all $t \in (0,\infty)$ that $f(t) = 4\pi^{-1}t^{ - 3 } ( t x - \sin( t x ) )$. Then
                     \begin{enumerate}[label=(\roman*)]
                      \item \label{item 3: integral} it holds  that $\int_{0}^\infty |f(t)|\, \d t=x^2<\infty$ and
                         \item \label{item 2: integral} it holds that $\int_{0}^\infty f(t)\, \d t=x|x|$.
                     \end{enumerate}
                 \end{athm}
                                 \end{tcolorbox}
                 \begin{aproof}
                     \argument{the fact that for all $y\in \R$ it holds that $|\sin (y)|\leq 1$}{that for all $y\in\R$,  $s\in (0,\infty)$, $r\in (s,\infty)$ it holds that
                     \begin{equation}\llabel{eqt1}
                         \int_{s}^r\biggl| \frac{ty-\sin(ty)}{t^3}\biggr|\, \d t\leq  \int_{s}^r \frac{|ty|+1}{t^3}\, \d t=|y|\biggl[ \int_{s}^r \frac{1}{t^2}\,\d t\biggr]+\biggl[ \int_{s}^r \frac{1}{t^3}\, \d t\biggr]=|y|\biggl(\frac{1}{s}-\frac{1}{r}\biggr)+\frac{1}{2s^2}-\frac{1}{2r^2}.
                     \end{equation}}
                     \argument{\lref{eqt1};}{for all $y\in \R$, $s\in (0,\infty)$ that 
                     \begin{equation}\llabel{eqt1'}
                         \limsup_{r\to \infty} \int_{s}^r\biggl| \frac{ty-\sin(ty)}{t^3}\biggr|\, \d t<\infty.
                     \end{equation}}
                     \argument{\lref{eqt1'};}
                     {that for every $s\in (0,\infty)$ there exists $I_s\colon \R\to\R$ which satisfies for all $y\in \R$ that
                     \begin{equation}\llabel{def: I}
                         I_s(y)=\int_{s}^\infty \frac{ty-\sin(ty)}{t^3}\, \d t.
                     \end{equation}}
                    \startnewargseq
                     \argument{\lref{def: I};the fact that for all $t\in (0,\infty)$, $y\in \R$ it holds that $\sin(ty)=\operatorname{sgn}(y) \sin (t|y|)$;the fact that for all $y\in \R$ it holds that $\operatorname{sgn}(y)|y|=y$}{that for all $y\in \R\backslash\{0\}$, $s\in (0,\infty)$ it holds that
                     \begin{equation}\llabel{eq1}
                     \begin{split}
                         I_s(y)&=\int_{s}^\infty \frac{ty-\operatorname{sgn}(y)\sin(t|y|)}{t^3}\, \d t=\int_{s}^\infty \biggl(\frac{ty|y|^2-\operatorname{sgn}(y)|y|^2
                         \sin(t|y|)}{t^3|y|^3}\biggr)|y|\, \d t\\
                         &=\int_{s}^\infty\biggl(\frac{ty|y|^2-y|y|
                         \sin(t|y|)}{t^3|y|^3}\biggr)|y|\, \d t=y|y| \biggl[\int_{s}^\infty\biggl(\frac{t|y|-
                         \sin(t|y|)}{t^3|y|^3}\biggr)|y|\,\d t\biggr]\\
                         &=y|y|\biggl[\int_{s|y|}^\infty\frac{t-
                         \sin(t)}{t^3}\,\d t\biggr].
                         \end{split}
                     \end{equation}} 
                     \argument{\lref{eq1};the fact that $\frac {\d(t-\sin(t))}{\d t}= 1-\cos(t)$;the fact that $\lim_{t\to\infty} \frac{t-\sin(t)}{2t^2}=0$;integral by part}{that for all $y\in \R\backslash\{0\}$, $s\in (0,\infty)$ it holds that
                     \begin{equation}\llabel{eq2}
                         I_s(y)=y|y|\biggl[\frac{s|y|-\sin(s|y|)}{2(s|y|)^2}-\int_{s|y|}^\infty\frac{1-\cos(t)}{-2 t^2}\, d t\biggr]=y|y|\biggl[\frac{s|y|-\sin(s|y|)}{2(s|y|)^2}+\frac 12\int_{s|y|}^\infty\frac{1-\cos(t)}{ t^2}\, d t\biggr].
                     \end{equation}}
                     \argument{integral by part;the fact that $\frac {\d(1-\cos(t))}{\d t}= \sin(t)$; the fact that $\lim_{t\to\infty} \frac{1-\cos(t)}{t}=0$}{that for all $y\in \R\backslash\{0\}$, $s\in (0,\infty)$ it holds that
                     \begin{equation}\llabel{eq3}
                         \int_{s|y|}^\infty\frac{1-\cos(t)}{ t^2}\,\ d t=\frac{1-\cos(s|y|)}{s|y|}+\int_{s|y|}^\infty\frac{\sin(t)}{t}\d t.
                     \end{equation}}
                     \argument{\lref{eq3};the fact that $\lim_{t\to 0}\frac{1-\cos(t)}{t}=0$;the Dirichlet integral}{that for all $y\in \R\backslash\{0\}$ it holds that
                     \begin{equation}\llabel{eq4}
                     \begin{split}
\lim_{s\to 0}\int_{s|y|}^\infty\frac{1-\cos(t)}{ t^2}\,\ d t=\biggl[\lim_{s\to 0}\frac{1-\cos(s|y|)}{s|y|}\bigg]+\biggl[\lim_{s\to 0}\int_{s|y|}^\infty\frac{\sin(t)}{t}\d t\biggr]=\int_{0}^\infty\frac{\sin(t)}{t}=\frac{\pi}{2}.
\end{split}
                     \end{equation} }
                    \argument{\lref{eq4};the fact that $\lim_{t\to 0}\frac{t-\sin(t)}{2t^2}=0$}{that for all $y\in \R\backslash\{0\}$ it holds that
                    \begin{equation}\llabel{eq5}
                    \begin{split}
                        &\lim_{s\to 0}\biggl[\frac{s|y|-\sin(s|y|)}{2(s|y|)^2}+\frac 12\int_{s|y|}^\infty\frac{1-\cos(t)}{ t^2}\,\ d t\biggr]\\
                        &= \lim_{s\to 0}\frac{s|y|-\sin(s|y|)}{2(s|y|)^2}+\frac{1}{2}\biggl[\lim_{s\to 0}\int_{s|y|}^\infty\frac{1-\cos(t)}{ t^2}\,\ d t\biggr]\\
                        &=\frac{\pi}{4}.
                        \end{split}
                    \end{equation}}
                    \argument{\lref{eq5};\lref{eq2};}{that for all $y\in \R\backslash\{0\}$ it holds that \llabel{arg2} $\lim_{s\to 0}I_s(y)=\frac{\pi y|y|}{4}$\dott}
                    \argument{\lref{arg2};the fact that for all $s\in (0,\infty)$ it holds that $I_s(0)=0$}{that for all $y\in \R$ it holds that
                    \begin{equation}\llabel{arg3}
                        \lim_{s\to 0}I_s(y)=\frac{\pi y|y|}{4}.
                    \end{equation}}
                    \argument{\lref{arg3};\lref{def: I}}{that for all $y\in \R$ it holds that
                    \begin{equation}\llabel{arg3'}
                        \lim_{s\to 0}\int_{s}^\infty\frac{ty-\sin(ty)}{t^3}\,\d t=\frac{\pi y|y|}{4}.
                    \end{equation}}
                    \argument{\lref{arg3'}; the fact that for all $t\in (0,\infty)$ it holds that $\sin(t)\leq t$}{that for all $y\in (0,\infty)$ it holds that
                    \begin{equation}\llabel{eq6}
                    \begin{split}
                        \int_{0}^\infty\frac{|ty-\sin(ty)|}{t^3}\, \d t=\lim_{s\to 0} \int_{s}^\infty\frac{|ty-\sin(ty)|}{t^3}\, \d t= \lim_{s\to 0}\int_{s}^\infty\frac{ty-\sin(ty)}{t^3}\, \d t =\frac{\pi y|y|}{4}=\frac{\pi y^2}{4}.
                        \end{split}
                    \end{equation}}
                    \argument{\lref{eq6};the fact that for all $t\in \R$ it holds that $\sin(t)=-\sin(-t)$}{that for all $y\in (0,\infty)$ it holds that
                    \begin{equation}\llabel{eq7}
                     \int_{0}^\infty\frac{|t(-y)-\sin(t(-y))|}{t^3}\, \d t = \int_{0}^\infty\frac{|ty-\sin(ty)|}{t^3}\, \d t =\frac{\pi y^2}{4}=\frac{\pi (-y)^2}{4}.
                                           \end{equation}}
                    \argument{\lref{eq7};\lref{eq6}}{that for all $y\in\R\backslash\{0\}$ it holds that 
                    \begin{equation}\llabel{eq8}
                         \int_{0}^\infty\frac{|ty-\sin(ty)|}{t^3}\, \d t=\frac{\pi y^2}{4}.
                    \end{equation}}
                    \argument{\lref{eq8}; the fact that $\int_{0}^\infty\frac{|0t-\sin(0t)|}{t^3}\, \d t=0$}{that for all $y\in \R$ it holds that
                    \begin{equation}\llabel{eq8'}
                          \int_{0}^\infty\frac{|ty-\sin(ty)|}{t^3}\, \d t=\frac{\pi y^2}{4}.
                    \end{equation}} 
                   \argument{\lref{eq8'};}{ \cref{item 3: integral}\dott}
                    \startnewargseq
                    \argument{\lref{arg3'};\cref{item 3: integral};}{that for all $y\in \R$ it holds that
                    \begin{equation}\llabel{eq9}
                       \int_{0}^\infty\frac{ty-\sin(ty)}{t^3}\,\d t= \lim_{s\to 0}\int_{s}^\infty\frac{ty-\sin(ty)}{t^3}\,\d t=\frac{\pi y|y|}{4}.
                    \end{equation}}
                    \argument{\lref{eq9};}{\cref{item 2: integral}\dott}
                 \end{aproof}

                 In the next elementary and well-known result, \cref{lem: tg3}, we present a characterization for square integrable
random variables with non-vanishing variance in terms of their characteristic functions. Only for completeness we include here a detailed proof for \cref{lem: tg3}.
                 \begin{tcolorbox}[colback=white!95!gray,
                  colframe=black,
                  boxrule=0.5pt,
                  sharp corners,
                  enhanced,
                  breakable,
                 ]
                 \begin{athm}{lemma}{lem: tg3}[\textcolor{red}{Characteristic functions of random variables with non-vanishing variance}]
                     Let $(\Omega,\cF,\P)$ be a probability space and let $Y\colon\Omega\to\R$ be a random variable with $\E[Y^2]<\infty$. Then the following two statements are equivalent:
                     \begin{enumerate}[label=(\roman*)]
                         \item \label{item 1: tg3}  It holds that $\operatorname{Var}( Y ) > 0$.
                         \item \label{item 2: tg3} There exists $\delta \in (0,\infty)$ such that for all $t \in (0,\delta)$ it holds that \begin{equation}
                         | \E[ \exp(tY\imag) ] | \leq 1 - \delta t^2 \leq \exp( - \delta t^2 ).
                         \end{equation}
                     \end{enumerate}
                 \end{athm}
                 \end{tcolorbox}
                 \begin{aproof}
                 Throughout this proof let $f\colon\R\to\R$ satisfy for all $t\in \R$ that
                 \begin{equation}\llabel{def: f}
                     f(t)=| \E[ \exp( t Y \imag ) ] |^2.
                 \end{equation}
                 \argument{\lref{def: f};the fact that for all $x\in \R$ it holds that $\exp(x\imag)=\cos(x)+\imag\sin(x)$}{that for all $t\in\R$ it holds that
                 \begin{equation}\llabel{eq1}
                 \begin{split}
                   f(t)&=| \E[ \exp( t Y \imag ) ] |^2=\bigl|\E[\cos(tY)+\imag \sin(tY)]\bigr|^2\\
                   &=\bigl|\E[\cos(tY)]+\imag\, \E[\sin(tY)]\bigr|^2=(\E[\cos(tY)])^2+(\E[\sin(tY)])^2.
                   \end{split}
                 \end{equation}}
                 \argument{\lref{eq1};the fact that $\E[|Y|]<\infty$; the dominated convergence theorem}{that for all $t\in \R$ it holds that
                 \begin{equation}\llabel{eq2}
                     f'(t)=- 2 \E[ \cos( t Y ) ] \E[ \sin( t Y ) Y ] + 2 \E[ \sin( t Y ) ] \E[ \cos( t Y ) Y ].
                 \end{equation}}
                 \argument{\lref{eq2};the assumption that $\E[Y^2]<\infty$; the dominated convergence theorem}{that for all $t\in \R$ it holds that
                 \begin{equation}\llabel{eq3}
                 \begin{split}
                    & f''(t)\\&=2 ( \E[ \sin( t Y ) Y ] )^2 - 2 \E[ \cos( t Y ) ] \E[ \cos( t Y ) Y^2 ] + 2 ( \E[ \cos( t Y ) Y ] )^2 - 2 \E[ \sin( t Y ) ] \E[ \sin( t Y ) Y^2 ].
                     \end{split}
                 \end{equation}}
                 \argument{\lref{eq1};}{that 
                 \begin{equation}\llabel{eq4}
                     f'(0)=0.
                 \end{equation}}
                 \argument{\lref{eq3};}{that 
                 \begin{equation}\llabel{eq5}
                     f''(0)=-2\E[Y^2]+2(\E[Y])^2=-2\operatorname{Var}(Y).
                 \end{equation}}
                 \argument{\lref{eq5};\lref{def: f};\lref{eq4};}{that
                 \begin{equation}\llabel{eq6}
                   \lim_{t\to 0} \frac{| \E[ \exp( t Y \imag ) ] |^2-1}{t^2}= \lim_{t\to 0} \frac{f(t)-1}{t^2}=  \lim_{t\to 0} \frac{f(t)-f(0)-f'(0) t}{2t^2}=\frac{f''(0)}{2}=-\operatorname{Var}(Y).
                 \end{equation}}
                 In the following we prove that (\cref{item 1: tg3}$\rightarrow$\cref{item 2: tg3}). In the following we thus assume that 
                 \begin{equation}\llabel{assume 1}
                    \operatorname{Var}(Y)>0.
                 \end{equation}
                 \startnewargseq
                 \argument{\lref{eq6};\lref{assume 1}}{that
                 \begin{equation}\llabel{eq6'}
                      \lim_{t\to 0}\frac{ |\E[\exp(tY\imag)]|^2-1}{t^2}<0.
                 \end{equation}}
                     \argument{\lref{eq6'};}{that there exist $\delta,\rho\in (0,\infty)$ which satisfy for all $t\in (0,\delta)$ that
                     \begin{equation}\llabel{eq9}
                         \frac{1-|\E[\exp(tY\imag)]|^2}{t^2}\geq \rho.
                     \end{equation}}
                     \startnewargseq
                     \argument{\lref{eq9};}{that for all $t\in (0,\delta)$ it holds that \llabel{eq10}
                         $|\E[\exp(tY\imag)]|^2\leq 1-\rho t^2$\dott
                    }
                     \argument{\lref{eq10};the fact that for all $x\in \R$ it holds that $1-x\leq e^{-x}$}{that for all $t \in ( 0, \min\{ \delta, \nicefrac{ \rho }{ 2 }\}  )$ it holds that
                     \begin{equation}\llabel{arg3}
                         \textstyle |\E[\exp(tY\imag)]|\leq  \sqrt{1-\rho t^2}\leq 1-\frac{\rho t^2}{2}\leq 1-\min\{\delta,\frac{\rho}{2}\} t^2\leq \exp\bigl( - \min\{ \delta, \nicefrac{\rho }{ 2 } \} t^2 \bigr).
                     \end{equation}}
                     \argument{\lref{arg3};}{that (\cref{item 1: tg3}$\rightarrow$\cref{item 2: tg3})\dott}
                      In the following we prove that (\cref{item 2: tg3}$\rightarrow$\cref{item 1: tg3}). In the following we thus assume that there exists $\delta \in (0,\infty)$ which satisfies for all $t \in (0,\delta)$ that
                 \begin{equation}\llabel{assume 2}
                     | \E[ \exp(tY\imag) ] | \leq 1 - \delta t^2 \leq \exp( - \delta t^2 ).
                 \end{equation}
                 \startnewargseq
                 \argument{\lref{eq6};\lref{assume 2};}{that
                 \begin{equation}\llabel{eq11}
                     \operatorname{Var}(Y)= \lim_{t\to 0}\frac{1- |\E[\exp(tY\imag)]|^2}{t^2}\geq  \lim_{t\to 0}\frac{1- \exp(-2\delta t^2)}{t^2}=2\delta>0.
                 \end{equation}}
                 \argument{\lref{eq11};}{that (\cref{item 2: tg3}$\rightarrow$\cref{item 1: tg3})\dott}
                 \end{aproof}
                 Let $s \colon \R \to \R$ be the signed square function (satisfying for all $x \in \R$ that $s(x) = x |x|)$, let $( \Omega, \cF, \P )$ be a probability space, let $\alpha \in (0,1)$, let $X_{ n,m } \colon \Omega \to \R$, $(n,m)\in \N^2$, be \iid\ random variables, and for every $M \in \N$ let $Z_M \colon \Omega \to \R$ satisfy $Z_M =\sum_{k=0}^\infty\alpha^{k}[\textstyle \frac 1M\sum_{m=1}^MX_{k+1,m}\bigr]$ (stationary momentum random variable) and let $z_M \in \R$ satisfy
\begin{equation}
  z_M = \E[ s( Z_M ) ] 
\end{equation}
(expectation of the signed square function evaluated at the stationary momentum random variable with mini-batch size $M$). 
In \cref{conclude: tg2} in \cref{lem: tg2} we prove that $z_M$ converges to zero as the mini-batch size $M$ tends to infinity and we also sharply specify the rate of convergence. More specifically, from \cref{conclude: tg2} in \cref{lem: tg2} we can conclude for every $r \in (0,\infty)$ that there exists $c \in \R$ such that for all $M \in \N$ it holds that
$| z_M | \leq c M^{ - r }$
(convergence of $z_M$ with convergence rate $r$ to zero) if and only if the rate of convergence $ r \leq \nicefrac{ 3 }{ 2 } $ is smaller or equal than $ \nicefrac{ 3 }{ 2 } $.
           \begin{tcolorbox}[colback=white!95!gray,
                  colframe=black,
                  boxrule=0.5pt,
                  sharp corners,
                  enhanced,
                  breakable,
                 ]
                 \begin{athm}{prop}{lem: tg2}
                     Let $(\Omega,\cF,\P)$ be a probability space, let $\alpha\in [0,1)$, let $X_{n,m}\colon \Omega\to\R$, $(n,m)\in \N^2$, be bounded \iid\ random variables, assume $\E[X_{1,1}]=0\neq \E[(X_{1,1})^3]$, and for every $M\in \N$ let $Z_M\colon \Omega\to\R$ satisfy $Z_M=\sum_{k=0}^\infty\alpha^{k}[\textstyle \frac 1M\sum_{m=1}^MX_{k+1,m}\bigr]$. Then there exists $\fC\in \R\backslash\{0\}$ such that 
                     \begin{equation}\label{conclude: tg2}
                        \textcolor{magenta}{ \limsup\nolimits_{M\to\infty} \bigl|M^{3/2}\E\bigl[Z_M|Z_M|\bigr]-\fC\bigr|=0}.
                     \end{equation}
                 \end{athm}
                 \end{tcolorbox}
                 \begin{aproof}
                     Throughout this proof for every $m\in \N$ let $Y_m\colon \Omega\to\R$ satisfy
    \begin{equation}\llabel{def: Y}
       \displaystyle Y_m=\sum_{k=0}^\infty\alpha^kX_{k+1,m}
    \end{equation} 
    and let $A,B\in \R$ satisfy $A=(1-\alpha^3)^{-1}\E[(X_{1,1})^3]$ and $B=(1-\alpha^2)^{-1}\E[(X_{1,1})^2]$.
    \argument{\cref{lem: integral};}{that for all $M\in \N$ it holds that
    \begin{equation}\llabel{eq1}
        MZ_M|Z_M|=\frac{4}{\pi}\int_{0}^\infty \frac{tM^{1/2}Z_M-\sin (tM^{1/2}Z_M)}{t^3}\,\d t.
    \end{equation}}
    \argument{the assumption that $X_{n,m}$, $(n,m)\in \N^2$ are bounded \iid\@;\cref{cor: upper order of Adam coefficient}}{that 
    \begin{equation}\llabel{eq2}
    \begin{split}
     \textstyle \sup\limits_{M\in\N} \E\biggl[\biggl|(1-\alpha)\sum\limits_{k=0}^\infty \alpha^k\biggl[\frac{1}{\sqrt{M}}\sum\limits_{m=1}^MX_{k+1,m}\biggr]\biggr|^{2}\biggr]&=\textstyle \sup\limits_{M\in\N}M \,\E\biggl[\biggl|(1-\alpha)\sum\limits_{k=0}^\infty \alpha^k\biggl[\frac{1}{M}\sum\limits_{m=1}^MX_{k+1,m}\biggr]\biggr|^{2}\biggr]\\
    &<\infty.
          \end{split}
    \end{equation}}
    \argument{\lref{eq2};the fact that for all $M\in \N$ it holds that $Z_M=\sum_{k=0}^\infty\alpha^{k}[\textstyle \frac 1M\sum_{m=1}^MX_{k+1,m}\bigr]$}{that \begin{equation}\llabel{eq3}
        \sup_{M\in \N}\E[M(Z_M)^2]<\infty.
    \end{equation}}
    \argument{\cref{lem: integral};\lref{eq3};}{that for all $M\in \N$ it holds that
    \begin{equation}\llabel{eq4}
    \begin{split}
        &\E\biggl[\frac{4}{\pi}\int_{0}^\infty \frac{|tM^{1/2}Z_M-\sin (tM^{1/2}Z_M)|}{t^3}\,\d t\biggr]\\
        &=\frac{4}{\pi} \,\E\biggl[\int_{0}^\infty \frac{|tM^{1/2}Z_M-\sin (tM^{1/2}Z_M)|}{t^3}\,\d t\biggr]=\E[(M^{1/2}Z_M)^2]<\infty.
        \end{split}
    \end{equation}}
    \argument{\lref{eq4};\lref{eq1}; the Fubini-Tonelli theorem}{that for all $M\in \N$ it holds that
    \begin{equation}\llabel{eq5}
        \E\bigl[MZ_M|Z_M|\bigr]=\frac{4}{\pi}\int_{0}^\infty \E\biggl[\frac{tM^{1/2}Z_M-\sin (tM^{1/2}Z_M)}{t^3}\biggr]\,\d t.
    \end{equation}}
    \argument{the fact that for all $n,m\in \N$ it holds that $\E[X_{n,m}]=0$; the fact that for all $M\in \N$ it holds that $Z_M=\sum_{k=0}^\infty\alpha^{k}[\textstyle \frac 1M\sum_{m=1}^MX_{k+1,m}\bigr]$}{that for all $t\in (0,\infty)$ it holds that
\begin{equation}\llabel{eq6}
    \E[Z_M]=\textstyle\sum_{k=0}^\infty\alpha^{k}[\textstyle \frac1M\sum_{m=1}^M\E[X_{k+1,m}]\bigr]=0.
\end{equation}}
\argument{\lref{eq5};\lref{eq6}}{that for all $M\in\N$ it holds that
\begin{equation}\llabel{eq7}
\begin{split}
   M \E\bigl[Z_M|Z_M|\bigr]&=\frac{-4}{\pi}\int_{0}^\infty \frac{\E[\sin (tM^{1/2}Z_M)]}{t^3}\,\d t=\frac{-4}{\pi}\int_{0}^\infty \frac{\E\bigl[\imagine(\exp(tM^{1/2}Z_M \imag))\bigr]}{t^3}\,\d t\\
   &=\frac{-4}{\pi}\int_{0}^\infty\frac{\imagine\bigl(\E[\exp(tM^{1/2}Z_M\bfi)]\bigr)}{t^3}\,\d t.
    \end{split}
\end{equation}}
\argument{\lref{def: Y}; the fact that for all $M\in \N$ it holds that $Z_M=\sum_{k=0}^\infty\alpha^{k}[\textstyle \allowbreak \frac1M\sum_{m=1}^M\allowbreak X_{k+1,m}\bigr]$}{that for all $M\in \N$ it holds that
\begin{equation}\llabel{eq8}
   \textstyle Z_M=\frac 1M\bigl[\sum_{m=1}^MY_m\bigr].
\end{equation}}
\argument{\lref{eq8};the fact that $Y_m$, $m\in\N$, are \iid\@}{that for all $t\in (0,\infty)$, $M\in\N$ it holds that
\begin{equation}\llabel{eq9}
\begin{split}
\textstyle\E[\exp(tM^{1/2}Z_M\bfi)]&=\textstyle\E\bigl[\exp\big(t\sqrt{M}\sum_{m=1}^M \frac{Y_m\imag}{M}\bigr)\bigr]=\prod_{m=1}^M\E\bigl[\exp\bigl(t\sqrt{M}(\frac{Y_m\imag}{M})\bigr)\bigr]\\
    &\textstyle=\bigl(\E[\exp(tM^{-1/2}Y_1\imag)]\bigr)^M.
    \end{split}
\end{equation}}
\argument{\lref{eq9};\lref{eq7};}{for all $M\in \N$ that
\begin{equation}\llabel{eq9'}
    M ^{3/2}\E\bigl[Z_M|Z_M|\bigr]=\frac{-4}{\pi}\int_{0}^\infty\frac{\sqrt{M}\imagine\bigl(\bigl(\E[\exp(tM^{-1/2}Y_1\imag)]\bigr)^M\bigr)}{t^3}\,\d t.
\end{equation}}
\argument{Taylor expansion;}{that there exists $\rho\in (0,\infty)$ such that for all $y\in \R$ it holds that
\begin{equation}\llabel{eq10}
    \biggl|\exp(y\imag )-\biggl(1+y\imag +\frac{y^2\imag^2}{2}+\frac{y^3\imag^3}{6}\biggr)\biggr|\leq \rho|y^4\imag^4|.
\end{equation}}
\argument{\lref{eq10};}{that there exists $\rho\in (0,\infty)$ such that for all $y\in \R$ it holds that
\begin{equation}\llabel{eq11}
    \biggl|\exp(y\imag )-\biggl(1+y\imag -\frac{y^2}{2}-\frac{y^3\imag}{6}\biggr)\biggr|\leq \rho|y^4|.
\end{equation}}
\argument{\lref{eq11};}{that there exists $\rho\in (0,\infty)$ such that for all $t\in(0,\infty)$, $M\in \N$ it holds that
\begin{equation}\llabel{eq12}
\begin{split}
   &\textstyle \Bigl|\E[\exp(tM^{-1/2} Y_1\imag)]-\E\bigl[1+tM^{-1/2} Y_1\imag -\frac{(tM^{-1/2} Y_1)^2}{2}-\frac{(tM^{-1/2} Y_1)^3\imag}{6}\bigr]\Bigr|\\
   &\leq \textstyle \E\Bigl[\bigl|\exp(tM^{-1/2} Y_1\imag)-\bigl(1+tM^{-1/2} Y_1\imag -\frac{(tM^{-1/2} Y_1)^2}{2}-\frac{(tM^{-1/2} Y_1)^3\imag}{6}\bigr)\bigr|\Bigr]\leq \E[\rho|tM^{-1/2} Y_1|^4].
   \end{split}
\end{equation}}
\argument{\lref{eq12};}{that there exists $\rho\in (0,\infty)$ such that for all $t\in(0,\infty)$, $M\in \N$ it holds that
\begin{equation}\llabel{eq13}
    \textstyle \Bigl|\E[\exp(tM^{-1/2} Y_1\imag)]-\E\bigl[1+tM^{-1/2} Y_1\imag -\frac{t^2M^{-1} (Y_1)^2}{2}-\frac{t^3M^{-3/2} (Y_1)^3\imag}{6}\bigr]\Bigr|\leq \rho^4t^4 M^{-2}\E[|Y_1|^4].
\end{equation}}
\argument{\lref{def: Y}; the fact that $X_{n,m}$, $(n,m)\in \N^2$, are bounded}{that there exists $\varrho\in (0,\infty)$ which satisfies
\begin{equation}\llabel{evd1}
    \P(|Y_1|\leq \varrho)=1.
\end{equation}}
\startnewargseq
\argument{\lref{def: Y};the fact that for all $n\in \N$ it holds that $\E[X_{n,1}]=0$;}{that
\begin{equation}\llabel{evd2}
    \textstyle\E[Y_1]=\sum_{k=0}^\infty\alpha^k\E[X_{k+1,1}]=0.
\end{equation}}
\argument{\lref{def: Y};\cref{item 2: non-negative moment} in \cref{lem: non-negative moment}}{that 
\begin{equation}\llabel{evd3}
    \E[(Y_1)^3]=\textstyle\E\bigl[\bigl(\sum_{k=0}^\infty\alpha^k X_{k+1,1}\bigr)^3\bigr]=(1-\alpha^3)^{-1}\E[(X_{1,1})^3]=A.
\end{equation}
}
\argument{\lref{def: Y};the fact that $X_{k,1}$, $k\in \N$, are centered \iid\@}{that
\begin{equation}\llabel{evd4}
    \begin{split}
    \E[(Y_1)^2]&=\textstyle\E\bigl[\bigl(\sum_{k=0}^\infty\alpha^k X_{k+1,1}\bigr)^2\bigr]=\sum_{k,l\in (\N_0)^2}\alpha^{k}\alpha^l\E[X_{k+1,1}X_{l+1,1}]\\
   &\textstyle=\sum_{k=0}^\infty\alpha^{2k}\E[(X_{k+1,1})^2]=(1-\alpha^2)^{-1}\E[(X_{1,1})^2]=B.
    \end{split}
\end{equation}}
\argument{\lref{eq13};\lref{evd1};\lref{evd2};\lref{evd3};\lref{evd4};}{that there eixist $\rho,\fC\in (0,\infty)$ such that for all $t\in (0,\infty)$, $M\in \N$ that 
\begin{equation}\llabel{eq14}
     \textstyle \Bigl|\E[\exp(tM^{-1/2} Y_1\imag)]-\bigl(1 -\frac{t^2M^{-1} B}{2}-\frac{t^3M^{-3/2} A\imag}{6}\bigr)\Bigr|\leq \rho^4t^4 M^{-2}\E[|Y_1|^4]\leq \fC t^4M^{-2}.
\end{equation}}
\argument{\lref{eq14};}{that for all $t\in (0,\infty)$ there exist $r_M^t\in \bbC$, $M\in\N$, which satisfy for all $M\in \N$ that
\begin{equation}\llabel{eq15}
    \limsup_{N\to\infty}(N^{3/2}|r_N^t|)=0 \qqandqq  \E[\exp(tM^{-1/2} Y_1\imag)]=1 -\frac{t^2M^{-1} B}{2}-\frac{t^3M^{-3/2} A\imag}{6}+r_M^t.
\end{equation}}
\startnewargseq
\argument{the assumption that $\E[X_{1,1}]=0\neq \E[(X_{1,1})^3]$;}{that \llabel{arg1} $\P(X_{1,1}=0)<1$\dott}
\argument{\lref{arg1};}{that \llabel{arg2} $\E[(X_{1,1})^2]>0$\dott}
\argument{\lref{arg2};the fact that $0\leq \alpha<1$}{that \llabel{arg3} $B>0$\dott}
\argument{\lref{arg3};\lref{eq15};\cref{lem: tg1};}{that for all $t\in (0,\infty)$ it holds that
\begin{equation}\llabel{eq16}
    \begin{split}
      & \textstyle \lim_{M\to\infty}\bigl[\sqrt{M}\imagine\bigl(\bigl(\E[\exp(tM^{-1/2} Y_1\imag)]\bigr)^M\bigr)\bigr]\\
      &=\textstyle
        \lim_{M\to\infty}\bigl[\sqrt{M}\imagine\bigl(\bigl(1 -\frac{t^2M^{-1} B}{2}-\frac{t^3M^{-3/2} A\imag}{6}+r_M^t\bigr)^M\bigr)\bigr]=-\frac{At^3\exp(-Bt^2/2)}{6}.
    \end{split}
\end{equation}}
\argument{\lref{eq16};}{for all $t\in (0,\infty)$ that
\begin{equation}\llabel{eq16'}
\lim_{M\to\infty}\frac{\sqrt{M}\imagine\bigl(\bigl(\E[\exp(tM^{-1/2} Y_1\imag)]\bigr)^M\bigr)}{t^3}=-\frac{A\exp(-Bt^2/2)}{6}.
\end{equation}}
\argument{the fact that for all $z\in \bbC$, $M\in \N$ it holds that $\overline{z^M}=\overline{z}^M$}{that for all $z\in \bbC$, $M\in\N$ it holds that
\begin{equation}\llabel{eqt1}
    2\imag\imagine(z^M)=z^M-\overline{z^M}=z^M-\overline{z}^M=(z-\overline{z})\sum_{j=0}^{M-1}z^{M-1-j}\overline{z}^j.
\end{equation}}
\argument{\lref{eqt1};}{for all $z\in \bbC$, $M\in \N$ it holds that
\begin{equation}\llabel{eqt2}
\begin{split}
    |2\imag\imagine(z^M)|&=\biggl|\textstyle (z-\overline{z})\sum\limits_{j=0}^{M-1}z^{M-1-j}\overline{z}^j\biggr|\leq |z-\overline{z}|\biggl[\sum\limits_{j=0}^{M-1}|z|^{M-1-j}|\overline{z}|^j\biggr]\\
    &\leq  \textstyle |z-\overline{z}|\biggl[\sum\limits_{j=0}^{M-1}|z|^{M-1}\biggr]=|z-\overline{z}| M|z|^{M-1}=|2\imagine(z)|M|z|^{M-1}.
    \end{split}
\end{equation}}
\argument{\lref{eqt2};}{for all $z\in \bbC$, $M\in \N$ that
\begin{equation}\llabel{eqt3}
    |\imagine(z^M)|\leq M|\imagine(z)||z|^{M-1}.
\end{equation}}
\argument{\cref{lem: tg3};the fact that $0=\E[Y_1]<\E[(Y_1)^2]$}{that there exist $\delta,\scrc\in (0,\infty)$ which satisfy for all $t\in (0,\infty)$ with $tM^{-1/2}<\delta$ that
\begin{equation}\llabel{eqt4}
    \bigl|\E[\exp(tM^{-1/2} Y_1\imag)]\bigr|\leq \exp(-\scrc (tM^{-1/2})^2)=\exp(-\scrc t^2M^{-1}).
\end{equation}}
\startnewargseq
\argument{\lref{eqt3};\lref{eqt4};the fact that $\E[Y_1]=0$; the fact that for all $M\in \N\backslash\{1\}$ it holds that $(M-1)M^{-1}\geq \frac 12$}{that for all $t\in (0,\infty)$, $M\in \N\backslash\{1\}$  with $t< \delta M^{1/2}$ it holds that
\begin{equation}\llabel{eq17}
\begin{split}
   & \bigl|\imagine\bigl(\bigl(\E[\exp(tM^{-1/2} Y_1\imag)]\bigr)^M\bigr)\bigr|\leq M\bigl|\imagine\bigl(\E[\exp(tM^{-1/2} Y_1\imag)]\bigr)\bigr| \bigl|\E[\exp(tM^{-1/2} Y_1\imag)]\bigr|^{M-1}\\
    &\leq M\bigl|\imagine\bigl(\E[\exp(tM^{-1/2} Y_1\imag)]\bigr)\bigr|(\exp(-\scrc t^2M^{-1}))^{M-1}\\
    &=M\bigl|\E[\sin(tM^{-1/2}Y_1)]\bigr|\exp(-\scrc t^2(M-1
    )M^{-1})\\
    &\leq M\bigl|\E[\sin(tM^{-1/2}Y_1)-tM^{-1/2}Y_1]\bigr|\exp(-\scrc t^2/2).
    \end{split}
\end{equation}}
\argument{the fact that $\lim_{t\to 0} \frac{|\sin(t)-t|}{|t|^3}=\frac 16$;the fact that $\lim_{t\to \infty} \frac{|\sin(t)-t|}{|t|^3}=0$; the fact that $\lim_{t\to -\infty} \frac{|\sin(t)-t|}{|t|^3}=0$;}{that there exists $\fC\in (0,\infty)$ which satisfies for all $t\in \R$ that
\begin{equation}\llabel{eq18}
    |\sin(t)-t|\leq \fC|t|^3.
\end{equation}}
\startnewargseq
\argument{\lref{evd1};\lref{eq17};\lref{eq18};}{that for all $t\in (0,\infty)$ it holds that
\begin{equation}\llabel{eq20.1}
\begin{split}
  & \sup_{M\in \N\backslash\{1\}}\frac{\sqrt{M}\bigl|\imagine\bigl(\bigl(\E[\exp(tM^{-1/2} Y_1\imag)]\bigr)^M\bigr)\bigr|\mathbbm 1_{(0,\delta\sqrt{M})}(t)}{t^3}\\
  &\leq \sup_{M\in \N\backslash\{1\}}\frac{M\sqrt{M}\bigl|\E[\sin(tM^{-1/2}Y_1)-tM^{-1/2}Y_1]\bigr|\exp(-\scrc t^2/2)\mathbbm 1_{(0,\delta\sqrt{M})}(t)}{t^3}\\
  &\leq \sup_{M\in \N\backslash\{1\}}\frac{M\sqrt{M}\fC t^3M^{-3/2} \E[|Y_1|^3]\exp(-\scrc t^2/2)\mathbbm 1_{(0,\delta\sqrt{M})}(t)}{t^3}\leq \fC \varrho ^3\exp(-\scrc t^2/2).
   \end{split}
\end{equation}}
\argument{\lref{eq16'}}{that for all $t\in (0,\infty)$ it holds that
\begin{equation}\llabel{eq20.2}
     \lim_{M\to\infty}\frac{\sqrt{M}\imagine\bigl(\bigl(\E[\exp(tM^{-1/2} Y_1\imag)]\bigr)^M\bigr)\mathbbm 1_{(0,\delta \sqrt{M})}(t)}{t^3}=-\frac{A\exp(-Bt^2/2)}{6}.
\end{equation}}
\argument{\lref{eq20.2};\lref{eq20.1};the dominated convergence theorem;the fact that 
\begin{equation}
    \int_{0}^\infty \fC\varrho \exp(-\scrc t^2/2)<\infty
\end{equation}}{that 
\begin{equation}\llabel{eq21}
\begin{split}
&\lim_{M\to\infty}\biggl[\int_{0}^{\delta\sqrt{M}}\frac{\sqrt{M}\imagine\bigl(\bigl(\E[\exp(tM^{-1/2}Y_1\imag)]\bigr)^M\bigr)}{t^3}\,\d t\biggr]\\
&=\lim_{M\to\infty}\biggl[\int_{0}^{\infty}\frac{\sqrt{M}\imagine\bigl(\bigl(\E[\exp(tM^{-1/2}Y_1\imag)]\bigr)^M\bigr)\mathbbm 1_{(0,\delta\sqrt{M})}(t)}{t^3}\,\d t\biggr]=\int_{0}^\infty \frac{-A\exp(-Bt^2/2)}{6}\,\d t.
\end{split}
\end{equation}}
\argument{\lref{eq9};the fact that for all $x\in \R$ it holds that $\imagine(\exp(x\imag))=\sin(x)$}{that
\begin{equation}\llabel{eq22}
\begin{split}
&\lim_{M\to\infty}\biggl[\int_{\delta\sqrt{M}}^{\infty}\frac{\sqrt{M}\bigl|\imagine\bigl(\bigl(\E[\exp(tM^{-1/2}Y_1\imag)]\bigr)^M\bigr)\bigr|}{t^3}\,\d t\biggr]\\
&=\lim_{M\to\infty}\biggl[\int_{\delta\sqrt{M}}^{\infty}\frac{\sqrt{M}\big|\imagine(\E[\exp(tM^{1/2}Z_M\bfi)])\bigr|}{t^3}\,\d t\biggr]\\
&=\lim_{M\to\infty}\biggl[\int_{\delta\sqrt{M}}^{\infty}\frac{\sqrt{M}\big|\E[\sin(tM^{1/2}Z_M)\bigr|}{t^3}\,\d t\biggr]\\
&\leq \lim_{M\to\infty}\biggl[\int_{\delta\sqrt{M}}^{\infty}\frac{\sqrt{M}}{t^3}\,\d t\biggr]=\lim_{M\to\infty}\biggl[\frac{\sqrt{M}}{2(\delta\sqrt{M})^2}\biggr]=0.
\end{split}
\end{equation}}
\argument{\lref{eq22};\lref{eq21};}{that
\begin{equation}\llabel{eq23}
\lim_{M\to\infty}\biggl[\int_{0}^{\infty}\frac{\sqrt{M}\imagine\bigl(\bigl(\E[\exp(tM^{-1/2}Y_1\imag)]\bigr)^M\bigr)}{t^3}\,\d t\biggr]=\int_{0}^\infty \frac{-A\exp(-Bt^2/2)}{6}\,\d t.
\end{equation}}
\argument{\lref{eq23};\lref{eq9'};the fact that $A\neq 0$; the fact that $B>0$}{that
\begin{equation}\llabel{eq24}
   \lim_{M\to\infty} \bigl( M^{3/2}\E\bigl[Z_M|Z_M|\bigr]\bigr)=\frac{4}{\pi}\int_{0}^\infty \frac{A\exp(-Bt^2/2)}{6}\,\d t\neq 0.
\end{equation}}
                 \end{aproof}

                 \subsection{Non-vanishing evaluations of the MUON vector field for non-symmetric data}\label{subsec: non vanishing}

                  In the next elementary result, \cref{lem: NS polynomial}, we establish in \cref{conclude: NS polynomial} polynomial representations without quadratic terms for the \NS\ functions. \cref{lem: NS polynomial} is used in our proof of \cref{non-zero vector field 2} below, which shows in \cref{conclude: non-zero vector field 2} in the case of centered data $\E[ X_{1,1} ] = 0$ that the \MUON\ vector field does not vanish at $0$.
 \begin{tcolorbox}[colback=white!95!gray,
                  colframe=black,
                  boxrule=0.5pt,
                  sharp corners,
                  enhanced,
                  breakable,
                 ]
                 \begin{athm}{lemma}{lem: NS polynomial}[\textcolor{red}{Polynomial representations without quadratic terms for the \NS\ functions}]
                    Let $K\in \N_0$ and let $\nscoe=(\nscoe_{i,j})_{(i,j)\in (\N_0)^2}\colon(\N_0)^2\to \R$ satisfy $\textcolor{magenta}{\min\{\nscoe_{0,0},\nscoe_{0,1}\}>0}$ and $\textcolor{magenta}{\#(\nscoe^{ - 1 }( \R\backslash\{0\} )) <\infty}$. Then there exists a polynomial $P\colon \R\to\R$ such that for all $\theta,\Theta\in\R$ with $ \textcolor{magenta}{\Theta=\theta(\nscoe_{0,0}|\theta|+\nscoe_{0,1})^{-1}}$ it holds that
\begin{equation}\label{conclude: NS polynomial}
    \textcolor{magenta}{\textstyle \TNewtonSchulzAlgorithm{\nscoe}{K}{\theta}=\bigl[\prod_{k=1}^K\nscoe_{k,0}\bigr] \Theta + \Theta^3 P( \Theta )}
\end{equation}
(cf.\ \cref{definition: NS}).
                 \end{athm}
                 \end{tcolorbox}
                 \begin{aproof}
                 Throughout this proof let $\psi\colon \R\to\R$ satisfy for all $\theta\in \R$ that
                 \begin{equation}\llabel{def: T}
        \psi(\theta)=\theta(\nscoe_{0,0}|\theta|+\nscoe_{0,1})^{-1}.
                 \end{equation}
                   In the following we prove that for all $k\in \N_0$ there exists a polynomial $P\colon \R\to\R$ such that for all $\theta\in\R$ it holds that
\begin{equation}\llabel{need to prove}
   \textstyle \TNewtonSchulzAlgorithm{\nscoe}{k}{\theta}=\bigl[\prod_{n=1}^k\nscoe_{n,0}\bigr]\theta(\nscoe_{0,0}|\theta|+\nscoe_{0,1})^{-1}+\theta^3(\nscoe_{0,0}|\theta|+\nscoe_{0,1})^{-3}P\bigl(\theta(\nscoe_{0,0}|\theta|+\nscoe_{0,1})^{-1}\bigr).
\end{equation}
We prove \lref{need to prove} by induction on $k\in\N_0$. For the base case $k=0$ note that \cref{T_B_D} shows that for all $\theta\in \R$ it holds that
\begin{equation}
\begin{split}
    \textstyle \TNewtonSchulzAlgorithm{\nscoe}{0}{\theta}=\theta(\nscoe_{0,0}|\theta|+\nscoe_{0,1})^{-1}\textstyle=\bigl[\prod_{n=1}^0\nscoe_{n,0}\bigr]\theta(\nscoe_{0,0}|\theta|+\nscoe_{0,1})^{-1}+\theta^3(\nscoe_{0,0}|\theta|+\nscoe_{0,1})^{-3}0.
    \end{split}
\end{equation}
This proves \lref{need to prove} in the case $k=0$. For the induction step we assume that there exist $k\in\N$ and a polynomial $P\colon \R\to\R$ which satisfy for all $\theta\in\R$ that
\begin{equation}\llabel{assume}
   \textstyle \TNewtonSchulzAlgorithm{\nscoe}{k-1}{\theta}=\bigl[\prod_{n=1}^{k-1}\nscoe_{n,0}\bigr]\theta(\nscoe_{0,0}|\theta|+\nscoe_{0,1})^{-1}+\theta^3(\nscoe_{0,0}|\theta|+\nscoe_{0,1})^{-3}P\bigl(\theta(\nscoe_{0,0}|\theta|+\nscoe_{0,1})^{-1}\bigr).
\end{equation}
\argument{\lref{def: T};\lref{assume};}{that for all $\theta\in \R$ it holds that
\begin{equation}\llabel{eq1}
    \textstyle \TNewtonSchulzAlgorithm{\nscoe}{k-1}{\theta}=\bigl[\prod_{n=1}^{k-1}\nscoe_{n,0}\bigr]\psi(\theta)+(\psi(\theta))^3P(\psi(\theta)).
\end{equation}}
\argument{\lref{eq1};\cref{T_B_D}}{that for all $\theta\in \R$ it holds that
\begin{equation}\llabel{eq2}
\begin{split}
\TNewtonSchulzAlgorithm{\nscoe}{k}{\theta}
&=\textstyle\sum_{i=0}^\infty\nscoe_{k,i} (\TNewtonSchulzAlgorithm{\nscoe}{k-1}{\theta})^{2i+1}=\textstyle\sum_{i=0}^\infty\nscoe_{k,i} \bigl(\bigl[\prod_{n=1}^{k-1}\nscoe_{n,0}\bigr]\psi(\theta)+(\psi(\theta))^3P(\psi(\theta))\bigr)^{2i+1}\\
&=\textstyle \nscoe_{k,0} \bigl(\bigl[\prod_{n=1}^{k-1}\nscoe_{n,0}\bigr]\psi(\theta)+(\psi(\theta))^3P(\psi(\theta))\bigr)\\
&+\textstyle\sum_{i=1}^\infty\nscoe_{k,i} \bigl(\bigl[\prod_{n=1}^{k-1}\nscoe_{n,0}\bigr]\psi(\theta)+(\psi(\theta))^3P(\psi(\theta))\bigr)^{2i+1}\\
&=\textstyle  \bigl[\prod_{n=1}^{k}\nscoe_{n,0}\bigr]\psi(\theta)+\nscoe_{k,0}(\psi(\theta))^3P(\psi(\theta))\\
&+(\psi(\theta))^3\textstyle\sum_{i=1}^\infty\nscoe_{k,i} (\psi(\theta))^{2i-2} \bigl(\bigl[\prod_{n=1}^{k-1}\nscoe_{n,0}\bigr]+(\psi(\theta))^2P(\psi(\theta))\bigr)^{2i+1}.
\end{split}
\end{equation}}
\argument{the assumption that $\#(\nscoe^{ - 1 }( \R\backslash\{0\} )) <\infty$}{that \llabel{argg1} $\#\{i\in\N\colon \nscoe_{k,i}\neq 0\}<\infty$\dott}
\argument{\lref{argg1};the fact that $P$ is polynomial}{that
\begin{equation}
    \textstyle\R\ni x\mapsto \nscoe_{k,0}P(x)+\textstyle\sum_{i=1}^\infty\nscoe_{k,i} x^{2i-2} \bigl(\bigl[\prod_{n=1}^{k-1}\nscoe_{n,0}\bigr]+x^2P(x)\bigr)^{2i+1}\in \R
\end{equation}
is \llabel{arg1} a polynomial\dott}
\argument{\lref{arg1};\lref{eq2}}{that there exists a polynomial $Q\colon \R\to\R$ such that for all $\theta\in \R$ it holds that
\begin{equation}\llabel{eq3}
    \TNewtonSchulzAlgorithm{\nscoe}{k}{\theta}=\textstyle\big[\prod_{n=1}^{k}\nscoe_{n,0}\big]\psi(\theta)+(\psi(\theta))^3Q(\psi(\theta)).
\end{equation}}
\argument{\lref{eq3};\lref{def: T}}{that there exists a polynomial $Q\colon \R\to\R$ such that for all $\theta\in \R$ it holds that
\begin{equation}\llabel{eq4}
   \textstyle \TNewtonSchulzAlgorithm{\nscoe}{k}{\theta}=\bigl[\prod_{n=1}^k\nscoe_{n,0}\bigr]\theta(\nscoe_{0,0}|\theta|+\nscoe_{0,1})^{-1}+\theta^3(\nscoe_{0,0}|\theta|+\nscoe_{0,1})^{-3}Q\bigl(\theta(\nscoe_{0,0}|\theta|+\nscoe_{0,1})^{-1}\bigr).
\end{equation}}
\argument{\lref{eq4}; \lref{assume}; induction}{\lref{need to prove}\dott}
\startnewargseq
\argument{\lref{need to prove};}{\cref{conclude: NS polynomial}\dott}
                 \end{aproof}
\begin{samepage}
 \begin{tcolorbox}[colback=white!95!gray,
                  colframe=black,
                  boxrule=0.5pt,
                  sharp corners,
                  enhanced,
                  breakable,
                 ]
\begin{athm}{prop}{non-zero vector field 2}[\textcolor{red}{Non-vanishing evaluations of the \MUON\ vector field for centered data}]
Let $(\Omega,\cF,\P)$ be a probability space, let $K\in\N_0$, $\alpha\in [0,1)$, $\lambda\in (0,\infty)$, let $\nscoe=(\nscoe_{i,j})_{(i,j)\in(\N_0)^2}\colon (\N_0)^2\to\R$ satisfy $\textcolor{magenta}{\min\{\nscoe_{0,0},\nscoe_{0,1}\}>0}$ and $\textcolor{magenta}{\#(\nscoe^{ - 1 }( \R\backslash\{0\} )) <\infty}$, assume for all $i\in\N\cap [0,K]$ that $\textcolor{magenta}{\nscoe_{i,0}>0}$, let $X_{n,m}\colon \Omega\to\R$, $(n,m)\in\N^2$, be bounded centered \iid\ random variables, and assume $\E[(X_{1,1})^3]\neq 0$.
Then 
\begin{equation}\label{conclude: non-zero vector field 2}
  \textcolor{magenta}{\liminf\nolimits_{M\to\infty} \bigl|M^{3/2}\E\bigl[ \TbigNewtonSchulzAlgorithm{\nscoe}{K}{\textstyle\textstyle-2\lambda(1-\alpha)\bigl[\sum_{k=0}^\infty\alpha^{k}[\textstyle \frac 1M\sum_{m=1}^MX_{k+1,m}\bigr]\bigr]}\bigr]\bigr|> 0}.
\end{equation}
(cf.\ \cref{definition: NS}).
\end{athm}
\end{tcolorbox}
\end{samepage}
\begin{aproof}
    Throughout this proof for every $n\in \N_0$, $M\in \N$ let $\bbX_n^M\colon \Omega\to\R$ satisfy
    \begin{equation}\llabel{def: bbX}
         \textstyle \bbX_n^M=\frac 1M \textstyle\sum_{m=1}^MX_{n+1,m}
    \end{equation}
  for every $M\in \N$ let $Z_M\colon \Omega\to\R$ satisfy
    \begin{equation}\llabel{def: Z}
      \textstyle  Z_M=-2\lambda(1-\alpha)\sum_{k=0}^\infty \alpha^k\bbX_k^M,
    \end{equation}
    and let $\scrc\in \R$ satisfy
    \begin{equation}\llabel{def: scrc}
       \textstyle \P(|X_{1,1}|\leq \frac{\scrc}{2\lambda})=1.
    \end{equation}
    \argument{\cref{lem: NS polynomial};}{that there exists a polynomial $P\colon \R\to\R$ such that for all $\theta,\Theta\in\R$ with $\Theta=\theta(\nscoe_{0,0}|\theta|+\nscoe_{0,1})^{-1}$ it holds that
\begin{equation}\llabel{eq1}
   \textstyle \TNewtonSchulzAlgorithm{\nscoe}{K}{\theta}=\bigl[\prod_{k=1}^K\nscoe_{k,0}\bigr] \Theta + \Theta^3 P( \Theta ).
\end{equation}}
\argument{\lref{eq1};}{that there exist $A\in \R$ and a polynomial $P\colon \R\to\R$ which satisfy for all $\theta,\Theta\in\R$ with $\Theta=\theta(\nscoe_{0,0}|\theta|+\nscoe_{0,1})^{-1}$ that
\begin{equation}\llabel{eq2}
\begin{split}
   \textstyle \TNewtonSchulzAlgorithm{\nscoe}{K}{\theta}=\bigl[\prod_{k=1}^K\nscoe_{k,0}\bigr] \Theta+\Theta^3[A+\Theta P(\Theta)]
   =\bigl[\prod_{k=1}^K\nscoe_{k,0}\bigr] \Theta+A\Theta^3+\Theta^4P(\Theta).
   \end{split}
\end{equation}}
\startnewargseq
\argument{\lref{eq2};the fact that for all $x\in \R$, $y\in [0,\infty)$ it holds that $x(y+(\nscoe_{0,1})^3)^{-1}=\frac{x}{(\nscoe_{0,1})^3}-\frac{xy}{(\nscoe_{0,1})^6}+\frac{(\nscoe_{0,1})^{-6}xy^2}{y+(\nscoe_{0,1})^3}$; the fact that for all $x\in \R$ it holds that $x(\nscoe_{0,0}|x|+\nscoe_{0,1})^{-1}=\frac{x}{\nscoe_{0,1}}-\frac{\nscoe_{0,0}x|x|}{(\nscoe_{0,1})^2}+\frac{(\nscoe_{0,0})^2x^3}{(\nscoe_{0,1})^3}-\frac{(\nscoe_{0,0})^3x|x|^3(\nscoe_{0,1})^{-3}}{\nscoe_{0,0}|x|+\nscoe_{0,1}}$}{that for all $\theta\in \R$ it holds that
\begin{equation}\llabel{eq3}
\begin{split}
       &\textstyle \TNewtonSchulzAlgorithm{\nscoe}{K}{\theta}\\
       &=\textstyle\frac{[\prod_{k=1}^K\nscoe_{k,0}]\theta}{\nscoe_{0,1}}-\frac{[\prod_{k=1}^K\nscoe_{k,0}]\nscoe_{0,0}\theta|\theta|}{(\nscoe_{0,1})^2}+\frac{[\prod_{k=1}^K\nscoe_{k,0}](\nscoe_{0,0})^2\theta^3}{(\nscoe_{0,1})^3}-\frac{[\prod_{k=1}^K\nscoe_{k,0}](\nscoe_{0,0})^3\theta|\theta|^3(\nscoe_{0,1})^{-3}}{\nscoe_{0,0}|\theta|+\nscoe_{0,1}}\\
       &\textstyle+\frac{A\theta^3}{(\nscoe_{0,1})^{3}}
       -\frac{A\theta^3(3(\nscoe_{0,1})^{2}\nscoe_{0,0}|\theta|+3\nscoe_{0,1}(\nscoe_{0,0})^2|\theta|^2+(\nscoe_{0,0})^3|\theta|^3)}{(\nscoe_{0,1})^{6}}\\
&\textstyle+\frac{A(\nscoe_{0,1})^{-6}\theta^3(3(\nscoe_{0,1})^2(\nscoe_{0,0})|\theta|+3\nscoe_{0,1}(\nscoe_{0,0})^2|\theta|^2+(\nscoe_{0,0})^3|\theta|^3)^2}{(\nscoe_{0,1})^3+3(\nscoe_{0,1})^2\nscoe_{0,0}|\theta|+3\nscoe_{0,1}(\nscoe_{0,0})^2|\theta|^2+(\nscoe_{0,0})^3|\theta|^3}+\frac{\theta^4P(\theta(\nscoe_{0,0}|\theta|+\nscoe_{0,1})^{-1})}{(\nscoe_{0,0}|\theta|+\nscoe_{0,1})^{4}}.
       \end{split}
\end{equation}}
\argument{the fact that $\nscoe_{0,0}>0$;the fact that $\nscoe_{0,1}>0$;}{that there exists $\fC\in (0,\infty)$ such that for all $\theta\in \R$ with $|\theta|\leq \scrc$ it holds that
\begin{equation}\llabel{eq4}
\begin{split}
&\biggl|\frac{A\theta^3(3(\nscoe_{0,1})^{2}\nscoe_{0,0}|\theta|+3\nscoe_{0,1}(\nscoe_{0,0})^2|\theta|^2+(\nscoe_{0,0})^3|\theta|^3)}{(\nscoe_{0,1})^{6}}\biggr|\\
&+\biggl|\frac{A(\nscoe_{0,1})^{-6}\theta^3(3(\nscoe_{0,1})^2(\nscoe_{0,0})|\theta|+3\nscoe_{0,1}(\nscoe_{0,0})^2|\theta|^2+(\nscoe_{0,0})^3|\theta|^3)^2}{(\nscoe_{0,1})^3+3(\nscoe_{0,1})^2\nscoe_{0,0}|\theta|+3\nscoe_{0,1}(\nscoe_{0,0})^2|\theta|^2+(\nscoe_{0,0})^3|\theta|^3}\biggr|\\
&\leq |\theta|^4\biggl|\frac{A(3(\nscoe_{0,1})^{2}\nscoe_{0,0}+3\nscoe_{0,1}(\nscoe_{0,0})^2|\theta|+(\nscoe_{0,0})^3|\theta|^2)}{(\nscoe_{0,1})^{6}}\biggr|\\
&+|\theta|^4\biggl|\frac{A(\nscoe_{0,1})^{-6}\theta(3(\nscoe_{0,1})^2(\nscoe_{0,0})+3\nscoe_{0,1}(\nscoe_{0,0})^2|\theta|+(\nscoe_{0,0})^3|\theta|^2)^2}{(\nscoe_{0,1})^3}\biggr|\\
&\leq \fC|\theta|^4.
\end{split}
\end{equation}}
\argument{the fact that $P$ is polynomial;}{that there exists $\fC\in (0,\infty)$ such that for all $x\in \R$ with $|x|\leq (\nscoe_{0,0})^{-1}$ it hold that
\begin{equation}\llabel{eq5}
   |P(x)|\leq \fC.
\end{equation}}
\argument{\lref{eq5};the fact that for all $\theta\in \R$ it holds that $|\theta(\nscoe_{0,0}|\theta|+\nscoe_{0,1})^{-1}|\leq (\nscoe_{0,0})^{-1}$}{that there exist $\fC_1,\fC_2\in (0,\infty)$ such that for all $\theta\in \R$ it holds that
\begin{equation}\llabel{eq6}
\frac{\theta^4|P(\theta(\nscoe_{0,0}|\theta|+\nscoe_{0,1})^{-1})|}{(\nscoe_{0,0}|\theta|+\nscoe_{0,1})^{4}}\leq \fC_1\theta^4(\nscoe_{0,1})^{-4}=\fC_2\theta^4.
\end{equation}}
\argument{\lref{eq6};\lref{eq3};\lref{eq4};the fact that for all $\theta\in \R$ it holds that
\begin{equation}
    \biggl|\frac{[\prod_{k=1}^K\nscoe_{k,0}](\nscoe_{0,0})^3\theta|\theta|^3(\nscoe_{0,1})^{-3}}{\nscoe_{0,0}|\theta|+\nscoe_{0,1}}\biggr|\leq \textstyle[\prod_{k=1}^K\nscoe_{k,0}](\nscoe_{0,0})^3(\nscoe_{0,1})^{-4}|\theta|^4
\end{equation}}{that there exists $\fC\in (0,\infty)$ which satisfies for all $\theta\in \R$ with $|\theta|\leq\scrc$ that
\begin{equation}\llabel{eq7}
    \biggl| \TNewtonSchulzAlgorithm{\nscoe}{K}{\theta}-\biggl(\frac{[\prod_{k=1}^K\nscoe_{k,0}]\theta}{\nscoe_{0,1}}-\frac{[\prod_{k=1}^K\nscoe_{k,0}]\nscoe_{0,0}\theta|\theta|}{(\nscoe_{0,1})^2}+\frac{[\prod_{k=1}^K\nscoe_{k,0}](\nscoe_{0,0})^2\theta^3}{(\nscoe_{0,1})^3}+\frac{A\theta^3}{(\nscoe_{0,1})^{3}}\biggr)\biggr|\leq \fC|\theta|^4.
\end{equation}}
\startnewargseq
\argument{\lref{def: bbX};\lref{def: scrc};the assumption that $X_{n,m}$, $(n,m)\in\N^2$, are \iid\@}{that for all $n\in \N_0$, $M\in \N$ it holds $\P$-a.s.\ that
\begin{equation}\llabel{eq9}
\begin{split}
\textstyle    |\bbX_n^M|&\textstyle\leq \frac 1M\sum_{m=1}^M|X_{n+1,m}|\leq \frac 1M\sum_{m=1}^M\frac{\scrc}{2\lambda}=\frac{\scrc}{2\lambda}.
\end{split}
\end{equation}}
\argument{\lref{eq9};\lref{def: Z};}{for all $M\in \N$ it holds $\P$-a.s.\ that
\begin{equation}\llabel{eq10}
   \textstyle |Z_M|\leq 2\lambda(1-\alpha)\sum_{k=0}^\infty\alpha^k|\bbX_k^M|\leq 2\lambda(1-\alpha)\sum_{k=0}^\infty\frac{\alpha^k\scrc}{2\lambda}=\scrc.
\end{equation}}
\argument{\lref{eq10};\lref{eq7};}{that for all $M\in \N$ it holds $\P$-a.s.\ that
\begin{equation}\llabel{eq11}
    \begin{split}
       & \textstyle \Bigl|\TNewtonSchulzAlgorithm{\nscoe}{K}{Z_M}
          -\Bigl(\frac{[\prod_{k=1}^K\nscoe_{k,0}]Z_M}{\nscoe_{0,1}}-\frac{[\prod_{k=1}^K\nscoe_{k,0}]\nscoe_{0,0}Z_M|Z_M|}{(\nscoe_{0,1})^2}+\frac{[\prod_{k=1}^K\nscoe_{k,0}](\nscoe_{0,0})^2(Z_M)^3}{(\nscoe_{0,1})^3}+\frac{A(Z_M)^3}{(\nscoe_{0,1})^{3}}\Bigr)\Bigr|\\
          &\leq \fC|Z_M|^4.
    \end{split}
\end{equation}}
\argument{\lref{eq11};}{that for all $M\in\N$ it holds that
\begin{equation}\llabel{eq12}
     \begin{split}
          &\textstyle\Bigl|\E[\TNewtonSchulzAlgorithm{\nscoe}{K}{Z_M}]\\
          &\quad \textstyle-\Bigl(\frac{[\prod_{k=1}^K\nscoe_{k,0}]\E[Z_M]}{\nscoe_{0,1}}-\frac{[\prod_{k=1}^K\nscoe_{k,0}]\nscoe_{0,0}\E[Z_M|Z_M|]}{(\nscoe_{0,1})^2}+\frac{[\prod_{k=1}^K\nscoe_{k,0}](\nscoe_{0,0})^2\E[(Z_M)^3]}{(\nscoe_{0,1})^3}+\frac{A\E[(Z_M)^3]}{(\nscoe_{0,1})^{3}}\Bigr)\Bigr|\\
           &\leq \fC\E[(Z_M)^4].
    \end{split}
\end{equation}}
\argument{\lref{def: bbX};\lref{def: Z};the fact that for all $n,m\in \N$ it holds that $\E[X_{n,m}]=0$}{that for all $M\in \N$ it holds that
\begin{equation}\llabel{eq13}
\begin{split}
    \textstyle \E[Z_M]=-2\lambda(1-\alpha)\sum\limits_{k=0}^\infty\alpha^k\E[\bbX_k^M]=-2\lambda(1-\alpha)\sum\limits_{k=0}^\infty\alpha^k\bigl[\frac 1M\sum_{m=1}^M\E[X_{k+1,m}]\bigr]=0.
    \end{split}
\end{equation}}
\argument{\lref{eq12};\lref{eq13};}{that for all $M\in\N$ it holds that
\begin{equation}\llabel{eq14}
    \begin{split}
          \textstyle\Bigl| \E[\TNewtonSchulzAlgorithm{\nscoe}{K}{Z_M}]
    +\Bigl(\frac{[\prod_{k=1}^K\nscoe_{k,0}]\nscoe_{0,0}\E[Z_M|Z_M|]}{(\nscoe_{0,1})^2}-\frac{[\prod_{k=1}^K\nscoe_{k,0}](\nscoe_{0,0})^2\E[(Z_M)^3]}{(\nscoe_{0,1})^3}-\frac{A\E[(Z_M)^3]}{(\nscoe_{0,1})^{3}}\Bigr)\Bigl|\leq \fC\,\E[(Z_M)^4].
    \end{split}
\end{equation}}
\argument{\lref{eq14};}{that for all $M\in\N$ it holds that
\begin{equation}\llabel{eq15}
    \begin{split}
\textstyle\bigl|\E[\TNewtonSchulzAlgorithm{\nscoe}{K}{Z_M}]\bigr|&\textstyle\geq \frac{[\prod_{k=1}^K\nscoe_{k,0}]\nscoe_{0,0}|\E[Z_M|Z_M|]|}{(\nscoe_{0,1})^2}-\frac{([\prod_{k=1}^K\nscoe_{k,0}](\nscoe_{0,0})^2+|A|)|\E[(Z_M)^3]|}{(\nscoe_{0,1})^{3}}-\fC\,\E[(Z_M)^4].
    \end{split}
\end{equation}}
\argument{\lref{def: Z};\cref{lem: tg2};the fact that $0\leq\alpha<1$}{that
\begin{equation}\llabel{evd1}
     \liminf_{M\to\infty}\bigl|M^{3/2}\E\bigl[Z_M|Z_M|\bigr]\bigr|>0.
\end{equation}}
\argument{\cref{cor: upper order of Adam coefficient};}{that for all $p\in \N\backslash\{1\}$, $r\in \N$ it holds that
\begin{equation}\llabel{eq18'}
    \textstyle \limsup\limits_{M\to\infty}\Bigl|M^{\nicefrac{(2r+1)}{2}}\,\E\bigl[\textstyle\bigl(2\lambda(1-\alpha)\sum_{k=0}^\infty \alpha^k[\frac{1}{M}\sum_{m=1}^MX_{k+1,m}]\bigr)^{2r+1}\bigr]\Bigr|=0.
\end{equation}
\begin{equation}\llabel{eq18}
   \text{and}\qquad \textstyle \limsup\limits_{M\to\infty}\Bigl(M^{\nicefrac{p}{2}}\,\E\Bigl[\textstyle\bigl|2\lambda(1-\alpha)\sum_{k=0}^\infty \alpha^k[\frac{1}{M}\sum_{m=1}^MX_{k+1,m}]\bigr|^{p}\Bigr]\Bigr)<\infty.
\end{equation}}
\argument{\lref{def: Z};\lref{eq18}}{that
\begin{equation}\llabel{eq19}
     \limsup_{M\to\infty}\bigl|M^{3/2}\E[(Z_M)^3]\bigr|=0\qqandqq  \limsup_{M\to\infty}\bigl(M^{2}\E[(Z_M)^4]\bigr)<\infty.
\end{equation}}
\argument{\lref{eq19};}{that
\begin{equation}\llabel{evd2}
    \limsup_{M\to\infty}\bigl|M^{3/2}\E\bigl[(Z_M)^3\bigr]\bigr|=0\qqandqq  \limsup_{M\to\infty}\bigl(M^{3/2}\E[(Z_M)^4]\bigr)=0.
\end{equation}}
\argument{\lref{evd2};\lref{eq15};\lref{evd1}}{that
\begin{equation}\llabel{eq20}
    \begin{split}
&\liminf_{M\to\infty}\bigl|M^{3/2}\E[\TNewtonSchulzAlgorithm{\nscoe}{K}{Z_M}]\bigr|\\&\geq \frac{\liminf_{M\to\infty}\big|M^{3/2}[\prod_{k=1}^K\nscoe_{k,0}]\nscoe_{0,0}\E\bigl[Z_M|Z_M|\bigr]\bigr|}{(\nscoe_{0,1})^2}\\
&-\frac{\limsup_{M\to\infty}\bigl|M^{3/2}([\prod_{k=1}^K\nscoe_{k,0}](\nscoe_{0,0})^2+|A|)\E[(Z_M)^3]\bigr|}{(\nscoe_{0,1})^3}
-\textstyle\limsup_{M\to\infty}\bigl(M^{3/2}\fC\E[(Z_M)^4]\bigr)\\
&>0.
    \end{split}
\end{equation}}
\argument{\lref{eq20};}{\cref{conclude: non-zero vector field 2}\dott}
\end{aproof}
\begin{tcolorbox}[colback=white!95!gray,
                  colframe=black,
                  boxrule=0.5pt,
                  sharp corners,
                  enhanced,
                  breakable,
                 ]
\begin{athm}{cor}{non-zero vector field 3}[\textcolor{red}{Non-vanishing evaluations of the \MUON\ vector field for non-symmetric data}]
Let $(\Omega,\cF,\P)$ be a probability space, let $K\in\N_0$, $\alpha\in [0,1)$, $\lambda\in (0,\infty)$, let $\nscoe=(\nscoe_{i,j})_{(i,j)\in(\N_0)^2}\colon (\N_0)^2\to\R$ satisfy $\textcolor{magenta}{\min\{\nscoe_{0,0},\nscoe_{0,1}\}>0}$ and $\textcolor{magenta}{\#(\nscoe^{ - 1 }( \R\backslash\{0\} )) <\infty}$, assume for all $i\in\N\cap [0,K]$ that $\textcolor{magenta}{\nscoe_{i,0}>0}$, let $X_{n,m}\colon \Omega\to\R$, $(n,m)\in\N^2$, be bounded \iid\ random variable, and assume $\E[(X_{1,1}-\E[X_{1,1}])^3]\neq 0$.
Then 
\begin{equation}\llabel{conclude}
\textcolor{magenta}{\liminf\nolimits_{M\to\infty}\bigl| M^{3/2}\E\bigl[ \TbigNewtonSchulzAlgorithm{\nscoe}{K}{\textstyle\textstyle-2\lambda(1-\alpha)\bigl[\sum_{k=0}^\infty\alpha^{k}[\textstyle \frac 1M\sum_{m=1}^M(X_{k+1,m}-\E[X_{1,1}])\bigr]\bigr]}\bigr]\bigr|\neq 0}
\end{equation}
(cf.\ \cref{definition: NS}).
\end{athm}
\end{tcolorbox}
\begin{aproof}
    \argument{\cref{non-zero vector field 2}}{\lref{conclude}\dott}
\end{aproof}
\subsection{Non-convergence of MUON}\label{subsec: muon non convergence}

In the next elementary result, \cref{NSLipschitz}, we show for every $d, \fd \in \N$, $\varepsilon\in  (0, \infty)$, $n \in \N_0$ and appropriate $a, b, c \in \R$ that the function $\R^{ d \times \fd } \ni A \mapsto \TNewtonSchulzAlgorithm{\nscoe}{n}{A} \in \R^{d\times\fd}$ (cf.\ \cref{definition: NS}) is globally Lipschitz continuous.
\begin{samepage}
\begin{tcolorbox}[colback=white!95!gray,
                  colframe=black,
                  boxrule=0.5pt,
                  sharp corners,
                  enhanced,
                  breakable,
                 ]
\begin{athm}{lemma}{NSLipschitz}[\textcolor{red}{Lipschitz continuity for the \NS\ functions}]
 Let $d,\fd\in \N$, $K\in\N_0$, let $\nscoe=(\nscoe_{i,j})_{(i,j)\in (\N_0)^2}\colon(\N_0)^2\to \R$ satisfy $\textcolor{magenta}{\min\{\nscoe_{0,0},\nscoe_{0,1}\}>0}$ and $\textcolor{magenta}{\#(\nscoe^{ - 1 }( \R\backslash\{0\} )) <\infty}$, and assume for all $i\in \N\cap[0,K]$, $x\in \R$ that $\textcolor{magenta}{\sum_{j=0}^\infty\nscoe_{i,j}x^{2j}>0}$.
Then for all $n\in \{0,1,\dots,K\}$ there exists $\fC\in \R$ such that for all $A,B\in \R^{d\times\fd}$ it holds that
\begin{equation}\llabel{conclude}
   \lVert \TNewtonSchulzAlgorithm{\nscoe}{n}{A}-\TNewtonSchulzAlgorithm{\nscoe}{n}{B}\rVert\leq \fC\|A-B\|
\end{equation}
(cf.\ \cref{definition: NS}).
\end{athm}
\end{tcolorbox}
\end{samepage}
\begin{aproof}
Throughout this proof for every $n\in \N_0$, $A \in \R^{d \times \fd}$ let
$
	\NewtonSchulzAlgorithm{A}{n}\in \R^{d \times \fd}
$
 satisfy 
\begin{equation}
\llabel{def: NS}
\begin{split} 
	\NewtonSchulzAlgorithm{A}{n}=\TNewtonSchulzAlgorithm{\nscoe}{n}{A}.
\end{split}
\end{equation}
    \argument{\cref{NS error};}{that there exists $\scrc\in (0,\infty)$ such that for all $n\in\{0,1,\dots,K\}$ it holds that
    \begin{equation}\llabel{eq1}
        \|\NewtonSchulzAlgorithm{A}{n}\|\leq \scrc.
    \end{equation}}
    \argument{\lref{eq1};the fact that for all $A\in \R^{d\times\fd}$ it holds that $\|\NewtonSchulzAlgorithm {A}{0}\|=(\nscoe_{0,0}\|A\|+\nscoe_{0,1})^{-1}\|A\|\leq (\nscoe_{0,0})^{-1}$}{that there exists $\scrc\in (0,\infty)$ which satisfies for all $n\in\{0,1,\dots,K\}$, $A\in\R^{d\times\fd}$  that
    \begin{equation}\llabel{def: scrc}
        \|\NewtonSchulzAlgorithm{A}{n}\|\leq \scrc.
    \end{equation}}
    \startnewargseq
    We prove \lref{conclude} by induction on $n\in \{0,1,\dots,K\}$. For the base case $n=0$ note that \lref{def: NS} shows that there exists $\fC\in (0,\infty)$ such that for all $A,B\in \R^{d\times\fd}$ it holds that
    \begin{equation}\llabel{base}
    \begin{split}
       \textstyle \|\NewtonSchulzAlgorithm{A}{0}-\NewtonSchulzAlgorithm{B}{0}\|&\textstyle=\bigl\|\frac{A}{\nscoe_{0,0}\|A\|+\nscoe_{0,1}}-\frac{B}{\nscoe_{0,0}\|B\|+\nscoe_{0,1}}\bigr\|\\
       &\textstyle\leq \bigl\|\frac{A}{\nscoe_{0,0}\|A\|+\nscoe_{0,1}}-\frac{B}{\nscoe_{0,0}\|A\|+\nscoe_{0,1}}\bigr\|+\bigl\|\frac{B}{\nscoe_{0,0}\|A\|+\nscoe_{0,1}}-\frac{B}{\nscoe_{0,0}\|B\|+\nscoe_{0,1}}\bigr\|\\
       &=\textstyle \bigl\|\frac{A-B}{\nscoe_{0,0}\|A\|+\nscoe_{0,1}}\bigr\|+\bigl\|\frac{\nscoe_{0,0}B(\|A\|-\|B\|)}{(\nscoe_{0,0}\|A\|+\nscoe_{0,1})(\nscoe_{0,0}\|B\|+\nscoe_{0,1})}\bigr\|\\
       &\leq\textstyle  \|A-B\|(\nscoe_{0,1})^{-1}+\frac{\nscoe_{0,0}\|B\|\|A-B\|}{(\nscoe_{0,0}\|A\|+\nscoe_{0,1})(\nscoe_{0,0}\|B\|+\nscoe_{0,1})}\\
       &\leq\|A-B\|(\nscoe_{0,1})^{-1}+\|A-B\|(\nscoe_{0,1})^{-1}\leq \fC\|A-B\|.
        \end{split}
    \end{equation}
    This establishes \lref{conclude} in the base case $n=0$. For the induction step we assume that there exist $n\in \N\cap(0,K]$, $\fC\in (0,\infty)$ which satisfies for all $A,B\in \R^{d\times\fd}$ that
    \begin{equation}\llabel{assume}
      \|  \NewtonSchulzAlgorithm{A}{n-1}-\NewtonSchulzAlgorithm{B}{n-1}\|\leq \fC\|A-B\|.
    \end{equation}
    In the following we prove that for all $i\in \N_0$ there exists $\scrC\in (0,\infty)$ such that for all $A,B\in\R^{d\times\fd}$ it holds that
    \begin{equation}\llabel{need to prove}
       \| (\NewtonSchulzAlgorithm{A}{n-1} (\NewtonSchulzAlgorithm{A}{n-1})^{\top})^i\NewtonSchulzAlgorithm{A}{n-1} -(\NewtonSchulzAlgorithm{B}{n-1} (\NewtonSchulzAlgorithm{B}{n-1})^{\top})^i\NewtonSchulzAlgorithm{B}{n-1} \|\leq \scrC\|A-B\|.
    \end{equation}
    \startnewargseq
    We prove \lref{need to prove} by induction on $i\in\N_0$. Note that \lref{assume} ensures \lref{need to prove} in the base case $i=0$. For the induction step we assume that there exist $i\in\N$, $\scrC\in (0,\infty)$ which satisfy for all $A,B\in \R^{d\times\fd}$ that
    \begin{equation}\llabel{assume 2}
        \| (\NewtonSchulzAlgorithm{A}{n-1} (\NewtonSchulzAlgorithm{A}{n-1})^{\top})^{i-1}\NewtonSchulzAlgorithm{A}{n-1} -(\NewtonSchulzAlgorithm{B}{n-1} (\NewtonSchulzAlgorithm{B}{n-1})^{\top})^{i-1}\NewtonSchulzAlgorithm{B}{n-1} \|\leq \scrC\|A-B\|.
    \end{equation}
    \argument{\lref{def: scrc};\lref{assume};\lref{assume 2};the Cauchy-Schwarz inequality}{that for all $A,B\in \R^{d\times\fd}$ it holds that
    \begin{equation}\llabel{evd1}
    \begin{split}
       & \|(\NewtonSchulzAlgorithm{A}{n-1} (\NewtonSchulzAlgorithm{A}{n-1})^{\top})^i\NewtonSchulzAlgorithm{A}{n-1} -(\NewtonSchulzAlgorithm{B}{n-1} (\NewtonSchulzAlgorithm{B}{n-1})^{\top})^i\NewtonSchulzAlgorithm{B}{n-1} \|\\
       &\leq \|(\NewtonSchulzAlgorithm{A}{n-1} (\NewtonSchulzAlgorithm{A}{n-1})^{\top})^{i-1}\NewtonSchulzAlgorithm{A}{n-1}(\NewtonSchulzAlgorithm{A}{n-1})^{\top}\NewtonSchulzAlgorithm{A}{n-1} -(\NewtonSchulzAlgorithm{B}{n-1} (\NewtonSchulzAlgorithm{B}{n-1})^{\top})^{i-1}\NewtonSchulzAlgorithm{B}{n-1}(\NewtonSchulzAlgorithm{A}{n-1})^{\top}\NewtonSchulzAlgorithm{A}{n-1} \|\\
       &+\|(\NewtonSchulzAlgorithm{B}{n-1} (\NewtonSchulzAlgorithm{B}{n-1})^{\top})^{i-1}\NewtonSchulzAlgorithm{B}{n-1}(\NewtonSchulzAlgorithm{A}{n-1})^{\top}\NewtonSchulzAlgorithm{A}{n-1} -(\NewtonSchulzAlgorithm{B}{n-1} (\NewtonSchulzAlgorithm{B}{n-1})^{\top})^{i-1}\NewtonSchulzAlgorithm{B}{n-1}(\NewtonSchulzAlgorithm{B}{n-1})^{\top}\NewtonSchulzAlgorithm{A}{n-1} \|\\
       &+\|(\NewtonSchulzAlgorithm{B}{n-1} (\NewtonSchulzAlgorithm{B}{n-1})^{\top})^{i-1}\NewtonSchulzAlgorithm{B}{n-1}(\NewtonSchulzAlgorithm{B}{n-1})^{\top}\NewtonSchulzAlgorithm{A}{n-1} -(\NewtonSchulzAlgorithm{B}{n-1} (\NewtonSchulzAlgorithm{B}{n-1})^{\top})^{i-1}\NewtonSchulzAlgorithm{B}{n-1}(\NewtonSchulzAlgorithm{B}{n-1})^{\top}\NewtonSchulzAlgorithm{B}{n-1} \|\\
       &\leq \scrC \scrc^2\|A-B\|+  \fC \scrc^{2i}\|A-B\|+\fC \scrc^{i+1}\|A-B\|\\
       &= (\scrC\scrc^2+2\fC \scrc^{2i})\|A-B\|.
       \end{split}
    \end{equation}}
     \argument{\lref{evd1};\lref{assume 2};induction}{\lref{need to prove}\dott}
     \startnewargseq
     \argument{the assumption that $\#(\nscoe^{ - 1 }( \R\backslash\{0\} )) \allowbreak<\infty$}{that \llabel{argg1} $\#\{i\in\N\colon \nscoe_{k,i}\neq 0\}\allowbreak<\infty$\dott}
    \argument{\lref{argg1};\cref{T_B_D};\lref{def: scrc};\lref{assume};\lref{need to prove}}{that there exist $M\in\N$ and $\scrC_i$, $i\in\{1,2,\dots,M+1\}$, such that for all $A,B\in \R^{d\times\fd}$ it holds that
    \begin{equation}\llabel{evd2}
        \begin{split}
          \|\NewtonSchulzAlgorithm{A}{n} -\NewtonSchulzAlgorithm{B}{n} \|&=\bigl\|\bigl(\textstyle\sum_{ i = 0 }^{ \infty } \nscoe_{ n, i } ( \NewtonSchulzAlgorithm{A}{n-1} (\NewtonSchulzAlgorithm{A}{n-1})^{\top} )^{ i} \NewtonSchulzAlgorithm{A}{n-1}\bigr)-\bigl(\textstyle\sum_{ i = 0 }^{ \infty } \nscoe_{ k, i } ( \NewtonSchulzAlgorithm{B}{n-1} (\NewtonSchulzAlgorithm{B}{n-1})^{\top} )^{ i} \NewtonSchulzAlgorithm{B}{n-1}\bigr)\bigr\|\\
          &=\bigl\|\bigl(\textstyle\sum_{ i = 0 }^{ M } \nscoe_{ n, i } ( \NewtonSchulzAlgorithm{A}{n-1} (\NewtonSchulzAlgorithm{A}{n-1})^{\top} )^{ i} \NewtonSchulzAlgorithm{A}{n-1}\bigr)-\bigl(\textstyle\sum_{ i = 0 }^{ M } \nscoe_{ k, i } ( \NewtonSchulzAlgorithm{B}{n-1} (\NewtonSchulzAlgorithm{B}{n-1})^{\top} )^{ i} \NewtonSchulzAlgorithm{B}{n-1}\bigr)\bigr\|\\
          &\leq \textstyle\sum_{ i = 0 }^{ M } |\nscoe_{ n, i }|\|(  \NewtonSchulzAlgorithm{A}{n-1} (\NewtonSchulzAlgorithm{A}{n-1})^{\top} )^{ i} \NewtonSchulzAlgorithm{A}{n-1}-  ( \NewtonSchulzAlgorithm{B}{n-1} (\NewtonSchulzAlgorithm{B}{n-1})^{\top} )^{ i} \NewtonSchulzAlgorithm{B}{n-1}\|\\
          &\leq \textstyle\sum_{i=0}^M\scrC_i \|A-B\|=\scrC_{M+1}\|A-B\|.
        \end{split}
    \end{equation}}
    \argument{\lref{evd2};\lref{assume};induction}{\lref{conclude}
    \dott}
\end{aproof}
\begin{samepage}
\begin{tcolorbox}[colback=white!95!gray,
                  colframe=black,
                  boxrule=0.5pt,
                  sharp corners,
                  enhanced,
                  breakable,
                 ]
\begin{athm}{theorem}{MUON non-convergence 1d big batch}[\textcolor{red}{Non-convergence of \MUON}]
     Let $(\Omega,\cF,\P)$ be a probability space, let $K\in\N_0$, $\alpha\in (0,1)$,  $\xi,\lambda\in (0,\infty)$, let $\nscoe=(\nscoe_{i,j})_{(i,j)\in (\N_0)^2}\colon(\N_0)^2\to \R$ satisfy $\textcolor{magenta}{\min\{\nscoe_{0,0},\nscoe_{0,1}\}>0}$ and $\textcolor{magenta}{\#(\nscoe^{ - 1 }( \R\backslash\{0\} )) <\infty}$,  assume for all $i\in \N\cap[0,K]$, $x\in \R$ that $\textcolor{magenta}{\sum_{j=0}^\infty\nscoe_{i,j}x^{2j}>0}$, let $\smalll=(\smalll(\theta,x))_{(\theta,x)\in \R\times \R} \allowbreak \colon \R\times \R\to\R$ satisfy for all  $\theta,x\in \R$ that $\smalll(\theta,x)=\lambda|\theta-x|^2$, let $X_{n,m}\colon \Omega\to\R$, $(n,m)\in \N^2$, be bounded \iid\ random variables, assume $\E[(X_{1,1}-\E[X_{1,1}])^3]\neq 0$, let $(\gamma_n)_{n\in \N}\subseteq (0,\infty)$ be non-increasing, assume $\limsup_{ n \to \infty } ( ( \gamma_{n+1} )^{ - 1 } \gamma_n ) < \alpha^{-1} < \sum_{ n=1 }^{ \infty } \gamma_n = \infty$, and for every $M\in \N$ let $\Theta^M \colon\N_0\times\Omega\to\R$ and $\bfm^M\colon \N_0\times\Omega\to\R$ be stochastic processes which satisfy for all $n\in \N$ that 
        \begin{equation}\llabel{def: bfm}
           \textcolor{magenta}{ \bfm_0^M=0},\qquad \textcolor{magenta}{\bfm_n^M= \alpha \bfm_{n-1}^M+(1-\alpha)\bigl[\textstyle \frac 1M \sum_{m=1}^M(\nabla_\theta\smalll)(\Theta_{n-1}^M,X_{n,m})\bigr]},
            \end{equation}
            \begin{equation}\llabel{def: Theta}
           \textcolor{magenta}{|\Theta_0^M|\leq \xi},\qquad \text{and}\qquad \textcolor{magenta}{\Theta_n^{M}=\Theta_{n-1}^{M}-\gamma_n\TNewtonSchulzAlgorithm{\nscoe}{K}{\bfm_n^{M}}}.
         \end{equation}
         Then there exists $\bfM\in \N$ such that for all $M\in \N\cap[\bfM,\infty)$ it holds that
         \begin{equation}\llabel{conclude}
           \textcolor{magenta}{\textstyle\limsup_{ n \to \infty }\E\bigl[ \min\{ 1, | \Theta_n^M - \E[X_{1,1}]|\} \bigr] > 0}.
           \end{equation}
  (cf.\ \cref{definition: NS}).
\end{athm}
\end{tcolorbox}
\end{samepage}
\begin{aproof}
\argument{the assumption that for all $i\in \N\cap[0,K]$, $x\in \R$ it holds that $\sum_{j=0}^\infty\nscoe_{i,j}x^{2j}>0$}{that for all \llabel{arg1} $i\in \N\cap[0,K]$ it holds that $\nscoe_{i,0}>0$\dott}
    \argument{\lref{arg1};\cref{non-zero vector field 3};}{that there exists $\bfM\in \N$ which satisfies for all $M\in \N\cap[\bfM,\infty)$ that
\begin{equation}\llabel{def: X}
  \E\Bigl[ \TbigNewtonSchulzAlgorithm{\nscoe}{K}{\textstyle\textstyle-2\lambda(1-\alpha)\bigl[\sum_{k=0}^\infty\alpha^{k}[\textstyle \frac 1M\sum_{m=1}^M(X_{k+1,m}-\E[X_{1,1}])\bigr]\bigr]}\Bigr]\neq 0.
\end{equation}}
\startnewargseq
In the following we prove that for all $M\in \N\cap[\bfM,\infty)$ it holds that
\begin{equation}\llabel{need to prove}
     \limsup_{ n \to \infty }\E\bigl[ \min\{ 1, | \Theta_n^M - \E[X_{1,1}]|\} \bigr] > 0.
\end{equation}
We prove \lref{need to prove} by contradiction. In the following we thus assume that there exists $M\in\N\cap[\bfM,\infty)$ which satisfies that
\begin{equation}\llabel{assume}
    \limsup_{ n \to \infty }\E\bigl[ \min\{ 1, | \Theta_n^M - \E[X_{1,1}]|\} \bigr] = 0.
\end{equation}
\argument{\lref{def: bfm};\lref{def: Theta};the assumption that for all $\theta,x\in \R$ it holds that $\smalll(\theta,x)=\lambda|\theta-x|^2$;\cref{theo: MUON bound}}{that there exists $\fC\in (0,\infty)$ such that
\begin{equation}\llabel{eq1}
  \P\bigl( \textstyle \sup_{n\in \N_0}|\Theta_n^M|\leq \fC\bigr)=1.
\end{equation}}
\argument{\lref{assume};}{that for all $\delta\in (0,\infty)$ it holds that
\begin{equation}\llabel{eq2}
    \textstyle\limsup_{n\to\infty}\P(|\Theta_n^M-\E[X_{1,1}]|\geq \delta)=0.
\end{equation}}
\argument{\lref{eq2};\lref{eq1};}{that
\begin{equation}\llabel{eq4}
     \textstyle \limsup_{ n \to \infty }\E\bigl[  | \Theta_n^M - \E[X_{1,1}]| \bigr] = 0.
\end{equation}}
\argument{\cref{T_B_D};\cref{NS error};\cref{NSLipschitz}}{that there exists $\fC\in (0,\infty)$ such that for all $x,y\in \R$ it holds that
\begin{equation}\llabel{eq5}
    |\TNewtonSchulzAlgorithm{\nscoe}{K}{x}|\leq \fC\qqandqq |\TNewtonSchulzAlgorithm{\nscoe}{K}{x}-\TNewtonSchulzAlgorithm{\nscoe}{K}{y}|\leq \fC|x-y|.
\end{equation}}
\argument{\lref{eq5};\lref{def: bfm};\lref{def: Theta};\lref{eq4};\cref{vector field}}{that
\begin{equation}\llabel{eq6}
   \E\Bigl[ \TbigNewtonSchulzAlgorithm{\nscoe}{K}{\textstyle\textstyle(1-\alpha)\sum_{k=0}^\infty\alpha^k\bigl[\frac {2\lambda}M\sum_{m=1}^M(\E[X_{1,1}]-X_{k+1,m})\bigr]}\Bigr]=0.
\end{equation}}
\argument{\lref{eq6};\lref{def: X}}{that
\begin{equation}\llabel{eq7}
\begin{split}
    0 &= \E\Bigl[ \TbigNewtonSchulzAlgorithm{\nscoe}{K}{\textstyle\textstyle(1-\alpha)\sum_{k=0}^\infty\alpha^k\bigl[\frac {2\lambda}M\sum_{m=1}^M(\E[X_{1,1}]-X_{k+1,m})\bigr]}\Bigr]\\
   &\textstyle=  \E\Bigl[ \TbigNewtonSchulzAlgorithm{\nscoe}{K}{\textstyle\textstyle-2\lambda(1-\alpha)\bigl[\sum_{k=0}^\infty\alpha^{k}[\textstyle \frac 1M\sum_{m=1}^M(X_{k+1,m}-\E[X_{1,1}])\bigr]\bigr]}\Bigr]\neq 0.
    \end{split}
\end{equation}}
This contradiction proves \lref{need to prove}\dott
\end{aproof}
\subsubsection*{Acknowledgements}
We thank Benno Kuckuck and Philippe von Wurstemberger for useful suggestions regarding the presentation of this article. This work has been partially supported by the National Natural Science Foundation of China (NSFC) under grant number W2531010. We also gratefully acknowledge the Cluster of Excellence EXC 2044/2-390685587, Mathematics Münster: Dynamics-Geometry-Structure funded by the Deutsche Forschungsgemeinschaft (DFG, German Research Foundation). Most of the specific formulations in the proofs of this work have been created using \cite{Bennoargumentcommand}. 
\bibliographystyle{acm}
\bibliography{bibfileMUONconvergenceI}

\end{document}